\documentclass[12pt,a4paper,twoside]{book}
\usepackage[french] {babel}
\usepackage{amssymb,amsmath,amsfonts,graphics}
\usepackage{stmaryrd}
\usepackage[all]{xy}

\def\rondI{\newbox\boxx{\hbox{$1$}}\hskip-11pt\raisebox{0pt}{$\displaystyle{\textrm{\Large$\bigcirc$}}$}}
\def\rondII{\newbox\boxx{\hbox{$2$}}\hskip-11pt\raisebox{0pt}{$\displaystyle{\textrm{\Large$\bigcirc$}}$}}
\def\rondIII{\newbox\boxx{\hbox{$3$}}\hskip-11pt\raisebox{0pt}{$\displaystyle{\textrm{\Large$\bigcirc$}}$}}
\def\rondIV{\newbox\boxx{\hbox{$4$}}\hskip-11pt\raisebox{0pt}{$\displaystyle{\textrm{\Large$\bigcirc$}}$}}
\def\rondV{\newbox\boxx{\hbox{$5$}}\hskip-11pt\raisebox{0pt}{$\displaystyle{\textrm{\Large$\bigcirc$}}$}}
\def\rondVI{\newbox\boxx{\hbox{$6$}}\hskip-11pt\raisebox{0pt}{$\displaystyle{\textrm{\Large$\bigcirc$}}$}}

\begin{document}

\newcommand{\hooklongrightarrow}{\lhook\joinrel\longrightarrow}

\thispagestyle{empty}
\begin{center}
\LARGE\textbf{IMAGES DIRECTES\\
 ET FONCTIONS $L$ \\
  EN COHOMOLOGIE RIGIDE}

\vskip30mm

Jean-Yves ETESSE
  \footnote{(CNRS - Institut de Math\'ematique, Universit\'e de Rennes 1, Campus de Beaulieu - 35042 RENNES Cedex France)\\
E-mail : Jean-Yves.Etesse@univ-rennes1.fr}
\end{center}
 
\cleardoublepage

Classification AMS: 11G, 11S40, 13B22, 13B35, 13B40, 14F20, 14F30, 14G10, 14G22, 14K, 14L05, 14L15.\\

Mots cl\'es: Alg\`ebres de Monsky-Washnitzer, ($F$)-cristaux, ($F$)-isocristaux convergents ou surconvergents, cohomologie cristalline, cohomologie rigide, cohomologie \'etale, cohomologie syntomique, fonctions $L$, sch\'ema ab\'elien ordinaire, groupe $p$-divisible ordinaire, cristal de Dieudonn\'e.

\newpage
\strut
\newpage

\noindent\textbf{R\'esum\'e}\\

Soient $k$ un corps parfait de caract\'eristique $p> 0$, $\mathcal{V}$ un anneau de valuation discr\`ete complet de corps r\'esiduel $k$ et corps des fractions $K$ de caract\'eristique 0, et $S$ un $k$-sch\'ema s\'epar\'e de type fini.\\

 Lorsque $S$ est lisse sur $k$, nous prouvons ici partiellement une conjecture de Berthelot portant sur la surconvergence des images directes du faisceau structural par un morphisme $f:X\rightarrow S$
propre et lisse; pour $k$ parfait et $\mathcal{V}$ mod\'er\'ement ramifi\'e de telles images directes sont toujours convergentes non seulement pour le faisceau structural mais aussi pour (presque) tous les $F$-isocristaux convergents. Plus g\'en\'eralement, nous prouvons cette surconvergence lorsque $f$ est relevable sur $\mathcal{V}$, ou que $X$ est une intersection compl\`ete relative dans des espaces projectifs sur $S$, et en prenant pour coefficients des isocristaux surconvergents quelconques. \\

 Nous appliquons ensuite ces r\'esultats aux fonctions $L$ avec pour coefficients de telles images directes avec structure de Frobenius: nous en d\'eduisons des propri\'et\'es de rationalit\'e ou de m\'eromorphie pour ces fonctions $L$ (conjecture de Dwork), et nous \'etudions leurs z\'eros et p\^oles unit\'es $p$-adiques (conjecture de Katz); un cas explicite concerne les sch\'emas ab\'eliens ordinaires.\\

 Un expos\'e plus pr\'ecis des r\'esultats par chapitres est fourni dans l'introduction.\\
  
\noindent\textbf{Abstract}\\

Let $k$ be a perfect field of characteristic $p>0$, $\mathcal{V}$ a complete discrete valuation ring with residue field $k$ and field of fractions $K$ of characteristic 0, and $S$ a separated $k$-scheme of finite type.\\

When $S$ is smooth over $k$, we partially prove here a conjecture of Berthelot about the overconvergence of the higher direct images of the structure sheaf under a proper smooth morphism $f:X\rightarrow S$; when $k$ is perfect and $\mathcal{V}$ is tamely ramified such direct images are always convergent, not only for the structure sheaf but also for (almost) every convergent $F$-isocrystals. More generally, we prove this overconvergence when $f$ is liftable over $\mathcal{V}$, or when $X$ is a relative complete intersection in some projective spaces over $S$, and taking as coefficients any overconvergent isocrystals. \\

We then apply these results to $L$-functions with coefficients such direct images with Frobenius structure: we derive rationality or meromorphy for these $L$-functions (Dwork's conjecture), and we study their $p$-adic unit zeroes and poles (Katz's conjecture) ; and explicit case concerns the ordinary abelian schemes.\\

A more precise presentation of results by chapters is given in the introduction.

\newpage 

\noindent\textbf{Sommaire}\\

  \begin{enumerate}
 	\item[0.] Introduction\dotfill 9
		\begin{enumerate}
		\item[0.1.] Rel\`evements et alg\`ebres de Monsky-Washnitzer
		\item[0.2.] Espaces rigides analytiques et images directes
		\item[0.3.] $F$-isocristaux convergents sur un sch\'ema lisse, et images directes
		\item[0.4.] Images directes de $F$-isocristaux surconvergents
		\item[0.5.] Cohomologie syntomique
		\item[0.6.] Fonctions $L$
		\end{enumerate}
	\item[I.] Rel\`evements et alg\`ebres de Monsky-Washnitzer\dotfill17
		\begin{enumerate}
		\item[1.] G\'en\'eralit\'es\dotfill17
		\item[2.] Des \'equivalences de cat\'egories\dotfill23
		\item[3.] Sch\'emas formels et rel\`evements de sch\'emas\dotfill 26
		\end{enumerate}
	\item[II.] Espaces rigides analytiques et images directes\dotfill 47
		\begin{enumerate}
		\item[1.] Changement de base pour un morphisme propre\dotfill 48
		\item[2.] Sorites sur les voisinages stricts\dotfill 58
		\item[3.] Images directes d'isocristaux\dotfill 70
			\begin{enumerate}
			\item[3.1.] Sections surconvergentes\dotfill 70
			\item[3.2.] D\'efinition des images directes\dotfill 73
			\item[3.3.] Changement de base\dotfill 74
			\item[3.4.] Surconvergence des images directes\dotfill 78
			\end{enumerate}
		\end{enumerate}
	\item[III.] $F$-isocristaux convergents sur un sch\'ema lisse, et images directes\dotfill 99
		\begin{enumerate}
		\item[1.] $F$-isocristaux convergents sur un sch\'ema affine et lisse\dotfill 99
			\begin{enumerate}
			\item[1.1.] Notations\dotfill 99
			\item[1.2.] Des \'equivalences de cat\'egories\dotfill 102
			\end{enumerate}
		\item[2.] $F$-isocristaux convergents sur un sch\'ema lisse formellement\\ relevable\dotfill 106
			\begin{enumerate}
			\item[2.1.] Espaces rigides associ\'es aux sch\'emas formels\dotfill 106
			\item[2.2.] $F$-isocristaux convergents et $\mathcal{O}_{\mathcal{X}_{K}}$-modules\dotfill 108
			\end{enumerate}
		\item[3.] Images directes de $F$-isocristaux convergents\dotfill 110
			\begin{enumerate}
			\item[3.1.] Rel\`evements de Teichm\" uller\dotfill 110
			\item[3.2.] Amplitude des $F$-isocristaux\dotfill 112
			\item[3.3.] Convergence des images directes\dotfill 113
			\item[3.4.] Fibres des $F$-isocristaux convergents\dotfill 126
			\item[3.5.] Cas fini \'etale\dotfill 127
			\end{enumerate}
		\end{enumerate}
	\item[IV.] Images directes de $F$-isocristaux surconvergents\dotfill 131
		\begin{enumerate}
		\item[1.] Frobenius\dotfill 131
		\item[2.] Cas relevable\dotfill 134
		\item[3.] Cas propre et lisse\dotfill 144
		\item[4.] Cas fini \'etale\dotfill 149
		\item[5.] Cas plongeable\dotfill 158
		\end{enumerate}
	\item[V.] Cohomologie syntomique\dotfill 161
		\begin{enumerate}
		\item[1.] Site syntomique\dotfill 161
		\item[2.] Cohomologie syntomique \`a supports compacts\dotfill 167
		\item[3.] Comparaison avec cohomologie \'etale et  cohomologie rigide\dotfill 172
		\end{enumerate}
	\item[VI.] Fonctions $L$\dotfill 183
		\begin{enumerate}
		\item[1.] Fonction $L$ des $F$-modules convergents ou surconvergents\dotfill 183
			\begin{enumerate}
			\item[1.1.] Rel\`evements de Teichm\" uller\dotfill 183
			\item[1.2.] $F$- modules convergents\dotfill 183
			\item[1.3.] $F$- modules surconvergents\dotfill 187
			\item[1.4.] La conjecture de Dwork pour les $F$- modules \\ surconvergents\dotfill 189
			\end{enumerate}
		\item[2.] Fonction $L$ des $F$-isocristaux convergents ( resp. $F$-cristaux)\dotfill 191
			\begin{enumerate}
			\item[2.1.] $F$-isocristaux convergents\dotfill 191
			\item[2.2.] $F$-cristaux\dotfill 195
			\end{enumerate}
		\item[3.] Fonction $L$ des $F$-isocristaux Dwork-surconvergents\dotfill 196 
			\begin{enumerate}
			\item[3.1.] $F$-isocristaux Dwork-surconvergents\dotfill 196
			\item[3.2.] Fonction $L$ des $F$-isocristaux Dwork-surconvergents\dotfill 198
			\end{enumerate}
		\item[4.] Fonction $L$ des $F$-isocristaux surconvergents\dotfill 199
			\begin{enumerate}
			\item[4.1.] D\'efinitions\dotfill 199
			\item[4.2.] La conjecture de Dwork\dotfill 202
			\item[4.3.] La conjecture de Katz\dotfill 202
			\end{enumerate}
		\item[5.] Sch\'emas ab\'eliens ordinaires\dotfill 205
			\begin{enumerate}
			\item[5.1.] $F$-cristaux ordinaires\dotfill 205
			\item[5.2.] Caract\'erisation des sch\'emas ab\'eliens ordinaires\dotfill 206
			\item[5.3.] Explicitation des fonctions $L_{\alpha}^{(r)}$ et conjecture de Dwork\dotfill 212
			\end{enumerate}
		\end{enumerate}
	\item[] Bibliographie\dotfill 215

 \end{enumerate}
 
\vskip 10mm
\chapter*{0.  Introduction}
\markboth{ \sc j.-y. etesse}{\sc 0.  Introduction}
\section*{0.1. Rel\`evements et alg\`ebres de Monsky-Washnitzer}

Le chapitre I, qui rassemble un certain nombre de r\'esultats g\'en\'eraux, est largement pr\'eparatoire pour les autres chapitres.\\

On d\'eveloppe au \S1 un formalisme qui aboutira au \S2 \`a pr\'eciser les liens entre le compl\'et\'e faible $A^{\dag}$ d'une $\mathcal{V}$-alg\`ebre de type fini $A$ et le compl\'et\'e $\hat{A}$, o\`u $\mathcal{V}$ est excellent: par exemple, si $A^{\dag}$ est r\'eduit, alors $A^{\dag}$ est int\'egralement ferm\'e dans $\hat{A}$ [th\'eo (2.2)], ce qui est l'analogue d'un th\'eor\`eme de Bosch, Dwork et Robba.\\

 On g\'en\'eralise, avec les th\'eor\`emes (2.3) et (2.4), des r\'esultats de [Et 4] en prouvant la pleine fid\'elit\'e du foncteur de la cat\'egorie des $A^{\dag}$-alg\`ebres \'etales vers la cat\'egorie des $\hat{A}$-alg\`ebres \'etales: ce foncteur est une \'equivalence de cat\'egories en se restreignant aux alg\`ebres finies \'etales, et un foncteur quasi-inverse est fourni par le passage \`a la fermeture int\'egrale [th\'eo (2.4)].\\
 
 Le \S3 rassemble des r\'esultats de rel\`evement de sch\'emas, de la caract\'eristique $p> 0$ \`a la caract\'eristique 0, qui nous serviront, aux chapitres II, III, IV, \`a \'etablir la convergence ou la surconvergence des images directes d'isocristaux dans les cas relevables: sur une base affine, les deux cas les plus importants seront celui d'un morphisme fini [th\'eo (3.4)], avec son utilisation ult\'erieure au IV pour relever le Frobenius via [cor (3.6)], et celui d'un morphisme projectif lisse [th\'eo (3.3)] et son corollaire [cor (3.3.7)] pour les intersections compl\`etes, ou plus particuli\`erement d'un morphisme fini \'etale [th\'eo (3.1)], ou fini \'etale galoisien [cor (3.7)].\\

\section*{0.2. Espaces rigides analytiques et images directes}

Dans le chapitre II on d\'efinit les images directes d'isocristaux surconvergents par un morphisme $f: X\rightarrow S$ et on prouve leur surconvergence dans le cas o\`u $f$ est propre et lisse relevable, ou $X$ est une intersection compl\`ete relative dans des espaces projectifs sur $S$.\\

L'outil essentiel pour ce dernier r\'esultat est le th\'eor\`eme de changement de base propre du \S1: on \'etablit celui-ci d'abord dans le cadre des sch\'emas formels [th\'eo (1.1)], puis dans le cadre des espaces rigides analytiques [th\'eo (1.2)] gr\^ace aux travaux de Bosch et L\" utkebohmert.\\

La notion de voisinage strict \'etant \'etroitement li\'ee \`a celle de surconvergence des isocristaux, on d\'eveloppe au \S2 quelques propri\'et\'es de ces voisinages stricts dans le cas plongeable: dans le cas cart\'esien, l'image inverse d'un syst\`eme fondamental de voisinages stricts est un syst\`eme fondamental de voisinages stricts [prop (2.2.3)]; le m\^eme r\'esultat vaut pour un morphisme fini et plat, ou fini \'etale, ou fini \'etale galoisien [prop (2.3.1)].\\

Apr\`es avoir rappel\'e au \S3 la d\'efinition des images directes d'isocristaux surconvergents donn\'ee par Berthelot dans une note non publi\'ee [B 5] (voir les articles [C-T], [LS] de Chiarellotto-Tsuzuki et Le Stum pour la publication des d\'etails), on \'etablit leur surconvergence pour un morphisme propre et lisse relevable, en m\^eme temps que le th\'eor\`eme de changement de base [th\'eo (3.4.4)]: un ingr\'edient essentiel est l'extension aux voisinages stricts du th\'eor\`eme de changement de base propre [th\'eo (3.3.2)].\\

Des r\'esultats de rel\`evement du chapitre I on d\'eduit alors la surconvergence des images directes d'un isocristal surconvergent par un morphisme $f: X\rightarrow S$ projectif et lisse, o\`u $S$ est lisse sur le corps de base $k$ et $X$ relevable en un $\mathcal{V}$-sch\'ema plat [th\'eo (3.4.8.2)] (resp $X$ est une intersection relative dans des espaces projectifs sur $S$ [cor (3.4.8.6)]). Une variante, \'etudi\'ee dans le th\'eor\`eme (3.4.9), ram\`ene la preuve de la surconvergence du cas projectif lisse au cas o\`u la base $S$ est affine et lisse sur $k$.\\

Ces th\'eor\`emes (3.4.8.2), (3.4.8.6) et (3.4.9) [resp. le th\'eo (3.4.4)] r\'esolvent partiellement une conjecture de Berthelot [B 2] sur la surconvergence des images directes: dans le chapitre IV nous en donnons une version avec structure de Frobenius.\\

Dans le cas o\`u la base $S$ est une courbe affine et lisse sur un corps alg\'ebriquement clos $k$ (resp. $S$ est une courbe lisse sur un corps parfait $k$) la conjecture a \'et\'e prouv\'ee pour le faisceau structural par Trihan [Tri] (resp. Matsuda et Trihan [M-T]); dans le cas relevable, la conjecture a aussi \'et\'e prouv\'ee ind\'ependamment par Tsuzuki par des voies diff\'erentes [Tsu 4].\\

\section*{0.3. $F$-isocristaux convergents sur un sch\'ema lisse, et images directes}

Avec le chapitre III on introduit une structure de Frobenius sur les isocristaux consid\'er\'es, tout en restant dans le cadre \guillemotleft convergent\guillemotright. Les r\'esultats concernent deux aspects de la th\'eorie: d'une part l'allure des $F$-isocristaux convergents sur un sch\'ema lisse aux \S1 et \S2, d'autre part la convergence des images directes par un morphisme propre et lisse $f: X\rightarrow S$ au \S3.\\

Au \S1 on obtient une description des $F$-isocristaux convergents sur un sch\'ema affine et lisse [cor (1.2.3)] qui est l'analogue de celle de type Monsky-Washnitzer obtenue par Berthelot dans le cas surconvergent [B 3]. Cette description s'\'etend au cas lisse relevable dans le \S2 [th\'eo (2.2.1)].\\

 L'un des ingr\'edients essentiels au \S3 est le passage par la cohomologie cristalline et l'utilisation des objets \guillemotleft consistants\guillemotright \   au sens de Berthelot-Ogus [B-O]: dans ce cas on est amen\'e \`a supposer le corps $k$ parfait, l'indice de ramification $e$ plus petit que $p-1$ et \`a se limiter aux cristaux plats pour pouvoir appliquer le th\'eor\`eme de changement de base en cohomologie cristalline.\\

Avec les restrictions pr\'ec\'edentes on obtient la convergence des images directes dans le cas propre et lisse, et le th\'eor\`eme de changement de base: si le morphisme propre et lisse $f: X\rightarrow S$ est relevable en un morphisme propre et lisse sur $\mathcal{V}$, on peut lever la restriction de platitude [th\'eo (3.3.1)]. En fait, si le morphisme $f: X\rightarrow S$ est projectif et lisse et que, ou bien $f$ est relevable, ou bien que $X$ est une intersection compl\`ete relative dans des espaces projectifs sur $S$, alors on peut lever toutes les restrictions pr\'ec\'edentes (i.e. sur l'indice de ramification, sur la platitude et la perfection de $k$) [th\'eo (3.3.2)]: la preuve n'utilise pas la cohomologie cristalline, mais repose sur le cas plongeable du [II, 3.4]. De plus les images directes commutent au passage aux fibres [prop (3.4.4)]. Dans le cas fini \'etale toutes les restrictions pr\'ec\'edentes sont lev\'ees [th\'eo (3.5.1)], et dans le cas galoisien la cohomologie convergente commute aux points fixes sous le groupe du rev\^etement [th\'eo (3.5.1)].\\

\section*{0.4. Images directes de $F$-isocristaux surconvergents}

Au chapitre IV on \'etudie les images directes des $F$-isocristaux dans le cadre \guillemotleft surconvergent\guillemotright. Tout le probl\`eme est d'obtenir de \guillemotleft bons\guillemotright \  rel\`evements du Frobenius.\\

Sur la base $S$, la question est r\'esolue par le diagramme (1.2.4), qu'il s'agit ensuite de \guillemotleft tirer en haut\guillemotright \ par le morphisme propre et lisse $f:X\rightarrow S.$ \\

Dans le cas o\`u $f$ est relevable on arrive \`a effectuer en parrall\`ele des rel\`evements du Frobenius et des bons choix de compactifications. D'o\`u le th\'eor\`eme de surconvergence des images directes dans le cas relevable [th\'eo (2.1)]: le point cl\'e est de montrer que le Frobenius, sur les images directes surconvergentes, est un isomorphisme; par le th\'eor\`eme de changement de base on est ramen\'e \`a la m\^eme propri\'et\'e dans le cas convergent vu au [III, (3.3.1)]. Dans le cas o\`u $f$ est projectif lisse relevable, ou $X$ intersection compl\`ete relative dans des espaces projectifs sur $S$, abord\'es au \S3, ou fini \'etale au \S4, on construit comme au II \S3 des foncteurs $R^{i}f_{rig\ast}$ sur la cat\'egorie des $F$-isocristaux surconvergents [th\'eos (3.1), (4.1)].\\

Lorsque $f$ est propre et lisse et ou bien g\'en\'eriquement projectif relevable, ou bien  g\'en\'eriquement projectif et $X$ intersection compl\`ete relative dans des espaces projectifs sur $S$, on utilise les r\'esultats du III: on est ainsi amen\'e \`a supposer $k$ parfait, $e\leqslant p-1$ et \`a se restreindre \`a des $F$-isocristaux plats. Gr\^ace aux propri\'et\'es de descente \'etale des $F$-isocristaux surconvergents de [Et 5] et de pleine fid\'elit\'e du foncteur d'oubli de la cat\'egorie surconvergente vers la cat\'egorie convergente de Kedlaya [Ked 2], on prouve encore la surconvergence des images directes [th\'eo (3.2)].\\

Dans le \S5, o\`u $f$ est suppos\'e seulement plongeable, on proc\`ede diff\'eremment: comme il n'y a plus de rel\`evements globaux du Frobenius comme pr\'ec\'edemment on utilise la fonctorialit\'e des images directes pour construire un morphisme de Frobenius; il reste alors \`a prouver que c'est un isomorphisme. Par un r\'esultat de Berthelot il suffit de voir que tel est le cas dans la cat\'egorie convergente: le th\'eor\`eme de changement de base et un r\'esultat de Bosch-G\"untzer-Remmert [B-G-R] nous ram\`ene \`a le v\'erifier aux points ferm\'es de $S$, pour lesquels l'assertion est connue [III, (3.3.1.18)]. Pour $f$ plongeable on a ainsi des foncteurs images directes sur la cat\'egorie des $F$-isocristaux surconvergents [th\'eo (5.2)].\\

\section*{0.5. Cohomologie syntomique}

La cohomologie syntomique introduite au chapitre V fait le lien entre la cohomologie \'etale et la cohomologie rigide, lien qui sera utilis\'e au chapitre VI pour r\'esoudre une conjecture de Katz sur les z\'eros et p\^oles unit\'es $p$-adiques des fonctions $L$.\\

Comme pour le topos \'etale [Mi, II, theo 3.10] on obtient au \S1 une description du topos syntomique [th\'eo (1.3)] calqu\'ee sur celle de [SGA 4,T 1, IV, th\'eo (9.5.4)], SGA 4 qui travaille en termes de sous-topos ouvert et du sous-topos ferm\'e compl\'ementaire. On en d\'eduit les suites exactes courtes usuelles de localisation [th\'eo (1.5)].\\

Au \S2 la cohomologie syntomique \`a supports compacts est d\'efinie: en particulier il faut s'assurer de l'ind\'ependance par rapport \`a la compactification choisie [prop (2.1)]. Les suites exactes courtes de localisation du \S1 fournissent alors la suite exacte longue de localisation en cohomologie syntomique \`a supports compacts [th\'eo (2.6)].\\

Au \S3 on \'etablit que les cohomologies syntomique et \'etale \`a supports compacts d'un faisceau lisse $\mathcal{F}$ co\" incident [th\'eo (3.2)]. De m\^eme la cohomologie rigide \`a supports compacts d'un $F$-isocristal surconvergent unit\'e $E$ associ\'e \`a $\mathcal{F}$ co\" incide avec une limite de la  cohomologie syntomique \`a supports compacts d'un $E_{n}^{m\mbox{-}cris}$ associ\'e \`a $\mathcal{F}$ [ th\'eo (3.3.13)(1) et (3.3.16)]. Comme il existe une suite exacte courte sur le site syntomique qui relie  $\mathcal{F}$ et $E_{n}^{m\mbox{-}cris}$ [th\'eo (3.3.12)], on en d\'eduit que la cohomologie \'etale \`a supports compacts de $\mathcal{F}$ s'identifie aux points fixes du Frobenius agissant sur la cohomologie rigide \`a supports compacts de $E$ [th\'eo(3.3.13) (2) et (3)].\\

\section*{0.6. Fonctions $L$}

Dans le chapitre VI on donne une d\'efinition unifi\'ee des fonctions $L$ des $F$-modules (resp. $F$-isocristaux) convergents ou surconvergents et des $F$-cristaux: elle redonne celle utilis\'ee en cohomologie cristalline par Katz [K 1] ou [Et 2], ou celle en cohomologie rigide de [E-LS 1], ou celle utilis\'ee par Wan [W 2]. Le but ici est d'obtenir des propri\'et\'es de m\'eromorphie ou de rationalit\'e de ces fonctions $L$ (conjecture de Dwork) et d'\'etudier leurs z\'eros et p\^oles unit\'es $p$-adiques ( conjecture de Katz); on explicite ces r\'esultats pour les sch\'emas ab\'eliens ordinaires au \S5.\\

Au \S1 les r\'esultats de Wan permettent d'\'etablir la m\'eromorphie attendue pour les fonctions $L$ de $F$-modules surconvergents.\\

Au \S2 la d\'efinition des fonctions $L$ des $F$-isocristaux convergents donne une forme globale des fonctions $L$ des $F$-modules convergents du \S1, et g\'en\'eralise les fonctions $L$ des $F$-cristaux du cas cristallin [Et 2], et les fonctions $L$ des $F$-isocristaux surconvergents du cas rigide [E-LS 1]: si $e\leqslant p-1$, on obtient la rationalit\'e des fonctions $L$ des images directes du faisceau structural par un morphisme propre et lisse [th\'eo (2.1.4)].\\

Au \S3 on introduit la notion de $F$-isocristaux Dwork-surconvergents: il s'agit des $F$-isocristaux convergents dont le Frobenius est surconvergent. Contrairement aux $F$-isocristaux surconvergents on n'y suppose pas la connexion surconvergente; mais la surconvergence du Frobenius suffit \`a assurer la m\'eromorphie $p$-adique des fonctions $L_{\alpha}^{(r)}$ [th\'eo (3.2.10)], ce qui prouve la conjecture de Dwork pour de tels coefficients.\\

Au \S4 les fonctions $L_{\alpha}^{(r)}$ d'un $F$-isocristal surconvergent $E$ sont d\'efinies comme \'etant celles du $F$-isocristal convergent $\mathcal{E}$ associ\'e: un $F$-isocristal surconvergent fournit un $F$-isocristal Dwork-surconvergent par oubli de la surconvergence de la connexion, ainsi, gr\^ace au \S3, la m\'eromorphie de $L_{\alpha}^{(r)}(E)$ est \'etablie [th\'eo (4.2.1)(i)]; la rationalit\'e de $L^{(r)}(E)$ r\'esulte de la cohomologie rigide [th\'eo (4.2.1)(ii)].\\

En utilisant les r\'esultats du IV et du V de surconvergence des images directes on en d\'eduit la m\'eromorphie de 
$$L_{\alpha}^{(r)}(R^{i}f_{rig\ast}(E))\  \ (4.3.2.3)$$
et la rationalit\'e de 
$$L^{(r)}(R^{i}f_{rig\ast}(E))\  \ (4.3.2.2)\ .$$

Toujours dans ce \S4 on aborde la preuve de la conjecture de Katz dans le cas surconvergent [th\'eo (4.3.1)]: on y l\`eve l'hypoth\`ese de prolongement \`a une compactification qui avait \'et\'e faite dans [E-LS 2; th\'eo 6.7], gr\^ace aux suites exactes en cohomologie syntomique \`a supports compacts du V [V, (3.3.13)(2)], tout en pr\'ecisant la preuve de Emerton et Kisin [(4.3.1.3)].\\

On explicite au \S5 les cons\'equences des r\'esultats pr\'ec\'edents pour les sch\'emas ab\'eliens ordinaires. Nous avons adopt\'e la d\'efinition d'ordinarit\'e d'un morphisme de Illusie [I$\ell$ 2], celle d'ordinarit\'e d'un groupe $p$-divisible de Raynaud [R 3], et celle d'ordinarit\'e d'un $F$-cristal de Katz et Deligne [K 2] et [De$\ell$ 2]. Si $f:X \rightarrow S$ est un sch\'ema ab\'elien, on note $G$ le groupe $p$-divisible associ\'e et $\mathbb{D}(G)$ son $F$-cristal de Dieudonn\'e: tout d'abord il est \'equivalent de dire que $X/S$ est ordinaire, ou que $G$ est ordinaire, ou que $\mathbb{D}(G)$ est ordinaire [th\'eo (5.2.9)]. Pour un sch\'ema ab\'elien ordinaire $f:X \rightarrow S$ on a des descriptions des fonctions $L$ suivantes
$$L_{\alpha}(R^{i}f_{cris\ast}(\mathcal{O}_{X/W}))\ , \quad   L_{\alpha}^{(r)}(R^{i}f_{cris\ast}(\mathcal{O}_{X/W}))$$
 en termes de fonctions $L$ usuelles associ\'ees aux cristaux de Dieudonn\'e des groupes $p$-divisibles $G^{\acute{e}t}$ et $G^{tm}$ provenant de $G$ [th\'eo (5.3.1)].\\

 \noindent\textit{Remerciements}. L'auteur tient ici \`a remercier P. Berthelot et B. Le Stum pour avoir mis \`a sa disposition leurs preprints [B 5] et [LS], ainsi que M.-F. Ch\'eriaux pour l'aide apport\'ee \`a la frappe du manuscrit.

\cleardoublepage

\vskip 10mm
\chapter*{I.  Rel\`evements et alg\`ebres de Monsky-Washnitzer}
\markboth{\sc j.-y. etesse}{\sc I.  Rel\`evements et alg\`ebres de Monsky-Washnitzer}
\section*{1. G\'en\'eralit\'es}

\textbf{1.0. Notations}\\

Tous les anneaux consid\'er\'es dans cet article sont (sauf mention du contraire) commutatifs et unitaires.
\vskip3mm
Soient $\mathcal{V}$ un anneau noeth\'erien, $I \not\subseteq \mathcal{V}$ un id\'eal, $A$ une $\mathcal{V}$-alg\`ebre telle que l'anneau $A$ soit noeth\'erien et $IA \neq A$. On note $\hat{A}$ le s\'epar\'e compl\'et\'e $I$-adique de $A$, $A_{n} = A/I ^{n+1} A$ et $A^{\dag} \subset \hat{A}$ le compl\'et\'e faible de $A$ au-dessus de la paire $(\mathcal{V}, I)$ [M-W, \S\  1] : on d\'esignera toujours par un indice $(\quad )_{0}$ la r\'eduction mod $I$ d'une $\mathcal{V}$-alg\`ebre ou d'un $\mathcal{V}$-morphisme.
\vskip3mm
Si $B$ est une $\mathcal{V}$-alg\`ebre, on dit que $B$ est faiblement compl\`ete de type fini (f.c.t.f. en abr\'eg\'e) si $B$ est la compl\'et\'ee faible d'une $\mathcal{V}$-alg\`ebre de type fini ; une telle alg\`ebre $B$ est appel\'ee ``w.c.f.g.''  dans la terminologie de [M-W, \S\  2].

\vskip 3mm
Consid\'erons la partie multiplicative $T = 1 + IA$ de $A$ ; notons $A_{T} = T^{-1} A$ et $(\tilde A, \tilde I)$ le hens\'elis\'e de $(A, IA)$ au sens de Raynaud [R 2, d\'ef 4 p 24] : rappelons qu'on a suppos\'e $IA \neq A$, si bien que $0 \not\in T$ ; sinon les anneaux $A_{T}, \tilde{A}, \hat{A}$ et $A^{\dag}$ seraient \'egaux \`a ${\{0\}}$.
\vskip 3mm
On rappelle [Et 4] que si $A$ est une $\mathcal{V}$-alg\`ebre de type fini, alors il existe des morphismes canoniques $A_{T} \rightarrow \tilde{A} \rightarrow A^{\dag} \rightarrow \hat{A}$ tous fid\`element plats et que tous ces anneaux ont m\^eme s\'epar\'e compl\'et\'e $I$-adique \'egal \`a $\hat{A}$.

\vskip 3mm
\noindent \textbf{Proposition (1.1)}. 
\textit{Soient $A, \mathcal{A}, \mathcal{B}$ des anneaux noeth\'eriens munis de morphismes}
$$ A \displaystyle \mathop{\longrightarrow}^{\varphi_{1}} \mathcal{A} \displaystyle \mathop{\longrightarrow}^{\varphi_{2}} \mathcal{B} \displaystyle \mathop{\longrightarrow}^{\varphi_{3}} \hat{A}$$ \textit{avec} $\varphi_{3}$ \textit{fid\`element plat. On suppose que} $Spec\ \hat{A} \longrightarrow Spec\  A$ \textit{est un morphisme normal (resp. r\'egulier)} [EGA IV, (6.8.1)]:
\textit{cette derni\`ere hypoth\`ese est v\'erifi\'ee si $A$ est excellent. Alors :}

\begin{itemize}
\item[(1)] \textit{Le morphisme
$$ h = Spec\ (\varphi _{2} \circ \varphi_{1}) : Spec\  \mathcal{B} \rightarrow Spec \ A$$
est normal (resp. r\'egulier).}
\item[(2)] \textit{Si de plus $\varphi_{2}$  est plat, alors
$$ f = Spec\ (\varphi_{2}) : Spec\ \mathcal{B} \rightarrow \ Spec\ \mathcal{A}$$ 
est un morphisme normal.}
\item[(3)] \textit{Si $\varphi_{2}$ est plat et $(\mathcal{A}, I \mathcal{A})$ est un couple hens\'elien tel  que $ \hat{\mathcal{A}}  \simeq \hat{A}$, alors
$$ f  : \mbox{Spec}\ \mathcal{B} \rightarrow \ Spec\ \mathcal{A}$$
est un morphisme normal \`a fibres g\'eom\'etriquement int\`egres.}
\item[(4)] \textit{Si $\varphi_{2}$ est plat, $(\mathcal{A}, I \mathcal{A})$ est un couple hens\'elien tel que $\hat{\mathcal{A}} \simeq \hat{A}$ et $\mathcal{A}$ est r\'eduit, alors $\mathcal{A}$ est int\'egralement ferm\'e dans $\mathcal{B}$ et dans $\hat{A}$.} 
\end{itemize}

\vskip 3mm
\noindent \textit{D\'emonstration}. Le (1) r\'esulte de [EGA IV, (6.5.4) (i) (resp. (6.5.2) (i))]. \par 
Pour le (2) notons $g = \mbox{Spec}\ \varphi_{1} : \mbox{Spec}\  \mathcal{A} \rightarrow \mbox{Spec}\  A$. Soit $\mathfrak{q} \in \mbox{Spec}\  \mathcal{A}$ et $k'$ une extension finie du corps r\'esiduel $k(\mathfrak{q})$ : il s'agit de montrer que

$$f^{-1} (\mathfrak{q})_{k'} = \mbox{Spec}\  (k' \otimes_{\mathcal{A}} \mathcal{B})$$

\noindent est normal [EGA IV, (6.8.1)]. Comme Spec $k'$ est normal et $h$ un morphisme normal, on sait par [EGA IV, (6.14.1)] que Spec $(k' \otimes_{A} \mathcal{B})$ est normal. \\
Consid\'erons les applications\\

$k' \otimes_{\mathcal{A}} \mathcal{B}   \displaystyle \mathop{\longrightarrow}^{\psi} k' \otimes_{A} \mathcal{B} \simeq  k' \otimes_{\mathcal{A}} (\mathcal{A}\  \otimes_{A} \mathcal{B}) \displaystyle \mathop{\longrightarrow}^{\varphi}  k' \otimes_{\mathcal{A}} \mathcal{B}$\par
$x \otimes b \quad  \longmapsto   \quad\quad\quad\quad  x \otimes (1_{\mathcal{A}} \otimes b)$\par
$\qquad \qquad \qquad \qquad \qquad x \otimes (a \otimes b) \qquad \longmapsto x \otimes (\varphi_{2}(a) . b)$ ; \\

\noindent clairement $\varphi \circ \psi = Id$. Pour $\mathfrak{p} \in \mbox{Spec}\ (k' \otimes_{\mathcal{A}} \mathcal{B})$, $\mathfrak{m} := \varphi^{-1}(\mathfrak{p})$ on a $\mathfrak{p} = \psi^{-1} (\varphi^{-1}(\mathfrak{p})) = \psi^{-1}(\mathfrak{m})$ et $\psi$ et $\varphi$ induisent des morphismes

$$(k' \otimes_{\mathcal{A}} \mathcal{B})_{\mathfrak{p}} \displaystyle \mathop{\longrightarrow}^{\psi'} (k' \otimes_{A} \mathcal{B})_{\mathfrak{m}} \displaystyle \mathop{\longrightarrow}^{\varphi'}  (k' \otimes_{\mathcal{A}} \mathcal{B})_{\mathfrak{p}}$$

\noindent dont le compos\'e est encore l'identit\'e.\\
Par hypoth\`ese $\mathcal{D} := (k' \otimes_{A} \mathcal{B})_{\mathfrak{m}}$ est int\'egralement clos, de corps des fractions not\'e $L$ ; en particulier $\mathcal{C} := (k' \otimes_{\mathcal{A}} \mathcal{B})_{\mathfrak{p}}$ est int\`egre : notons $K$ son corps des fractions. Il s'agit de montrer que $\mathcal{C}$ est int\'egralement clos : soit $x \in K$ un \'el\'ement entier sur $\mathcal{C}$, annul\'e par le polyn\^ome $R(X) = X^n + \displaystyle \mathop{\Sigma}_{i=0}^{n-1} a_{i}\ X^i$, $a_{i} \in \mathcal{C}$. L'injection $\psi '$ induit une injection :

$$\psi '' : K \hookrightarrow L\  ;$$

\noindent puisque $\mathcal{D}$ est int\'egralement clos on a $\psi''(x) \in \mathcal{D}$ et $ x_{1} := \varphi' (\psi''(x)) \in \mathcal{C}$ est racine de $R(X)$ : en effet $\psi''$ induit $\tilde \psi'' : K[K] \rightarrow L[X]$, $R(X) \mapsto \tilde R(X) \in \mathcal{D}[X]$ et $ \varphi'$ induit $\tilde \varphi' : \mathcal{D}[X] \rightarrow \mathcal{C}[X]$, $\tilde R(X) \mapsto R(X)$.  D'o\`u $R(X) = (X - x_{1}) R_{1}(X)$ avec $R_{1}(X) \in \mathcal{C}[X]$. Si $x = x_{1}$ on a termin\'e, sinon $x$ est racine de $R_{1}(X)$ et on it\`ere : finalement $x \in \mathcal{C}$, donc $\mathcal{C}$ est int\'egralement clos, i.e. $f^{-1}(\mathfrak{q})_{k'}$ est normal.\\

Pour le (3), on sait par le (2) que le morphisme $g : \mbox{Spec}\  \hat{\mathcal{A}} = \mbox{Spec}\  \hat{A} \rightarrow
\mbox{Spec}\  \mathcal{A}$ est normal, donc ses fibres sont g\'eom\'etriquement int\`egres par un th\'eor\`eme de Raynaud [R 2, th\'eo 3, p.126]. Or $\hat{A}$ est une $\mathcal{B}$-alg\`ebre fid\`element plate, donc les fibres de $f$ sont aussi g\'eom\'etriquement int\`egres.\\

Pour le (4), le th\'eor\`eme de Raynaud [loc. cit.] prouve que $\mathcal{A}$ est int\'egralement ferm\'e dans $\hat{A}$ : comme $\varphi_{3}$ est injectif, $\mathcal{A}$ est aussi int\'egralement ferm\'e dans $\mathcal{B}$. $\square$\\

Pour la commodit\'e des r\'ef\'erences nous avons rassembl\'e ci-apr\`es quelques lemmes qui r\'esultent des EGA.

\vskip 3mm
\noindent \textbf{Lemme (1.2)}.
\textit{Soient $\mathcal{A}$ et $\mathcal{B}$ deux anneaux noeth\'eriens tels que $\mathcal{B}$ soit une $\mathcal{A}$-alg\`ebre. Les propri\'et\'es (i) et (ii) ci-apr\`es sont \'equivalentes :}

\begin{itemize}
\item[(i)] \textit{$\mbox{Spec}\  \mathcal{B} \rightarrow \mbox{Spec}\  \mathcal{A}$ est un morphisme r\'egulier.}

\item[(ii)] \textit{Pour tout $\mathfrak{q} \in \mbox{Spec}\  \mathcal{A}$ et tout  $\mathfrak{p} \in \mbox{Spec}\  \mathcal{B}$ au-dessus de $\mathfrak{q}$, le morphisme $\mathcal{A}_{\mathfrak{q}} \rightarrow \mathcal{B}_{\mathfrak{p}}$ est formellement lisse pour les topologies pr\'eadiques respectives (i.e. d\'efinies par $\mathfrak{q}\  \mathcal{A}_{\mathfrak{q}}$ et $\mathfrak{p}\  \mathcal{B}_{\mathfrak{p}}$ respectivement).  } 
\end{itemize}

\vskip 3mm
\noindent \textit{D\'emonstration}. En utilisant la d\'efinition d'un morphisme r\'egulier [EGA IV, (6.8.1)], l'\'equivalence r\'esulte de [EGA $O_{IV}$, (22.5.8) et (19.7.1)]. $\square$

\vskip 3mm
\noindent \textbf{Lemme (1.3)}.
\textit{Soient $\mathcal{A}$ et $\mathcal{B}$ deux anneaux noeth\'eriens tels que $\mathcal{B}$ soit une $\mathcal{A}$-alg\`ebre. Si  $Spec\  \mathcal{B} \rightarrow Spec\  \mathcal{A}$ est formellement lisse pour les topologies discr\`etes, alors c'est un morphisme r\'egulier.} 

\vskip 3mm
\noindent \textit{D\'emonstration}.  Si $\mathcal{B}$ est une $\mathcal{A}$-alg\`ebre formellement lisse pour les topologies discr\`etes  [EGA $O_{IV}$, (19.3.1)], alors pour tout $\mathfrak{q} \in \mbox{Spec}\ \mathcal{A}$ et tout $\mathfrak{p} \in \mbox{Spec}\  \mathcal{B}$ au-dessus de $\mathfrak{q}$, le morphisme

$$\mathcal{A}_{\mathfrak{q}} \rightarrow \mathcal{B}_{\mathfrak{p}}$$ 
est formellement lisse pour les topologies discr\`etes [EGA $O_{IV}, (19.3.5) (iv)]$, donc aussi pour les topologies pr\'eadiques sur $\mathcal{A}_{\mathfrak{q}}$ et $\mathcal{B}_{\mathfrak{p}}$ [EGA $O_{IV}$, (19.3.8)] ; d'o\`u la conclusion par le lemme (1.2). $\square$

\vskip 3mm
\noindent \textbf{Lemme (1.4)}.
\textit{Soient $\mathcal{A}$ un anneau, $\mathcal{B}$ une $\mathcal{A}$-alg\`ebre et $\mathcal{J} \subset \mathcal{A}$ un id\'eal.\\
Si $Spec\  \mathcal{B} \rightarrow Spec\ \mathcal{A}$ est formellement lisse pour les topologies discr\`etes, alors c'est un morphisme formellement lisse pour les topologies $\mathcal{J}$-adiques sur $\mathcal{A}$ et $\mathcal{B}$.}

\vskip 3mm
\noindent \textit{D\'emonstration}. Via [EGA $O_{IV}$, (19.3.8)]. $\square$

\vskip 3mm
\noindent \textbf{Lemme (1.5)}.
\textit{Soient $\mathcal{A}$ un anneau et $S \subset \mathcal{A}$ une partie multiplicative. Alors}
 \begin{itemize}
\item[(i)] \textit{Le morphisme $f : X = Spec\ (S^{-1} \mathcal{A}) \rightarrow Spec\  \mathcal{A} = Y$ est formellement \'etale pour les topologies discr\`etes. }
\item[(ii)] \textit{En particulier $f$ est r\'egulier si $\mathcal{A}$ est noeth\'erien.} 
\item[(iii)] \textit{De plus : \\
$\ast$ si $y \in Y \backslash\  f(X)$, alors $f^{-1}(y) = \phi$\\
$\ast$ si $y \in f(X)$, alors $f$ induit un isomophisme}
$$ f^{-1}(y) \displaystyle \mathop{\longrightarrow}^{\sim} Spec\  k(y).$$ 
\end{itemize}

\vskip 3mm
\noindent \textit{D\'emonstration}. La premi\`ere assertion n'est autre que [EGA $O_{IV}$, (19.10.3) (ii)] ; la deuxi\`eme en r\'esulte gr\^ace au lemme (1.3). La derni\`ere assertion est cons\'equence de [Bour, AC II, \S\  2, n$^\circ\  5$, prop 11]. $\square$

\vskip 3mm
\noindent \textbf{Lemme (1.6)}.
\textit{Soient $\mathcal{A}$ et $\mathcal{B}$ deux anneaux noeth\'eriens tels que $\mathcal{B}$ soit une $\mathcal{A}$-alg\`ebre;  notons
$$f : Spec\ \mathcal{B} \rightarrow Spec\ \mathcal{A}$$
le morphisme canonique. Alors :}
\begin{enumerate}
\item[(1)] \textit{On a les implications :}
\begin{itemize}
\item[(i)] \textit{$f$ est r\'eduit et $\mathcal{A}$  r\'eduit $\Rightarrow \mathcal{B}$ r\'eduit.}
\item[(ii)] \textit{$f$ normal et $\mathcal{A}$ normal $\Rightarrow \mathcal{B}$ normal.}
\item[(iii)]\textit{$f$ r\'egulier et $\mathcal{A}$ r\'egulier $\Rightarrow \mathcal{B}$ r\'egulier.}
\end{itemize}
\item[(2)] \textit{Si $f$ est fid\`element plat, on a les implications :}
\begin{itemize}
\item[(i)] \textit{$\mathcal{B}$ r\'eduit $\Rightarrow \mathcal{A}$ r\'eduit.}
\item[(ii)] \textit{$\mathcal{B}$ normal $\Rightarrow \mathcal{A}$ normal.}
\item[(iii)] \textit{$\mathcal{B}$ r\'egulier $\Rightarrow \mathcal{A}$ r\'egulier.}
\end{itemize}
\end{enumerate}

\vskip 3mm
\noindent \textit{D\'emonstration}.\\
(1) (i) On utilise la d\'efinition d'un morphisme r\'eduit [EGA IV, (6.8.1)] et la caract\'erisation des sch\'emas noeth\'eriens r\'eduits de [EGA IV, (5.8.5)] via les propri\'et\'es $\ll R_{0}$ et $S_{1}\gg$ : comme $\mathcal{A}$ v\'erifie $\ll R_{0}$ et $S_{1}\gg$, il en est de m\^eme de $\mathcal{B}$ via [EGA IV, (6.5.3) (ii) et (6.4.2)], donc $\mathcal{B}$ est r\'eduit.\\

\noindent Pour (ii) (resp. (iii)) on utilise [EGA IV, (6.5.4) (ii)] (resp. [EGA IV, (6.5.2) (ii)]).\\

\noindent (2) (i) Le fait qu'un sch\'ema $X$ est r\'eduit s'exprime via les anneaux locaux $\mathcal{O}_{X,x}$ [EGA $0_{I}$, (4.1.4)] : le (i) r\'esulte de [EGA IV, (2.1.13)].\\

\noindent Le (ii), c'est [EGA IV, (6.5.4) (i)] et le (iii) c'est [EGA IV, (6.5.2) (i)]. $\square$\\

La proposition suivante g\'en\'eralise des assertions de [Et 4, prop 2, prop 11].

\vskip 3mm
\noindent \textbf{Proposition (1.7)}.
\textit{Soient $A$ une $\mathcal{V}$-alg\`ebre (resp. une $\mathcal{V}$-alg\`ebre de type fini) telle que $A$ soit un anneau noeth\'erien. Notons $\mathcal{B}$ l'un des anneaux $\tilde{A}, \hat{A}$ (resp. $A^{\dag}$). Alors }

\begin{enumerate}
\item[(1)] \textit{$A_{T}$ et $\mathcal{B}$ sont des anneaux de Zariski.}
\item[(2)] \textit{On a les implications :}
\begin{itemize}
\item [(i)] \textit{$A$ r\'eduit $\Rightarrow A_{T}$ r\'eduit $\Leftrightarrow \tilde{A}$ r\'eduit.}
\item [(ii)] \textit{$A$ normal $\Rightarrow A_{T}$ normal $\Leftrightarrow \tilde{A}$ normal.}
\item [(iii)] \textit{$A$ r\'egulier $\Rightarrow \hat{A}$ r\'egulier (resp. $\Rightarrow A^{\dag}$ r\'egulier)\\
$\Rightarrow \tilde{A}$ r\'egulier $\Rightarrow A_{T}$ r\'egulier.}
\item [(iv)] \textit{$A$ int\'egralement clos $\Rightarrow A_{T}$ int\'egralement clos.}
\end{itemize}
\item[(3)] \textit{Si $f_{T} : Spec\  \hat{A} \rightarrow Spec\  A_{T}$ est le morphisme canonique et $f$ le compos\'e
$f : Spec\ \hat{A}  \displaystyle \mathop{\longrightarrow}_{f_{T}} Spec\ A_{T} \longrightarrow Spec\ A$, on a l'implication\\
$f$ r\'eduit (resp. normal ; resp. r\'egulier)\\
$\Rightarrow f_{T}$ r\'eduit (resp. normal ; resp. r\'egulier).}
\item[(4)]ÊÊ
\begin{itemize}
\item[(i)] \textit{Si $f$ est r\'eduit on a les implications :\\
 $A$ r\'eduit $\Rightarrow A_{T}$ r\'eduit $\Leftrightarrow
\mathcal{B}$ r\'eduit.}
\item[(ii)] \textit{Si $f$ est normal on a les implications :\\ $A$ normal $\Rightarrow A_{T}$ normal $\Leftrightarrow \mathcal{B}$ normal.}
\item[(iii)] \textit{Si $f$ est r\'egulier on a les implications :\\ $A$ r\'egulier $\Rightarrow A_{T}$ r\'egulier $\Leftrightarrow \mathcal{B}$ r\'egulier.}
\item[(iv)] \textit{Si $f$ est normal on a les implications : \\
$A$ int\'egralement clos $\Rightarrow A_{T}$ int\'egralement clos $\Leftrightarrow 
\mathcal{B}$ int\'egralement clos.}
\end{itemize}
\end{enumerate}

\newpage
\noindent \textit{D\'emonstration}. 
\begin{enumerate}
\item [(1)] $A_{T}$ est noeth\'erien [Bour, AC II,  \S\ 2, $\mbox{n}^\circ\  4$, cor 2 de prop 10],
de  m\^eme que $\tilde{A}$ [R 2, p 125] et $\hat{A}$ [Bour, AC III,  \S\ 3, $\mbox{n}^\circ\  4$, prop 8]. De plus si $A$ est de type fini sur $\mathcal{V}$, alors $A^{\dag}$ est noeth\'erien [M-W, theo 2.1]. 
De plus par [Et 4, \S\  1] on a les inclusions $IA_{T} \subset \mbox{Rad}\  A_{T}$, $I \tilde{A} \subset \mbox{Rad}\  \tilde{A}$, $IA^{\dag} \subset \mbox{Rad}\  A^{\dag}$, 
$I \hat{A} \subset \mbox{Rad}\  \hat{A}$ ; donc les anneaux $A_{T}$, $\tilde{A}$, $A^{\dag}$ et $\hat{A}$ sont de Zariski.
\item [(2)] Gr\^ace aux lemmes (1.5) et (1.3), l'implication
$$A\  \mbox{r\'eduit (resp. normal ; resp. r\'egulier)} \Rightarrow A_{T}\   \mbox{de m\^eme},$$
r\'esulte du [lemme (1.6), (1)] et l'\'equivalence
$$A_{T}\  \mbox{r\'eduit (resp. normal)} \Leftrightarrow \tilde{A}\  \mbox{de m\^eme}, $$
c'est [R 2, p 125].\\
Si $A$ est r\'egulier, alors $\hat{A}$ l'est [EGA O$_{IV}$, (17.3.8.1)] : par fid\`ele platitude de $\hat{A}$ sur $A^{\dag}$, $\tilde{A}$ et $A_{T}$, ces derniers sont aussi r\'eguliers [EGA IV, (6.5.2)(i)].\\
Si $A$ est int\'egralement clos, rappelons que $0 \notin T$, donc $A_{T}$ est int\`egre de m\^eme corps des fractions que $A$, par suite $A_{T}$ est int\'egralement clos [Bour, AC V, \S\ 1, $\mbox{n}^\circ\  5$, prop 16].
\item[(3)] R\'esulte du [lemme (1.5) (iii)] et de [EGA IV, (6.8.1)].
 \item[(4)] D'apr\`es le (2) et le (3) on est ramen\'e \`a prouver le (4) en rempla\c{c}ant $A$ par $A_{T}$ et $f$ par le morphisme fid\`element plat $f_{T}$.
\end{enumerate}

\vskip 3mm
Les assertions (i) \`a (iii) sont fournies par le lemme (1.6).

\vskip 3mm
Pour le (iv) supposons d'abord $\mathcal{B}$ int\'egralement clos : par fid\`ele platitude de $\mathcal{B}$ sur $A_{T}$ on en d\'eduit que $\mbox{Spec}\ A_{T}$ est connexe, normal [EGA IV, (2.1.13)] et noeth\'erien, donc  $A_{T}$ est int\`egre [EGA I, (4.5.6)] et int\'egralement clos par [Bour, AC II, \S\  3, $\mbox{n}^\circ\  3$, cor 4 du th\'eo 1 et AC V, \S\  1, $\mbox{n}^\circ\  2$, cor de prop 8].

\vskip 3mm
R\'eciproquement supposons $\mathcal{A} := A_{T}$ int\'egralement clos : par fid\`ele platitude de $\hat{A}$ sur $\mathcal{B}$ il nous suffit de montrer que $\hat{A} = \hat{\mathcal{A}}$ est int\'egralement clos. Puisque $\mathcal{A}$ est int\`egre, l'id\'eal $I \mathcal{A}$ est sans torsion, par suite 
$\mathcal{A} = \displaystyle\mathop{\cap}_{\mathfrak{p}} \mathcal{A}_{\mathfrak{p}}$ et $I \mathcal{A} = \displaystyle \mathop{\cap}_{\mathfrak{p}} I \mathcal{A}_{\mathfrak{p}}$ o\`u $\mathfrak{p}$ parcourt l'ensemble $M$ des id\'eaux maximaux de $\mathcal{A}$ [Bour, AC II, \S\  3, $\mbox{n}^\circ\  3$, cor 4 du th\'eo 1] ; comme $\mathcal{A}/I \mathcal{A} \simeq A/I A \neq \{0 \}$ par hypoth\`ese, il existe donc $\mathfrak{p} \in M$ tel que $I \mathcal{A}_{\mathfrak{p}} \displaystyle \mathop{\subset}_{\neq} \mathcal{A}_{\mathfrak{p}}$. Soit $\mathfrak{p} \in M$ tel que $I \mathcal{A}_{\mathfrak{p}}  \displaystyle \mathop{\subset}_{\neq} \mathcal{A}_{\mathfrak{p}}$ : l'id\'eal $I \mathcal{A}_{\mathfrak{p}}$ est contenu dans le seul id\'eal maximal $\mathfrak{p}  \mathcal{A}_{\mathfrak{p}}$ de $ \mathcal{A}_{\mathfrak{p}}$.
Or l'inclusion $\mathcal{A} \hookrightarrow \mathcal{A}_{\mathfrak{p}}$ donne l'inclusion $\varphi : \hat{\mathcal{A}} \hookrightarrow \widehat{(\mathcal{A}_{\mathfrak{p}})} := \displaystyle \mathop{lim}_{\leftarrow \atop{n}}\ \mathcal{A}_{\mathfrak{p}}/I^n\ \mathcal{A}_{\mathfrak{p}}$, o\`u
$\widehat{(\mathcal{A}_{\mathfrak{p}})}$ est local d'id\'eal maximal $\mathfrak{p} \widehat{(\mathcal{A}_{\mathfrak{p}})}$ [Bour, AC III, \S\  3, $\mbox{n}^\circ\  4$, prop 8 (ii)] : le sch\'ema $Spec\  \widehat{(\mathcal{A}_{\mathfrak{p}})} $ est connexe et son image par le morphisme dominant

$$Spec\ \varphi : Spec\  \widehat{(\mathcal{A}_{\mathfrak{p}})} \rightarrow Spec\ \hat{\mathcal{A}}$$

\noindent est un connexe dense, donc $Spec\ \hat{\mathcal{A}}$ est connexe ; comme $\hat{\mathcal{A}}$ est noeth\'erien (cf (1)) et normal (cf (4) (ii)), il en r\'esulte que $\hat{\mathcal{A}}$ est int\`egre [EGA I, (4.5.6)] et int\'egralement clos. $\square$

\vskip6mm
\section*{2. Des \'equivalences de cat\'egories}

\vskip 3mm
\noindent \textbf{Th\'eor\`eme (2.1)}.
\textit{Soit $A$ une $\mathcal{V}$-alg\`ebre telle que l'anneau $A$ soit noeth\'erien. On suppose que le morphisme canonique $Spec\  \hat{A} \rightarrow Spec\ A$ est normal (vrai par exemple si $A$ est excellent). Alors}
\begin{enumerate}
\item[(1)] \textit{Le morphisme canonique
$$f : Spec\ \hat{A} \rightarrow Spec\ \tilde{A} $$
est normal, fid\`element plat, \`a fibres g\'eom\'etriquement int\`egres.}
\item[(2)] 
\begin{itemize}
\item[(i)] \textit{Si $\tilde{A}$ est r\'eduit alors $\tilde{A}$ est int\'egralement ferm\'e dans $\hat{A}$}.
\item[(ii)] \textit{Le} (i) \textit{est v\'erifi\'e si $A$ est r\'eduit.}
\item[(iii)] \textit{On a les \'equivalences :\\
$\tilde{A}$ int\`egre (resp. int\'egralement clos)\\
$\Leftrightarrow \hat{A}$ int\`egre (resp. int\'egralement clos).}
\end{itemize}
\end{enumerate}

\vskip 3mm
\noindent \textit{D\'emonstration}. Pour le (1), $\hat{A}$ est le s\'epar\'e compl\'et\'e $I$-adique de $\tilde{A}$, et $\tilde{A}$ est de Zariski [prop (1.7) (1)], donc $f$ est fid\`element plat [Bour, AC III, \S\  3, n$^\circ\ 5$, prop 9]. L'hypoth\`ese entra\^{\i}ne alors que $f$ est normal \`a fibres g\'eom\'etriquement int\`egres via [prop (1.1) (3)], car $(\tilde{A}, \tilde{I})$ est un couple hens\'elien [R 2, th\'eo 3, p126].\\

Le (2) (i) r\'esulte de [prop (1.1) (4)], car $(\tilde{A}, \tilde{I})$ est un couple hens\'elien, et le (2) (ii) provient de [prop (1.7) (2) (i)].\\

Dans le (2) (iii) l'assertion dans le cas ``int\'egralement clos'' est prouv\'ee dans [prop (1.7), (4) (iv)]. Pour le cas ``int\`egre'', comme $\tilde{A} \hookrightarrow \hat{A}$ est fid\`element plat, il suffit de prouver que si $\tilde{A}$ est int\`egre, alors $\hat{A}$ l'est : supposons $\tilde{A}$ int\`egre, alors $\hat{A}$ est r\'eduit [prop (1.7) (4) (i)] ; or les fibres de $f$ sont g\'eom\'etriquement int\`egres, ainsi  $\mbox{Spec}\  \hat{A}$ est irr\'eductible [EGA IV, (2.3.5) (iii)] et r\'eduit, donc $\hat{A}$ est int\`egre. $\square$

\vskip 3mm
\noindent \textbf{Th\'eor\`eme (2.2)}.
\textit{Soit $A$ une $\mathcal{V}$-alg\`ebre de type fini ; on suppose que le morphisme canonique $Spec\  \hat{A} \rightarrow Spec\  A$ est normal (vrai par exemple si $\mathcal{V}$ est excellent). Alors :}

\begin{enumerate}
\item[(1)] \textit{Les morphismes canoniques
$$g : Spec\ \hat{A} \longrightarrow Spec\ A^\dag\ ,\  h : Spec\ A^\dag \longrightarrow Spec\ \tilde{A} $$
sont normaux, fid\`element plats, \`a fibres g\'eom\'etriquement int\`egres}.
\item[(2)] (i) \textit{Si $\tilde{A}$ est r\'eduit, $\tilde{A}$ est int\'egralement ferm\'e dans $A^\dag$ et dans $\hat{A}$.} \\
(ii) \textit{Si $A^\dag$ est r\'eduit, $A^\dag$ (resp. $A^\dag_K$) est int\'egralement ferm\'e dans $ \hat{A}$ (resp. dans $\widehat{A_K}$).} \\
(iii) \textit{Les hypoth\`eses (i) et (ii) sont satisfaites si $A$ est r\'eduit.} \\
(iv) \textit{On a les \'equivalences : } \\
\textit{$\tilde{A}$ int\`egre (resp. int\'egralement clos)\\
 $ \Updownarrow$ \\
${A}^\dag$ int\`egre (resp. int\'egralement clos)\\
$\Updownarrow$ \\
$\hat{A}$ int\`egre (resp. int\'egralement clos)}.
\end{enumerate}

\vskip 3mm
\noindent \textit{D\'emonstration}.
\begin{enumerate}
\item [(1)] La preuve pour $g$ est identique \`a celle du th\'eor\`eme pr\'ec\'edent en remarquant que ($A^\dag, I A^\dag)$ est un couple hens\'elien [Et 4, th\'eo 3]. \\
 Le cas de $h$ a \'et\'e trait\'e dans [prop (1.1) (3)].
\item [(2)] Le (i) et le (ii) ont \'et\'e vus dans le (4) de [prop (1.1)]. \\
Le (iii) provient du (4) (i) de [prop (1.7)].\\
Dans le (iv) l'assertion dans le cas ``int\'egralement clos'' est prouv\'ee dans [prop (1.7) (4) (iv)]. Pour le cas ``int\`egre'' les inclusions $\tilde{A}\ \subset\ A^\dag\ \subset \hat{A}$ ram\`enent la preuve au (2) (iii) du th\'eor\`eme (2.1). $\square$
\end{enumerate}

 \vskip 3mm
 \noindent \textbf{Remarque (2.2.1)}.
 \noindent Lorsque $\mathcal{V}$ est un anneau de valuation discr\`ete complet et $A$ une $\mathcal{V}$-alg\`ebre de type fini, $A$ est excellent [EGA IV, (7.8.3)]. Si $A^\dag$ est r\'eduit le (2) (ii) du th\'eor\`eme (2.2) prouve que $A^\dag$ est int\'egralement ferm\'e dans $\hat{A}$ : on retrouve ainsi l'analogue du th\'eor\`eme 2 de Bosch, Dwork, Robba de [Bo-Dw-R] lorsque la valuation de $K$ de [loc. cit.] est discr\`ete ; cf. aussi [Bo].\\
 
Les th\'eor\`emes (2.3) et (2.4) qui suivent am\'eliorent la th\'eor\`eme 15 de [Et 4].

\vskip 3mm
\noindent \textbf{Th\'eor\`eme (2.3)}.
\textit{Soit $A$ une $\mathcal{A}$-alg\`ebre telle que l'anneau $A$ soit noeth\'erien. On suppose que le morphisme canonique $Spec\  \hat{A} \rightarrow Spec\ A$ est normal (vrai si $A$ est excellent). Alors}
\begin{enumerate}
\item[(1)] \textit{Le foncteur $\mathcal{E}$ de la cat\'egorie des $\tilde{A}$-sch\'emas \'etales dans la cat\'egorie des $\hat{A}$-sch\'emas \'etales d\'efini par}

$$f : Spec\ \hat{A} \rightarrow Spec\ \tilde{A} $$
\textit{est pleinement fid\`ele.}
\item[(2)]  (i) \textit{Le foncteur $\mathcal{F} : B \mapsto B\otimes_{\tilde{A}} \hat{A}$ de la cat\'egorie des $\tilde{A}$-alg\`ebres finies \'etales dans la cat\'egorie des $\hat{A}$-alg\`ebres finies \'etales d\'efini par $f$ est une \'equivalence de cat\'egories}.\\
(ii) \textit{Si $\tilde{A}$ est r\'eduit (c'est le cas par exemple si $A$ est r\'eduit), le foncteur $\mathcal{G}$ qui \`a une $\hat{A}$-alg\`ebre finie \'etale $C$ associe la fermeture int\'egrale de $\tilde{A}$ dans $C$ est un foncteur quasi-inverse de $\mathcal{F}$}.
\end{enumerate}

\vskip 3mm
\noindent \textit{D\'emonstration}. Le (1) et le (2) (i) se d\'emontrent comme le (1) du th\'eor\`eme 15 de [Et 4].\\

Pour (2) (ii) il suffit, gr\^ace au fait que $\mathcal{F}$ est une \'equivalence de cat\'egories, de prouver que si $B$ est une $\tilde{A}$-alg\`ebre finie \'etale, alors $B$ est la fermeture int\'egrale de $\tilde{A}$ dans $B\otimes_{\tilde{A}} \hat{A} \simeq \hat{B}.$\\

\quad Supposons $\tilde{A}$ r\'eduit et soit $B$ une $\tilde{A}$-alg\`ebre finie \'etale : alors $(B, IB)$ est un couple hens\'elien [R 2, prop 2 (1) p 124], $B$ est r\'eduit par le (1) (i) du [lemme (1.6)] et $\mbox{Spec}\  \hat{B} \rightarrow \mbox{Spec}\  B$ est un morphisme normal [EGA IV, (6.8.3) (iii)] ; par le (4) de la [prop (1.1)] on en d\'eduit que $B$ est int\'egralement ferm\'e dans $\hat{B}$. Comme $B$ est fini sur $\tilde{A}$, $B$ est bien la fermeture int\'egrale de $\tilde{A}$ dans $\hat{B}.\   \square$

\vskip 3mm
\noindent \textbf{Th\'eor\`eme (2.4)}.
\textit{Soit $A$ une $\mathcal{V}$-alg\`ebre de type fini ; on suppose que le morphisme canonique  $Spec\  \hat{A} \rightarrow Spec\ A$ est normal (vrai par exemple si $\mathcal{V}$ est excellent) et on d\'esigne par $(\mathcal{A, B})$ l'un des couples $(\tilde{A}, A^{\dag}), (A^{\dag}, \hat{A})$. Alors}

\begin{enumerate}
\item[(1)] \textit{Le foncteur $\mathcal{E}$ de la cat\'egorie des $\mathcal{A}$-sch\'emas \'etales dans la cat\'egorie des $\mathcal{B}$-sch\'emas \'etales d\'efini par}
                       $$f : Spec\ \mathcal{B} \rightarrow Spec\ \mathcal{A} $$
\textit{est pleinement fid\`ele.}

\item[(2)]  (i) \textit{Le foncteur : $\mathcal{F} : B \mapsto B\otimes_{\mathcal{A}} \mathcal{B}$ de la cat\'egorie des $\mathcal{A}$-alg\`ebres finies \'etales dans la cat\'egorie des $\mathcal{B}$-alg\`ebres finies \'etales d\'efini par $f$ est une \'equivalence de cat\'egories}.\\
(ii) \textit{Si $\mathcal{A}$ est r\'eduit (c'est le cas par exemple si $A$ est r\'eduit), le foncteur $\mathcal{G}$ qui \`a une $\mathcal{B}$-alg\`ebre finie \'etale $C$ associe la fermeture int\'egrale de $\mathcal{A}$ dans $C$ est un foncteur quasi-inverse de $\mathcal{F}$}.
\end{enumerate}

\vskip 3mm
\noindent \textit{D\'emonstration}.
\begin{enumerate}
\item [(1)] Ici $f$ est fid\`element plat et quasi-compact, donc universellement submersif [EGA I, (3.9.4) (ii) et (7.3.5)] ; de plus les fibres de $f$ sont g\'eom\'e\-triquement int\`egres [th\'eo (2.2)]. En vertu de [SGA 1, IX, cor 3.4] le foncteur $\mathcal{E}$ d\'efini par $f$ est pleinement fid\`ele.
\item [(2)] Le (i) se montre comme [Et 4, th\'eo 15, 2 (i) et (ii)]. La preuve du (ii) est la m\^eme que celle du (2) (ii) du [th\'eor\`eme (2.3)] si l'on rappelle que toute $A^{\dag}$-alg\`ebre finie $B$ est  ``f.c.t.f '' [Et 4, prop 1] et que $(B, IB)$ est un couple hens\'elien [loc. cit, th\'eo 3]. $\square$
\end{enumerate}
 
 \vskip 3mm
 \noindent \textbf{Remarque (2.4.1)}.
 \noindent La partie (1) des th\'eor\`emes (2.3) et (2.4) pr\'ec\'edents est une g\'en\'eralisation de [EGA IV, (18.9.5)].

 \vskip6mm
\section*{3. Sch\'emas formels et rel\`evements de sch\'emas}

\vskip 3mm
\noindent \textbf{Th\'eor\`eme (3.1)}.
\textit{Soit $A$ une $\mathcal{V}$-alg\`ebre normale de caract\'eristique z\'ero telle que $A$ soit un anneau noeth\'erien ; on suppose $A$ excellent et $0 \notin T : = 1 + IA$. On note $A_{0} = A/IA$. Alors}
\begin{enumerate}
\item[(1)] \textit{Si $\varphi : S' \rightarrow Spec\ A_{0} = : S$ est un morphisme fini \'etale, il existe un morphisme fini}

$$\psi : Spec\ B \rightarrow Spec\ A $$
\textit{relevant $\varphi$, o\`u $B$ est normal (resp. $B$ est int\'egralement clos si $A$ l'est) et tel que}

$$\psi_{T} : Spec\ B_{T} \rightarrow Spec\ A_{T}$$
\textit{soit un rel\`evement fini \'etale de $\varphi$.}

\item[(2)]   \textit{De plus il existe $g_{0} \in A_{0}$ et $g \in A$ relevant $g_{0}$ tels que }

$$\psi_{g} : Spec\ B_{g} \rightarrow Spec\ A_{g}$$
\textit{soit un rel\`evement fini \'etale de}

$$ \phi_{g_{0}}  : S'_{g_{0}} \rightarrow Spec\ (A_{0_{g_{0}}}).$$

\item[(3)] \textit{Supposons de plus $\mathcal{V}$ excellent, $A$ de type fini sur $\mathcal{V}$ et fixons une pr\'esentation}

$$A = \mathcal{V} [t_{1}, ...,t_{n}] / (f_{1},...,f_{r}).$$
\end{enumerate}

\textit{Notons $P$ la fermeture projective de $Spec\ A$ dans $\mathbb{P}^n_{\mathcal{V}}$, 
$P'$ le normalis\'e de $P$ et $P''$ la fermeture int\'egrale de $P$ dans l'anneau $\mathcal{R}\  (Spec\  B)$ des fonctions rationnelles sur $Spec\ B$. Les morphismes structuraux $P'' \rightarrow P'$ et $P' \rightarrow P$ sont finis, leur compos\'e $\theta : P'' \rightarrow P$ aussi, et on a des carr\'es cart\'esiens}

$$\xymatrix{
Spec\  B\  \ar@{^{(}->}[r] \ar[d]_{\psi} & P'' \ar[d]^{\theta}\\
Spec\  A\  \ar@{^{(}->}[r] &  P
}
\qquad \qquad
\xymatrix{Spec\  B\  \ar@{^{(}->}[r] \ar[d]_{\psi} & P'' \ar[d]^{\theta'}\\
Spec\  A\  \ar@{^{(}->}[r] &  P'
}
$$

\noindent \textit{o\`u les fl\`eches horizontales sont des immersions ouvertes fournissant par passage aux s\'epar\'es compl\'et\'es des carr\'es cart\'esiens}

\begin{center}
$\begin{array}{ccc}
\mathcal{S'} := Spf \hat{B}  & \overset{j'}{\hooklongrightarrow} & \hat{P}'' =: \overline{\mathcal{S}'}\\
\hat{\psi} \downarrow & &\downarrow \hat{\theta}\\
\mathcal{S} := Spf  \hat{A} & \underset{\tilde{j}}{\hooklongrightarrow}  & \hat{P} =: \tilde{\mathcal{S}}
\end{array}$
\qquad\qquad
$\begin{array}{ccc}
Spf\  \hat{B}  & \hooklongrightarrow & \hat{P}'' = \overline{\mathcal{S}'}\\
\hat{\psi} \downarrow & &\downarrow \hat{\theta'}\\
Spf\  \hat{A} & \underset{j}{\hooklongrightarrow} & \hat{P}' = \overline{\mathcal{S}}\\
\end{array}$
\end{center}

\noindent \textit{o\`u $\hat{\psi}$ est un rel\`evement fini \'etale de $\varphi, \hat{\theta}$ et $\widehat{\theta}'$ sont finis,  
$\widehat{P'}$ est normal et les fl\`eches horizontales sont des immersions ouvertes.} \\

\noindent \textit{De plus $\mathcal{O}_{\widehat{P''}}$ est la fermeture int\'egrale de $\mathcal{O}_{\widehat{P}}$ 
(resp. de $\mathcal{O}_{\widehat{P'}})$ dans $\mathcal{O}_{\mathcal{S}'}$ [EGA II, (6.3.2)].}\\

\noindent \textit{Enfin, $\hat{\theta}$ (resp.$\widehat{\theta}')$ est plat si et seulement si sa r\'eduction module I est plate.}\\

\noindent Donnons tout de suite un corollaire et sa d\'emonstration avant de prouver le th\'eor\`eme.

\vskip 3mm
\noindent \textbf{Corollaire (3.2)}.
 \textit{Soit $\mathcal{V}$ un anneau de valuation discr\`ete complet de caract\'eristique z\'ero, d'id\'eal maximal $I$ et de corps r\'esiduel $k$.  Si $A_{0}$ est une $k$-alg\`ebre lisse et $\varphi : Spec\ B_{0} \rightarrow Spec\ A_{0}$ est un morphisme fini \'etale, alors il existe une $\mathcal{V}$-alg\`ebre lisse $A$ et un morphisme fini $\psi : Spec\ B \rightarrow Spec\ A$ relevant $\varphi$ et satisfaisant aux propri\'et\'es du th\'eor\`eme (3.1).}
 
\vskip 3mm
\noindent \textit{D\'emonstration du corollaire (3.2)}. L'existence d'une $\mathcal{V}$-alg\`ebre lisse $A$ relevant $A_{0}$ r\'esulte du th\'eor\`eme 6 de Elkik [E$\ell$]. Le $\mathcal{V}$ du corollaire est r\'egulier [EGA $O_{IV}$, (17.1.4) (ii)] et excellent [EGA IV, (7.8.3)] : comme $A$ est lisse sur $\mathcal{V}$, $A$ est r\'egulier [EGA IV, (17.5.8)] et excellent [EGA IV, (7.8.3)]. Il suffit alors d'appliquer le th\'eor\`eme (3.1).
$\square$
\vskip 3mm
\noindent \textit{D\'emonstration du th\'eor\`eme (3.1)}.\\

\noindent \textbf{(1) et (2)} . D'apr\`es [EGA IV, (18.3.2)] il existe une $\hat{A}$-alg\`ebre finie \'etale $C$ telle que $Spec\ C \rightarrow Spec\  \hat{A}$ rel\`eve $\varphi$. Puisque $A$ est noeth\'erien normal on peut d\'ecomposer $Spec\ A$ en somme de ses composantes connexes $\displaystyle \mathop{\coprod_{i}} Spec\ A_{i}$, avec $A_{i}$ int\'egralement clos, et $\displaystyle \mathop{\coprod_{i}} Spec\ \hat{A}_{i}$ est une d\'ecomposition de $Spec\  \hat{A}$ en somme de ses composantes connexes, avec $\hat{A}_{i}$ int\'egralement clos [prop (1.7) (4) (iv)] : $C$ est aussi normal noeth\'erien et on le d\'ecompose de m\^eme. Comme $Spec\ C \rightarrow Spec\ \hat{A}$ est fini et plat, on est ramen\'e au cas o\`u ce morphisme est surjectif avec $C$ et $\hat{A}$ int\'egralement clos.\\

Soient $L$ le corps des fractions de $\hat{A}$ et $L_{1}$ celui de $C$ : d'apr\`es [EGA II, (6.1.8)]  $L \hookrightarrow L_{1}$ est une extension finie de corps de caract\'eristique nulle (donc l'extension est s\'eparable) et $C$ est la fermeture int\'egrale de $\hat{A}$ dans $L_{1}$. Par le th\'eor\`eme de l'\'el\'ement primitif il existe $x \in L_{1}$ tel que $L_{1} = L[x] : x$ est s\'eparable sur $L$ [Bour, A V, prop 6 p 38], son polyn\^ome minimal $f(X) \in L[X]$ est s\'eparable [Bour, A V, prop 5 p 38], donc $f \wedge f' = 1$ dans $L[X]$ [Bour, A V, prop 3 p 36]  ; appliquant B\'ezout dans $L[X]$, il existe $g_{1} \in \hat{A}$ tel que $f$, $f' \in (\hat{A})_{g_{1}} [X]$ et tel qu'il existe $u, v \in (\hat{A})_{g_{1}} [X]$ v\'erifiant l'identit\'e $u f + v f' = 1$ dans $(\hat{A})_{g_{1}} [X]$. Ainsi le morphisme canonique 

$$\mu : (\hat{A})_{g_{1}} \longrightarrow   (\hat{A})_{g_{1}} [X]\  /\  (f)$$

\noindent est fini \'etale [Mi, I, 3.4] et s'ins\`ere dans le carr\'e cocart\'esien\\

$$
\xymatrix{
(\hat{A})_{g_{1}} \ar[r]^{\mu} \ar[d] & (\hat{A})_{g_{1}}[X]/(f) = : D \ar[d]\\
Frac(\hat{A})_{g_{1}} = L\quad  \ar@{^{(}->}[r]   & \quad L_{1} = L[X]/(f) = L\otimes_{(\hat{A})_{g_{1}}} D\ ;
}
$$

\noindent par suite $\mu$ est fid\`element plat puisqu'il est injectif. La platitude de $\mu$ fournit l'injectivit\'e de $D \hookrightarrow L_{1}$, donc $D$ est int\`egre (de corps des fractions $L_{1}$) et normal [EGA IV, (17.5.7)], donc int\'egralement clos ; donc $D$ est la fermeture int\'egrale de $(\hat{A})_{g_{1}}$ dans $L_{1}$. Par changement de base, $\mu$ fournit le morphisme fini \'etale fid\`element plat

$$\rho : \widehat{(\hat{A})_{g_{1}}} \hookrightarrow \widehat{(\hat{A})_{g_{1}}} [X]\  /\  (\hat{f}),$$
o\`u $\hat{f}$ est l'image de $f$ dans $\widehat{(\hat{A})_{g_{1}}} [X].$\\

Soit $g_{2} \in A$ un rel\`evement de ($g_{1}$  mod $I$) $\in A_{0}$ : comme $\widehat{(\hat{A})_{g_{1}}}$ est une $A$-alg\`ebre formellement \'etale par les topologies $I$-adiques il existe un unique $A$- morphisme (en fait un $\hat{A}$-morphisme)

$$\nu :  \widehat{(\hat{A})_{g_{1}}} \longrightarrow  (\widehat{A_{g_{2}}})$$

\noindent relevant l'identit\'e de  $A_{g_{2}}\ /\ IA_{g_{2}}$ et $\nu$ est un isomorphisme [EGA $O_{I}$, (6.6.21)]. Notons $f_{1}(X) \in A_{g_{2}} [X]$ un polyn\^ome unitaire relevant 

$$
f(X) \textrm{mod} I \in (\hat{A})_{g_{1}}\ [X]\ /\   I (\hat{A})_{g_{1}} [X] \simeq A_{g_{2}} [X]\ /\ IA_{g_{2}} [X]\  ; 
$$

\noindent alors, si $\hat{f}_{1}$ d\'esigne l'image de $f_{1}$ dans $(\widehat{A_{g_{2}}}) [X]$, $(\widehat{A_{g_{2}}}) [X]\ /\ (\hat{f}_{1})$ est fini et plat sur $(\widehat{A_{g_{2}}})$ car $\hat{f_{1}}$ est unitaire [Mi, I, 2.6 (a)]. Comme $\rho$ est fini \'etale, que $D$ et $\widehat{A_{g_{2}}} [X]\ /\ (\hat{f_{1}})$ ont m\^eme r\'eduction mod $I$ et que $\widehat{(\hat{A})_{g_{1}}}$ s'identifie \`a $\widehat{A_{g_{2}}}$ via $\nu$, il existe un unique $\widehat{A_{g_{2}}}$-morphisme

$$\widehat{A_{g_{2}}} [X]\ /\ (\nu(\hat{f})) \rightarrow \widehat{A_{g_{2}}} [X]\ /\ (\hat{f_{1}})$$

\noindent qui est un isomorphisme [EGA $O_{I}$, (6.6.21)] relevant l'identit\'e de $D/ID$.\\
Puisqu'on a des injections

$$(\hat{A})_{g_{1}} \hookrightarrow (\hat{A})_{g_{1},T} \hookrightarrow \widehat{((\hat{A)}_{g_{1},T})} = \widehat{((\hat{A})_{g_{1}})},$$

$$A_{g_{2}} \hookrightarrow A_{g_{2}T} \hookrightarrow \widehat{(A_{g_{2},T}}) = \widehat{A_{g_{2}}}$$
et que $f$ et $f_{1}$ sont unitaires on en d\'eduit que

$$\nu(\hat{f}) = \widehat{f_{1}},$$
i.e. que dans l'\'ecriture $L_{1} = L[X]\ /\ (f)$ on peut supposer $f = f_{1} \in A_{g_{2}} [X]$ : ainsi il existe $g_{3} \in A$ tel que $f, f' \in A_{g_{3}} [X]$ et tel qu'il existe $u, v \in A_{g_{3}} [X]$ v\'erifiant l'identit\'e $u f + v f' = 1$ dans $A_{g_{3}} [X]$. En particulier le morphisme canonique

$$\eta : A_{g_{3}} \rightarrow A_{g_{3}} [X]\ /\ (f)$$

\noindent est fini \'etale et s'ins\`ere dans le diagramme commutatif

$$
\xymatrix{
A_{g_{3}}[X]/(f) \ar[r]^{\tau}  & (\hat{A})_{g_{3}}[X]/(f) \ar[r] & L_{1}& \\
A_{g_{3}} \ar[u] ^{\eta}\ar@{^{(}->}[r] & (\hat{A})_{g_{3}} \ar[u] \ar@{^{(}->}[r] & L \ar[u] &\  ;
}
$$

\noindent par suite $\eta$ est injectif, donc fid\`element plat.\ Par le m\^eme raisonnement que ci-dessus pour $D = (\hat{A})_{g_{1}} [X]\ /\ (f)$ on montre que $(\hat{A})_{g_{3}} [X]\ /\ (f)$ est int\`egre et int\'egralement clos de corps des fractions $L_{1}$. La platitude de $\eta$ fournit l'injectivit\'e de $\tau$ : donc $A_{g_{3}} [X]\ /\ (f)$ est int\`egre de corps des fractions not\'e $K_{1}$ ; on notera $K$ le corps des fractions de $A$. On sait par [EGA IV, (18.10.12)] que $K_{1}$ est une extension finie \'etale de $K$ et que $A_{g_{3}} [X]\ /\ (f)$ est la fermeture int\'egrale de $A_{g_{3}}$ dans $K_{1}$ : comme $f$ est irr\'eductible dans $L[X]$, il l'est aussi dans $K[X] \subset L[X]$, d'o\`u $K_{1} = K[X]\ /\ (f)$.\\

Comme $A$ est excellent, il est universellement japonais [EGA IV, (7.8.3) (vi)] : la fermeture int\'egrale $B$ de $A$ dans $K_{1} = K[X]\ /\ (f)$ est une $A$-alg\`ebre finie, $Spec\ B \rightarrow Spec\ A$ est surjectif et $B_{g_{3}} = A_{g_{3}} [X]\ /\ (f)$ ; $B$ est aussi la fermeture int\'egrale de $A$ dans $A_{g_{3}} [X]\ /\ (f)$, $B_{T}$ est la fermeture int\'egrale de $A_{T}$ dans  $A_{g_{3,T}} [X]\ /\ (f)$, et le corps des fractions de $B$ est $K_{1}$ (donc $B$ est int\'egralement clos). Comme $Spec\ \hat{A} \rightarrow Spec\ A$ est un morphisme r\'egulier [EGA IV, (7.8.3) (v)], donc normal, la fermeture int\'egrale de $\hat{A}$ dans $(\hat{A})_{g_{3}} [X]\ /\ (f)$ est \'egale \`a $\hat{B} = B \otimes_{A} \hat{A}$ [EGA IV, (6.14.4)] : or $(\hat{A})_{g_{3}} [X]\ /\ (f)$ est int\'egralement clos ; donc $\hat{B}$ est la fermeture int\'egrale de $\hat{A}$ dans $L_{1}$, d'o\`u $\hat{B} = C$. On en d\'eduit que $B$ et $C$ ont m\^eme r\'eduction mod $I$, d'o\`u l'existence du morphisme fini

$$\psi : Spec\ B \rightarrow Spec\ A$$

\noindent relevant $\varphi$ et il suffit de prendre $g = g_{3}, g_{0} = g$ mod $I$.\\

De plus $Spec\  \hat{A} = Spec\  \hat{A}_{T} \rightarrow Spec\ A_{T}$ est un morphisme normal puisque $A_{T}$ est excellent : par suite la fermeture int\'egrale de $\hat{A} = \widehat{A_{T}}$ dans $(\hat{A})_{g_{3}} [X]\ /\ (f)$, qu'on sait \^etre \'egale \`a $C = \hat{B}$, est aussi \'egale \`a $B_{T} \otimes_{{A_{T}}} \hat{A} = \widehat{B_{T}} = \hat{B}$ [EGA IV, (6.14.4)]. Par passage aux s\'epar\'es compl\'et\'es $\psi$ induit

$$
\hat{\psi} = Spec\ \hat{B} \rightarrow\ Spec\ \hat{A}
$$ 

\noindent qui s'identifie \`a notre morphisme fini \'etale

$$
Spec\ C \rightarrow Spec\ \hat{A}\  ;
$$

\noindent par fid\`ele platitude de $\hat{A}$ sur $A_{T}$ le morphisme

$$\psi_{T} = Spec\  B_{T} \rightarrow Spec\  A_{T}$$

\noindent est donc fini \'etale et c'est clairement un rel\`evement de $\varphi : S' \rightarrow Spec\  A_{0}$.\\

\noindent \textbf{(3)} On se ram\`ene \`a $A$ et $B$ int\'egralement clos comme en (1) dont on reprend les notations. Le sch\'ema $P$ est int\`egre, car on a suppos\'e $A$ int\`egre [EGA I, (6.10.5)], d'o\`u $\mathcal{R}(P) = \mathcal{R}\ (Spec\ A) = \mbox{Frac}\ A = K$ est un corps [EGA I, (8.1.5)]. Le sch\'ema $P$ est excellent, donc pour chaque ouvert $U = Spec\ R \subset P$, $R$ est japonais [EGA IV, (7.8.3)] : ainsi la fermeture int\'egrale $P'$ (resp. $P''$) de $P$ dans $\mathcal{R}\  (Spec\ A) = K$ (resp. dans $\mathcal{R}\  (Spec\ B) = K_{1} = K[X]\ /\ (f))$ est un $P$-sch\'ema fini et on a $\mathcal{R}(P') = K$ (resp. $\mathcal{R}(P'') = K_{1})$ [EGA II, (6.3.7)] : \'evidemment $P''$ est aussi la fermeture int\'egrale de $P'$ dans $K_{1}$ et $P'' \rightarrow P'$ est fini. De plus $P'$ est int\'egralement clos car noeth\'erien normal et int\`egre [EGA II, (6.3.8)] ; comme $A$ est int\'egralement clos, $P'$ est aussi la fermeture int\'egrale de $P$ dans $Spec\ A$ : on a donc une immersion ouverte

$$Spec\ A \hookrightarrow P' $$
[R 2, cor 2, p 42]. Par [EGA II, Rq entre (6.3.4) et (6.3.5)] ou [EGA IV, (6.14.4)] les carr\'es

$$\xymatrix{
Spec\  B\  \ar@{^{(}->}[r] \ar[d]_{\psi} & P'' \ar[d]^{\theta}\\
Spec\  A\  \ar@{^{(}->}[r] &  P
}
\qquad \qquad
\xymatrix{Spec\  B\  \ar@{^{(}->}[r] \ar[d]_{\psi} & P'' \ar[d]^{\theta'}\\
Spec\  A\  \ar@{^{(}->}[r] &  P'
}
$$

\noindent sont cart\'esiens.\\

On conclut la d\'emonstration du th\'eor\`eme (3.1) par passage aux s\'epar\'es compl\'et\'es : que $\mathcal{O}_{\widehat{P''}}$ soit la fermeture int\'egrale de $\mathcal{O}_{\hat{P}}$ (resp. $\mathcal{O}_{\widehat{P'}})$ dans $\mathcal{O}_{\mathcal{S}'}$ r\'esulte de [EGA IV, (6.14.4)] compte tenu du fait que $\hat{P} \rightarrow P$ est un morphisme normal, car $P$ est excellent. De m\^eme $\widehat{P'} \rightarrow P'$ est normal : ainsi $\widehat{P'}$ est normal car $P'$ l'est [prop (1.7)].\\

Il nous reste \`a montrer que si $\hat{\theta}$ mod $I$ (resp. $\widehat{\theta}'$ mod $I$) est plat, alors $\hat{\theta}$  (resp. $\widehat{\theta}'$) est plat : faisons la d\'emonstration pour $\hat{\theta}$. Comme ci-dessus on peut supposer $P$ int\`egre et se limiter \`a un ouvert affine $V = Spec\ R$ de $P$ : alors l'ensemble $U$ des points de $V$ tels que la restriction de $\theta$ \`a $\theta^{-1}(U)$ soit plate est un ouvert non vide de $V$ [EGA IV, (11.1.1), (2.4.6), (6.9.1)] et $U$ contient la r\'eduction $V_{0}$ de $V$ modulo $I$ par hypoth\`ese. En posant : $\tilde{U} : = U \times_{V} Spec\ \hat{R}$ et $\tilde{P''} : = \tilde{U} \times_{P} P''$, soit $\tilde{\theta} : \tilde{P''} \rightarrow \tilde{U}$ l'image inverse de $\theta$ par le changement de base $\tilde{U} \rightarrow P : \tilde{\theta}$ est plat. Or $\tilde{U}$ est un ouvert de $Spec\ \hat{R}$ qui contient $V_{0}$, donc $\tilde{U} = Spec\  \hat{R}$ via [EGA IV, (18.5.4.3)] car $(Spec\ \hat{R}, V_{0})$ est un couple hens\'elien [Et 4, th\'eo 3]. Par passage aux compl\'et\'es formels, $\hat{\theta} : \widehat{P''} \rightarrow \hat{P}$ est plat. $\square$

\vskip 3mm
\noindent \textbf{Th\'eor\`eme (3.3)}.
\textit{Soient $\mathcal{V}$ un anneau excellent normal, $I \subset \mathcal{V}$ un id\'eal et  $\mathcal{V}_{0} = \mathcal{V}/I$  tel que $(\mathcal{V}, \mathcal{V}_{0})$ soit un couple hens\'elien au sens de [EGA IV, (18.5.5)] ; on suppose $\mathcal{V}_{0}$ normal. Soient $S_{0} = Spec\ A_{0}$ un $\mathcal{V}_{0}$-sch\'ema affine et lisse, $A = \mathcal{V} [t_{1},...,t_{d}] / J$  une $\mathcal{V}$-alg\`ebre lisse relevant $A_{0}$ et dont on a fix\'e une pr\'esentation et $S = Spec\ A$: on note $\hat{A}$ le s\'epar\'e compl\'et\'e de $A$ pour la topologie $I$-adique, $A^{\dag}$ son compl\'et\'e faible, $\tilde{A}$ l'hens\'elis\'e de $A$ au sens de Raynaud et $\hat{S}=Spec \ \hat{A},\ S^{\dag}=Spec\ A^{\dag},\  \tilde{S}= Spec\ \tilde{A}$. On d\'esigne par $\overline{S}$ l'adh\'erence sch\'ematique de $S$ dans $\mathbb{P}^d_{\mathcal{V}}$, et par $\mathcal{S}$ (resp. $\overline{\mathcal{S}})$ le compl\'et\'e formel $I$-adique de $S$ (resp. de $\overline{S}).$}\\

\noindent \textit{Soit $f : X_{0} \rightarrow S_{0}$ un $\mathcal{V}_{0}$-morphisme projectif. Alors}\\
\begin{itemize}
\item[(3.3.1)] \textit{Il existe un carr\'e cart\'esien}
$$
\begin{array}{c}
\xymatrix{
 X \ar@{^{(}->}[r] \ar[d]_{h} & \overline{X} \ar[d]^{\overline{h}}\\
 S  \ar@{^{(}->}[r] &  \overline{S}
}
\end{array}
\leqno{(3.3.1.1)}
$$
 \textit{dans lequel $\overline{h}$ est projectif, $h$ est un rel\`evement projectif de $f$ et les fl\`eches horizontales sont des immersions ouvertes.}\\

\item[(3.3.2)] \textit{Consid\'erons le diagramme commutatif \`a carr\'es cart\'esiens d\'eduit de (3.3.1.1)}

$$
\begin{array}{c}
\xymatrix{
X_{\hat{S}}\ar[r] \ar[d]_{h_{\hat{S}}}&X_{S^{\dag}}\ar[r] \ar[d]_{h_{S^{\dag}}}&X_{\tilde{S}}\ar[r] \ar[d]_{h_{\tilde{S}}}& X \ar@{^{(}->}[r] \ar[d]_{h} & \overline{X} \ar[d]^{\overline{h}}&{}\\
\hat{S}\ar[r]& S^{\dag}\ar[r]& \tilde{S}\ar[r]&S  \ar@{^{(}->}[r] &  \overline{S}& 
}
\end{array}
\leqno{(3.3.2.1)}
$$
 \textit{dans lequel $h,h_{\tilde{S}}, h_{S^{\dag}}, h_{\hat{S}}$ sont des rel\`evements projectifs de $f$.\\
\vskip 3mm
Alors on a \'equivalence entre les propri\'et\'es suivantes:
	\begin {itemize} 
	\item[(i)] $X$ est plat sur $\mathcal{V}$ et $f$ est plat.
	\item [(ii)] $h_{\hat{S}}$ est plat.
	\item[(iii)] $h_{S^{\dag}}$ est plat.
	\item[(iv)] $h_{\tilde{S}}$ est plat.
	\end{itemize}
	}
\vskip 3mm
\item[(3.3.3)]
\textit{Alors on a \'equivalence entre les propri\'et\'es suivantes:
	\begin {itemize} 
	\item[(i)]$X$ est plat sur $\mathcal{V}$ et $f$ est lisse.
	\item [(ii)] $h_{\hat{S}}$ est lisse.
	\item[(iii)] $h_{S^{\dag}}$ est lisse.
	\item[(iv)] $h_{\tilde{S}}$ est lisse.
	\end{itemize}
	}
\vskip 3mm
 \item[(3.3.4)] \textit{Le carr\'e cart\'esien (3.3.1.1) fournit par passage aux compl\'et\'es formels  un carr\'e cart\'esien de $\mathcal{V}$-sch\'emas formels}

$$
\begin{array}{c}
\xymatrix{
\mathcal{X}\ \ar@{^{(}->}[r] \ar[d]_{\hat{h}} & \overline{\mathcal{X}} \ar[d]^{\hat{\overline{h}}}\\
\mathcal{S}\   \ar@{^{(}->}[r] &  \overline{\mathcal{S}}
}
\end{array}
\leqno{(3.3.4.1)}
$$

 \textit{dans lequel $\hat{h}$ est un rel\`evement projectif de $f$, $\hat{\overline{h}}$ est projectif et les fl\`eches horizontales sont des immersions ouvertes.\\
 De plus on a \'equivalence entre les propri\'et\'es suivantes:
	\begin {itemize} 
	\item[(i)]$X$ est plat sur $\mathcal{V}$ et $f$ est plat (resp. $f$ est lisse).
	\item [(ii)] $\hat{h}$ est plat (resp. $\hat{h}$ est lisse).
	\end{itemize}}
\end{itemize}

\vskip 3mm
\noindent \textit{D\'emonstration}. Quitte \`a d\'ecomposer les sch\'emas normaux noeth\'eriens $Spec\ \mathcal{V}$ et $S_{0}$ en somme de leurs composantes connexes on supposera dans toute la suite que  $Spec\ \mathcal{V}$ et $S_{0}$ sont connexes, donc int\'egralement clos.\\
\textit{Pour (3.3.1)}. Le morphisme projectif $f$ se factorise en $f : X_{0} \displaystyle\mathop{\hookrightarrow}^{i_{0}} \mathbb{P}^n_{A_{0}} \displaystyle\mathop{\rightarrow}^{s_{0}} S_{0} = Spec\ A_{0}$ o\`u $i_{0}$ est une immersion ferm\'ee et $s_{0}$ est le morphisme canonique. Soient $(x_{0},x_{1},..., x_{n})$ les coordonn\'ees projectives sur $\mathbb{P}^n_{A_{0}}$ (resp. sur $\mathbb{P}^n_{S}$): alors $X_{0}$ est isomorphe \`a $Proj(A_{0}[x_{0},x_{1},..., x_{n}]/\mathcal{J}^{0})$ pour un certain id\'eal homog\`ene $\mathcal{J}^{0}$ de l'anneau noeth\'erien $A_{0}[x_{0}, x_{1}, ..., x_{n}]$: $\mathcal{J}^{0}$ est engendr\'e par un nombre fini de polyn\^omes homog\`enes $f^1_{0},...,f^r_{0}$. Pour $\alpha\in\lbrace1,...,r\rbrace$ relevons $f^{\alpha}_{0}$ en un polyn\^ome homog\`ene $f^{\alpha} \in A[x_{0},x_{1},...,x_{n}]$ de m\^eme degr\'e en relevant coefficient par coefficient de $A_{0}$ \`a $A$, et soit $\mathcal{J}$ l'id\'eal (homog\`ene) engendr\'e par $f^1,...,f^r$ 
$$
\mathcal{J}=(f^1,...,f^r) \subset A[x_{0},x_{1},...,x_{n}]\ .
$$
\noindent D\'esignons par
$$i : X=Proj(A[x_{0},x_{1},...,x_{n}]/ \mathcal{J}) \hookrightarrow \mathbb{P}^n_{S}$$
l'immersion ferm\'ee et par $p_{S} : \mathbb{P}^n_{S} \rightarrow S = Spec\ A$ la projection canonique. Alors le morphisme compos\'e

$$h = p_{S} \circ  i : X \rightarrow S $$
est un rel\`evement projectif de $f$. Notons $p_{\overline{S}} : \mathbb{P}^n_{\overline{S}} \rightarrow \overline{S}$ la projection canonique, $\overline{X}$ la fermeture int\'egrale de $\mathbb{P}^n_{\overline{S}}$ dans $X$, $\overline{i} : \overline{X} \hookrightarrow \mathbb{P}^n_{\overline{S}}$ l'immersion ferm\'ee et $\overline{h} = p_{\overline{S}} \circ \overline{i} : \overline{h}$ est projectif. On dispose ainsi d'un carr\'e cart\'esien\\

$$
\begin{array}{c}
\xymatrix{
X\ \ar@{^{(}->}[r]_{j_{\overline{Y}}} \ar[d]_{h} & \overline{X} \ar[d]^{\overline{h}}\\
S\   \ar@{^{(}->}[r]_{j_{\overline{S}}} &  \overline{S}
}
\end{array}
\leqno{(3.3.1.1)}
$$

\noindent o\`u les fl\`eches horizontales sont des immersions ouvertes.\\

\noindent\textit{Pour (3.3.2).} On a le lemme suivant:\\ 

\noindent \textbf{Lemme(3.3.2.2)}. \textit{On a l'\'equivalence:}\\
$$X \ est\ plat\ sur\ \mathcal{V} \Longleftrightarrow X_{\hat{S}}\ est\ plat\ sur\ \mathcal{V}\ .  $$

\noindent\textit{D\'emonstration du lemme.}  Notons $v$ le morphisme compos\'e $X \displaystyle\mathop{\rightarrow}^{h}  S \rightarrow Spec\ \mathcal{V}$ et $w$ le compos\'e de $v$ avec le morphisme plat $X_{\hat{S}}\rightarrow X$. Si $v$ est plat, $w$ l'est. Supposons que $w$ soit plat: d'apr\`es le crit\`ere de platitude par fibres [EGA\ IV,(11.3.10)] $v$ est plat sur sa fibre sp\'eciale $X_{0}$, donc au-dessus de $\mathcal{V}_{0}$. Or par le th\'eor\`eme de platitude g\'en\'erique [EGA IV, (6.9.1)] il existe un ouvert non vide $V$ de $Spec\ \mathcal{V}$ au-dessus duquel $v$ est lisse; comme cet ouvert $V$ contient $Spec\ \mathcal{V}_{0}$ d'apr\`es le raisonnement fait ci-dessus, il r\'esulte de [EGA IV, (18.5.4.3)] que $V= Spec\ \mathcal{V}$ puisque $(\mathcal{V}, \mathcal{V}_{0})$ est un couple hens\'elien. D'o\`u le lemme. $\square$\\

On montre de la m\^eme fa\c con le lemme suivant:\\

\noindent \textbf{Lemme(3.3.2.3)}. \textit{Avec $\mathcal{V}$ comme dans le th\'eor\`eme (3.3), soit B une $\mathcal{V}$-alg\`ebre de type fini, de s\'epar\'e compl\'et\'e $I$-adique not\'e $\hat{B}$. Alors on a l'\'equivalence:}\\
$$B \ est\ plat\ sur\ \mathcal{V} \Longleftrightarrow \hat{B}\ est\ plat\ sur\ \mathcal{V} \ . $$\\

Par fid\`ele platitude des morphismes $\hat{S}\rightarrow S^{\dag}$ et $S^{\dag}\rightarrow \tilde{S}$ les propri\'et\'es $(ii),(iii), (iv)$ sont \'equivalentes.\\ 

Supposons $(ii)$ v\'erifi\'e: alors $X_{\hat{S}}$ est plat sur $\mathcal{V}$, donc par le lemme (3.3.2.2) $X$ est plat sur $\mathcal{V}$ et $(i)$ est clair. Il nous reste \`a prouver que $(i)\Rightarrow (ii)$. \\
Supposons la propri\'et\'e $(i)$ v\'erifi\'ee. Puisque $X_{\hat{S}}$ est plat sur $\mathcal{V}$ (car $X$ est plat sur $\mathcal{V}$) et $f$ est plat, le crit\`ere de platitude par fibres [EGA\ IV,(11.3.10)] prouve que $h_{\hat{S}}$ est plat en tous les points au-dessus de $S_{0}$. Or $Spec\ \hat{A}$ est connexe car $S_{0} = Spec\ A_{0}$ l'est [Et 4, cor 2 du th\'eo 3] ; comme $A$ est normal et que le morphisme $A \rightarrow \hat{A}$ est normal, car r\'egulier [EGA IV, (7.8.3) (v)], il r\'esulte de la [prop (1.7) 4 (ii)] que $\hat{A}$ est normal, donc que $\hat{A}$ est int\'egralement clos. Par le th\'eor\`eme de platitude g\'en\'erique [EGA IV, (6.9.1)] il existe un ouvert non vide $V$ de $Spec\ \hat{A}$ tel que la restriction $h_{V} :h^{-1}_{\hat{S}}(V) \rightarrow V$ de $h_{\hat{S}}$ soit plate. Or $V$ contient $Spec\ A_{0}$ d'apr\`es le raisonnement fait ci-dessus, donc $V = Spec\ \hat{A} =: \hat{S}$ puisque  $(\hat{S}, S_{0})$ est un couple hens\'elien [EGA IV, (18.5.4.3)]. \\

\textit{Pour (3.3.3).} Par fid\`ele platitude des morphismes $\hat{S}\rightarrow S^{\dag}$ et $S^{\dag}\rightarrow \tilde{S}$ les propri\'et\'es $(ii),(iii), (iv)$ sont \'equivalentes. Comme $(ii)\Rightarrow (i)$ est clair, il nous reste \`a prouver que $(i)\Rightarrow (ii)$. \\
Supposons la propri\'et\'e $(i)$ v\'erifi\'ee. D'apr\`es (3.3.2) $h_{\hat{S}}$ est plat. Appliquons [EGA IV, (12.2.4)] au morphisme projectif et plat $h_{\hat{S}}$: l'ouvert $W$ des $s \in \hat{S}$ o\`u $h_{s} : (X_{\hat{S}})_{s}  \rightarrow s$ est lisse contient $Spec\  A_{0}$ puisque $f$ est lisse, donc $W = \hat{S}$ l\`a encore en utilisant le caract\`ere hens\'elien du couple $(\hat{S}, S_{0})$. Ainsi on a obtenu un rel\`evement projectif et lisse $h_{\hat{S}} : X_{\hat{S}} \rightarrow \hat{S}$ de $f$. \\

 \textit{Pour (3.3.4).} Il suffit de prendre le compl\'et\'e formel du carr\'e (3.3.1.1): on obtient un carr\'e cart\'esien de $\mathcal{V}$-sch\'emas formels\\
$$
\begin{array}{c}
\xymatrix{
\mathcal{X}\ \ar@{^{(}->}[r]^{j_{\overline{Y}}} \ar[d]_{\hat{h}} & \overline{\mathcal{X}} \ar[d]^{\hat{\overline{h}}} &\\
\mathcal{S}\   \ar@{^{(}->}[r]_{j_{\overline{S}}} &  \overline{\mathcal{S}} &  ,\ 
}
\end{array}
\leqno{(3.3.4.1)}
$$

\noindent dans lequel $\hat{h}$ et $\hat{\overline{h}}$ sont projectifs. Puisque $\hat{A}$ est un anneau de Zariski et que $\hat{h}$ est le compl\'et\'e de $h, h_{\hat{S}}$, il  r\'esulte de [Bour, AC\ III, \S5; \no4, prop 2, et \no2, th\'eo 1] que la platitude (resp. la lissit\'e) de $\hat{h}$ \'equivaut \`a celle de $h_{\hat{S}}$ ; d'o\`u l'\'equivalence des propri\'et\'es $(i)$ et $(ii)$ . $\square$\\

\noindent \textbf{Remarques (3.3.5)}. Supposons que $\mathcal{V}$ est un anneau de valuation discr\`ete complet, $I$ son id\'eal maximal, $k$ son corps r\'esiduel et $\pi$ une uniformisante.\\
(i) L'hypoth\`ese de platitude de $X$ sur $\mathcal{V}$ dans (3.3.2), (3.3.3) et (3.3.4) \'equivaut \`a $\mathcal{O}_{X}$ sans $\pi$-torsion.\\
(ii) Lorsque $S= Spec\ \mathcal{V}$, un exemple de Serre [S 1] prouve qu'il existe des cas o\`u $X$ n'est pas plat sur $\mathcal{V}$: dans ce cas $h$ n'est pas plat. Nous allons voir ci-dessous en (3.3.7) une condition suffisante de platitude de $h$, celle o\`u le morphisme projectif $f$ identifie $X$ \`a une intersection compl\`ete dans un espace projectif (cf d\'efinition (3.3.6)).\\

En nous inspirant de la d\'efinition donn\'ee par Deligne des intersections compl\`etes dans un fibr\'e projectif [SGA 7, II, exp.XI, 1.4] on adoptera la d\'efinition suivante:\\

\noindent \textbf{D\'efinition (3.3.6)}.
\textit{Soit $f:X\rightarrow S$ un morphisme de sch\'emas. On dira que $X$ est une intersection compl\`ete relativement \`a $S$ dans des espaces projectifs sur $S$ si il existe un recouvrement de $S$ par des ouverts de Zariski $S_{\alpha},\ S= \underset{\alpha}{\bigcup} \ S_{\alpha}$ et, en d\'esignant par $f_{\alpha}:X_{\alpha}\rightarrow S_{\alpha}$ le morphisme d\'eduit de $f$ par le changement de base $S_{\alpha}\rightarrow S$, il existe, pour chaque $\alpha$, un couple $(n_{\alpha}, r_{\alpha})\in \mathbb{N}^{\ast}\times \mathbb{N}$ et une $S_{\alpha}$-immersion ferm\'ee $i_{\alpha}$ qui factorise $f_{\alpha}$ }
$$
\begin{array}{c}
\xymatrix{
X_{\alpha}\ \ar@{^{(}->}[r]^{i_{\alpha}} \ar[rd]_{f_{\alpha}} & \mathbb{P}^{n_{\alpha}}_{S_{\alpha}}\ar[d]^{p_{\alpha}}  \\
&S_{\alpha}\ 
}
\end{array}
\leqno{(3.3.6.1)}
$$
\textit{o\`u $p_{\alpha}$ d\'esigne la projection canonique, telle que:\\
\begin{enumerate}
\item[(i)] fibre \`a fibre, $i_{\alpha}$ est de codimension $r_{\alpha}$, i.e.  pour tout $s= Spec\ k(s)\in S_{\alpha}$, si $i_{\alpha, s}: \ X_{\alpha, s}\hookrightarrow \mathbb{P}^{n_{\alpha}}_{s}$ d\'esigne la fibre de $i_{\alpha}$ au-dessus de $s$, on a $r_{\alpha}$= codim ($X_{\alpha, s}, \mathbb{P}^{n_{\alpha}}_{s}$) 
\item[(ii)] il existe un recouvrement de $S_{\alpha}$ par des ouverts de Zariski $S_{\alpha, \beta},\ S_{\alpha}= \underset{\beta}{\bigcup} \ S_{\alpha, \beta}$ tels que, en d\'esignant par 
$$
\begin{array}{c}
\xymatrix{
X_{\alpha, \beta}\ \ar@{^{(}->}[r]^{i_{\alpha, \beta}} \ar[rd]_{f_{\alpha, \beta}} & \mathbb{P}^{n_{\alpha}}_{S_{\alpha, \beta}}\ar[d]^{p_{\alpha, \beta}}  \\
&S_{\alpha, \beta}\ 
}
\end{array}
\leqno{(3.3.6.2)}
$$
 le diagramme d\'eduit de (3.3.6.1) par le changement de base $S_{\alpha, \beta}\rightarrow S_{\alpha}$, chaque $X_{\alpha, \beta}=Proj(\mathcal{O}_{S_{\alpha, \beta}}[x_{0},x_{1},...,x_{n_{\alpha}}]/ \mathcal{J}_{\alpha, \beta})$ est l'intersection de $r_{\alpha}$ hypersurfaces (de certains degr\'es) au sens sch\'ematique, i.e.  l'id\'eal $\mathcal{J}_{\alpha, \beta}$ est engendr\'e par $r_{\alpha}$ \'el\'ements homog\`enes. 
\end{enumerate}
On dit alors que $X_{\alpha}$ est une intersection compl\`ete relativement \`a $S_{\alpha}$ dans $\mathbb{P}^{n_{\alpha}}_{S_{\alpha}}$ de codimension $r_{\alpha}$.
}\\

Nous allons \'etablir le corollaire suivant du th\'eor\`eme (3.3):\\

\noindent \textbf{Corollaire (3.3.7)}.
\textit{Soient $\mathcal{V}$ un anneau local normal d'id\'eal maximal $I \subset \mathcal{V}$, complet pour la topologie $I$-adique, de corps r\'esiduel $k = \mathcal{V}/I$. Soient $S_{0}$ un $k$-sch\'ema lisse et s\'epar\'e, $f : X_{0} \rightarrow S_{0}$ un $k$-morphisme projectif et lisse tel que $X_{0}$ est une intersection compl\`ete, relativement \`a $S_{0}$, dans des espaces projectifs sur $S_{0}$, $S_{\alpha,\beta} = Spec\ A_{0}$ un $k$-sch\'ema affine et lisse tel qu'en (3.3.6.2), $A = \mathcal{V} [t_{1},...,t_{d}] / J$  une $\mathcal{V}$-alg\`ebre lisse relevant $A_{0}$, dont on a fix\'e une pr\'esentation, et $S = Spec\ A$: on note $\hat{A}$ le s\'epar\'e compl\'et\'e de $A$ pour la topologie $I$-adique, $A^{\dag}$ son compl\'et\'e faible, $\tilde{A}$ l'hens\'elis\'e de $A$ au sens de Raynaud et $\hat{S}=Spec \ \hat{A},\ S^{\dag}=Spec\ A^{\dag},\  \tilde{S}= Spec\ \tilde{A}$. On d\'esigne par $\overline{S}$ l'adh\'erence sch\'ematique de $S$ dans $\mathbb{P}^d_{\mathcal{V}}$, et par $\mathcal{S}$ (resp. $\overline{\mathcal{S}})$ le compl\'et\'e formel $I$-adique de $S$ (resp. de $\overline{S}$).}\\
\noindent \textit{ Alors le $X$ construit en (3.3.1) est plat sur $\mathcal{V}$ et les morphismes $h_{\tilde{S}}, h_{S^{\dag}}, h_{\hat{S}}$ de (3.3.2) sont des rel\`evements projectifs et lisses du morphisme $f_{\alpha,\beta}: X_{\alpha,\beta}\rightarrow S_{\alpha,\beta}$ d\'eduit de $f$ par le changement de base $S_{\alpha,\beta}\rightarrow S_{0}$. De plus on dispose d'un diagramme tel que (3.3.4.1) dans lequel $\hat{h}$ est un rel\`evement projectif et lisse de $f$.}\\

\noindent \textit{D\'emonstration}. L'anneau $\mathcal{V}$ est excellent [EGA IV, (7.8.3)(iii)]. Quitte \`a faire le changement de base $S_{\alpha,\beta}\rightarrow S_{0}$ on supposera dans la suite de la d\'emonstration que $f_{\alpha,\beta}= f$. Quitte \`a d\'ecomposer les sch\'emas normaux noeth\'eriens $Spec\ \mathcal{V}$, $S_{0}$ et $X_{0}$ en somme de leurs composantes connexes on supposera dans toute la suite que  $Spec\ \mathcal{V}$, $S_{0}$ et $X_{0}$ sont connexes, donc int\'egralement clos. En utilisant alors les notations de la preuve de (3.3), avec $X_{0}= Proj(A_{0}[x_{0},x_{1},..., x_{n}]/\mathcal{J}^{0})$, on dispose du diagramme commutatif
$$
\begin{array}{c}
\xymatrix{
X_{0}\ \ar@{^{(}->}[r]^{i_{0}} \ar[rd]_{f} & \mathbb{P}^{n}_{S_{0}}\ar[d]^{p_{S_{0}}} & \\
&S_{0}& . 
}
\end{array}
$$
Comme $f$ est plat et $S_{0}$ connexe, pour tout $s \in S_{0}$, la dimension de $X_{0,s}= f^{-1}(s)$ est constante [H 2, III, cor 9.10], donc la codimension Codim ($X_{0, s}, \mathbb{P}^{n}_{k(s)}$) aussi: notons $r$ cette derni\`ere. Puisque $X_{0}$ est une intersection compl\`ete, relativement \`a $S_{0}$, dans $\mathbb{P}^{n}_{S_{0}}$, l'id\'eal $\mathcal{J}^{0}$ poss\`ede une famille g\'en\'eratrice constitu\'ee de $r$ \'el\'ements homog\`enes: notons ceux-ci $f^1_{0},...,f^r_{0}$, que l'on rel\`eve comme dans la preuve de (3.3) en $f^1,...,f^r$; d'o\`u un diagramme tel que (3.3.1.1). Notons $\mathbb{A}^{n}_{S, (j)}, \ j \in \llbracket0,n\rrbracket$, chacun des espaces affines qui recouvrent $\mathbb{P}^{n}_{S}$ et $X_{(j)}= Spec\ (A[x_{0},..., x_{n}]/(x_{j}-1, f^1,...,f^r))$ le sch\'ema induit sur $\mathbb{A}^{n}_{S, (j)}$ par $X$. Il s'agit de montrer que, pour tout $j$, $X_{(j)}$ est plat sur $S$. Pour $ j \in \llbracket0,n \rrbracket$, notons $X_{0,(j)}$ (resp $\mathbb{A}^{n}_{S_{0}, (j)}$) la r\'eduction de $X_{(j)}$ (resp $\mathbb{A}^{n}_{S, (j)}$) mod $I$ et $f_{0,(j)}^{\ell}(x_{0},..., x_{j-1}, x_{j+1},..., x_{n})= f_{0}^{\ell}(x_{0},..., x_{j-1}, x_{j}=1, x_{j+1},..., x_{n}), \ell \in \llbracket1,r \rrbracket$. Pour prouver la platitude de $X_{(j)}$ sur $S$, il suffit d'apr\`es [Mi, I, Rk 2.6 (d)] de prouver que $(f_{0,(j)}^{1},..., f_{0,(j)}^{r})$ est une suite r\'eguli\`ere, ce qui \'equivaut ici [$EGA \ 0_{IV}$, (15.2.2) et (15.2.3)] \`a prouver qu'elle est quasi-r\'eguli\`ere. Or $i_{0}$ est une immersion ferm\'ee r\'eguli\`ere d'apr\`es [EGA IV, (17.12.1)], i.e., pour tout $j \in \llbracket0,n \rrbracket$, l'id\'eal $\mathcal{J}^{0}_{(j)}:=(f_{0,(j)}^{1},..., f_{0,(j)}^{r})$ est r\'egulier [EGA IV, (16.9.2)]: donc, d'apr\`es [$EGA \ 0_{IV}$, (15.2.2)] et [EGA IV, (16.1.2.2)], l'homomorphisme $(ii)$ de [EGA IV, (16.9.3)] est bijectif. Comme $i_{0}$ est r\'eguli\`ere, donc quasi-r\'eguli\`ere, le faisceau conormal $\mathcal{N}_{X_{0}/\mathbb{P}^{n}_{S_{0}}}$ est localement libre [EGA IV, (16.9.8)], i.e. pour tout $j \in \llbracket0,n \rrbracket$ $\mathcal{J}^{0}_{(j)}/(\mathcal{J}^{0}_{(j)})^{2}$ est localement libre: or ici, fibre \`a fibre au-dessus de chaque point $s \in S_{0}$, il est de rang \'egal \`a rg ($\mathcal{J}^{0}_{(j)}/(\mathcal{J}^{0}_{(j)})^{2}$)= Codim ($X_{0,s}, \mathbb{P}^{n}_{k(s)}$)= $r$. Donc, pour tout $j \in \llbracket0,n \rrbracket$, $\mathcal{J}^{0}_{(j)}/(\mathcal{J}^{0}_{(j)})^{2}$ est localement libre de rang $r$. Or les images de $f_{0,(j)}^{1},..., f_{0,(j)}^{r}$ dans $\mathcal{J}^{0}_{(j)}/(\mathcal{J}^{0}_{(j)})^{2}$ l'engendrent en tant que $\mathcal{O}_{X_{0,(j)}}/\mathcal{J}^{0}_{(j)}$-module: \'etant au nombre de $r$ c'en est une base [Bour, AC II, \S3, cor 5 du th\'eo 1]; ainsi [EGA IV, (16.9.3)] la suite $(f_{0,(j)}^{1},..., f_{0,(j)}^{r})$ est quasi-r\'eguli\`ere. Ceci ach\`eve la preuve du corollaire. $\square$ \\

Dans le th\'eor\`eme qui suit on particularise les hypoth\`eses faites sur l'anneau $\mathcal{V}$ dans le th\'eor\`eme (3.1) ce qui permet d'\'etendre celui-ci du cas ``fini \'etale'' au cas ``fini":

\vskip 3mm
\noindent \textbf{Th\'eor\`eme (3.4)}.
\textit{Soient $\mathcal{V}$ un anneau excellent normal de caract\'eristique z\'ero, $I \subset \mathcal{V}$ un id\'eal et $\mathcal{V}_{0} = \mathcal{V}/I$ tel que $(\mathcal{V, V}_{0})$ soit un couple hens\'elien au sens de} [EGA IV, (18.5.5)] \textit{. Soient $A_{0}$ et $C_{0}$ deux $\mathcal{V}_{0}$-alg\`ebres lisses et}

$$\varphi_{0} : A_{0} \rightarrow C_{0} $$
\textit{un $\mathcal{V}_{0}$-morphisme fini (resp. fini et plat ; resp. fini et fid\`element plat ; resp. fini \'etale) ; fixons deux $\mathcal{V}$-alg\`ebres lisses $A$ et $C$ relevant respectivement $A_{0}$ et $C_{0}$ et notons $\hat{A},\  \hat{C}$ leurs s\'epar\'es compl\'et\'es $I$-adiques.}\\

\noindent \textit{Alors}

\begin{enumerate}
\item[(1)] \textit{Il existe un $\mathcal{V}$-morphisme}

                       $$\varphi : \hat{A} \rightarrow \hat{C} $$
\textit{fini (resp. fini et plat ; resp. fini et fid\`element plat ; resp. fini \'etale) relevant $\varphi_{0}.$}

\item[(2)]  (i) \textit{Il existe une $\mathcal{V}$-alg\`ebre de type fini $B$ normale (resp. int\'egralement close si $A$ l'est) relevant $C_{0}$, un $\mathcal{V}$-morphisme fini} 

$$\psi : A \rightarrow B $$
\textit{relevant $\varphi_{0}$ et un $\mathcal{V}$-isomorphisme $\hat{C} \simeq \hat{B}$ s'ins\'erant dans un triangle commutatif}

$$\xymatrix{
\hat{A} \ar[r]^{\hat{\psi}} \ar[rd]_{\varphi} & \hat{B} \ar[d]^{\simeq} &\\
& \hat{C} &  .
}$$

(ii) \textit{De plus les $\mathcal{V}$-morphismes}

$$\psi_{T} : A_{T} \rightarrow B_{T}\  \mbox{et}\  \psi^{\dag} = \psi \otimes_{A} A^{\dag} : A^{\dag} \rightarrow B^{\dag} \simeq B \otimes_{A} A^{\dag}$$ 
\textit{sont finis (resp. finis et plats ; resp. finis et fid\`element plats ; resp. finis \'etales) et rel\`event $\varphi_{0}$ ; les morphismes $\widehat{\psi}_{T}$ et $\widehat{\psi}^{\dag}$ s'identifient \`a $\hat{\psi}$.}

\item[(3)] \textit{Fixons de plus une pr\'esentation de la $\mathcal{V}$-alg\`ebre $A$}

$$A = \mathcal{V} [t_{1},...,t_{n}] / (f_{1},...,f_{r})$$
\textit{et reprenons les notations du (3) du th\'eor\`eme (3.1). On a les m\^emes carr\'es cart\'esiens, mais o\`u cette fois $\hat{\psi}$ est un rel\`evement fini (resp. fini et plat ; resp. fini et fid\`element plat ; resp. fini \'etale) de $\varphi_{0}$ et o\`u $\theta, \theta', \hat{\theta}, \widehat{\theta'}$ sont finis. De plus $\hat{\theta}\  (resp.\  \widehat{\theta'})$ est plat si et seulement si sa r\'eduction modulo I est plate.}
\end{enumerate}

\vskip 3mm
\noindent \textit{D\'emonstration de (3.4)}.

\begin{enumerate}
\item[(1)]  L'existence de $A$ et $C$ r\'esulte du th\'eor\`eme 6 de Elkik [E$\ell$] et celle de $\varphi$ se d\'eduit de la lissit\'e formelle de $\hat{A}$ sur $\mathcal{V}$ :  donc $\varphi$ est un morphisme fini [Bour, AC III, $\S$ 2, n$\circ$ 11, prop 14].\\
On conclut comme dans le th\'eor\`eme 17 de [Et 4].

\item[(2)]  et (3) D\'ecomposant les sch\'emas normaux $Spec\  A$, $Spec\  C$ en somme de leurs composantes connexes, on peut supposer $A$ et $C$ int\'egralement clos : alors $\hat{A}$ et $\hat{C}$ sont aussi  int\'egralement clos [prop (1.7) (4) (iv)], car $A$ et $C$ sont excellents. Le th\'eor\`eme de platitude g\'en\'erique [EGA IV, (6.9.1)] prouve alors la surjectivit\'e du morphisme fini $Spec(\varphi) : Spec\ \hat{C} \rightarrow Spec\ \hat{A}$.\\

Soient $L$ le corps des fractions de $\hat{A}$ et $L_{1}$ celui de $\hat{C}$ : la suite de la d\'emonstration est alors identique  celle du th\'eor\`eme (3.1).\ $\square$
\end{enumerate}

\vskip 3mm
\noindent \textbf{Corollaire (3.5)}.
\textit{Avec $\mathcal{V}$ comme dans le th\'eor\`eme (3.4) fixons une $\mathcal{V}$-alg\`ebre lisse $A$. Notons $\mathcal{C}^\dag_{f}$ (resp. $\mathcal{C}^\dag_{fp}\  ;\  resp.  \mathcal{C}^\dag_{ffp}\  ;\  resp. \mathcal{C}^\dag_{f \acute {e} t}) $ la cat\'egorie des $A^\dag$-alg\`ebres finies $B$ (resp. finies et plates ; resp. finies et fid\`element plates ; resp. finies \'etales) telles que $B$ soit formellement lisse sur $\mathcal{V}$ pour la topologie $I$-adique : on notera $\hat{\mathcal{C}}$ (resp. $\tilde{\mathcal{C}}\  ;\  resp. \stackrel{\circ}{\mathcal{C}})$ les cat\'egories analogues obtenues en  rempla\c cant $A^\dag$ par $\hat{A}$ (resp. par $\tilde{A}$\  ;\  resp. par $A_{0}$) ; ici on a omis les indices ``$f$'', ``$fp$'' .... Alors}

\begin{enumerate}
\item[(1)] \textit{Le foncteur}

$$\mathcal{F} : \mathcal{B} \mapsto \mathcal{B}\  \otimes_{A^{\dag}} \hat{A} \  (resp.\  \mathcal{F} :  \mathcal{B} \mapsto \mathcal{B}\  \otimes_{\tilde{A}} A^{\dag}) $$
 \textit{est une \'equivalence de cat\'egories de la cat\'egorie  $\mathcal{C}^{\dag}_{f}$\ (resp. $\tilde{C_{f}})$  dans la cat\'egorie $\hat{\mathcal{C}_{f}}\  (resp.\  \mathcal{C}^{\dag}_{f})$ : on a les m\^emes r\'esultats avec les indices $fp, ffp$ et $f \acute {e} t$.} 
 
 \item[(2)] \textit{Un foncteur quasi-inverse de $\mathcal{F}$ est fourni par le foncteur $\mathcal{G}$ qui \`a une $\hat{A}$-alg\`ebre (resp. $A^{\dag}$-alg\`ebre) finie $\mathcal{D}$ associe la fermeture int\'egrale de $A^{\dag}$ (resp.\ $\tilde{A})$ dans $\mathcal{D}.$ }
 
 \item[(3)] (i) \textit{Le foncteur}
 $$\tilde{\mathcal{H}} (resp.\ \mathcal{H}^{\dag}\  ;\  resp.\ \mathcal{H}^{\wedge})\  \mathcal{B} \mapsto \mathcal{B}\ /\ I \mathcal{B}$$
 \textit{est un foncteur plein et essentiellement surjectif de la cat\'egorie $\tilde{\mathcal{C}_{f}}\ (resp.\ \hfill\break \mathcal{C}^{\dag}_{f}\  ; \ resp.\  \hat{\mathcal{C}}_{f})$ dans la cat\'egorie $\stackrel{\circ} {\mathcal{C}}_{f}$ : on a les m\^emes r\'esultats avec les indices $fp$ et $ffp.$\\
 (ii) Le foncteur}
 
 $$\mathcal{H}_{\acute {e} t} : \mathcal{B} \mapsto \mathcal{B}\ /\ I \mathcal{B}$$  
 \textit{est une \'equivalence de cat\'egories de la cat\'egorie $\tilde{\mathcal{C}}_{f \acute {e} t}$\  (resp. $\mathcal{C}^{\dag}_{f \acute {e} t}\ ; \ resp.\hfill\break  \hat{\mathcal{C}}_{f \acute {e} t})$ dans la cat\'egorie  $\stackrel{\circ} {\mathcal{C}}_{f \acute {e} t}.$}
 \end{enumerate}
 
 \vskip 3mm
\noindent \textit{D\'emonstration}.\\
Le (3) (ii) est ici pour m\'emoire car montr\'e dans [th\'eo 2.4].

\begin{enumerate}
\item[(1)] et (2). On fait la d\'emonstration pour l'indice ``$f$'', et $\mathcal{F} : \mathcal{C}^{\dag}_{f} \rightarrow \hat{\mathcal{C}}_{f}$, les autres \'etant analogues.\\

\qquad Le foncteur $\mathcal{F}$ est fid\`ele d'apr\`es [EGA IV (2.2.16)].\\

\qquad Prouvons  que $\mathcal{F}$ est essentiellement surjectif. Soit : $\varphi : \hat{A} \rightarrow C$ un objet de $\hat{\mathcal{C}_{f}}$ ; la r\'eduction mod $I$, $C_{0}$, de $C$ est une $\mathcal{V}_{0}$-alg\`ebre lisse que l'on rel\`eve par le th\'eor\`eme de Elkik en une $\mathcal{V}$-alg\`ebre lisse $D$ : par lissit\'e formelle de $C$ sur $\mathcal{V}$ il existe un $\mathcal{V}$-morphisme $\hat{\theta}_{1} : C \rightarrow \hat{D}$ et c'est un isomorphisme ; comme le diagramme suivant commute :

 $$\xymatrix{
& C \ar[dl]_{\hat{\theta}_{1}} \\
\hat{D} & \hat{A} \ar@{-->}[l] \ar[u]_{\varphi}\\
& \mathcal{V} \ar[lu] \ar[u]
}
$$

\noindent on peut voir $\hat{\theta}_{1}$ comme un $\hat{A}$-isomorphisme au moyen de la fl\`eche en pointill\'e $\hat{\theta}_{1} \circ \varphi$. Le th\'eor\`eme (3.4) fournit alors une $A$-alg\`ebre finie $B$, $\psi : A \rightarrow B$, et un $\hat{A}$-isomorphisme

$$\hat{\theta}_{2} : \hat{D} \tilde{\rightarrow} \hat{B}$$
s'ins\'erant dans le diagramme commutatif

$$\xymatrix{
\hat{A} \ar[r]^{\varphi} \ar[d]_{\hat{\psi}} \ar@{-->}[rd] & C \ar[d]^{\hat{\theta}_{1}}_{\simeq} \\
\hat{B} & \hat{D} \ar[l]^{\hat{\theta}_{2}}_{\simeq} 
}
$$

\noindent donc $\varphi\   \tilde{\rightarrow}\ \mathcal{F} (\psi^{\dag} : A^{\dag} \rightarrow B^{\dag})$ et $\mathcal{F}$ est essentiellement surjectif.\\

\quad Prouvons que $\mathcal{F}$ est plein. Soient $\varphi : \hat{A} \rightarrow C$ et $\varphi' : \hat{A} \rightarrow C'$ deux objets de $\hat{\mathcal{C}_{f}}$ et $\lambda : C \rightarrow C'$ un $\hat{A}$-morphisme. On vient de montrer qu'il existe des morphismes finis $\psi : A \rightarrow B, \psi' : A \rightarrow B'$ s'ins\'erant dans le diagramme commutatif

$$
\xymatrix{\hat{A} \ar@/^/[rrd]^{\varphi} \ar@/_/[rdd]_{\hat{\psi'}} \ar[rd]^{\hat{\psi}} \\
&  \hat{B} \ar[d]_{\lambda'} & C \ar[l]_{\simeq}^{\theta} \ar[d]^{\lambda}\\
& \hat{B'} & C' \ar[l]_{\simeq}^{\theta'}
}
$$

o\`u l'on a pos\'e $\lambda' = \theta' \circ \lambda \circ \theta^{-1}$. Il s'agit de trouver un $A^{\dag}$-morphisme $\lambda^{\dag} : B^{\dag} \rightarrow B'^{\dag}$  induisant $\lambda'$ par passage aux compl\'et\'es.\\
Comme dans le th\'eor\`eme (3.4) on se ram\`ene au cas o\`u $A, \hat{A}, C$ et $C'$ sont int\'egralement clos. On a vu dans la d\'emonstration du th\'eor\`eme (3.4) (semblable \`a celle de (3.1)) qu'il existe $a \in A$ et $f, g \in A_{a}[X]$ tels que $\psi_{a}$ et $\psi'_{a}$ soient finis \'etales et s'identifient aux morphismes canoniques :

$$\psi_{a} : A_{a} \rightarrow B_{a} = A_{a}[X]\ /\ (f) $$

$$\psi'_{a} : A_{a} \rightarrow B'_{a} = A_{a}[X]\ /\ (g) ;$$

$\hat{\psi}_{a}, \hat{\psi'_{a}}$ sont finis \'etales et s'identifient \`a

$$(\hat{\psi})_{a} : (\hat{A})_{a} \rightarrow (\hat{B})_{a} = (\hat{A})_{a} [X]\ /\ (f)$$

$$(\hat{\psi'})_{a} : (\hat{A})_{a} \rightarrow (\hat{B'})_{a} = (\hat{A})_{a} [X]\ /\ (g)$$
et font commuter le triangle dont les fl\`eches sont finies \'etales

$$\xymatrix{
\widehat{(A_{a})} =  \widehat{((\hat{A})_{a})} \ar[r]^{\widehat{((\hat{\psi})_{a})}}  \ar[rd]_{\widehat{((\hat{\psi'})_{a})}} & \widehat{((\hat{B})_{a})} =   \ar[d]^{\widehat{(\lambda'_{a})}}    \widehat{(B_{a})} &\\
& \widehat{((\hat{B'})_{a})} = \widehat{(B'_{a})} & .
}
 $$

Par l'\'equivalence de (3) (ii) il existe un morphisme $\rho$ fini \'etale, qui rel\`eve $(\widehat{\lambda'_{a}})$ et fait commuter le triangle

$$\xymatrix{
(A_{a})^\dag \ar[r]^{(\psi_{a})^\dag} \ar[rd]_{(\psi'_{a})^{\dag}} & (B_{a})^{\dag} \ar[d]^{\rho} &\\
& (B'_{a})^{\dag} & .
}
$$

D'apr\`es la proposition 2 de [Et 5] on a les \'egalit\'es\\
$$A^{\dag} = (A_{a})^{\dag}\ \cap \hat{A},\quad B^{\dag} = (B_{a})^{\dag}\ \cap \hat{B}\ \quad , \quad 
B'^{\dag} = (B'_{a})^{\dag} \cap \hat{B'}\  ;$$ 

\noindent comme $\rho$ et $\lambda'$ induisent tous deux $(\widehat{\lambda'_{a}})$, la restriction de $\rho$ et $\lambda'$ \`a $B^{\dag} = (B_{a})^{\dag} \cap \hat{B}$ fournit le $A^{\dag}$-morphisme cherch\'e

$$ \lambda^{\dag} : B^{\dag} \rightarrow B'^{\dag}\ .$$
\item[(2)] Se d\'emontre comme le (2) (ii) du th\'eor\`eme (2.4).
\item [(3)] (i) Il suffit de le prouver pour $\mathcal{H}^{\wedge}$.\\

\quad On prouve que $\mathcal{H}^{\wedge}$ est essentiellement surjectif par une m\'ethode analogue \`a celle utilis\'ee pour $\mathcal{F}$ dans le (1) : \'etant donn\'ee une $A_{0}$-alg\`ebre finie $B_{0}$, on trouve une $A$-alg\`ebre finie $B$ relevant $B_{0}$ et donc $\hat{B}$ est une $\hat{A}$-alg\`ebre finie r\'epondant \`a la question.\\

\quad Pour montrer que $\mathcal{H}^{\wedge}$ est plein on part d'un $A_{0}$-morphisme $\lambda_{0} : B_{0} \rightarrow B'_{0}$ ; avec $B$ et $B'$ comme ci-dessus et lissit\'e formelle de $\hat{B}$ sur $\mathcal{V}$ on en d\'eduit un $\mathcal{V}$-morphisme (en fait un $\hat{A}$-morphisme) $\lambda : \hat{B} \rightarrow \hat{B'}$ qui rel\`eve $\lambda_{0}$.\ $\square$
\end{enumerate}

 \vskip 3mm
\noindent \textbf{Corollaire (3.6)}.
\textit{Avec les hypoth\`eses et notations du (3) du th\'eor\`eme (3.4), supposons de plus donn\'e un carr\'e cart\'esien de $\mathcal{V}$-sch\'emas formels}

$$
\xymatrix{
\mathcal{X}\ \ar@{^{(}->}[r] \ar[d]_{h} & \overline{\mathcal{X}} \ar[d]^{\overline{h}}\\
\mathcal{S}\   \ar@{^{(}->}[r]_{j} &  \overline{\mathcal{S}}
}
$$

\noindent \textit{o\`u $\overline{h}$ est propre. Alors le $\mathcal{V}$-sch\'ema formel $\mathcal{X}'$ d\'efini par le produit fibr\'e}

$$
\xymatrix{
\mathcal{X'}\ \ar[r] \ar[d]_{h'} & \mathcal{X} \ar[d]^{h}\\
\mathcal{S'}\   \ar [r]_{\hat{\psi}} & \mathcal{S}
}
$$

\noindent \textit{admet une compactification $\mathcal{\overline{X'}}$ d\'efinie par le cube \`a faces cart\'esiennes}

$$
\xymatrix{
&\mathcal{X} \ar@{^{(}->}[rr] \ar@{.>}[dd]_(.4){h} |\hole && \overline{\mathcal{X}} \ar[dd]^{\overline{h}} \\
 \mathcal{X}'  \ar[rr] \ar[dd]_{h'} \ar[ur] && \overline{\mathcal{X}'} \ar[dd] \ar[ur]_{\theta'}\\
&\mathcal{S} \ar@{^{(}.>}[rr]^{j\quad } |\hole && \overline{\mathcal{S}}  \\
 \mathcal{S}'  \ar@{^{(}->}[rr]_{ j'} \ar@{.>}[ur]^{\hat{\psi}} && \overline{\mathcal{S}' } \ar[ur]_{\hat{\theta}} \\
}
$$

\noindent \textit{o\`u $\hat{\theta}$ est fini.\\
De plus on a les m\^emes r\'esultats en rempla\c cant $\overline{\mathcal{S}}$ par $\tilde{\mathcal{S}}$.}

\vskip 3mm
\noindent \textit{D\'emonstration}. Comme $j'$ est une immersion ouverte il en est de m\^eme de $\mathcal{X}' \rightarrow \overline{\mathcal{X'}}$. De plus $\overline{\mathcal{X'}} \rightarrow  \overline{\mathcal{X}}$ est fini puisque $\hat{\theta}$ l'est [th\'eo (3.4) (3)] ; d'o\`u le corollaire. $\square$

 \vskip 3mm
\noindent \textbf{Corollaire (3.7)}.
\textit{Avec les hypoth\`eses et notations du th\'eor\`eme (3.4), supposons de plus $\varphi_{0}$ fini \'etale galoisien du groupe $G$. Alors}

\begin{enumerate}
\item[(1)] \textit{ $\hat{\psi} : \mathcal{S}' = Spf \ \hat{B} \rightarrow \mathcal{S} = Spf\ \hat{A}$ est fini \'etale galoisien de groupe $G$.}
\item[(2)] \textit{On a un carr\'e cart\'esien}

$$
\xymatrix{
\mathcal{S}'_{0}\ \ar@{^{(}->}[r]^{j'_{0}} \ar[d]_{h_{0}} & \overline{S'_{0}} \ar[d]^{\overline{h}_{0}}\\
S_{0}\   \ar@{^{(}->}[r]_{j_{0}} &  \overline{\mathcal{S}}_{0}
}
$$

\textit{o\`u $h_{0} = Spec\ \varphi_{0} = \hat{\psi}\  mod\   I, \overline{S_{0}}$ et $\overline{S'_{0}}$ sont propres sur $\mathcal{V}_{0}$ et normaux, $\overline{h_{0}}$ est fini, et $j_{0}, j'_{0}$ sont des immersions ouvertes dominantes. De plus $G$ agit sur $\overline{S'_{0}}$ par $\overline{S}_{0}$-automorphismes et on a un isomorphisme}
 
 $$\mathcal{O}_{\overline{S}_{0}}\  \tilde{\longrightarrow}\  (\overline{h}_{0 \ast}(\mathcal{O}_{\overline{S'_{0}}} ))^G. $$
\end{enumerate}

\vskip 3mm
\noindent \textit{D\'emonstration}.\\
 Le (1) est classique.\\
Dans le (2) on prend pour $\overline{S}_{0}$ le normalis\'e de $\widehat{P'}$ mod $I$, et pour $\overline{S'_{0}}$ la fermeture int\'egrale de $\overline{S}_{0}$ dans $S'_{0}$ ; d'o\`u le carr\'e cart\'esien ci-dessus, et un carr\'e commutatif

$$
\xymatrix{
j_{0 \ast}  h_{0 \ast}  \mathcal{O}_{{S}'_{0}} =   \overline{h}_{0 \ast} j'_{0 \ast} \mathcal{O}_{{S}'_{0}}  &  \overline{h}_{0 \ast} \mathcal{O}_{\overline{S'_{0}}} \ar@{_{(}->}[l] \\
j_{0 \ast}\  \mathcal{O}_{s_{0}} \ar[u]^{\varphi_{0}}& \quad \mathcal{O}_{\overline{S_{0}}} \ar@{_{(}->}[l]^{j_{\overline{S}}}  \ar[u]
}
$$

\noindent \`a fl\`eches horizontales injectives.\\

Pour l'action de $G$ on peut supposer $S_{0}$ connexe, donc int\'egralement clos puisqu'il est lisse sur $\mathcal{V}_{0}$ : alors $\overline{S_{0}}$ est int\`egre et normal, donc $\mathcal{O}_{\overline{S}_{0}}$ est un faisceau d'anneaux int\'egralement clos, de corps des fractions celui de $\mathcal{O}_{S_{0}}$.\\

Consid\'erons $g \in G$ et une section $x$ de $\overline{h}_{0 \ast}  \mathcal{O}_{\overline{S'_{0}}}$ : alors $g(x)$ est une section de $j_{0 \ast}  h_{0 \ast}  \mathcal{O}_{{S}'_{0}}$ qui est enti\`ere sur $\mathcal{O}_{\overline{S}_{0}}$, donc $g(x)$ est en fait une section de $\overline{h}_{0 \ast}\  \mathcal{O}_{\overline{S'_{0}}}$ car $\overline{S'}_{0}$ est la fermeture int\'egrale de $\overline{S}_{0}$ dans $S'_{0}$.\\

Si l'on suppose de plus $x$ fixe par $G$, alors $x$ est une section de $j_{0 \ast}  \mathcal{O}_{{S}_{0}}$ car $\varphi_{0}$ est galoisien de groupe $G$ ; or $x$ est enti\`ere sur $\mathcal{O}_{\overline{S}_{0}}$, donc $x$ est une section de $\mathcal{O}_{\overline{S}_{0}}$ puisque $\mathcal{O}_{\overline{S}_{0}}$ est int\'egralement clos. D'o\`u le corollaire. $\square$ 


\cleardoublepage


\vskip 10mm
\chapter*{II.  Espaces rigides analytiques et images directes}
\markboth{\sc j.-y. etesse}{\sc II.  Espaces rigides analytiques et images directes}

\section*{0. Notations}
 Pour les notions sur les espaces rigides analytiques nous renvoyons le lecteur \`a [B 3] et [B-G-R].
Sauf mention contraire dans tout ce II on d\'esigne par $\mathcal{V}$ un anneau de valuation discr\`ete complet, de corps r\'esiduel $k = \mathcal{V}/\mathfrak{m}$ de caract\'eristique $p > 0$, de corps des fractions $K$ de caract\'eristique $0$, d'uniformisante $\pi$ et d'indice de ramification $e$.\\

On suppose donn\'e un entier $a \in \mathbb{N}^{\ast}$ et on d\'esigne par $C(k)$ un anneau de Cohen de $k$ de caract\'eristique $0$ [Bour, AC IX, $\S$ 2, $\no3$, prop 5] : $C(k)$ est un anneau de valuation discr\`ete complet d'id\'eal maximal $p\ C(k)$ [EGA $O_{IV}$, 19.8.5] et on note $K_{0}$ son corps des fractions, $K_{0} =$ Frac $(C(k))$. Il existe une injection fid\`element plate $C(k) \hookrightarrow \mathcal{V}$ qui fait de $\mathcal{V}$ un $C(k)$-module libre de rang $e$ [EGA $O_{IV}$, 19.8.6, 19.8.8] et [Bour, AC IX, $\S$ 2, $\no1$, prop 2]. On fixe un rel\`evement $\sigma : \mathcal{V} \rightarrow \mathcal{V}$ de la puissance $a^{i \grave{e}me}$ du Frobenius absolu de $k$, tel que $\sigma(\pi) = \pi$ comme dans  [Et 5, I, 1.1] (un tel $\sigma$ existe d'apr\`es loc. cit.); on notera encore $\sigma$ l'extension naturelle de $\sigma$ \`a $K$: lorsque $k$ est parfait $C(k)$ est isomorphe \`a l'anneau $W(k)$ des vecteurs de Witt de $k$ et $\sigma$ est un automorphisme de $K$. Si $k \hookrightarrow k'$ est une extension de corps de caract\'eristique $p > 0$, $\mathcal{V'} : = \mathcal{V} \otimes_{C(k)} C(k'), K' = \operatorname{Frac} (\mathcal{V'})$, on peut relever la puissance $a^{i \grave{e}me}$ du Frobenius absolu de $k'$ en un morphisme $\sigma' : K' \rightarrow K'$ au-dessus de $\sigma : K \rightarrow K$ [Et 5, I, 1.1].\\
\newpage
\section*{1. Changement de base pour un morphisme propre}

 Le th\'eor\`eme suivant est un pr\'ealable pour les th\'eor\`emes de changement de base en g\'eom\'etrie rigide.
 \vskip 3mm
\noindent \textbf{Th\'eor\`eme (1.1)}.
\textit{Soient $\mathcal{V}$ un anneau noeth\'erien et $I \not\subseteq \mathcal{V}$ un id\'eal ; on suppose $\mathcal{V}$ s\'epar\'e et complet pour la topologie $I$-adique.
Soit }

$$
\xymatrix{
\mathcal{X'}\ \ar[r]^{v} \ar[d]_{g} &\mathcal{X} \ar[d]^{f}\\
\mathcal{S'}\   \ar[r]^{u}&  \mathcal{S}
}
$$

\noindent \textit{un carr\'e cart\'esien de $\mathcal{V}$-sch\'emas formels (pour la topologie $I$-adique) de type fini, avec $f$ propre. }

\begin{enumerate}
\item[(1)] \textit{Soit $\mathcal{F}$ un $\mathcal{O}_{\mathcal{X}}$-module coh\'erent. Alors, pour tout entier $i \geqslant 0$, $R^i f_{\ast}(\mathcal{F})$ est un $\mathcal{O}_{\mathcal{S}}$-module coh\'erent et le morphisme de changement de base
$$u^{\ast}  R^{i}  f_{\ast}(\mathcal{F}) \longrightarrow R^{i}  g_{\ast} v^{\ast} (\mathcal{F})$$
est un isomorphisme.}

\item[(2)] \textit{Soient $\mathcal{E}^{\bullet}_{\mathcal{X}}$ un complexe born\'e de $\mathcal{O}_{\mathcal{S}}$-modules plats, \`a composantes des $\mathcal{O}_{\mathcal{X}}$-modules coh\'erents et  $\mathcal{E}^{\bullet}_{\mathcal{X'}} = \mathcal{E}^{\bullet}_{\mathcal{X}} \otimes_{\mathcal{O}_{\mathcal{S}}} \mathcal{O}_{\mathcal{S}'}$. On suppose $\mathcal{S}$ plat sur $\mathcal{V}$ et $u$ plat. Alors, pour tout entier $i \geqslant 0$, $R^{i} f_{\ast} \mathcal{(E}^{\bullet}_{\mathcal{X}})$ est un $\mathcal{O}_{\mathcal{S}}$-module coh\'erent et le morphisme de changement de base
$$u^{\ast}  R^{i}  f_{\ast}(\mathcal{E}^{\bullet}_{\mathcal{X}})\  \tilde{\longrightarrow}\  R^{i}  g_{\ast}  (\mathcal{E}^{\bullet}_{\mathcal{X'}}) $$
est un isomorphisme, qui s'interpr\`ete aussi comme un isomorphisme
$$\displaystyle \mathop{\lim}_{\leftarrow\atop{n}} u^{\ast}_{n}  R^{i} f_{n \ast}  (\mathcal{E}^{\bullet}_{X_{n}}) \tilde{\longrightarrow}\  \displaystyle \mathop{\lim}_{\leftarrow\atop{n}} R^{i} g_{n \ast} (\mathcal{E}^{\bullet}_{X'_{n}})\ ;$$
(cf. notations plus bas).}
\end{enumerate}

\vskip 3mm
\noindent \textit{D\'emonstration}.
\begin{enumerate}
\item[(1)] La coh\'erence de $R^{i}  f_{\ast}(\mathcal{F})$ sur $\mathcal{O}_{\mathcal{S}}$ est rappel\'ee pour m\'emoire [EGA III, (3.4.2)]. Notons $\mathcal{F}_{n} = \mathcal{F} / I^{n+1} \mathcal{F}$  et $\varphi_{n}$ le morphisme canonique $\varphi_{n} : R^{i} f_{\ast}(\mathcal{F}) \rightarrow R^{i} f_{\ast}(\mathcal{F}_{n})$.  Posons $C_{n} = \mbox{Coker}\ \varphi_{n},\  \mathcal{H} = R^{i} f_{\ast}(\mathcal{F}), \mathcal{H}_{n} = R^{i} f_{\ast}(\mathcal{F}_{n})$ et $\mathcal{H}'_{n} = \mathcal{H} / I^{n+1} \mathcal{H}$. Comme $I^{n+1} \mathcal{H} \subset \mbox{Ker}\  \varphi_{n}$ on a une surjection

$$\mathcal{H}'_{n}\    \twoheadrightarrow \mathcal{H} / \mbox{Ker}\  \varphi_{n}$$

de noyau not\'e $K_{n}$. Dans les suites exactes\\

(1.1.1) $\qquad 0 \longrightarrow K_{n} \longrightarrow \mathcal{H}'_{n} \longrightarrow \mathcal{H} / \mbox{Ker}\ \varphi_{n} \longrightarrow 0$

(1.1.2) $\qquad 0 \longrightarrow \mathcal{H} / \mbox{Ker}\ \varphi_{n} \longrightarrow \mathcal{H}_{n}
\longrightarrow C_{n} \longrightarrow 0$,\\

\noindent les syst\`emes projectifs $(\mathcal{H}'_{n})_{n}$ et $(\mathcal{H} / \mbox{Ker}\ \varphi_{n})_{n}$ v\'erifient la condition de Mittag-Leffler (not\'ee M-L) car les fl\`eches de transition sont surjectives. De plus $(\mathcal{H}_{n})_{n}$ v\'erifie M-L d'apr\`es [EGA III, (3.4.3)], donc $(C_{n})_{n}$ aussi [EGA ${O_{III}}$, (13.2.1)]. D'o\`u l'exactitude de la suite\\

(1.1.3) $\qquad  0 \longrightarrow \displaystyle \mathop{\lim}_{\leftarrow\atop{n}}\ \mathcal{H} / \mbox{Ker}\ \varphi_{n} \longrightarrow \displaystyle \mathop{\lim}_{\leftarrow\atop{n}}\ \mathcal{H}_{n}\
 \longrightarrow \displaystyle \mathop{\lim}_{\leftarrow\atop{n}}\  C_{n} \longrightarrow 0\ .$
 
 D'autre part on a un isomorphisme [EGA III, (3.4.3)]
 
 $$R^{i} f_{\ast}(\mathcal{F}) \tilde{\longrightarrow}  \displaystyle \mathop{\lim}_{\leftarrow\atop{n}}\  \mathcal{H}_{n}\  ;$$
 comme $R^{i} f_{\ast}(\mathcal{F})$ est un $\mathcal{O}_{\mathcal{S}}$-module coh\'erent [EGA III, (3.4.2)], il est s\'epar\'e et complet pour la topologie $I$-adique, donc on a aussi un isomorphisme
 
 $$R^{i} f_{\ast}(\mathcal{F}) \tilde{\longrightarrow}  \displaystyle \mathop{\lim}_{\leftarrow\atop{n}}\  \mathcal{H}'_{n}. $$ 
 Ainsi le morphisme compos\'e
 
 $$\displaystyle \mathop{\lim}_{\leftarrow\atop{n}}\ \mathcal{H}'_{n} \longrightarrow \displaystyle \mathop{\lim}_{\leftarrow\atop{n}} \mathcal{H} /  \mbox{Ker}\ \varphi_{n} \hookrightarrow   \displaystyle \mathop{\lim}_{\leftarrow\atop{n}}\ \mathcal{H}_{n}$$
 est un isomorphisme, donc chacune des fl\`eches est un isomorphisme (car la seconde est injective par (1.1.3)). Par suite $\displaystyle \mathop{\lim}_{\leftarrow\atop{n}}\ C_{n} = 0$ ; comme $(C_{n})_{n}$ v\'erifie ML, on en d\'eduit que le pro-objet $\ll (C_{n})_{n} \gg$ associ\'e est le pro-objet nul [G 1, 195-03, \S 2], i.e. pour tout $n$, il existe $n' \geqslant n$ tel que $C_{n'} \rightarrow C_{n}$ soit la fl\`eche nulle. On a donc un isomorphisme de pro-objets
 $$\ll  \mathcal{H} / \mbox{Ker}\ \varphi_{n})_{n} \gg \tilde{\longrightarrow} \ll (\mathcal{H}_{n})_{n} \gg .$$
 De la m\^eme fa\c con $\ll (K_{n})_{n} \gg$ est le pro-objet nul, d'o\`u un isomorphisme de pro-objets
 
 $$\ll (\mathcal{H'}_{n})_{n} \gg \tilde{\longrightarrow} \ll  (\mathcal{H} /  \mbox{Ker}\ \varphi_{n})_{n} \gg,$$
 et par composition avec le pr\'ec\'edent, un isomorphisme de pro-objets\\
 
 (1.1.4) $\qquad  \ll (\mathcal{H'}_{n})_{n} \gg \tilde{\longrightarrow} \ll (\mathcal{H}_{n})_{n} \gg.$\\
 
 En notant $\mathcal{C}$ (resp. $\mathcal{C}'$\ ; resp. Pro $\mathcal{C}$) la cat\'egorie des $\mathcal{O}_{\mathcal{S}}$-modules (resp. des $\mathcal{O}_{\mathcal{S}'}$-modules ; resp. des pro-objets de $\mathcal{C})$ le foncteur $u^{\ast} : \mathcal{C} \rightarrow \mathcal{C}'$ s'\'etend en un foncteur $\tilde{{u}^{\ast}} : Pro\  \mathcal{C} \rightarrow \mathcal{C}'$. Comme $u^{\ast}(\mathcal{H})$ est un $\mathcal{O}_{\mathcal{S}'}$-module coh\'erent il est s\'epar\'e et complet pour la topologie $I$-adique, d'o\`u un isomorphisme
 
  $$u^{\ast}(\mathcal{H}) \tilde{\longrightarrow} \displaystyle \mathop{\lim}_{\leftarrow\atop{n}}\ u^{\ast}(\mathcal{H}'_{n})\  ;$$
 or $(\mathcal{H}_{n})_{n}$ et $(\mathcal{H'}_{n})_{n}$\ v\'erifient M.L, donc d'apr\`es [G 1, 195-05, fin du \S\  2] et (1.1.4) on a des isomorphismes

$$u^{\ast}(\mathcal{H})  \tilde{\longrightarrow}\  \displaystyle \mathop{\lim}_{\leftarrow\atop{n}}\  u^{\ast}(\mathcal{H}'_{n})$$
$\qquad  \qquad  \qquad \qquad  \qquad  \qquad \qquad \simeq \tilde{u^{\ast}} ( \ll (\mathcal{H'}_{n})_{n} \gg) $ [G1, 195]

$\qquad  \qquad  \qquad \qquad  \qquad  \qquad \qquad  \simeq \tilde{u^{\ast}} ( \ll (\mathcal{H}_{n})_{n} \gg) $ (1.1.4)

$\qquad  \qquad  \qquad \qquad  \qquad  \qquad \qquad  \simeq \displaystyle \mathop{\lim}_{\leftarrow\atop{n}}\ u^{\ast}  (\mathcal{H}_{n})$ [G 1, 195]

$\qquad  \qquad  \qquad \qquad  \qquad  \qquad \qquad  \simeq \displaystyle \mathop{\lim}_{\leftarrow\atop{n}}\ u^{\ast}_{n}  (\mathcal{H}_{n})$

\noindent o\`u l'on note \\

$$
\xymatrix{
X'_{n}\ \ar[r]^{v_{n}} \ar[d]_{g_{n}} & X_{n} \ar[d]^{f_{n}}\\
S'_{n}\   \ar[r]_{u_{n}} &  S_{n}
}
$$

\noindent le carr\'e cart\'esien d\'eduit de celui de la proposition par r\'eduction mod $I^{n+1}$.\\

Or le th\'eor\`eme de changement de base pour un morphisme propre [SGA 4, T3, XII, th\'eo 5.1] fournit un isomorphisme

$$u^{\ast}_{n}(\mathcal{H}_{n}) =  u^{\ast}_{n}(R^{i} f_{n \ast}(\mathcal{F}_{n})) \tilde{\longrightarrow}\ R^{i} g_{n \ast} (v^\ast_{n}(\mathcal{F}_{n})),$$
d'o\`u 

$$u^{\ast}(R^{i} f_{\ast}(\mathcal{F}))  \tilde{\longrightarrow}\  \displaystyle \mathop{\lim}_{\leftarrow\atop{n}}\ R^{i} g_{n \ast} (v^\ast_{n}(\mathcal{F}_{n})) = \displaystyle \mathop{\lim}_{\leftarrow\atop{n}}\ R^{i} g_{\ast} (v^\ast_{n}(\mathcal{F}_{n}))$$

$\qquad\qquad\ \  \quad \qquad \quad \tilde{\longrightarrow}\ R^{i} g_{\ast} (v^{\ast}(\mathcal{F})) $ [EGA III, (3.4.3)].\\

D'o\`u le (1) du th\' eor\`eme.\\

 \item[(2)] On note $\mathcal{V}_{n} = \mathcal{V}/I^{n+1}$, $\mathcal{E}^{\bullet}_{X_{n}} = \mathcal{E}^{\bullet}_{\mathcal{X}} / I^{n+1} \mathcal{E}^{\bullet}_{\mathcal{X}}$ et $\mathcal{E}^{\bullet}_{X'_{n}} = \mathcal{E}^{\bullet}_{\mathcal{X}'} / I^{n+1}\ \mathcal{E}^{\bullet}_{\mathcal{X}'}$.\\

 \noindent D'apr\`es [SGA 4, XVII, th\'eo 4.3.1] on a un isomorphisme dans la cat\'egorie d\'eriv\'ee des $\mathcal{V}_{n-1}$-modules sur $S'_{n-1}$ :

$$\mathcal{V}_{n-1}\  \displaystyle \mathop{\otimes}^{\mathbb{L}}{_{_{\mathcal{V}_{n}}}} \  \mathbb{R}  f_{{n}^{\ast}}  (\mathcal{E}^{\bullet}_{{X}_{n}})\  \tilde{\longrightarrow}\ \mathbb{R} f_{n-1^{\ast}} (\mathcal{V}_{n-1} \displaystyle \mathop{\otimes}^{\mathbb{L}}{_{\mathcal{V}_{n}}}\  \mathcal{E}^{\bullet}_{{X}_{n}}).$$

D'apr\`es les hypoth\`eses, $\mathcal{E}^{\bullet}_{{X}_{n}} $ est \`a composantes plates sur $\mathcal{V}_{n}\  (\mathcal{S}$ est plat sur $\mathcal{V}$) d'o\`u un isomorphisme

$$\mathcal{V}_{n-1}\  \displaystyle \mathop{\otimes}^{\mathbb{L}}{_{\mathcal{V}_{n}}}\ \mathcal{E}^{\bullet}_{{X}_{n}}\ \tilde{\longrightarrow}\ \mathcal{E}^{\bullet}_{{X}_{n-1}}.$$
Or on a un isomorphisme

$$\mathcal{V}_{n-1}  \displaystyle \mathop{\otimes}^{\mathbb{L}}{_{\mathcal{V}_{n}}}\ \mathbb{R} f_{{n}^{\ast}}  (\mathcal{E}^{\bullet}_{{X}_{n}})\ \overset{\sim} {\longrightarrow}\ \mathcal{O}_{{S}_{n-1}}\  \displaystyle \mathop{\otimes}^{\mathbb{L}}{_{\mathcal{O}_{{S}_{n}}}}\ \mathbb{R}  f_{{n}^{\ast}}\ (\mathcal{E}^{\bullet}_{{X}_{n}})\  ,$$
d'o\`u un isomorphisme\\

 (1.1.5) $\qquad \mathbb{R}  f_{{n}^{\ast}} (\mathcal{E}^{\bullet}_{{X}_{n}})\ \displaystyle \mathop{\otimes}^{\mathbb{L}}{_{\mathcal{O}_{{S}_{n}}}}\ \mathcal{O}_{{S}_{n-1}}\ \tilde{\longrightarrow}\ \mathbb{R}  f_{{n-1}^{\ast}} (\mathcal{E}^{\bullet}_{{X}_{n-1}}).$\\

 Comme $R^{i+j} f_{{n}^{\ast}}\ (\mathcal{E}^{\bullet}_{{X}_{n}})$ est l'aboutissement d'une suite spectrale de terme $E^{i,j}_{1}$ donn\'e par
 
 $$E^{i,j}_{1} = R^{j}  f_{{n}^{\ast}} (\mathcal{E}^i_{{X}_{n}})$$
 et que $E^{i,j}_{1}$ est coh\'erent sur $\mathcal{O}_{{S}_{n}}$ puisque $f_{n}$ est propre, on en d\'eduit que $R^{i+j} f_{{n}^{\ast}}\ (\mathcal{E}^{\bullet}_{{X}_{n}})\ $ est coh\'erent sur $\mathcal{O}_{{S}_{n}}$, donc gr\^ace \`a (1.1.5) et [B-O, page B-4], que $(R\  f_{{n}^{\ast}}\ (\mathcal{E}^{\bullet}_{{X}_{n}}))_{n}$ est un objet ``consistant" au sens de loc. cit : par suite [B-O, page B-7] on a des isomorphismes canoniques de $\mathcal{O}_{\mathcal{S}}$-modules coh\'erents\\
 
 (1.1.6) $\qquad H^{i} (\mathbb{R}  \displaystyle \mathop{\lim}_{\leftarrow\atop{n}}\ \mathbb{R}f_{{n}^{\ast}} (\mathcal{E}^{\bullet}_{{X}_{n}}))\  \tilde{\longrightarrow} \displaystyle \mathop{\lim}_{\leftarrow\atop{n}}\ R^{i} f_{{n}^{\ast}} (\mathcal{E}^{\bullet}_{{X}_{n}}).$\\
 
 Dans la d\'emonstration du (1) on a rappel\'e l'isomorphisme (pour tout $\mathcal{O}_{\mathcal{X}}$-module coh\'erent $\mathcal{F})$\\
 
 $$R^{i} f_{\ast} (\mathcal{F})\  \tilde{\longrightarrow}\ \displaystyle \mathop{\lim}_{\leftarrow\atop{n}}\ R^{i} f_{{n}^{\ast}}(\mathcal{F}_{n})$$
 
 $\qquad  \qquad  \qquad \qquad  \qquad  \qquad \quad  \simeq \quad H^{i} (\mathbb{R} \displaystyle \mathop{\lim}_{\leftarrow\atop{n}}\ \mathbb{R}  f_{{n}^{\ast}}(\mathcal{F}_{n}))$ \\
 
 [B-O, page B-7], \\
 d'o\`u un isomorphisme\\
 
 (1.1.7) $\qquad \mathbb{R} f_{\ast} (\mathcal{E}^{\bullet}_{{\mathcal{X}}})\   \tilde{\longrightarrow}\ \mathbb{R} \displaystyle \mathop{\lim}_{\leftarrow\atop{n}}\ \mathbb{R} f_{{n}^{\ast}}(\mathcal{E}^{\bullet}_{{\mathcal{X}}_{n}}),$\\
 i.e. compte tenu de (1.1.6)\\
 
 (1.1.8) $\qquad R^{i} f_{\ast} (\mathcal{E}^{\bullet}_{{\mathcal{X}}})\   \tilde{\longrightarrow}\  \displaystyle \mathop{\lim}_{\leftarrow\atop{n}}\ R^{i} f_{{n}^{\ast}}(\mathcal{E}^{\bullet}_{{X}_{n}})\ ,$\\ 
 et c'est un $\mathcal{O}_{\mathcal{S}}$-module coh\'erent.\\
 
 \noindent Puisque $u_{n}$ est plat, le th\'eor\`eme de changement de base [SGA 4, XVII, th\'eo (4.3.1)] fournit un isomorphisme
 
 $$u^{\ast}_{n}  \mathbb{R}  f_{{n}^{\ast}}(\mathcal{E}^{\bullet}_{{X}_{n}})  \tilde{\longrightarrow}\ \mathbb{R} g_{{n}^{\ast}} (\mathcal{E}^{\bullet}_{{X}'_{n}}),$$
 et par application de $\mathbb{R} \displaystyle \mathop{\lim}_{\leftarrow\atop{n}}\ $ un isomorphisme\\
 
  (1.1.9) $\qquad u ^{\ast} \mathbb{R} f_{\ast} (\mathcal{E}^{\bullet}_{{\mathcal{X}}})\   \tilde{\longrightarrow}\   \mathbb{R} g_{\ast} (\mathcal{E}^{\bullet}_{{\mathcal{X}'}}).$\\
  
  \noindent En passant \`a la cohomologie on obtient les isomorphismes annonc\'es au (2) du th\'eor\`eme (1.1).\  $\square$
 \end{enumerate}

\vskip 3mm
\noindent \textbf{Th\'eor\`eme (1.2)}.
\textit{Soient $\mathcal{V}$ un anneau de valuation discr\`ete s\'epar\'e et complet pour la topologie $\mathfrak{m}$-adique, o\`u $\mathfrak{m}$ est son id\'eal maximal, $k = \mathcal{V}/ \mathfrak{m}$ son corps r\'esiduel suppos\'e de caract\'eristique $p>0$, $K$ son corps des fractions de caract\'eristique $0$.}\\

\noindent \textit{Soit}

$$
\xymatrix{
X'\ \ar[r]^{v} \ar[d]_{g} & X \ar[d]^{f}\\
S'\   \ar[r]_{u} &  S
}
$$

\noindent \textit{un carr\'e cart\'esien de $K$-espaces analytiques rigides, avec $f$ propre.}\\

 \noindent (1.2.1) \textit{Soit $E$ un $\mathcal{O}_{X}$-module coh\'erent. Alors, pour tout entier $i \geqslant 0$},
 
 \begin{enumerate}
 \item[(1)] \textit{$R^{i} f_{\ast}(E)$ est un $\mathcal{O}_{S}$-module coh\'erent.}
\item[(2)] \textit{Le morphisme de changement de base}

$$u ^{\ast} (R^{i} f_{\ast} (E)) \longrightarrow R^{i} g_{\ast} (v^{\ast}(E))$$
\textit{est un isomorphisme.}
\end{enumerate}

\noindent (1.2.2)  \textit{Soient $E^{\bullet}$ un complexe born\'e de $\mathcal{O}_{S}$-modules plats, \`a composantes des $\mathcal{O}_{X}$-modules coh\'erents et $E'^{\bullet} = E^{\bullet}    \otimes_{\mathcal{O}_{\mathcal{S}}} \mathcal{O}_{S'} $. On suppose u plat. Alors pour tout entier $i \geqslant 0$}\\

\begin{enumerate}
\item[(1)] \textit{$R^i{} f_{\ast}(E^{\bullet})$ est un $\mathcal{O}_{S}$-module coh\'erent.}
\item[(2)] \textit{Le morphisme de changement de base
$$u^{\ast} R^{i} f_{\ast}(E^{\bullet})\ \longrightarrow R^{i} g_{\ast}(E'^{\bullet})$$
est un isomorphisme. }
\end{enumerate}

\vskip 10mm
\noindent \textit{D\'emonstration du th\'eor\`eme (1.2)}.\\

\noindent \textbf{1\`ere \'etape}. Supposons d'abord d\'emontr\'e le th\'eor\`eme dans le cas o\`u les $K$-espaces analytiques rigides sont tous quasi-compacts et quasi-s\'epar\'es.

Prouvons (1.2.1) dans ce cas. L'assertion (1) est locale sur $S$ puisque $R^{i} f_{\ast}(E)$ est le faisceau associ\'e au pr\'efaisceau

$$V \mapsto\ H^i(f^{-1} (V), E)$$

\noindent o\`u $V$ parcourt les ouverts de $S$ [SGA 4, V, prop 5.1].

Soient $\psi : V \hooklongrightarrow S$ un ouvert affino¬\"{\i}de de $S$ et $W$ d\'efini par le carr\'e cart\'esien

$$
\begin{array}{c}
\xymatrix{
W\ \ar[r]^{\theta} \ar[d]_{f'} & X \ar[d]^{f}  &\\ 
V\   \ar[r]_{\psi} &  S  & \quad;
}
\end{array}
\leqno(1.2.1.1)  
$$

\noindent d'apr\`es [loc. cit.] on a alors un isomorphisme canonique

$$\psi^\ast  R^{i}  f_{\ast}(E)\   \cong\ R^{i}  f'_{\ast}(\theta^{\ast}(E)).$$

\noindent Or ici $V$ est quasi-compact, quasi-s\'epar\'e, et $W$ aussi car $f'$ est propre ; d'o\`u l'assertion (1) via le cas quasi-compact, quasi-s\'epar\'e.\\

Pour le (2) on reprend le carr\'e cart\'esien (1.2.1.1) : soient 

$$u' : V' = V\times_{S}  S' \longrightarrow\ V $$
$$\mbox{et}\ u'' : V'' \hookrightarrow V'\  \mbox{un ouvert affino¬\"{\i}} \mbox{de de}\  V',$$

\noindent et on consid\`ere le diagramme commutatif suivant

$$
\xymatrix{& X'\ar@{.>}[dd]^(.6){g} |\hole\ar@{=} [rr] & & X' \ar@{.>}[dd]_(.6){g} |\hole \ar[rr]^{v}& &
X \ar[dd]^f \\
W'' \ar[rr]^(.8){v''} \ar[dd]_{g''} \ar[ur]^{\theta''} & & W' \ar[rr]^(.8){v'} \ar[dd]_(.7){g'} \ar[ur]^{\theta'}  & & W \ar[dd]_(.7){f'} \ar[ur]_{\theta} & \\
& S' \ar@{==}[rr]  |\hole & & S' \ar@{.>}[rr]^(.8){u} |\hole & &  S  \\
V'' \ar@{^{(}->}[rr]^{u''} \ar@{.>}[ur]^{\psi''} & & V' \ar[rr]^{u'} \ar@{.>}[ur]_{\psi'} & & V \ar[ur]_{\psi} 
}
$$

\vskip 2mm
\noindent dont les faces verticales sont cart\'esiennes.\\

Vu le caract\`ere local de l'isomorphisme (2) cherch\'e, il suffit de montrer que l'on a un isomorphisme

$$\psi''^{\ast} u^{\ast} R^{i} f_{\ast}(E) \displaystyle \mathop{\longrightarrow}^\sim \psi''^{\ast} R^{i} g_{\ast}\  v^{\ast}(E).$$

\noindent Or ici $V$ et $V''$ sont quasi-compacts, quasi-s\'epar\'es, donc $W$ et $W''$ aussi car $f'$ et $g''$ sont propres et on peut appliquer l'isomorphisme de changement de base du cas quasi-compact, quasi-s\'epar\'e pour le carr\'e cart\'esien

$$
\xymatrix{
W'\ \ar[r]^{v'\circ v''} \ar[d]_{g''} & W \ar[d]^{f'} &\\
V''\   \ar[r]_{u' \circ u''}  & V & \quad;
}
$$

\noindent d'o\`u une suite d'isomorphismes

$$\psi''^{\ast}  u^{\ast} R^{i} f_{\ast}(E) =  (u' u'')^{\ast} \psi^{\ast} R^{i} f_{\ast}(E)$$

$\qquad  \qquad  \qquad \qquad  \qquad  \qquad \qquad = (u' u'')^{\ast}  R^{i} f'_{\ast}   (\theta^{\ast}(E))$ 

$\qquad  \qquad  \qquad \qquad  \qquad  \qquad \qquad  \simeq R^{i} g''_{\ast} ((v' v'')^{\ast}(\theta^{\ast}(E)))$

$\qquad  \qquad  \qquad \qquad  \qquad  \qquad \qquad  = R^{i} g''_{\ast} (\theta''^{\ast}(v^{\ast}(E))) $

$\qquad  \qquad  \qquad \qquad  \qquad  \qquad \qquad  = \psi''^{\ast} R^{i} g_{\ast} (v^{\ast}(E)).$

\noindent D'o\`u le (2).

Pour prouver (1.2.2) \`a partir du cas quasi-compact, quasi-s\'epar\'e la d\'emarche est analogue.\\

\noindent \textbf{2\`eme \'etape}. Prouvons le th\'eor\`eme dans le cas quasi-compact, quasi-s\'epar\'e.

(1.2.1.) Le (1) est un th\'eor\`eme de L¬\"utkebohmert [L¬\"u, theo 2.7].\\
Pour le (2) on adopte les notations de [Bo-L¬\"u 1, d\'emonstration de 4.1] : d'apr\`es [loc. cit.] il existe des mod\`eles formels $\stackrel{\circ}{\mathcal{X}}$, $\mathcal{S}$, $\stackrel{\circ}{\mathcal{S}'}$ de $X, S, S'$ respectivement et des \'eclatements admissibles

$$\tau_{\mathcal{X}} : \mathcal{X} \rightarrow\  \stackrel{\circ}{\mathcal{X}} \quad , \quad  \tau_{\mathcal{S}'} : \mathcal{S}'  \rightarrow\  \stackrel{\circ}{\mathcal{S}'} $$
et des morphismes

$$\varphi : \mathcal{X} \rightarrow \mathcal{S} \quad , \quad \theta : \mathcal{S}' \rightarrow \mathcal{S} $$
tels que

$$\varphi_{\mbox{rig}} = f\ \circ\  \tau_{\mathcal{X} \mbox{rig}} \quad \mbox{et} \quad  \theta_{\mbox{rig}} = u\ \circ\  \tau_{\mathcal{S}' \mbox{rig}}.$$

Remarquons que tous les sch\'emas formels pr\'ec\'edents sont admissibles au sens de [Bo - L¬\"u 1], donc sont plats sur $\mathcal{V}$ [Bo - L¬\"u 1, \S\ 1].

\noindent On dispose donc d'un carr\'e cart\'esien

$$
\xymatrix{
\mathcal{X}' \ar[r]^{\theta'} \ar[d]_{\varphi'} & \mathcal{X} \ar[d]^{\varphi} &\\
\mathcal{S'}  \ar[r]_{\theta} &  \mathcal{S} & \quad ,
}
$$

\noindent qui est un mod\`ele formel du carr\'e du th\'eor\`eme : de plus il existe un $\mathcal{O}_{\mathcal{X}}$-module coh\'erent $\mathcal{F}$ tel que $\mathcal{F}_{\mbox{rig}} = E$ [L¬\"u, 2.2], ou [Bo-L¬\"u 1, 5.6] et $\varphi$ est un morphisme propre [L¬\"u, 2.6]. D'apr\`es le th\'eor\`eme (1.1) le morphisme de changement de base

$$\theta^{\ast} R^{i} \varphi_{\ast}(\mathcal{F}) \rightarrow R^{i} \varphi '_{\ast}\ \theta '^{\ast}(\mathcal{F}) $$
est un isomorphisme ; par passage aux fibres g\'en\'eriques le (1.2.1) en r\'esulte.\\

(1.2.2) Par d\'efinition [Bo-L¬\"u 1, \S\  5] un $\mathcal{O}_{X}$-module coh\'erent $M$ est plat sur $S$, s'il existe un mod\`ele formel $h : \mathcal{X} \rightarrow \mathcal{S}$ de $f$ et un $\mathcal{O}_{\mathcal{X}}$-module coh\'erent $\mathcal{M}$ tels que $\mathcal{M}$ est plat sur $\mathcal{S}$ aux points de $\mathcal{X}_{K}$ (i.e. tels que $\mathcal{M}$ est rig-plat sur $\mathcal{S}$) : l'existence de $h$ r\'esulte de [Bo-L¬\"u 1, theo 4.1] et celle de $\mathcal{M}$ r\'esulte de [L¬\"u, lemma 2.2] o\`u l'on peut \'evidemment supposer $\mathcal{M}$ sans $\pi$-torsion.

\vskip 3mm
\noindent \textbf{Lemme (1.2.2.1)}.
\textit{Avec les notations pr\'ec\'edentes, les propri\'et\'es suivantes sont \'equivalentes : }

\begin{enumerate}
\item[(i)] \textit{$M$ est plat sur $S$.}
\item[(ii)] \textit{$\mathcal{M}$ est rig-plat sur $\mathcal{S}$.}
\item[(iii)] \textit{$\mathcal{M}_{K}$ est un $\mathcal{O}_{{\mathcal{S}_{K}}}$-module plat.}
\end{enumerate}
\textit{De plus le $\mathcal{M}$ du (ii) peut \^etre suppos\'e sans $\pi$-torsion : si $\mathcal{S} = Spf\ \mathcal{V}$ le (ii) revient \`a dire que $\mathcal{M}$ est plat sur $\mathcal{V}$ ($\mathcal{V}$ = anneau de valuation discr\`ete complet).}

\vskip 3mm
\noindent \textit{D\'emonstration du lemme}. Par d\'efinition (i) et (ii) sont \'equivalentes. Prouvons l'\'equivalence de (ii) et (iii). Puisque la notion de platitude est locale sur $X$ et $S$ [Bo-L¬\"u 1, \S\ 5], on peut supposer $\mathcal{X}$ et $\mathcal{S}$ affines, $\mathcal{X} = Sp f \hat{A}$ et $\mathcal{S} = Sp f \hat{B}$ o\`u $A$ et $B$ sont des $\mathcal{V}$-alg\`ebres de type fini. On note encore $\mathcal{M} = \Gamma(\mathcal{X}, \mathcal{M})$. Ici les ``rig-points'' au sens de [Bo-L¬\"u 1, \S\ 3] sont les rel\`evements de Teichm¬\"uller $\hat{\tau}(x)$ au sens de [Et 6, 2.1] :

$$
\xymatrix{
0 \ar[r] & \mathfrak{p}_{x}  \ar[r] \ar@{^{(}->}[d] & \hat{A} \ar[r]^{\hat{\tau}(x)} \ar@{^{(}->}[d] & \mathcal{V}(x) \ar[r] \ar[d] & 0 &\\ 
0 \ar[r] & \mathfrak{q}_{x} \ar[r] & \hat{A}_{K} \ar[r]_{\hat{\tau}_{k}(x)} & K(x) \ar[r]   & 0 & \quad .\\
}
$$

\noindent D'apr\`es [Bo-L¬\"u 1], dire que $\mathcal{M}$ est rig-plat sur $\hat{B}$ c'est dire que
$\mathcal{M}_{\mathfrak{p}_{x}}$ est plat sur $\hat{B}$ : or d'apr\`es [Bour, AC II, $\S$\ 2, \no 5, prop 11(iii)] on a un isomorphisme 

$$\mathcal{M}_{\mathfrak{p}_{x}} \tilde{\longrightarrow} (\mathcal{M}_{{K})_{\mathfrak{q}_{x}}}\  ,$$
d'o\`u l'\'equivalence du lemme. $\square$\\

\noindent Alors par [Bo-L¬\"u 2, theo 4.1 et cor 5.9] et quitte \`a faire un \'eclatement formel de $\mathcal{S}$ on peut supposer $\mathcal{M}$ plat sur $\mathcal{O}_{\mathcal{S}}$.\\

 Revenons \`a la preuve de (1.2.2) : le complexe $E^\bullet$ provient d'un complexe $\mathcal{E}^\bullet$ de $\mathcal{O}_{\mathcal{S}}$-modules, \`a composantes $\mathcal{E}^i$ des $\mathcal{O}_{\mathcal{X}}$-modules coh\'erents, et on vient de voir que, quitte \`a \'eclater formellement $\mathcal{S}$, on peut supposer les $\mathcal{E}^i$ plats sur $\mathcal{O}_{\mathcal{S}}$.\\

Il suffit alors d'appliquer le (2) du th\'eor\`eme (1.1) et de passer aux fibres g\'en\'eriques pour obtenir (1.2.2).\\
 
Nous rassemblons pour m\'emoire dans la proposition suivante quelques propri\'et\'es des immersions.

\vskip 3mm
\noindent \textbf{Proposition (1.3)}. \textit{Soient $\mathcal{V}$ et $K$ comme en (1.2). Alors}
\begin{enumerate}
\item[(1.3.1)]  \textit{Toute immersion (resp. immersion ouverte, resp. immersion ferm\'ee) $\alpha : X\hookrightarrow Y$ de $K$-espaces rigides analytiques quasi-compacts et quasi-s\'epar\'es admet un mod\`ele formel $\beta : \mathcal{X} \hookrightarrow \mathcal{Y}$ au sens de [Bo-L¬\"u 2, cor 5.10] qui est une immersion (resp. une immersion ouverte, resp. une immersion ferm\'ee) : en particulier $\mathcal{X}$ et $\mathcal{Y}$ sont plats sur $\mathcal{V}$.}

\item[(1.3.2)]  \textit{Pour tout $K$-espace analytique rigide $X$, le faisceau d'anneaux $\mathcal{O}_{X}$ est coh\'erent.}

\item[(1.3.3)] \textit{Soient $\alpha : X \hookrightarrow Y$ une immersion ferm\'ee de $K$-espaces analytiques rigides et $M$ un $\mathcal{O}_{X}$-module coh\'erent. Alors, pour tout entier $i > 0$ on a  }

$$R^{i} \alpha_{\ast}(M) = 0 $$

\noindent \textit{et le morphisme canonique}

$$\alpha^{\ast}\alpha_{\ast}(M) \rightarrow M $$

\noindent \textit{est un isomorphisme.}\\
\noindent \textit{De plus le morphisme $\alpha_{\ast}$ commute \`a tout changement de base plat $u : Y' \rightarrow Y$}.
\end{enumerate}

\vskip 3mm
\noindent \textit{D\'emonstration}. Le (1.3.1) n'est autre que [Bo-L¬\"u 2, cor 5.10].\\

Pour le (1.3.2) la propri\'et\'e est locale sur $X$ [B-G-R, 9.4.3] : on peut donc supposer $X$ affino¬\"{\i}de ; par suite $X$ est quasi-compact et quasi-s\'epar\'e et admet un mod\`ele formel $\mathcal{X}$.  Il suffit de prouver la coh\'erence du faisceau d'anneaux $\mathcal{O}_{\mathcal{X}}$ :  comme $\mathcal{X}$ est un sch\'ema formel de type fini sur l'anneau noeth\'erien  $\mathcal{V}$, le faisceau $\mathcal{O}_{\mathcal{X}}$ est un faisceau coh\'erent d'anneaux par [EGA I, (10.11.2)].\\

Pour (1.3.3), l'immersion ferm\'ee $\alpha$ est propre [B-G-R, 9.5.3, prop 2, 9.6.2 prop 5] d'o\`u l'isomorphisme canonique

$$\alpha^{\ast} \alpha_{\ast}(M) \overset{\sim}{\rightarrow} M $$

\noindent gr\^ace au th\'eor\`eme (1.2) ; de m\^eme pour la commutation de $\alpha_{\ast}$ aux changements de base plats.\\

L'assertion $R^i\ \alpha_{\ast}(M) = 0 $ pour $i > 0$ est locale sur $Y$ : on peut donc supposer $Y$ affino¬\"{\i}de, donc quasi-compact, quasi-s\'epar\'e et de m\^eme pour $X$ puisque $\alpha$ est propre. On prend alors un mod\`ele formel $\beta : \mathcal{X} \hookrightarrow \mathcal{Y}$ de $\alpha$ et un $\mathcal{O}_{\mathcal{X}}$-module coh\'erent $\mathcal{M}$ tel que $\mathcal{M}_{K} = M$ [L¬\"u, lemma 2.2]. Soient $\beta_{n} : \mathcal{X}_{n} \hookrightarrow \mathcal{Y}_{n}$ la r\'eduction de $\beta$ mod $\mathfrak{m}^{n+1}$ et $\mathcal{M}_{n} = \mathcal{M} /\mathfrak{m}^{n+1}\ \mathcal{M}$ ; alors

$$ \beta^{\ast}  \beta_{\ast}(\mathcal{M}) /\mathfrak{m}^{n+1} = \beta^{\ast}_{n}\ \beta_{n^{\ast}}(\mathcal{M}_{n})$$

\noindent et puisque $\mathcal{M}$ et $\beta_{\ast}(\mathcal{M})$ sont coh\'erents [th\'eo (1.1) (1)] on a des isomorphismes

$$\mathcal{M} \overset{\sim}{\rightarrow} \displaystyle \mathop{\mbox{lim}}_{\longleftarrow \atop{n}}\  \mathcal{M}_{n} \quad \mbox{et} \quad  \beta^{\ast}  \beta_{\ast}(\mathcal{M}) \overset{\sim}{\rightarrow}  \displaystyle \mathop{\mbox{lim}}_{\longleftarrow \atop{n}}\   \beta^{\ast}_{n}\ \beta_{n^{\ast}}(\mathcal{M}_{n}).$$

\noindent Or le morphisme canonique

$$\beta^{\ast}_{n}\ \beta_{n^{\ast}}(\mathcal{M}_{n}) \longrightarrow \mathcal{M}_{n}$$

\noindent est un isomorphisme d'apr\`es [SGA 4, VIII, 5.7], d'o\`u une nouvelle d\'emonstration de l'isomorphisme canonique

$$\alpha^{\ast} \alpha_{\ast}(M) \overset{\sim}{\rightarrow} M $$

\noindent en prenant ci-dessus la limite sur $n$ et en passant aux fibres g\'en\'eriques. L'\'egalit\'e $R^{i} \alpha_{\ast}(M) = 0 $ pour $i > 0$ s'obtient en appliquant [SGA 4, VIII, 5.6]. $\square$

\vskip 3mm
\section*{2. Sorites sur les voisinages stricts}

\textbf{2.0.} Rappelons la d\'efinition de ``voisinage strict'' dans un espace rigide analytique [G-K 2, 2.22].\\

Si $U$ est un ouvert admissible d'un espace rigide analytique $W$, un ouvert admissible $V \subset W$ est appel\'e voisinage strict de $U$ dans $W$ si $\{V, W \backslash U\}$ est un recouvrement admissible de $W$. Cette d\'efinition redonne celle de [B 3, (1.2.1)] dans le cas des tubes.\\
\newpage
\textbf{2.1.} Consid\'erons un diagramme commutatif

$$
\begin{array}{c}
\xymatrix{
X \ar@{^{(}->}[r]^{i_{\mathcal{X}}} \ar[d]_{f}  \ar[d] & \mathcal{X} \ar[d]^{h} \ar[r]^{\psi} \ar @{} [dr] |{\square} & \mathcal{Y}\ar[d]^{\overline{h}}\\
S \ar@{^{(}->} [r]_{i_{\mathcal{S}}} & \mathcal{S} \ar[r]_{\varphi} & \mathcal{T}
} 
\end{array}
\leqno{(2.1.1)}
$$

\noindent dans lequel le carr\'e de droite est un carr\'e cart\'esien de $\mathcal{V}$-sch\'emas formels s\'epar\'es plats de type fini, $f$ est un morphisme de $k$-sch\'emas, $i_{\mathcal{X}}$, $i_{\mathcal{Y}}
 = \psi  \circ i_{\mathcal{X}}$, $i_{\mathcal{S}}$ et $i_{\mathcal{T}} = \varphi  \circ i_{\mathcal{S}}$ sont des immersions.\\
 
 \noindent On note $h' :\ ]X[_{\mathcal{X}}\  \longrightarrow\  ]S[_{\mathcal{S}} ,~\overline{h}' :\  ]X[_{\mathcal{Y}}\ 
\longrightarrow\  ]S[_{\mathcal{T}},~ÊÊÊ \psi' :\ ]X[_{\mathcal{X}}\  \longrightarrow\  ]X[_{\mathcal{Y}},\\\varphi' :\ ]S[_{\mathcal{S}}\ \longrightarrow\  ]S[_{\mathcal{T}}$ les morphismes induits respectivement par $h,\  \overline{h},\  \psi,\\ 
 \varphi$\ [B 3, (1.1.11) (i)].

\vskip 3mm
\noindent \textbf{Proposition (2.1.2)}. \textit{ Avec les notations pr\'ec\'edentes, on a un diagramme commutatif \`a carr\'es cart\'esiens}

$$
\xymatrix{
]X[_{\mathcal{X}} \ar[r]^{\psi'} \ar[d]_{h'} & ]X[_{\mathcal{Y}} \ar[d]^{\overline{h}'}&\\
]S[_{\mathcal{S}} \ar[r]_{\varphi'}\ar@{^{(}->}[d]& ]S[_{\mathcal{T}}\ar@{^{(}->}[d]&\\
\mathcal{S}_{K} \ar[r]_{\varphi_{K}} & \mathcal{T}_{K}&.
}
$$

\vskip 3mm
\noindent \textit{D\'emonstration}. Pour le carr\'e du bas, \c ca r\'esulte de la d\'efinition des tubes et du morphisme de sp\'ecialisation. Pour le carr\'e du haut, il s'agit de v\'erifier que $]X[_{\mathcal{X}}$ satisfait la propri\'et\'e universelle du produit fibr\'e. Soit $Z$ un espace rigide analytique s'ins\'erant dans un diagramme commutatif

$$\xymatrix{
Z \ar[r]^{u}  \ar[d]  &  ]X[_{\mathcal{Y}} \ar[d] ^{\overline{h'}} \ar@{^{(}->}[r]  & \mathcal{Y}_{K} \ar[dd]^{\overline{h}_{K}} &\\
]S[_{\mathcal{S}} \ar[r]_{\varphi'} \ar@{^{(}->}[rd] &  ]S[ _{\mathcal{T}}\ar@{^{(}->}[rd] & &\\
& \mathcal{S}_{K} \ar[r]_{\varphi_{K}} & \mathcal{T}_{K} & \quad ;
}
$$

\noindent par propri\'et\'e universelle de $\mathcal{X}_{K} = \mathcal{S}_{K}  \times_{\mathcal{T}_{K}} \mathcal{Y}_{K}$ on en d\'eduit une fl\`eche $Z \displaystyle \mathop{\longrightarrow}^{v} \mathcal{X}_{K}$ qui s'ins\`ere dans le diagramme commutatif

$$
\xymatrix{
& X \ar@{=}[r]  \ar@{^{(}->}[d]^{i_{\mathcal{X}}} &  X  \ar@{^{(}->}[d]^{i_{\mathcal{Y}}} \\
& \mathcal{X} \ar[r]^{\psi}  & \mathcal{Y} \\
Z \ar[r]^{v}  \ar@/_4pc/[drr]_{u} & \mathcal{X}_{K} \ar[u]_{sp} \ar[r]^{\psi_{K}} & \mathcal{Y}_{K} \ar[u]_{sp} \\
& ]X[_{\mathcal{X}} \ar[r]^{\psi'} \ar@{^{(}->}[u] & ]X[_{\mathcal{Y}} \ar@{^{(}->}[u]
}
$$

\noindent o\`u $sp$ sont les morphismes de sp\'ecialisation. Puisque $sp (]X[_{\mathcal{Y}}) = X \displaystyle \mathop{\hookrightarrow}_{i_{\mathcal{Y}}} \mathcal{Y}$ et compte tenu de la commutativit\'e du diagramme pr\'ec\'edent le morphisme $v$ se factorise par $]X[_{\mathcal{X}}$ en

$$
\xymatrix{
Z \ar[r]^{v}  \ar[dr] & \mathcal{X}_{K} \\
& ]X[_{\mathcal{X}} \ar@{^{(}->}[u] &\quad ;
}
$$

\noindent d'o\`u la proposition. $\square$\\

$ \mathbf{2.2.}$  Consid\'erons \`a pr\'esent un diagramme commutatif

$$
\begin{array}{c}
$$\xymatrix{
X \ar@{^{(}->} [r]^{j_{X_{1}}} \ar [dr]_{f}& X_{1} \ar@{^{(}->}[r]^{j_{Y}} \ar[d]_{f_{1}}  \ar @{} [dr] |{\square} & Y \ar@{^{(}->}[r]^{i_{Y}} \ar[d]^{\overline{f}} & \mathcal{Y} \ar[d] ^{\overline{h}}\\
&S  \ar@{^{(}->}[r]_{j_{T}} & T \ar@{^{(}->}[r]_{i_{T}} & \mathcal{T} \ar[r]^{\rho}
 & \mathcal{W}\\
}
\end{array}
\leqno{(2.2.1)}
$$

\noindent dans lequel le carr\'e de gauche est cart\'esien $f, \ f_{1}$ et $\overline{f}$ sont des morphismes de $k$-sch\'emas, $\overline{h}$ et $\rho$ sont des morphismes de $\mathcal{V}$-sch\'emas formels s\'epar\'es plats de type fini, $ j_{X_{1}}, \  j_{Y}$ et $j_{T}$ sont des immersions ouvertes, $i_{Y}$ et $i_{T}$ sont des immersions ferm\'ees. Notons $X_{2}=\overline{h}^{-1}(S),\ Y_{2}=\overline{h}^{-1}(T)\ \mbox{et}\ f_{2}: X_{2}\rightarrow S,\  \overline{f}_{2}: Y_{2}\rightarrow T$\ les morphismes induits par $\overline{h}$. Soient $Y_{0}$ et $T_{0}$ les r\'eductions modulo $\pi$ de $\mathcal{Y}$ et $\mathcal{T}$ : les immersions ferm\'ees $i_{Y}$ et $i_{T}$ se factorisent respectivement via les immersions ferm\'ees $i_{2Y_{0}} : Y \hookrightarrow Y_{0}, \ i_{Y_{0}}: Y \hookrightarrow Y_{2}\hookrightarrow Y_{0}$ et $i_{T_{0}} : T \hookrightarrow T_{0}$. On d\'esigne par $\overline{h}_{X} :\  ]X[_{\mathcal{Y}} \ \longrightarrow\ ]S[_{\mathcal{T}},~ \overline{h}_{Y} :\   ]Y[_{\mathcal{Y}}\ \longrightarrow\ ]T[_{\mathcal{T}}$ les morphismes induits par $\overline{h}_{K} :\  \mathcal{Y}_{K}\ \longrightarrow\ \mathcal{T}_{K}$
[B 3, (1.1.11) (i)], $j'_{Y} :\ ]X[_{\mathcal{Y}}\ \longrightarrow\ ]Y[_{\mathcal{Y}},~i'_{Y_{0}} :\ ]Y[_{\mathcal{Y}}\ \longrightarrow]Y_{0}[_{\mathcal{Y}}\  = \mathcal{Y_{K}}, \ i'_{2Y_{0}}: ]Y_{2}[\rightarrow \mathcal{Y}_{K}$ ceux induits par l'identit\'e de $\mathcal{Y}_{K}$ et $j'_{T}\ :\ ]S[_{\mathcal{T}}\ \longrightarrow\ ]T[_{\mathcal{T}},~i'_{T_{0}} :\ ]T[_{\mathcal{T}}\ \longrightarrow\ ]T_{0}[_{\mathcal{T}}\ =\ \mathcal{T}_{K}$ ceux induits par l'identit\'e de $\mathcal{T}_{K}$. Si $V$ est un voisinage strict de $]S[_{\mathcal{T}}$ dans $]T[_{\mathcal{T}}$, alors  $W := \overline{h}^{-1}_{Y}(V)$ est un voisinage strict de $ ]X_{1}[_{\mathcal{Y}}$ (donc de $ ]X[_{\mathcal{Y}}$)dans 
$ ]Y[_{\mathcal{Y}}$ [B 3, (1.2.7)] et on note  $h_{V} := \overline{h}_{Y \mid W} : W \rightarrow V$.

\vskip 3mm
\noindent \textbf{Proposition (2.2.2)}.
\textit{Sous les hypoth\`eses (2.2) on a:}
\begin{enumerate}
	\item[(2.2.2.1)] \textit{Supposons que $\overline{f}^{-1}(S) = X$, alors le diagramme (2.2.1) induit un 	diagramme commutatif}

		$$
		\begin{array}{c}
		\xymatrix{
	]X[_{\mathcal{Y}} \ar@{^{(}->}[r]^{j'_{Y}} \ar[d]_{\overline{h}_{X}}  \ar@{} [dr] |{\square} & 		]Y[_{\mathcal{Y}}
	\ar@{^{(}->}[r]^{i'_{Y_{0}}} \ar[d]^{\overline{h}_{Y}} & \mathcal{Y}_{K}  \ar[d]^{\overline{h}_{K}} \\
	]S[_{\mathcal{T}} \ar@{^{(}->} [r]_{j'_{T}} & ]T[_{\mathcal{T}} \ar@{^{(}->} [r]_{i'_{T_{0}}} & 	\mathcal{T}_{K}\ ,
		} 
		\end{array}
		$$

	\textit{dans lequel le carr\'e de gauche est cart\'esien et les fl\`eches horizontales sont des 	immersions ouvertes .}\\

	\item[(2.2.2.2)] \textit{Supposons que $\overline{h}^{-1}(T) = Y $ et $\overline{f}^{-1}(S) = X$. Alors:\\
	Les deux carr\'es du diagramme pr\'ec\'edent sont cart\'esiens.\\
	Si de plus $V$ d\'ecrit un syst\`eme fondamental de voisinages stricts de $]S[_{\mathcal{T}}$ dans 	$]T[_{\mathcal{T}}$, alors $\overline{h}^{-1}_{Y}(V)$ d\'ecrit un syst\`eme fondamental de voisinages 	stricts de $]X[_{\mathcal{Y}}$ dans $]Y[_{\mathcal{Y}}$. }\\
\end{enumerate}

\vskip 3mm
\noindent \textit{D\'emonstration}.\\
 L'existence du diagramme commutatif dans (2.2.2.1) r\'esulte de (2.2.1) et [B 3, (1.1.11) (i)]  via la d\'efinition des tubes [B 3, (1.1.1)].\\
 Puisque $\overline{f}^{-1}(S) = X$ le diagramme (2.2.1) se d\'ecompose en

$$
\begin{array}{c}
$$\xymatrix{
X \ar@{^{(}->}[r]^{j_{Y}} \ar[d] \ar @{} [dr] |{\square} & Y \ar@{^{(}->}[r]^{i_{Y}} \ar[d]& \mathcal{Y} \ar@{=}[d] \\
 \overline{h}^{-1}(S) \ar[r] \ar[d]  \ar @{} [dr] |{\square} &  \overline{h}^{-1}(T) \ar[r] \ar[d] \ar @{} [dr] |{\square} & \mathcal{Y} \ar[d] ^{\overline{h}}\\
S  \ar@{^{(}->}[r]_{j_{T}} & T \ar@{^{(}->}[r]_{i_{T}} & \mathcal{T} \ar[r]^{\rho}
 & \mathcal{W}&.\\
}
\end{array}
\leqno{(2.2.2.3)}
$$
On est donc ramen\'e \`a \'etudier s\'epar\'ement le cas o\`u $\overline{h}^{-1}(T) = Y $ et $\overline{f}^{-1}(S) = X$ et celui o\`u $\overline{h}= id$: on montrera d'abord (2.2.2.2) et pour (2.2.2.1) il nous suffira de traiter le cas o\`u $h=id$.\\

 \textit{Pour (2.2.2.2)}. Les carr\'es du diagramme (2.2.2.1) sont cart\'esiens d'apr\`es la d\'efinition des tubes [B 3, (1.1.1)] et le fait que  $\overline{h}^{-1}(T) = Y $ et $\overline{f}^{-1}(S) = X$.\\

Pour la deuxi\`eme assertion de (2.2.2.2), on va utiliser les voisinages standarts $V_{\underline{\eta}, \underline{\lambda}}= \bigcup_{n\in\mathbb{N}}V_{\eta_{n}, \lambda_{n}}$ de Berthelot [B 3, (1.2.4)]: rappelons au passage que si $\lambda_{n}< \lambda'_{n}$ alors $V_{\eta_{n}, \lambda'_{n}}\subset V_{\eta_{n}, \lambda_{n}}$.\\

Si $V$ est un voisinage strict de $]S[_{\mathcal{T}}$ dans $]T[_{\mathcal{T}}$, alors $\overline{h}^{-1}_{Y}(V)$ est un voisinage strict de $]X[_{\mathcal{Y}}$ dans $]Y[_{\mathcal{Y}}$ [B 3, (1.2.7)].  Soit $W$ un voisinage strict de $]X[_{\mathcal{Y}}$ dans $]Y[_{\mathcal{Y}}$ ; d'apr\`es [B 3, (1.2.2)] on se ram\`ene au cas o\`u $]Y[_{\mathcal{Y}}$ est affino¬\"{\i}de, et avec les notations de [loc. cit.] il existe $\lambda_{0} < 1$ tel que, pour $\lambda_{0} \leqslant \lambda < 1$, on ait $U_{\lambda} \subset W$. Avec les notations de [B 3, (1.2.4) (i)] si $V_{\underline{\eta},\underline{\lambda}}(S)$ parcourt un syst\`eme fondamental de voisinages stricts de $]S[_{\mathcal{T}}$ dans $]T[_{\mathcal{T}}$,  alors il existe $\eta_{n}$ et $\lambda_{n}$ assez proches de 1 tels que $(\overline{h}_{Y})^{-1}\ (V_{\eta_{n}, \lambda_{n}}(S))\ \subset\ U_{\lambda}$ car les \'equations locales de $Y$ (resp de $Z := Y\  \backslash\  X$) sont obtenues par image inverse par $\overline{h}$ des \'equations locales de $T$ (resp de $T\ \backslash\ S$) (cf. aussi [LS, prop 3.2.8]).\\

\textit{Pour (2.2.2.1)}. Puisque $]X[_{\mathcal{Y}}$ et $]Y[_{\mathcal{Y}}$ sont des ouverts de $\mathcal{Y}_{K}$ [B 3, (1.1.2)] il en r\'esulte que $j'_{Y}$ et $i'_{Y_{0}}$ sont des immersions ouvertes ; de m\^eme pour $j'_{T}$ et $i'_{T_{0}}.$\\
  Pour montrer que le carr\'e de gauche du diagramme (2.2.2.1) est cart\'esien il nous suffit de traiter le cas o\`u $h=id$. Comme $\overline{f}^{-1}(S)=X$, les \'equations locales de $Z=Y\setminus X$ sont obtenues par image inverse par $\overline{f}$ des \'equations locales de $T\setminus S$: la d\'efinition des tubes [B 3, (1.1.1)] fournit alors  l'\'egalit\'e $]X[_{\mathcal{Y}}=]S[_{\mathcal{Y}}\ \bigcap\ ]Y[_{\mathcal{Y}}$, d'o\`u le carr\'e cart\'esien de (2.2.2.1). \ $\square$\\

\noindent \textbf{Proposition (2.2.3)}.
\textit{Sous les hypoth\`eses de (2.2) on a:}
\begin{enumerate}
\item[(2.2.3.1)]  \textit{Supposons $\overline{h}^{-1}(T)=Y, \overline{f}^{-1}(S)=\overline{h}^{-1}(S)=X$ et $\overline{h}$ est propre. Alors $\overline{h}_{K}, \overline{h}_{Y}$ et $\overline{h}_{X}$ sont propres. Si $V$ est un voisinage strict de $]S[_{\mathcal{T}}$ dans $]T[_{\mathcal{T}}$ et $W = \overline{h}^{-1}_{Y}(V)$, alors $h_{V} = \overline{h}_{Y \mid_{W}} : W \rightarrow V$ est propre.}
\item[(2.2.3.2)] \textit{Supposons $\overline{h}$ lisse sur un voisinage de $X$ dans $\mathcal{Y}$. Alors}
	\begin{enumerate}
	\item[(i)] \textit{$\overline{h}_{X}$ est lisse et quel que soit $V$ un voisinage strict de $]S[_{\mathcal{T}}$ dans $]T[_{\mathcal{T}}$  il existe un voisinage strict $W$ de $]X[_{\mathcal{Y}}$ dans $]Y[_{\mathcal{Y}}$ tel que $\overline{h}_{K}$ induise un morphisme lisse $h_{V}: W \rightarrow V$. De plus $h_{V}(W)$ est un ouvert admissible de $V$ et de $\mathcal{T}_{K}$, et $\Omega^{i}_{W / V}$ est un $\mathcal{O}_{W}$-module coh\'erent et localement libre. }
	\item[(ii)] \textit{Si l'on suppose aussi que $\overline{h}^{-1}(T)=Y$, et $ \overline{h}^{-1}(S)=X$, alors il existe un voisinage strict $V$ de $]S[_{\mathcal{T}}$ dans $]T[_{\mathcal{T}}$ tel qu'en posant $W= \overline{h}^{-1}_{Y}(V)$ le morphisme $\overline{h}_{K}$ induise un morphisme lisse \\
	$h_{V}:W \rightarrow V$ }
	\item[(iii)] \textit{Si en outre $g : \mathcal{T} \rightarrow Spf\ \mathcal{V}$ est lisse sur un voisinage de $S$ dans $\mathcal{T}$, alors on peut prendre le $V$ du (i)  lisse sur $K$ et ainsi $\Omega^{1}_{V /K}$ est localement libre de type fini sur le faisceau coh\'erent d'anneaux $\mathcal{O}_{V}$.Ê}
	\item[(iv)] \textit{Supposons $f$ surjectif, $\overline{h}^{-1}(T)=Y, \overline{h}^{-1}(S)=X$ et pour $V$ et $W$ comme en (ii) posons $V' = h_{V}(W)$.\\
Si $(W_{\lambda})_{\lambda}$ est un syst\`eme fondamental de voisinages stricts de $]X[_{\mathcal{Y}}$ dans $] Y [_{\mathcal{Y}}$ avec $W_{\lambda} \subset W$, alors $(h_{V}(W_{\lambda}))_{\lambda}$ est un syst\`eme fondamental de voisinages stricts de $] S [_{\mathcal{T}}$ dans $V'$.  }
	\item[(v)] \textit{ Sous les hypoth\`eses (iv) supposons $\overline{h}$ propre. Alors $h_{V}$ induit 	un morphisme propre lisse et surjectif}
$$h_{V'} : W\ \longrightarrow\ V'\ .$$
	\end{enumerate}
\end{enumerate}

\vskip 3mm
\noindent \textit{D\'emonstration}. \\
\textit{Prouvons (2.2.3.1)}. Sous nos hypoth\`eses les deux carr\'es de (2.2.2.1) sont cart\'esiens [cf (2.2.2.2)]. D'apr\`es L¬\"utkebohmert [L¬\"u, theo 3.1] le morphisme propre $\overline{h}$ induit un morphisme propre d'espaces analytiques rigides $\overline{h}_{K} : \mathcal{Y}_{K} \rightarrow \mathcal{T}_{K} $. Comme la notion de morphisme propre est stable par changement de base en g\'eom\'etrie rigide [B-G-R, fin de 9.6.2, p 396], on en d\'eduit que $\overline{h}_{Y},\overline{h}_{X}$ et $h_{V}$  sont propres.\\

\noindent \textit{Pour (2.2.3.2)}.

\begin{enumerate}

\item[(i)] L'ensemble $W'$ des points de $] Y [_{\mathcal{Y}}$ o\`u le morphisme $\overline{h}_{K}$ est lisse est un voisinage strict de $]X[_{\mathcal{Y}}$ dans $] Y [_{\mathcal{Y}}$ [B 3, (2.2.1)]. Si $V$ est un voisinage strict quelconque de $]S[_{\mathcal{T}}$ dans $]T[_{\mathcal{T}}, W= \overline{h}_{Y}^{-1}(V)\cap W'$  est un voisinage strict de $]X[_{\mathcal{Y}}$ dans $]Y[_{\mathcal{Y}}$ [B 3, (1.2.7) et (1.2.10)] et $\overline{h}_{K}$ induit donc un morphisme lisse $h_{V}: W \rightarrow V$; en particulier $\overline{h}_{X}: ]X[_{\mathcal{Y}} \rightarrow ]S[_{\mathcal{T}} \subset V$ est lisse.\\
Ê\quad Puisque $h_{V}$ est plat il est ouvert pour la topologie rigide [Bo-L¬\"u 2, 5.11], donc $h_{V}(W)$ est un ouvert admissible de $V$, donc de $\mathcal{T}_{K}$ car $V$ est un ouvert admissible de $\mathcal{T}_{K}$. La lissit\'e de $h_{V}$ prouve que $\Omega^i_{W / V}$ est un $\mathcal{O}_{W}$-module localement libre de type fini, et comme $\mathcal{O}_{W}$ est un faisceau coh\'erent d'anneaux [prop (1.3)], il r\'esulte de [EGA $O_{I}$, (5.4.1)] que $\Omega^i_{W / V}$ est un $\mathcal{O}_{W}$-module coh\'erent.\\

\item[(ii)] Si $V$ d\'ecrit un syst\`eme fondamental de voisinages stricts de $]S[_{\mathcal{T}}$ dans $]T[_{\mathcal{T}}$, alors $\overline{h}_{Y}^{-1}(V)$ d\'ecrit un syst\`eme fondamental de voisinages stricts de $]X[_{\mathcal{Y}}$ dans $]Y[_{\mathcal{Y}}$ [(2.2.2.2)]: ainsi (ii) r\'esulte de (i).\\

\item[(iii)] Supposons en outre $g : \mathcal{T} \rightarrow Spf \mathcal{V}$ lisse sur un voisinage de $S$ dans $\mathcal{T}$ ; alors l'ensemble des points de $] T [_{\mathcal{T}}$ o\`u $g_{K}$ est lisse est un voisinage strict de $] S [_{\mathcal{T}}$ dans $] T [_{\mathcal{T}}$ [B 3, (2.2.1)] : quitte \`a restreindre le $V$ du(i) [B 3, (1.2.10)] on peut supposer que la restriction de $g_{K}$ \`a $V$ est lisse. Ainsi $\Omega^1_{V / K}$ sera un $\mathcal{O}_{V}$-module localement libre de type fini, donc coh\'erent sur $\mathcal{O}_{V}$  [(1.3)]. D'o\`u (iii).\\

\item[(iv) et (v)] Puisque $f$ est surjectif et $\overline{h}$ plat au voisinage de $X$, le morphisme

$$\overline{h}_{X} : ] X [_{\mathcal{Y}}\  \longrightarrow\  ] S [_{\mathcal{T}} $$

induit par $\overline{h}$ est surjectif [B 3, (1.1.12)] et lisse puisque c'est aussi la restriction de $h_{V}$. On a vu en (i) que $V' = h_{V}(W)$ est un ouvert admissible de $V$ et $h_{V}(W)$ contient $] S [_{\mathcal{T}}$ par la surjectivit\'e de $\overline{h}_{X}$. Dans le diagramme commutatif

$$
\begin{array}{c}
 \xymatrix{
]X[_{\mathcal{Y}} \ar@{^{(}->}[r] \ar@{->>}[d]_{\overline{h}_{X}}   & W = h^{-1}_{V}(V') = h^{-1}_{V}(V)   \ar@{->>}[d]^{h_{V'}}  \ar@{=} [r] & W  \ar[d]^{h_{V}}\\
]S[_{\mathcal{T}} \ar@{^{(}->} [r] &V' := h_{V}(W) \ar@{^{(}->} [r] & V
} 
\end{array}
\leqno{(2.2.3.2.1)}
$$

les carr\'es sont cart\'esiens d'apr\`es (2.2.2.2)et $h_{V'}$, $\overline{h}_{X}$ sont lisses et surjectifs.\\

 \quad Soit $W_{\lambda}$ un voisinage strict de $]X[_{\mathcal{Y}}$ dans $]Y[_{\mathcal{Y}}$ avec $W_{\lambda} \subset  W$, o\`u $W$ est d\'efini ci-dessus : le morphisme plat $h_{V}$ envoie le recouvrement admissible  $\{ W_{\lambda} ;  W\  \backslash\  ]X[_{\mathcal{Y}} \}$ de $W$ sur le recouvrement admissible $\{ h_{V}(W_{\lambda}) ; \hfill\break h_{V}(W\  \backslash\   ]X[_{\mathcal{Y}}) \}$ de $V' = h_{V}(W)$. Par la surjectivit\'e de $h_{V'}$ et le fait que $h^{-1}_{V} (] S [_{\mathcal{T}}) =\  ] X [_{\mathcal{Y}}$ on a $h_{V}(W)\ \backslash\  ] X [_{\mathcal{Y}}) = V'\ \backslash\ ] S [_{\mathcal{T}} $ ; par suite $\ \{ h_{V}(W_{\lambda})) ; V'\ \backslash\ ] SÊ[_{\mathcal{T}} \}$ est un recouvrement admissible de $V'$, i.e. $h_{V}(W_{\lambda})$ est un voisinage strict de $] S [_{\mathcal{T}}$ dans $V'$.\\

Ê \quad Supposons maintenant que $(W_{\lambda})_{\lambda}$ d\'ecrive un syst\`eme fondamental de voisinages stricts de $]X[_{\mathcal{\mathcal{Y}}}$ dans $]Y[_{\mathcal{Y}}$ avec $W_{\lambda} \subset W$. Soit $V''$ un voisinage strict de $] S [_{\mathcal{T}}$ dans $V'$ : montrons que $h^{-1}_{V'}(V'')$ est un voisinage strict de $]X[_{\mathcal{Y}}$ dans $W$. D'abord $\{ V''; V' \backslash\ ] S [_{\mathcal{T}} \}$ est un recouvrement admissible de $V'$, donc par image inverse $\{ h^{-1}_{V'}(V'') ; h^{-1}_{V'}(V' \backslash ] S [_{\mathcal{T}}) \}$ est un recouvrement admissible de
 $h^{-1}_{V'}(V') = W$ ; en utilisant encore l'\'egalit\'e $h^{-1}_{V'} (] S [_{\mathcal{T}}) = ] X [_{\mathcal{Y}}$ on en d\'eduit que $h^{-1}_{V'}(V' \backslash\ ] S [_{\mathcal{T}}) = W \backslash\ ] X [_{\mathcal{Y}}$ : donc  $h^{-1}_{V'}(V'')$ est un voisinage strict de $] X [_{\mathcal{Y}}$ dans $W$. Ainsi il existe $\mu$ tel que $W_{\mu} \subset h^{-1}_{V'}(V'')$, d'o\`u 
 
 $$h_{V'}(W_{\mu}) \subset h_{V'}\  h^{-1}_{V'}(V'') = V''  ; $$
 
 par suite $(h_{V}(W_{\lambda}))_{\lambda}$ est bien un syst\`eme fondamental de voisinages stricts de $] S [_{\mathcal{T}}$ dans $V'$.\\
 
 \quad Si de plus $\overline{h}$ est propre alors $h_{V'}$ est de surcro\^\i t propre puisque (2.2.3.2.1) est \`a carr\'es cart\'esients et $h_{V}$ est propre. $\square$
 
 \end{enumerate}

 \vskip5mm
 
\textbf{2.3.} Soient $S = Spec\ A_{0}$ un $k$-sch\'ema lisse et $f : X = Spec\ B_{0}\ \rightarrow\ S$ un $k$-morphisme fini. D\'esignons par $A = \mathcal{V} [t_{1},..., t_{n}]  / (f_{1},..., f_{r})$ une $\mathcal{V}$-alg\`ebre lisse relevant $A_{0}$, $\mathcal{S} = Spf\ \hat{A}$, $\tilde{\mathcal{S}} = \hat{P}$ et $\overline{\mathcal{S}} = \widehat{P'}$ comme dans le th\'eor\`eme (3.4) du I : on sait [loc. cit.] que $\tilde{\mathcal{S}}$ et $\overline{\mathcal{S}}$ sont propres sur $\mathcal{V}$, que $\overline{\mathcal{S}}$ est normal, qu'il existe une $A$-alg\`ebre finie $B$ et un carr\'e cart\'esien de $\mathcal{V}$-sch\'emas formels

$$
 \xymatrix{
Spf \hat{B} = \mathcal{X} \ar@{^{(}->}[r] \ar[d]_{h}  & \overline{\mathcal{X}} = \widehat{P''_{1}}  \ar[d]^{\overline{h}}  \\
Spf \hat{A} = \mathcal{S}  \ar@{^{(}->} [r]_{j}  &\overline{\mathcal{S}} = \widehat{P'}
} 
$$

 \noindent o\`u $P''_{1}$ est la fermeture int\'egrale de $P'$ dans Spec $B$, avec $h$ fini, $\overline{h}$ fini et $j$ une immersion ouverte. On note $\overline{f} :\overline{X}\ \rightarrow\ \overline{S}$ la r\'eduction de $\overline{h} : \overline{\mathcal{X}} \rightarrow \overline{\mathcal{S}}$ sur $k = \mathcal{V} / \mathfrak{m}$.\\
 
Rappelons [I, th\'eo 3.4] qu'il existe $a \in A$ et $f(t) \in A_{a}[t]$ tels que $B$ est la fermeture int\'egrale de $A$ dans $A_{a} [t]  /  (f)$. Fixons d'autre part une pr\'esentation de la $\mathcal{V}$-alg\`ebre $B$

$$B\ \simeq\ \mathcal{V} [t'_{1} ,..., t'_{n'}] / (g_{1} ,..., g_{s}).	$$

\noindent Soient $\overline{\mathcal{Y}}$ le compl\'et\'e formel de la fermeture projective $P''_{2}$ de Spec $B$ dans $\mathbb{P}^{n'}_{\mathcal{V}}$ et $\overline{Y}$ sa r\'eduction sur $k$.\\

\quad Comme $P'$ est le normalis\'e de $P$ on a un triangle commutatif

$$
\xymatrix{
&   \overline{\mathcal{S}}:=\hat{P'} \ar[d]^{v} \\
\mathcal{S} \ar@{^{(}->}[ur]^{j} \ar@{^{(}->}[r]_{\tilde{j}} & \tilde{\mathcal{S}}:=\hat{P}\\
}
$$

\noindent o\`u $v$ est fini et $\tilde{j}$ une immersion ouverte. Un syst\`eme fondamental de voisinages stricts de $\mathcal{S}_{X}$ dans $\tilde{\mathcal{S}}_{K}$ est fourni par les intersections $\tilde{V}_{\lambda}$ de $(\mbox{Spec}\ A)^{an}_{K}$ avec les boules $B(0, \lambda^{\dag}) \subset \mathbb{A}^{n}_{K}$ pour $\lambda \rightarrow 1^+$ et $\tilde{V}_{\lambda} = \mbox{Spm}\ A_{\lambda}$, $A^{\dag}_{K} = \displaystyle \mathop{\lim}_{\rightarrow \atop{\lambda}} A_{\lambda}$  [B 3, (2.5.1)]. Puisque $v$ est propre, et \'etale au voisinage de $\mathcal{S}$, il existe $\lambda_{0}  > 1$ tel que tout $\lambda, 1 < \lambda \leqslant \lambda_{0}$, $v$ induise un isomorphisme entre $\tilde{V}_{\lambda}$ et un voisinage strict $V_{\lambda}$ de $\mathcal{S}_{K}$ dans $\overline{S}_{K}$ [B 3, (1.3.5)]  : on identifiera $V_{\lambda}$ et $\tilde{V}_{\lambda}$ dans la suite. \\

Notons $P''_{3}$ l'adh\'erence sch\'ematique de $\mbox{Spec}\  B$ plong\'e diagonalement dans $ P''_{1}\ \times_{\mathcal{V}} P''_{2}$, $\overline{\mathcal{Z}} = \widehat{P''_{3}}$ le compl\'et\'e formel de $P''_{3}$ et $\overline{Z}$ sa r\'eduction mod $\mathfrak{m}$. On a un diagramme commutatif

$$
\xymatrix{
&   \overline{Y} \ar@{^{(}->}[r] & \overline{\mathcal{Y}}  \\
X \ar@{^{(}->}[ur] \ar@{^{(}->}[r]  \ar@{^{(}->}[rd] & \overline{Z} \ar@{^{(}->}[r] \ar[u]_{v_{2}} \ar[d]^{v_{1}}& \overline{\mathcal{Z}} \ar[u]_{u_{2}} \ar[d]^{u_{1}}\\
& \overline{X} \ar@{^{(}->}[r] & \overline{\mathcal{X}}
}
$$

\noindent o\`u les $u_{i}, v_{i}$ sont propres et les $u_{i}$ sont \'etales au voisinage de $X$. D'apr\`es [B 3, (1.3.5)] $u_{1K}$ induit un isomorphisme  entre un voisinage strict de $] X [_{\overline{\mathcal{Z}}}$ dans $\overline{\mathcal{Z}}_{K}$ et un voisinage strict de $] X [_{\overline{\mathcal{X}}}\ \simeq\  ] X [_{\mathcal{X}}\  =\  \mathcal{X}_{K}$ dans $\overline{\mathcal{X}}_{K}$ et par suite un isomorphisme entre des syst\`emes fondamentaux de tels voisinages stricts. De m\^eme $u_{2K}$ induit un isomorphisme entre un syst\`eme fondamental de voisinages stricts de $] X [_{\overline{\mathcal{Z}}}$ dans $\overline{\mathcal{Z}}_{K}$ et un syst\`eme fondamental $(W'_{\lambda'})_{\lambda'} = ( \mbox{Spm}\ B_{\lambda'})_{\lambda'}$ de voisinages stricts de $\mathcal{X}_{K}$ dans $\overline{\mathcal{Y}}_{K}$. Par composition il en r\'esulte pour $\lambda'\rightarrow 1^+$ un isomorphisme entre les $W'_{\lambda'} = \mbox{Spm}\ B_{\lambda'}$ et un syst\`eme fondamental de voisinages stricts $(W''_{\lambda''})$ de $\mathcal{X}_{K}$ dans $\overline{\mathcal{X}}_{K}$ identifi\'es ci-apr\`es. Pour $\lambda > 1$, il existe donc $\lambda' > 1$ et des immersions ouvertes

$$\mbox{Spm}\ \hat{B}_{K} = \mathcal{X}_{K} \hookrightarrow \mbox{Spm}\ B_{\lambda'} \displaystyle \ \mathop{\hooklongrightarrow}^{j'_{\lambda \lambda'}}\  \overline{h}^{-1}_{K}(V_{\lambda}) =: W_{\lambda} .$$

\newpage
\noindent \textbf{Proposition (2.3.1)}.
\textit{Avec les notations de (2.3) on a : }

\begin{enumerate}
\item[(1)] \textit{Si $(V_{\lambda})_{\lambda}$ est un syst\`eme fondamental de voisinages stricts de $\mathcal{S}_{K}$ dans $\overline{\mathcal{S}}_{K}$, alors $(W_{\lambda})_{\lambda} :=  (\overline{h}^{-1}_{K}(V_{\lambda}))_{\lambda} $ est un syst\`eme fondamental de voisinages stricts de $\mathcal{X}_{K}$ dans $\overline{\mathcal{X}}_{K}$. }

\item[(2)] \textit{Supposons de plus $f$ fini et plat (resp. fini et fid\`element plat, resp.fini \'etale, resp. fini \'etale galoisien de groupe G) et $V_{\lambda} = \mbox{Spm}\ A_{\lambda}$. Alors il existe $\lambda_{0} > 1$ tel que pour tout $\lambda, 1 < \lambda \leqslant \lambda_{0}$, et $W_{\lambda} := \overline{h}^{-1}_{K}(V_{\lambda}) $, le morphisme induit par $\overline{h}_{K}$ }

$$h_{\lambda}	:= \overline{h}_{K \mid W_{\lambda}} : W_{\lambda} \longrightarrow V_{\lambda}$$

\textit{soit fini et plat (resp. fini et fid\`element plat, resp. fini \'etale, resp. fini \'etale galoisien de groupe G), avec $V_{\lambda}$ lisse sur $K$ et $\Omega^1_{V_{\lambda}/K}$ localement libre de type fini sur le faisceau coh\'erent d'anneaux $\mathcal{O}_{V_{\lambda}}$.}

\end{enumerate}

\vskip 3mm
\noindent \textit{D\'emonstration}.  On utilise les notations du (2.3).

\begin{enumerate}

\item[(1)] On a d\'ej\`a prouv\'e le (1) dans la proposition (2.1.2) : on en donne ici une autre d\'emonstration. Il suffit de faire la d\'emonstration dans le cas $V_{\lambda} = \mbox{Spm}\ A_{\lambda}$ et on peut supposer $A, B, P'$ et $P''_{1}$ int\'egralement clos.\\
Notons $W_{\lambda} = \overline{h}^{-1}_{K}(V_{\lambda})$ et $h_{\lambda} := \overline{h}_{K \mid W_{\lambda}} : W_{\lambda} \rightarrow V_{\lambda}$. Comme $P''_{1}$ est int\`egre les immersions ouvertes 

$$\mathcal{X}_{K} \hookrightarrow W_{\lambda} \hookrightarrow \overline{\mathcal{X}}_{K} $$

sont dominantes.\\
On a vu \`a la fin de la preuve du (1) de la proposition (2.2.1) que pour $\lambda$ assez proche de 1, $V_{\lambda}$ est lisse sur $K$ et $\Omega^1_{V_{\lambda}/K}$ localement libre de type fini sur l'anneau coh\'erent $\mathcal{O}_{V_{\lambda}}$, d'o\`u la fin du (2).\\
Comme $B_{\lambda'}$ est form\'e d'\'el\'ements entiers sur $A^{\dag}_{K}$, il existe $\mu$, $1 < \mu \leqslant \lambda$, tel que $B_{\lambda'}$ soit fini sur $A_{\mu}$. Par suite on a un diagramme commutatif  \`a carr\'e cart\'esien

$$
\xymatrix{W'_{\lambda'} = Spm\  B_{\lambda'}  \ar@/^/[rrd]^{j'_{\lambda \lambda'}} \ar@/_/[rdd]_{h'} \ar@{.>}[rd]_{j'_{\mu \lambda'}} \\
& W_{\mu} \ar[d]^{h_{\mu}} \ar@{^{(}->}[r]_{\alpha'_{\lambda \mu}} & W_{\lambda}  \ar[d]^{h_{\lambda}}\\
& V_{\mu}=Spm\  A_{\mu} \ar@{^{(}->}[r]_{\alpha_{\lambda \mu}} & V_{\lambda}=Spm\  A_{\lambda}
}
$$

avec une factorisation $j'_{\lambda \lambda'} = \alpha'_{\lambda \mu}\  \circ\  j'_{\mu \lambda'}$ et $h'$ fini : les immersions ouvertes $j'_{\lambda \lambda'}$ et $\alpha'_{\lambda \mu}$ admettent des mod\`eles formels qui sont des immersions ouvertes entre sch\'emas formels noeth\'eriens [Bo-L¬\"u 2, cor 5.10] ; en passant par un mod\`ele formel de $j'_{\mu \lambda'}$, il en r\'esulte que $j'_{\mu \lambda'}$ est aussi une immersion ouverte et $j'_{\mu \lambda'}$ est dominante car $j'_{\lambda \lambda'}$ et $\alpha'_{\lambda \mu}$ le sont. D'autre part il existe des mod\`eles formels propres de $h'$ et $h_{\mu}$ [L¬\"u, 2.6] : par suite $j'_{\mu \lambda'}$ est propre ; or $j'_{\mu \lambda'}$ est une immersion ouverte dominante, donc $j'_{\mu \lambda'}$ est un isomorphisme. \\

D'o\`u le (1).\\

\item[(2)] L\`a encore on peut supposer $A$ et $B$ int\'egralement clos : on traitera \`a part le cas fini \'etale galoisien.\\

\ \quad Notons $\varphi : A \rightarrow B$ le morphisme fini tel que Spec $\varphi : \mbox{Spec}\  B \rightarrow \mbox{Spec}\  A$ rel\`eve $f$ [I, th\'eo 3.4], et $\hat{\varphi} : \hat{A} \rightarrow \hat{B}$ (resp. $\varphi^{\dag}
 : A^{\dag} \rightarrow B^{\dag})$ le morphisme induit sur les s\'epar\'es compl\'et\'es (resp. sur les compl\'et\'es faibles). D'apr\`es le [I, th\'eo 3.4], $\varphi^{\dag}$ et $\hat{\varphi}$ sont finis et plats (resp. finis et fid\`element plats, resp. finis \'etales) si et seulement si $f$ l'est. Avec les notations du (2.3) on a :

 $$
 A^{\dag}_{K} = \displaystyle{\lim_{\rightarrow \atop{> \atop{\lambda \rightarrow 1}}}}\  A_{\lambda}, B^{\dag}_{K} = \displaystyle{\lim_{\rightarrow \atop{> \atop{\lambda \rightarrow 1}}}} B_{\lambda}\  \textrm{et}\   \varphi^{\dag}_{K} : A^{\dag}_{K} \rightarrow B^{\dag}_{K}
 $$
 
\noindent est la limite inductive des $\varphi_{\lambda} : A_{\lambda} \rightarrow B_{\lambda}$ avec 
 
 $$h_{\lambda} = \mbox{Spm}\ (\varphi_{\lambda}) : W_{\lambda} = \mbox{Spm}\ B_{\lambda}\ \rightarrow V_{\lambda} = \mbox{Spm}\ A_{\lambda}. $$

 Si $f$ est fini et plat (resp. ...) alors $\varphi^{\dag}_{K}$ l'est et pour $\lambda
$ assez proche de 1, $\varphi_{\lambda}$ l'est aussi par [EGA IV, 11.2.6, 8.10.5, 17.7.8], de m\^eme  pour $h_{\lambda}$. D'o\`u le (2) hormis la cas galoisien.\\

\quad Consid\'erons \`a pr\'esent le cas galoisien.\\

\quad Puisque $f$ est galoisien il est surjectif; ainsi \\
$$h : \mathcal{X} = \mbox{Spf}\ \hat{B}\ \rightarrow\ \mathcal{S} = \mbox{Spf}\ \hat{A}$$
 est fini \'etale surjectif et galoisien de groupe $G$, d'o\`u en particulier une injection : $\hat{A} \hookrightarrow \hat{B}$. Par suite $h_{K} : \mathcal{X}_{K} \rightarrow \mathcal{S}_{K}$ est fini \'etale surjectif et galoisien de groupe $G$.\\

\quad Remarquons ensuite que puisque $B$ est la fermeture int\'egrale de $A$ dans $A_{a}[t] / (f)$ et que Spec $A
 \rightarrow \mbox{Spec}\ A^{\dag} $ est un morphisme normal [I, prop (1.1)], il r\'esulte de [EGA IV, (6.14.4)] que $B^{\dag} = B \otimes_{A} A^{\dag} $ est la fermeture int\'egrale de $A^{\dag}$ dans $(B^{\dag})_{a} = (A^{\dag})_{a}\  [t] / (f) $ : par suite $B^{\dag}_{K}$ est la fermeture int\'egrale de $A^{\dag}_{K}$ dans $(A^{\dag})_{a,K}\  [t] / (f)$. L'anneau $A^{\dag}$ est r\'eduit par [I, prop (1.6)], car $A$ est r\'eduit, et $A^{\dag} \rightarrow B^{\dag}$ est fini \'etale car $\hat{A} \rightarrow \hat{B}$ l'est : donc $B^{\dag}$ est r\'eduit car $A^{\dag}$ est r\'eduit [I, lemme (1.5)]. Ainsi $B^{\dag}$ est int\'egralement ferm\'e dans $\hat{B}$ [I, th\'eo (2.2) (2) ii] ; d'o\`u $B^{\dag}_{K}$ est la fermeture int\'egrale de $A^{\dag}_{K}$ dans $\hat{B}_{K}$.\\
 
 \quad On a vu ci-dessus que, pour $\lambda$ suffisamment proche de 1, 
 
 $$h_{\lambda} : W_{\lambda} = \mbox{Spm}\ B_{\lambda} \rightarrow V_{\lambda} = \mbox{Spm}\ A_{\lambda} $$
 
 est fini \'etale. Compte tenu du diagramme commutatif

  $$
 \xymatrix{
 \hat{B}_{K} &    & B_{\lambda}  \ar@{_{(}->} [ll]  \\
 \hat{A}_{K} \ar@{^{(}->} [u] & A^{\dag}_{K} \ar@{_{(}->} [l] & A_{\lambda} \ar@{_{(}->} [l] \ar[u]
  }
 $$

 \noindent la fl\`eche $A_{\lambda} \rightarrow B_{\lambda} $ induite par $h_{\lambda}$ est injective, donc $h_{\lambda}$ est surjectif. Choisissons un ensemble fini $\{x_{i} \}$ de g\'en\'erateurs de $B_{\lambda}$ sur $A_{\lambda}$.  Comme chaque $x_{i}$ est entier sur $A_{\lambda}$, les \'el\'ements $g_{\hat{B}_{K}}(x_{i}) \in \hat{B}_{K}$, pour $g$ d\'ecrivant $G$ et $g_{\hat{B}_{K}} : \hat{B}_{K} \rightarrow \hat{B}_{K}$ induit par $g$, sont aussi entiers sur $A_{\lambda} \subset A^{\dag}_{K}$, donc a fortiori sur $B^{\dag}_{K} = \displaystyle \mathop{\lim}_{\rightarrow \atop{\mu}} B_{\mu}$. Or on a vu que $B^{\dag}_{K}$ est int\'egralement ferm\'e dans $\hat{B}_{K}$ : il existe donc $\lambda'$, $1 < \lambda' \leqslant \lambda$ tel que pour tout $i$ et tout $g \in G$ on ait $g_{\hat{B}_{K}}(x_{i}) \in B_{\lambda'}$. Ainsi l'action de $G$ s'\'etend de $\mathcal{X}_{K}$ \`a $W_{\lambda'} = \mbox{Spm}\ B_{\lambda'}$ : en effet $g \in G$ d\'efinit un morphisme $g_{\lambda \lambda'}: W_{\lambda'} \rightarrow W_{\lambda}$ s'ins\'erant dans le diagramme commutatif \`a carr\'e cart\'esien

$$
\xymatrix{W_{\lambda'} \ar@/^/[rrd]^{g_{\lambda \lambda'}} \ar@/_/[rdd]_{h_{\lambda'}} \ar@{.>}[rd]_{g_{\lambda'}} \\
& W_{\lambda'} \ar[d]^{h_{\lambda'}} \ar@{^{(}->}[r] & W_{\lambda}  \ar[d]^{h_{\lambda}}\\
& V_{\lambda'} \ar@{^{(}->}[r]_{\alpha_{\lambda \lambda'}} & \ V_{\lambda}\ ;}
$$

  \noindent d'o\`u la factorisation de $g_{\lambda \lambda'}$ par $W_{\lambda'}$.\\
  
  Montrons que le morphisme fini \'etale surjectif
  
  $$h_{\lambda'} : W_{\lambda'} \rightarrow V_{\lambda'}  $$
  
  est galoisien de groupe $G$. On a
  
  $$(B_{\lambda'})^G \subset (\hat{B}_{K})^G = \hat{A}_{K} ;	$$
  
  d'o\`u
  
  $$A_{\lambda'} \subset (B_{\lambda'})^G \subset B_{\lambda'} \cap \hat{A}_{K} $$
  
  et on dispose d'un carr\'e commutatif \\
  
$$
\xymatrix{
A_{\lambda'} \ar@{^{(}->}[r] \ar[d] & \hat{A}_{K} \ar[d]\\
B_{\lambda'} \ar@{^{(}->}[r]& \hat{B}_{K}
}
$$
  
  \noindent avec $A_{\lambda'} \rightarrow B_{\lambda'}$ fid\`element plat et $\hat{B}_{K} = \hat{A}_{K} \otimes_{A_{\lambda'}} B_{\lambda'}$: d'apr\`es [Et  5 prop 2] on en d\'eduit 
  
  $$A_{\lambda'} = B_{\lambda' } \cap \hat{A}_{K} = (B_{\lambda'})^G , $$
  
  d'o\`u la proposition (2.3.1). $\square$
 \end{enumerate}

 \vskip 3mm
 \section*{3. Images directes d'isocristaux}
 
\subsection*{3.1 Sections surconvergentes}
 
 On suppose donn\'e un diagramme commutatif tel que (2.2.1). Pour un voisinage strict $W$ (resp. un couple de voisinages stricts $W'\subset W$) de $]X[_{\mathcal{Y}}$ dans $]Y[_{\mathcal{Y}}$ on note $\alpha_{W}$ (resp.$\alpha_{W W'})$ l'immersion ouverte de $W$ dans $]Y[_{\mathcal{Y}}$(resp. de $W'$ dans $W$). Si $\mathcal{A}$ est un faisceau d'anneaux sur $W$ et $E$ un $\mathcal{A}$-module, on pose [B 3,(2.1.1.1)]:\\
 
 \noindent (3.1.1) $\qquad\qquad j^{\dag}_{W}E:=\displaystyle\mathop{\mbox{lim}}_{\longrightarrow\atop{W'\subset W}} \alpha^{}_{W W'^{\ast}} \alpha^{\ast}_{W W'} E$ ,\\
 \noindent la limite \'etant prise sur les voisinages $W'\subset W$.\\
 De m\^eme [B 3,(2.1.1.3)] :\\
 
 \noindent (3.1.2) $\qquad\qquad j^{\dag}_{Y} E:= \alpha^{}_{W ^{\ast}} j^{\dag}_{W} E$ .\\
 
 \noindent(3.1.3) Si $V$ est un voisinage strict de $]S[_{\mathcal{T}}$ dans $]T[_{\mathcal{T}}$, alors $W = \overline{h}^{\ -1}_{K} (V)\ \cap\ ]Y[_{\mathcal{Y}}\ =\overline{h}^{\ -1}_{Y} (V)\ $ est un voisinage strict de $]X[_{\mathcal{Y}}$ dans $]Y[_{\mathcal{Y}}$ [B 3, (1.2.7)] et on note $h_{V}$ la restriction de $\overline{h}^{}_{Y} $ \`a $W$, et $R^{i} \overline{h}_{K^{\ast}} j^{\dag}_{Y} E:= R^{i}Ê\overline{h}_{Y^{\ast}} j^{\dag}_{Y} E.$\\
 
 \vskip 3mm
 \noindent \textbf{Proposition (3.1.4)}.
\textit {Avec les hypoth\`eses et notations de (3.1.3) supposons que}  $\overline{h}_{Y} :\  ]Y[_{\mathcal{Y}}{\longrightarrow}]T[_{\mathcal{T}} $  \textit{soit quasi-compact et quasi-s\'epar\'e; soit $E$ un faisceau ab\'elien sur $W$.}
 \begin{enumerate}
 \item[(a)] \textit{Supposons que $\overline{h}^{\ -1}(T) = Y$ et $\overline{h}^{\ -1}(S) = X$; alors, pour tout entier $ i \geqslant 0 $, on a des isomorphismes canoniques}
   \begin{enumerate}
   \item [(3.1.4.1)] $ \qquad\qquad R^{i}h_{V^{\ast}}(j^{\dag}_{W}E ) \overset{\sim}{\longrightarrow} j^{\dag}_{V}R^{i}h_{V^{\ast}} (E ).  $
   \item [(3.1.4.2)] $ \qquad\qquad R^{i}\overline{h}_{K^{\ast}}(j^{\dag}_{Y}E ) \overset{\sim}{\longrightarrow} j^{\dag}_{T}R^{i}h_{V^{\ast}} (E ).  $\\
   
   \textit{Si de plus $\overline{h}$ est une immersion ferm\'ee, alors}
  \item [(3.1.4.3)] $ \qquad\qquad R^{i}h_{V^{\ast}}(j^{\dag}_{W}E )\ =\ 0$ \textit{pour $ i \geqslant 1 $}\\
  
   \textit{et le morphisme canonique}\\
   \item [(3.1.4.4)] $ \qquad\qquad \overline{h}_{K}^{\ast}\overline{h}_{K^{\ast}}j^{\dag}_{Y}E  \overset{\sim}{\longrightarrow} j^{\dag}_{Y} E $\\
      \item[]\textit{est un isomorphisme.}
    \end{enumerate}
  \item[(b)] \textit{Si l'on ne suppose plus que $\overline{h}^{\ -1}(T) = Y$ et $\overline{h}^{\ -1}(S) = X$, alors, pour tout entier $ i \geqslant 0 $, on a des isomorphismes canoniques}
   \begin{enumerate}
   \item [(3.1.4.5)] $ \qquad\qquad R^{i}h_{V^{\ast}}(j^{\dag}_{W}E ) \overset{\sim}{\longrightarrow} j^{\dag}_{V}R^{i}h_{V^{\ast}} (j^{\dag}_{W}E ).$
   \item [(3.1.4.6)] $ \qquad\qquad R^{i}\overline{h}_{K^{\ast}}(j^{\dag}_{Y}E ) \overset{\sim}{\longrightarrow} j^{\dag}_{T}R^{i}h_{V^{\ast}} (j^{\dag}_{W}E ) .$
    \end{enumerate}

 \end{enumerate}
  \vskip15mm
 \noindent\textit{D\'emonstration}
 \begin{enumerate}
   \item[(a)] Les deux foncteurs
 $$\mathcal{F}: E \longmapsto h_{V^{\ast}}j^{\dag}_{W}E$$ et $$\mathcal{G}: E \longmapsto j^{\dag}_{V}h_{V^{\ast}}E$$ de la cat\'egorie $\mathcal{C}$ des faisceaux ab\'eliens sur $W$ dans la cat\'egorie  des faisceaux ab\'eliens sur $V$ sont exacts \`a gauche [B 3,(1.1.3)(iii)]. Comme la cat\'egorie ab\'elienne $C$ admet suffisamment d'injectifs, les foncteurs d\'eriv\'es droits $R^{i}\mathcal{F}$ et  $R^{i}\mathcal{G}$ existent et  $ (R^{i} \mathcal{F})_{i}$ , $(R^{i} \mathcal{G})_{i}$ sont des $ \delta$-foncteurs universels:   puisque $j^{\dag}_{W}$ et $j^{\dag}_{V}$ sont exacts [loc. cit.],  et que $\alpha_{W}:  W \hooklongrightarrow ]Y[_{\mathcal{Y}}$  (resp. $\alpha_{V} : V \hooklongrightarrow ]T[_{\mathcal{T}})$   est exact sur la cat\'egorie des  $j^{\dag}_{W}\mathbb{Z}$-modules (resp. des  $j^{\dag}_{V}\mathbb{Z}$-modules) [B 3, d\'em. de (2.1.3)], on est ramen\'e pour le (a) \`a prouver que $\mathcal{F}= \mathcal{G}$.\\
 
 Comme $\overline{h}_{Y}$ est quasi-compact et quasi-s\'epar\'e, il en est de m\^eme par changement de base pour $h_{V}$, donc $h_{V^{\ast}}$ commute aux limites inductives filtrantes: de plus les hypoth\`eses entra\^inent que si $V'$ d\'ecrit un syst\`eme fondamental de voisinages stricts de $]S[_{\mathcal{T}}$ dans $]T[_{\mathcal{T}}$, alors $h^{-1}_{V}(V')= W'$ d\'ecrit un syst\`eme fondamental de voisinages stricts de $]X[_{\mathcal{Y}}$ dans $]Y[_{\mathcal{Y}}$ [(2.2.2.2)]; d'o\`u
   \begin{eqnarray*}
   h_{V^{\ast}}(j^{\dag}_{W}E )&=& \displaystyle\mathop{\mbox{lim}}_{\longrightarrow\atop{W'\subset W}} h_{V^{\ast}} \alpha^{}_{WW'^{\ast}} \alpha^{-1}_{W W'} E\\
     &=&\displaystyle\mathop{\mbox{lim}}_{\longrightarrow\atop{V' \subset V}} \alpha^{}_{VV'^{\ast}} h_{V'^{\ast}} \alpha^{-1}_{W W'} E\\
     &=& \displaystyle\mathop{\mbox{lim}}_{\longrightarrow\atop{V'\subset V}} \alpha^{}_{VV'^{\ast}}\alpha^{-1}_{VV'} h_{V^{\ast}}E\\
     &=&  j^{\dag}_{V}h_{V^{\ast}} (E ), 
      \end{eqnarray*}
ce qui prouve (3.1.4.1) et (3.1.4.2).\\

Si $ \overline{h}$ est une immersion ferm\'ee, alors $ \overline {h}_{K }$ en est une aussi [B 3, (0.2.4)(iv)], de m\^eme que $h_{V}$ par changement de base: en particulier $h_{V}$ est quasi-compact et quasi-s\'epar\'e. Ainsi (3.1.4.3) r\'esulte de (3.1.4.2) et (1.3.3). Pour (3.1.4.4) on utilise la suite d'isomorphismes
 \begin{eqnarray*}
   \overline{h}_{K}^{\ast} \overline{h}_{K^{\ast}}j^{\dag}_{Y}E & \overset{\sim}{\longrightarrow}&\overline{h}_{K}^{\ast}j^{\dag}_{T} h_{V^{\ast}}E \qquad\qquad (3.1.4.2)    \\
    & \overset{\sim}{\longleftarrow} & j_{Y}^{\dag} h^{\ast}_{V} h_{V^{\ast}}(E)\  \quad\qquad [B 3, (2.1.4.8)]
     \\
     & \overset{\sim}{\longrightarrow} & j_{Y}^{\dag}(E) \qquad\qquad\qquad (1.3.3). 
   \end{eqnarray*}\\
   
     \item[(b)]L'ouvert W = $\overline{h}_{Y}^{\ -1}(V)$ est un voisinage strict de $]X_{1}[_{\mathcal{Y}}$ (donc de $]X[_{\mathcal{Y}}$) dans $]Y[_{\mathcal{Y}}$; la d\'efinition de $j^{\dag}_{W}(E)$ fait intervenir une limite inductive sur les voisinages stricts $W^{'}$ de $]X[_{\mathcal{Y}}$ dans $]Y[_{\mathcal{Y}}$: si cette fois la limite inductive est prise sur les voisinages stricts $W^{'}_{1}$de $]X_{1}[_{\mathcal{Y}}$ dans $]Y[_{\mathcal{Y}}$ nous noterons $j^{\dag}_{W_{1}}(E)$ le r\'esultat, et on a [B 3,(2.1.7)]
        $$j^{\dag}_{W}(E) = j^{\dag}_{W} \circ j^{\dag}_{W_{1}}(E) = j^{\dag}_{W_{1}}\circ j^{\dag}_{W}(E).$$
 D'o\`u , en appliquant (3.1.4.1):
  \begin{eqnarray*}
   R^{i}h_{V^{\ast}}(j^{\dag}_{W}E )&=& R^{i}h_{V^{\ast}}j^{\dag}_{W_{1}} (j^{\dag}_{W}E )\\
     &=&j^{\dag}_{V} R^{i}h_{V^{\ast}}(j^{\dag}_{W}E);
     \end{eqnarray*} 
 de m\^eme pour (3.1.4.6). $\square$
 \end{enumerate}

\subsection*{3.2  D\'efinition des images directes}
   
  \noindent \textbf{(3.2.1)} On suppose fix\'e un diagramme commutatif tel que (2.2.1) et on fait l'hypoth\`ese suppl\'ementaire que $\overline{h}$ (resp $\rho$) est lisse sur un voisinage de $X$ dans $\mathcal{Y}$ (resp de $S$ dans $\mathcal{T}$)\\ 
$$\xymatrix{
X \ar@{^{(}->}[r]^{j_{Y}} \ar[d]_{f} & Y \ar@{^{(}->}[r]^{i_{Y}} \ar[d]^{\overline{f}} & \mathcal{Y} \ar[d] ^{\overline{h}}\\
S  \ar@{^{(}->}[r]_{j_{T}} & T \ar@{^{(}->}[r]_{i_{T}} & \mathcal{T} \ar[r]^{\rho} & \mathcal{W}.\\
}
$$

Soient $E\in Isoc^{\dag} ((X,Y)/ \mathcal{W})$, $W$ un voisinage strict de $]X[_{\mathcal{Y}}$ dans $]Y[_{\mathcal{Y}}$ et $E_{W}$ un $\mathcal{O}_{W}$-module coh\'erent tel que $j_{Y}^{\dag}E_{W}=: E_{\mathcal{Y}}$ soit une r\'ealisation de $E$ [B 3, (2.3.2)]; on notera $Isoc^{\dag} ((X,Y)/ \mathcal{W})_{plat}$ la sous-cat\'egorie pleine de $Isoc^{\dag} ((X,Y)/ \mathcal{W})$ form\'ee des $E$ tels qu'il existe un $E_{W}$ qui soit un $\mathcal{O}_{W}$-module coh\'erent et plat: lorsque $\mathcal{W}= Spf\mathcal{V}$ cette condition est automatiquement v\'erifi\'ee, i.e. on a  [B 3, (2.2.3)(ii)]:\
$Isoc^{\dag} ((X,Y)/ Spf\mathcal{V})_{plat}=Isoc^{\dag} ((X,Y)/ Spf\mathcal{V})=: Isoc^{\dag} ((X,Y)/ K)$. 
Lorsque $Y$ est propre sur $\mathcal{W}$ on notera $Isoc^{\dag} ((X,Y)/ \mathcal{W})_{plat}=Isoc^{\dag} (X/ \mathcal{W})_{plat}$ [B 3, (2.3.6)]: si de plus $\mathcal{W}= Spf\mathcal{V}$ on a
 $Isoc^{\dag} (X/ Spf\mathcal{V})_{plat}=Isoc^{\dag} (X/ Spf\mathcal{V})=: Isoc^{\dag} (X/ K)$.\\
 
Pour $E\in Isoc^{\dag} ((X,Y)/ \mathcal{W})$ Berthelot a d\'efini dans [B 5,(3.1.11)] les images directes en cohomologie rigide (cf. aussi [LS, (7.4)] et [C-T, 10]) par la formule
 \begin{eqnarray*}
 (3.2.1.1)\qquad\qquad\mathbb{R}\overline{f}_{rig^{\ast}}((X,Y)/ \mathcal{T};E)&:=& \mathbb{R}\overline{h}_{K^{\ast}}(j^{\dag}_{Y}E_{W}\otimes_{\mathcal{O}_{]Y[_{\mathcal{Y}}}}\Omega^{^{\bullet}}_{{]Y[_{\mathcal{Y}}}/ \mathcal{T}_{K}} ),\\
     &:=& \mathbb{R}\overline{h}_{Y^{\ast}}(j^{\dag}_{Y}E_{W}\otimes_{\mathcal{O}_{]Y[_{\mathcal{Y}}}}\Omega^{^{\bullet}}_{{]Y[_{\mathcal{Y}}}/ ]T[_{\mathcal{T}}} );
  \end{eqnarray*}
  \noindent et la cohomologie de ces complexes est ind\'ependante du $\mathcal{Y}$ choisi [B 5, (3.1.2)] [LS, 7.4.2].\\
 Lorsque $X = Y$, alors $\overline{f}:X \rightarrow T$ et on obtient la cohomologie convergente:\\
 
  \noindent(3.2.1.2)$\qquad\qquad\mathbb{R}\overline{f}_{conv^{\ast}}(X,Y/ \mathcal{T};E):= \mathbb{R}\overline{f}_{rig^{\ast}}((X,X)/ \mathcal{T};E).$\\
  
  \noindent\textbf{(3.2.2)} Sous les hypoth\`eses (3.2.1) supposons de plus $\overline{h}$ propre: ainsi $Y$ est une compactification $\overline{X}=Y$ de $X$ au-dessus de $\mathcal{T}$. Berthelot d\'efinit alors $\mathbb{R}f_{rig^{\ast}}(X,/ \mathcal{T};E)$ par la formule [B 5, (3.2.3)] ( cf. aussi [LS, 8.2])\\
   
\noindent$ (3.2.2.1)\qquad\qquad\mathbb{R}f_{rig^{\ast}}(X/ \mathcal{T};E):=\mathbb{R}\overline{f}_{rig^{\ast}}((X,\overline{X})/ \mathcal{T};E);$\\

  \noindent et la cohomologie de ce complexe est ind\'ependante du $\overline{X}$ choisi [B 5,(3.2.2)] [LS, 8.2.1] [CT, 10.5.3].\\

     \subsection*{3.3 Changement de base}
  
  \noindent\textbf{(3.3.1)} Avec les notations de (2.2) on consid\`ere un parall\'el\'epip\`ede commutatif\\  
  
  $$
\begin{array}{c}
\xymatrix{
& X \ar@{.>}[dd]^(.3){f} |\hole \  \ar@{^{(}->}[rr]^{j_{Y}} & & Y\ar@{.>}[dd]^(.3){\overline{f}} |\hole \ \ar@{^{(}->}[rr]^{i_{Y}} & & \mathcal{Y} \ar[dd]^{\overline{h}} \\
X^{'} \ \ar@{^{(}->}[rr]^(.8){j_{Y'}} \ar[dd]_{f^{'}} \ar[ur]^{\varphi'} & & Y' \ \ar@{^{(}->}[rr]^(.8){i_{Y'}} \ar[dd]_(.7){\overline{f}'} \ar[ur]^{{\overline{\varphi}}'} & & \mathcal{Y}' \ar[dd]_(.7){\overline{h}'} \ar[ur]^{\overline{g}'} \\
& S\ \ar@{^{(}.>}[rr]^(.7){j_{T}} |\hole & &T\ \ar@{^{(}.>}[rr]^(.8){i_{T}}  |\hole & & \mathcal{T}\ar[rr]^{\rho} && \mathcal{W}\\
S^{'} \  \ar@{^{(}->}[rr]_{j_{T'}} \ar@{.>}[ur]^{\varphi}& & T' \ \ar@{^{(}->}[rr]_{i_{T'}} \ar@{.>}[ur]_{\overline{\varphi}} & & \mathcal{T'} \ar[ur]_{\overline{g}}\ar[rr]_{\rho '}&&\mathcal{W'}\ar[ur]_{\theta} 
}
\end{array}
\leqno{(3.3.1.1)}
  $$
  \noindent dans lequel le cube de gauche est form\'e de $k$-sch\'emas s\'epar\'es de type fini, les morphismes $\overline{g}$, $\overline{h}$, $\overline{g}'$, $\overline{h}'$, $\rho $, $\rho '  $, $\theta$ sont des morphismes de $\mathcal{V}$-sch\'emas formels ($\mathcal{V}$-sch\'emas formels que l'on supposera s\'epar\'es, plats et de type fini), les $j$ (resp. les $i$) sont des immersions ouvertes (resp. ferm\'ees). On suppose $\theta$ lisse et $X'=X\times_{S}S'$, $Y'=Y\times_{T}T'$, et $\mathcal{Y}'=\mathcal{Y}\times_{\mathcal{T}}\mathcal{T'}$.\\
  On v\'erifie alors facilement que l'on a\\
  
\noindent (3.3.1.2)
 $ \qquad\qquad  ]X'[_{\mathcal{Y'}}=]X[_{\mathcal{Y}}\times_{]S[_{\mathcal{T}}}]S'[_{\mathcal{T'}},  
     \qquad\qquad  ]Y'[_{\mathcal{Y'}}=]Y[_{\mathcal{Y}}\times_{]T[_{\mathcal{T}}}]T'[_{\mathcal{T'}}.$\\
     
       Soient $V$ un voisinage strict de $]S[_{T}$ dans $]T[_{T}$ et $V'$ un voisinage strict de $]S'[_{T'}$ dans $]T'[_{T'}$ tel que $V'\subset{\overline{g}_{K}^{\ -1}}(V)\cap]T'[_{\mathcal{T'}}$; alors $W:={\overline{h}_{K}^{\ -1}}(V)\cap]Y[_{\mathcal{Y}}$ est un voisinage strict de $]X[_{Y}$ dans $]Y[_{Y}$ et $W':={\overline{h}_{K}^{ \prime-1}}(V')\cap]Y'[_{\mathcal{Y'}}$ est un voisinage strict de $]X'[_{Y'}$ dans $]Y'[_{Y'}$ tel que $W'\subset{\overline{g}_{K}^{  \prime -1}}(W)\cap]Y'[_{\mathcal{Y'}}$.\\
 Notons \\
     
          $h_{V}: W\rightarrow V, \ g_{V}: V'\rightarrow V, \ h'_{V'}: W'\rightarrow V', \ \mbox{et} \ g'_{W}: W'\rightarrow W $\\
          
\noindent les morphismes induits respectivement par $\overline{h}_{K},\ \overline{g}_{K},\ \overline{h}'_{K},\  \mbox{et}\  \overline {g}'_{K}.$\\
  Ainsi on dispose d'un cube commutatif
  
  $$
\begin{array}{c}
\xymatrix{
& W \ \ar@{.>} [dd]^(.3){h_{V}} |\hole  \ar@{^{(}->} [rr] & & {\mathcal{Y}_{K}} \ar [dd]^{\overline{h}_{K}} \\
W' \ \ar@{^{(}->} [rr] \ar [dd]_{h'_{V'}} \ar[ur]^{g'_{W}} & & {\mathcal{Y}'_{K}} \ar[dd]_(.7){\overline{h}'_{K}} \ar [ur]^{{\overline{g}'_{K}}} \\
& V \ \ar@{^{(}.>}[rr] |\hole & &{\mathcal{T}_{K}}\\
V'  \ \ar@{^{(}->}[rr] \ar@{.>}[ur]^{g_{V}}& & \mathcal{T}'_{K}  \ar[ur]_{\overline{g}_{K}} 
}
\end{array}
\leqno{(3.3.1.3)}
  $$
  
  \noindent dans lequel  $W'=W\times_{V}V'$ et $\mathcal{Y}'_{K}=\mathcal{Y}_{K}\times_{\mathcal{T}_{K}}\mathcal{T'}_{K}$. \\
  
  \noindent\textbf{Lemme (3.3.1.4)}.
   \textit{Avec les hypoth\`eses et notations de (3.3.1) supposons de plus que $\overline{g}$ soit plat (resp. que $\overline{g}$ soit lisse sur un voisinage de $S'$ dans $\mathcal{T'}$). Alors il existe un voisinage strict $V'$ de $]S'[_{\mathcal{T'}}$ dans $]T'[_{\mathcal{T'}}$ tel que $g_{V}$ soit plat (resp. que $g_{V}$ soit lisse).}\\
   
   \noindent\textit{D\'emonstration}. Dans le cas plat  c'est clair puisque $\overline{g}_{K}$est plat; dans le cas lisse, c'est le lemme (2.2.1) de [B 3]. $\square$\\
   
   \noindent\textbf{D\'efinition(3.3.1.5)}
   \textit{Soit} $\mathcal{E}$   \textit{un} $j_{T}^{\dag}\mathcal{O}_{V}$-\textit{module; son image inverse surconvergente est d\'efinie par la formule}
   
   $$(\varphi, \overline{\varphi}, \overline{g})^{\dag}(\mathcal{E}):= j^{\dag}_{T'}(\overline{g}_{K}^{\ast}\mathcal{E});$$
   \textit{lorsque} $\mathcal{E}$ \textit{est une r\'ealisation d'un isocristal E $\in Isoc^{\dag}((S,T)/\mathcal{W})$ on \'ecrira aussi}
   $$(\varphi, \overline{\varphi})^{\ast}(\mathcal{E})=(\varphi, \overline{\varphi}, \overline{g})^{\dag}(\mathcal{E})= j^{\dag}_{T'}(\overline{g}_{K}^{\ast}\mathcal{E}).$$\\

   \noindent\textbf{Th\'eor\`eme (3.3.2)}
   \textit{Sous les hypoth\`eses (3.3.1) supposons que $\overline{h}^{\ -1}(T)=Y, \overline{h}^{\ -1}(S)=X$ et $h_{V}$ propre (cette derni\`ere hypoth\`ese est v\'erifi\'ee si $\overline h$ est propre).}
\begin{enumerate}
      \item[(3.3.2.1)] \textit{
      Soit $E_{W}$ un $\mathcal{O}_{W}$-module coh\'erent. Alors, pour tout entier $i\geqslant 0$, on a:
      				    }
         \begin{enumerate}
         \item[(1)]$R^{i}\overline{h}_{K^{\ast}}(j^{\dag}_{Y}E_{W} )$ \textit{est un $j_{T}^{\dag}\mathcal{O}_{]T[_{\mathcal{T}}}$-module coh\'erent et on a un isomorphisme}
         
         $ \qquad\qquad R^{i}\overline{h}_{K^{\ast}}(j^{\dag}_{Y}E_{W} ) \overset{\sim}{\longrightarrow} j^{\dag}_{T}R^{i}h_{V^{\ast}} (E_{W} ).  $
         
         \item[(2)]
   		\begin{enumerate}
   			\item[(i)]\textit{On a des isomorphismes de changement de base au sens surconvergent}				\begin{eqnarray*}
				(\varphi, \overline{\varphi}, \overline{g})^{\dag}(R^{i}\overline{h}_{K^{\ast}}j^{\dag}_{Y}E_{W} )&\simeq & R^{i}\overline{h}'_{K^{\ast}}(\overline{g}'^{\ast}_{K}j^{\dag}_{Y}E_{W} )\\
		&\simeq &	 R^{i}\overline{h}'_{K^{\ast}}(j^{\dag}_{Y'}g'^{\ast}_{W}E_{W} )\\
		&\simeq & j^{\dag}_{T'}R^{i}\overline{h}'_{V'^{\ast}}(g'^{\ast}_{W}E_{W} ).
				\end{eqnarray*}   			
			\item[(ii)]\textit{ Si de plus $\overline{\varphi}^{\ -1}(S)= S'$, les isomorphismes pr\'ec\'edents deviennent}
				\begin{eqnarray*}
				\overline{g}_{K}^{\ast}R^{i}\overline{h}_{K^{\ast}}j^{\dag}_{Y}E_{W} &\simeq & R^{i}\overline{h}'_{K^{\ast}}(j^{\dag}_{Y'}g'^{\ast}_{W}E_{W} )\\		
						&\simeq & j^{\dag}_{T'}R^{i}\overline{h}'_{V'^{\ast}}(g'^{\ast}_{W}E_{W} ),
				\end{eqnarray*} 
				\textit{et l'on a des isomorphismes
				\begin{eqnarray*}
				\overline{g}_{V}^{\ast}R^{i}h_{V^{\ast}}j^{\dag}_{W}E_{W} &\simeq & R^{i}h'_{V'^{\ast}}(g'^{\ast}_{W}j_{W}^{\dag}E_{W} )\\		
						&\simeq & R^{i}h'_{V'^{\ast}}j_{W'}^{\dag}g'^{\ast}_{W}E_{W} .
				\end{eqnarray*} 
				}
   		  \end{enumerate}
   	     \end{enumerate}
      \item[(3.3.2.2)]
      \textit{
      Soit $E_{W}^{^{\bullet}}$ un complexe born\'e de $\mathcal{O}_{V}$-modules plats, \`a composantes des $\mathcal{O}_{W}$-module coh\'erents et 
      $$  E_{W'}^{^{\bullet}}=E_{W}^{^{\bullet}} \otimes_{\mathcal{O}_{V}} \mathcal{O}_{V'} = g_{W}^{\prime \ast}E_{W}^{^{\bullet}}.$$ 
      Alors, pour tout entier $i\geqslant 0$, on a:
      				    }
         \begin{enumerate}
         \item[(1)]$R^{i}\overline{h}_{K^{\ast}}(j^{\dag}_{Y}E_{W}^{^{\bullet}} )$ \textit{est un $j_{T}^{\dag}\mathcal{O}_{]T[_{\mathcal{T}}}$-module coh\'erent et on a un isomorphisme}
         
         $ \qquad\qquad R^{i}\overline{h}_{K^{\ast}}(j^{\dag}_{Y}E_{W}^{^{\bullet}} ) \overset{\sim}{\longrightarrow} j^{\dag}_{T}R^{i}h_{V^{\ast}} E_{W}^{^{\bullet}}.  $   
               
         \item[(2)] \textit {Supposons de plus $g_{V}$ plat, alors}
   		\begin{enumerate}
   			\item[(i)]\textit{On a des isomorphismes de changement de base au sens surconvergent}
				\begin{eqnarray*}
				(\varphi, \overline{\varphi}, \overline{g})^{\dag}(R^{i}\overline{h}_{K^{\ast}}j^{\dag}_{Y}E_{W}^{^{\bullet}} )&\simeq &  R^{i}\overline{h}'_{K^{\ast}}(j^{\dag}_{Y'}g'^{\ast}_{W}E_{W}^{^{\bullet}} )\\
				&\simeq & j^{\dag}_{T'}R^{i}\overline{h}'_{V'^{\ast}}(E_{W'}^{^{\bullet}} ).
				\end{eqnarray*}   			
			\item[(ii)]\textit{ Si de plus $\overline{\varphi}^{\ -1}(S)= S'$, les isomorphismes pr\'ec\'edents deviennent}
				\begin{eqnarray*}
				\overline{g}_{K}^{\ast}R^{i}\overline{h}_{K^{\ast}}j^{\dag}_{Y}E_{W}^{^{\bullet}} &\simeq & R^{i}\overline{h}'_{K^{\ast}}(j^{\dag}_{Y'}E_{W'}^{^{\bullet}})\\		
						&\simeq & j^{\dag}_{T'}R^{i}h'_{V'^{\ast}}(E_{W'}^{^{\bullet}} ),
				\end{eqnarray*} 
				\textit{et l'on a des isomorphismes				
				\begin{eqnarray*}
				\overline{g}_{V}^{\ast}R^{i}h_{V^{\ast}}j^{\dag}_{W}E_{W}^{^{\bullet}} &\simeq & R^{i}h'_{V'^{\ast}}g'^{\ast}_{W}j_{W}^{\dag}E_{W}^{^{\bullet}} \\		
						&\simeq & R^{i}h'_{V'^{\ast}}j_{W'}^{\dag}g'^{\ast}_{W}E_{W}^{^{\bullet}} .
				\end{eqnarray*} 
				}
   		  \end{enumerate}
   	     \end{enumerate}
 \end{enumerate}
 
 \noindent\textit{D\'emonstration}\\
 \textit{Le (1) de (3.3.2.1)} r\'esulte de (3.1.4.2) et (1.2.1).\\
 \textit{Pour le (i) de (3.3.2.1)(2)} on consid\`ere le diagramme commutatif
 
   $$
\begin{array}{c}
\xymatrix{
S\  \ar@{^{(}->}[r]^{j_{T}} \ar @{} [dr] |{\square} & T\\
S'_{1} \ \ar@{^{(}->} [r] _{j_{1T'}} \ar [u] ^{\varphi _{1}} & T' \ar[u] _{\overline{\varphi}} \ar@{=}[d]\\
S' \ \ar@{^{(}->} [r] _{j_{T'}} \ar[u] ^{j} & T'
		} 
\end{array}
$$

\noindent dans lequel le carr\'e du haut est cart\'esien et $\varphi= \varphi_{1}\circ j$. On a alors une suite d'isomorphismes

$$ \begin{array}{lll@{\qquad}l}
(\varphi, \overline{\varphi}, \overline{g})^{\ \dag}(R^{i}\overline{h}_{K^{\ast}}j^{\dag}_{Y}E_{W} ) &\simeq & j_{T'}^{\dag} \overline{g}_{K}^{ \ast}R^{i}\overline{h}_{K^{\ast}}j^{\dag}_{Y}E_{W} & [(3.3.1.5)]\\
		&\overset{\sim}{\rightarrow}  & j^{\dag}_{T'} j^{\dag}_{1T'}g^{ \ast}_{V}R^{i}h_{V^{\ast}}E_{W} & [\mbox{B}\ 3,(2.1.4.8)]\\
		&\overset{\sim}{\rightarrow}  & j^{\dag}_{T'} g_{V}^{\ast} R^{i}\overline{h}_{V^{\ast}}E_{W} & [\mbox{B}\ 3, (2.1.7)]\\
		& \overset{\sim}{\rightarrow}  & j^{\dag}_{T'}  R^{i} {h'}_{{V'}^{ \ast}}{g'}_{W}^{\ast}E_{W} & [\mbox{Th\'eo} (1.2.1)]\\
		& \overset{\sim}{\leftarrow}  & R^{i}\overline{h}'_{K^{\ast}}j^{\dag}_{Y'}{g'}_{W}^{\ast}E_{W} &  [(3.1.4.2)]\\
		& \overset{\sim}{\leftarrow}  & R^{i}\overline{h}_{K^{\ast}} {\overline{g}'}_{K}^{\ast} j^{\dag}_{Y}E_{W} & [\mbox{B}\ 3,(2.1.4.8)].
\end{array}$$

\textit{Pour le (ii) de (3.3.2.1)(2)} il suffit de remarquer que les hypoth\`eses impliquent que
$(\varphi, \overline{\varphi}, \overline{g})^{\ \dag}(\mathcal{E})=\overline{g}_K^{\ast}(\mathcal{E})$ pour tout faisceau ab\'elien $\mathcal{E}$ sur $\mathcal{T}_K$; la derni\`ere assertion r\'esulte de (3.1.4.1), [B 3,(2.1.4.7)] et (1.2.1) comme ci-dessus.
\\
\\
\textit{Pour (3.3.2.2)} on proc\`ede de m\^eme en utilisant cette fois le (1.2.2) du th\'eor\`eme (1.2). $\square$ \\
\\
\textit{Remarque (3.3.3)} En fait, dans le (3.3.2.1) (2) (ii) du th\'eor\`eme pr\'ec\'edent, si l'on ne suppose plus l'existence de $\overline{g}$, mais que l'on suppose toujours l'existence du carr\'e cart\'esien
$$
\begin{array}{c}
\xymatrix{
W'\ \ar@{->}[r]^{g'_W} \ar@{->}[d]_{h'_{V'}} & W \ar@{->}[d]^{h_V}\\
V'\ \ar@{->}[r]_{g_V} &V,
		}
\end{array}
$$
on obtient, pour tout $\mathcal{O}_{W}$-module coh\'erent $E_{W}$, un isomorphisme de changement de base
$$\begin{array}{ccc}
 g_{V}^{\ast}R^{i}h_{V^{\ast}}j^{\dag}_{W}E_{W} &\overset{\sim}{\rightarrow} & R^{i}h'_{V'^{\ast}}g'^{\ast}_{W}j_{W}^{\dag}E_{W} \\		
						&\stackrel{\sim}{\rightarrow} & R^{i}h'_{V'^{\ast}}j_{W'}^{\dag}g'^{\ast}_{W}E_{W} .
\end{array} 
\leqno{(3.3.3.1)}$$
De m\^eme, pour le (3.3.2.2)(2)(ii) du th\'eor\`eme, on a un isomorphisme
$$
 g_{V}^{\ast}R^{i}h_{V^{\ast}}{j_{W}}^{\dag}(E_{W}^{^{\bullet}}) \overset{\sim}{\rightarrow} R^{i}h'_{V'^{\ast}}j^{\dag}_{W'}(E_{W'}^{^{\bullet}}).
\leqno{(3.3.3.2)}
$$\\

\subsection*{3.4 Surconvergence des images directes.}

Nous allons consid\'erer dans le prochain th\'eor\`eme l'une des trois situations suivantes.\\

\textbf{(3.4.1) \textit{Dans le premier cas}}, nous consid\'erons un diagramme commutatif satisfaisant aux hypoth\`eses de (2.2)
$$
\begin{array}{c}
\xymatrix{
X\ \ar@{^{(}->}[r]^{j_{Y}} \ar @{}[dr] |{\square} \ar[d]_{f} & Y \ \ar@{^{(}->}[r]^{i_{Y}} \ar[d]^{\overline{f}} \ar @{}[dr] |{\square}& \mathcal{Y} \ar[d] ^{\overline{h}}\\
S \ \ar@{^{(}->}[r]_{j_{T}} & T\ \ar@{^{(}->}[r]_{i_{T}} & \mathcal{T} \ar[r]_{\rho} & \mathcal{W},\\
}
\end{array}
\leqno{(3.4.1.1)}
$$
et un diagramme commutatif
$$
\begin{array}{c}
\xymatrix{
S\ \ar@{^{(}->}[r]^{j_{T}} & T\ \ar@{^{(}->}[r]^{i_{T}}  & \mathcal{T} \ar[r] ^\rho&\mathcal{W}\\
S'\  \ar@{^{(}->}[r]_{j_{T'}} \ar@{->}[u]^{\varphi} & T'\ \ar@{^{(}->}[r]_{i_{T'}} \ar@{->}[u]^{\overline{\varphi}}& \mathcal{T}' \ar[r]_{\rho'} \ar@{->}[u]_{\overline{g}}& \mathcal{W}' \ar@{->}[u]_{\theta}\\
}
\end{array}
\leqno{(3.4.1.2)}
$$
\noindent tel qu'en prenant l'image inverse de (3.4.1.1) par (3.4.1.2) on obtienne un parall\'el\'epip\`ede commutatif tel que (3.3.1.1). On suppose de plus les carr\'es de (3.4.1.1) cart\'esiens ($\overline{h}^{\ -1}(T)=Y\ ,\ \overline{h}^{\ -1}(S)=X$), $ \overline{h}$ propre, $\overline{h}$ lisse sur un voisinage de X dans $\mathcal{Y}$, $\theta$ lisse, $\rho$ (resp $ \rho' $) lisse sur un voisinage de S dans $\mathcal{T}$ (resp de $S'$ dans $\mathcal{T}'$). On suppose \'egalement satisfaite l'une des deux hypoth\`eses suivantes: $\overline{g}$ est lisse sur un voisinage de $S' $ dans $\mathcal{T}'$, ou $\overline{g}$ est plat.\\

\textbf{(3.4.2) \textit{Dans le deuxi\`eme cas}}, nous consid\'erons un diagramme commutatif tel que (3.4.1.1) et un diagramme commutatif
$$
\begin{array}{c}
\xymatrix{
S\ \ar@{^{(}->}[r]^{j_{T}} & T\ \ar@{^{(}->}[r]^{i_{T}}  & \mathcal{T} \ar[r] ^\rho&\mathcal{W}\\
S'\  \ar@{^{(}->}[r]_{j_{T'}} \ar@{->}[u]^{\varphi} & T'\ \ar@{^{(}->}[r]_{i_{T'}} \ar@{->}[u]^{\overline{\varphi}}& \mathcal{T}' \ar[r]_{\rho'} & \mathcal{W}' \ar@{->}[u]_{\theta}\\
}
\end{array}
\leqno{(3.4.2.1)}
$$
satisfaisant aux m\^emes propri\'et\'es que (3.4.1.2) except\'e l'existence de $\overline{g}$, mais en supposant $ \rho\  propre$, et nous noterons
$$\xymatrix{X'=X \times_{S} S'\quad \ar@{^{(}->}[r]^{j_{Y'}}&\quad Y'=Y \times _{T}T'}$$
et
  $$
\begin{array}{c}
\xymatrix{
&X\ \ar@{.>} [dd]^(.3){f} |\hole  \ar@{^{(}->} [rr]^{j_{Y}} & & Y \ar [dd]^(.5){\overline{f}} \\
X' \ \ar@{^{(}->} [rr]^(.7){j_{Y'}} \ar [dd]_{f'} \ar[ur]^{\varphi'} & & Y' \ar[dd]_(.7){\overline{f}'} \ar [ur]^{\overline{\varphi}'} \\
& S \ \ar@{^{(}.>}[rr]^(.7){j_{T}} |\hole & &{T}\\
S'  \ \ar@{^{(}->}[rr]_{j_{T'}} \ar@{.>}[ur]_{\varphi}& & {T}' \ar[ur]_{\overline{\varphi}} 
}
\end{array}
\leqno{(3.4.2.2)}
  $$
l'image inverse par ($\varphi$, $\overline{\varphi}$) de (3.4.1.1).\\

\textbf{(3.4.3) \textit{Dans le troisi\`eme cas}}, qui g\'en\'eralise le premier, nous consid\'erons des diagrammes commutatifs tels que (3.4.1.1) et (3.4.2.1) mais sans supposer $\rho$ propre. Par contre nous supposons de plus l'existence d'un diagramme commutatif
$$
\begin{array}{c}
\xymatrix{
X' \ \ar@{^{(}->}[r]^{j_{Y'}} \ar @{}[dr] |{\square} \ar[d]_{f'} & Y' \ \ar@{^{(}->}[r]^{i_{Y'}} \ar[d]^{\overline{f}'} \ar @{}[dr] |{\square}& \mathcal{Y'} \ar[d] ^{\overline{h}'}\\
S' \ \ar@{^{(}->}[r]_{j_{T'}} & T'\ \ar@{^{(}->}[r]_{i_{T'}} & \mathcal{T}' \ar[r]^{\rho'} & \mathcal{W}'\\
}
\end{array}
\leqno{(3.4.3.1)}
$$
satisfaisant aux m\^emes hypoth\`eses que (3.4.1.1), et dans lequel $X', Y' $ satisfont aux propri\'et\'es de (3.4.2.2).\\

Dans le cas relevable le th\'eor\`eme suivant r\'esout une conjecture de Berthelot [B 2,(4.3)] et g\'en\'eralise le th\'eor\`eme 5 de loc. cit.\\
\\ 
\noindent\textbf{Th\'eor\`eme (3.4.4)}
\textit{Pour tout entier i $\geqslant$ 0, on a:
\begin{enumerate}
   \item[(3.4.4.1)] Sous les hypoth\`eses (3.4.3), on a:
      \begin{enumerate}
         \item[(i)]$\overline{f}$ induit un foncteur\\
         \\
                           		      $R^{i}\overline{f}_{rig \ast}((X,Y)/ \mathcal{T};-):{Isoc}^{\dag}((X,Y)/				\mathcal{W})_{plat} \longrightarrow {Isoc}^{\dag}((S,T)/W).$\\
		\item[(ii)] il existe un morphisme de changement de base\\
		\\
				${(\varphi,\overline{\varphi})}^{\ast}R^{i}\overline{f}_{rig \ast}((X,Y)/				\mathcal{T};E) \longrightarrow R^{i}\overline{f'}_{rig \ast}((X',Y')/				\mathcal{T}';{(\varphi',\overline{\varphi}')}^{\ast}(E))$\\
				\\
					et celui-ci est un isomorphisme dans ${Isoc}^{\dag}((S',T')/ \mathcal{W}')$.
	\end{enumerate}
     \item[(3.4.4.2)] Sous les hypoth\`eses (3.4.2) on a :
       \begin{enumerate}
       	\item[(i)] f induit un foncteur\\
	\\
		$R^{i}f_{rig \ast}(X/\mathcal{T};-):{Isoc}^{\dag}(X/\mathcal{W})_{plat}  					\longrightarrow {Isoc}^{\dag}(S/\mathcal{W}).$\\
		\item[(ii)]L'isomorphisme de changement de base de (3.4.4.1)(ii) existe et devient\\
		\\
		${(\varphi,\overline{\varphi})}^{\ast}R^{i}f_{rig \ast}(X/				\mathcal{T};E) \overset{\sim}{\longrightarrow} R^{i}\overline{f'}_{rig \ast}((X',Y')/				\mathcal{T}';{(\varphi',\overline{\varphi}')}^{\ast}(E)).$\\
	\item[(iii)] Si de plus $\rho'$ est propre, alors l'isomorphisme de (ii) pr\'ec\'edent devient un isomorphisme dans ${Isoc}^{\dag}(S'/\mathcal{W}')$:\\
	\\
	${\varphi}^{\ast}R^{i}f_{rig \ast}(X/\mathcal{T};E) \overset{\sim}	{\longrightarrow} 	R^{i}{f}'_{rig \ast}(X'/\mathcal{T}';{\varphi'}^{\ast}(E)).$\\
	\item[(iv)] Si $S'=T'$, l'isomorphisme du (ii) pr\'ec\'edent devient un isomorphisme dans $ Isoc(S'/ \mathcal{W}'):$\\
	\\
	${\varphi}^{\ast}j_{T}^{\ast}R^{i}f_{rig \ast}(X/\mathcal{T};E) \overset{\sim}		{\longrightarrow} R^{i}{f}'_{conv \ast}(X'/\mathcal{T}'; {\varphi'}^{\ast}(\hat{E}))$\\
	\\
	o\`u $\hat{E}$ $\in$ $Isoc(X/\mathcal{W})$ est l'isocristal convergent associ\'e \`a E $\in$ ${Isoc}^{\dag}(X/\mathcal{W})$ par le foncteur d'oubli ${Isoc}^{\dag}(X/\mathcal{W})\ \rightarrow\ {Isoc}(X/\mathcal{W}).$\\
	En particulier si $S'=T'=S$, on a un isomorphisme \\
	\\
	$j_{T}^{\ast}R^{i}f_{rig \ast}(X/\mathcal{T};E) \overset{\sim}					{\longrightarrow} R^{i}{f}_{conv \ast}(X/\mathcal{T};\hat{E})$\\
	\end{enumerate}
    \item[(3.4.4.3)] Sous les hypoth\`eses (3.4.3) avec $S=T\  et \ S'=T'$ on a :
    	\begin{enumerate}
	    \item[(i)] f induit un foncteur\\
	    \\
	    $R^{i}f_{conv \ast}: Isoc(X/\mathcal{W})_{plat} \longrightarrow \ Isoc(S/\mathcal{W}).$\\
	    \item[(ii)] Il existe un isomorphisme de changement de base dans $Isoc(S'/ \mathcal{W}')$\\
	    \\
	    ${\varphi}^{\ast}R^{i}f_{conv \ast}(X/\mathcal{T};\mathcal{E}) \overset{\sim}				      {\longrightarrow} 	R^{i}{f}'_{conv \ast}(X'/\mathcal{T}';{\varphi'}^{\ast}		(\mathcal{E})).$
	 \end{enumerate}
 \end{enumerate}
	}
	
	\vskip10mm
Avant de donner la preuve du th\'eor\`eme, faisons quelques remarques:\\

\noindent\textit{$\underline{Remarques \ (3.4.4.4)}$}:\\

\begin{enumerate}
   \item[(i)] Sous les hypoth\`eses de (3.4.4.2)(iii), et en supposant de plus que $\varphi 	$ est l'identit\'e de S et $\theta$ l'identit\'e de $\mathcal{W}$, l'isomorphisme de 	changement de base prouve que $R^{i}f_{rig \ast}(X/\mathcal{T};E)$ est ind\'ependant du sch\'ema formel $\mathcal{T}$ dans lequel S est plong\'e (avec 		bien s\^ur $\mathcal{T}$ propre sur $\mathcal{W}$, $\rho$ lisse sur un voisinage 	de S dans $\mathcal{T}$ et $\overline{h}$ v\'erifiant (3.4.4)).\\
         La m\^eme remarque s'applique \`a $R^{i}f_{conv \ast}(X/\mathcal{T};\mathcal	{E}).$\\
         
   \item[(ii)] D'apr\`es [I,(3.3)] les hypoth\`eses de (3.4.4.2)(iii) sont v\'erifi\'ees pour $	\mathcal{W}=Spf \mathcal{V}, S$ affine et lisse sur $k$ et certains morphismes $f$ projectifs et lisses.\\
   
   \item[(iii)] En conjuguant (i) et (ii) nous en d\'eduirons plus loin [Th\'eor\`eme (3.4.8.2)] que les constructions se recollent pour certains morphismes $f$ projectifs lisses et $S$ un $k$-sch\'ema lisse ( plus n\'ecessairement affine).\\
   
   \item[(iv)] Lorsque $\mathcal{W}= Spf\mathcal{V}$ on peut enlever l'indice ''plat'' dans (3.4.4.1)(i), (3.4.4.2)(i) et (3.4.4.3)(i) [cf (3.2.1)].\\

\end{enumerate}
\newpage

\noindent\textit{D\'emonstration de (3.4.4)}.\\

Le (3.4.4.3) est cons\'equence directe de (3.4.4.1).\\

La preuve de (3.4.4.1) et (3.4.4.2) va se faire en six \'etapes que l'on pr\'ecise ici:\\

Dans les quatres premi\`eres \'etapes on va faire la preuve de (3.4.4.1) sous l'hypoth\`ese plus particuli\`ere (3.4.1) (i.e. existence de $\overline{g}$):
\begin{enumerate}
  \item[-]\textit{\'etape $\  \rondI$}: montrer que $R^{i}\overline{f}_{rig \ast}((X,Y)/\mathcal{T};E)$ est un $j_{T}^{\dag}\mathcal{O}_{]T[_{\mathcal{T}}}$-module coh\'erent.
  \item[-]\textit{\'etape $\  \rondII$}: montrer l'isomorphisme de changement de base (3.4.4.1)(ii) 	lorsque $\overline{\varphi}^{-1}(S)=S'$.
  \item[-]\textit{\'etape $\  \rondIII$}: achever la preuve de (3.4.4.1)(i).
  \item[-]\textit{\'etape $\  \rondIV$}: achever la preuve de l'isomorphisme de changement de base (3.4.4.1)(ii).
  
  \item[-]\textit {\'etape $\  \rondV$}: on prouve (3.4.4.1) sous les hypoth\`eses (3.4.3).
  \item[-]\textit {\'etape $\  \rondVI$}: on prouve (3.4.4.2).\\
\end{enumerate}
Pla\c cons-nous d'abord sous les hypoth\`eses (3.4.1).\\
On se donne donc un diagramme tel que (3.3.1.1) avec $\overline{h}^{-1}(T)=Y$, $\overline{h}^{-1}(S)=X$: on note $S'_{1}=\overline{\varphi}^{-1}(S) \ \mbox{et}\ S' \xrightarrow[]{j}S_{1}^{\prime}\xrightarrow[]{\varphi_{1}}S $ la factorisation de $\varphi$ o\`u $j$ est une immersion ouverte; on en d\'eduit un diagramme commutatif \`a carr\'es verticaux cart\'esiens\\
\\
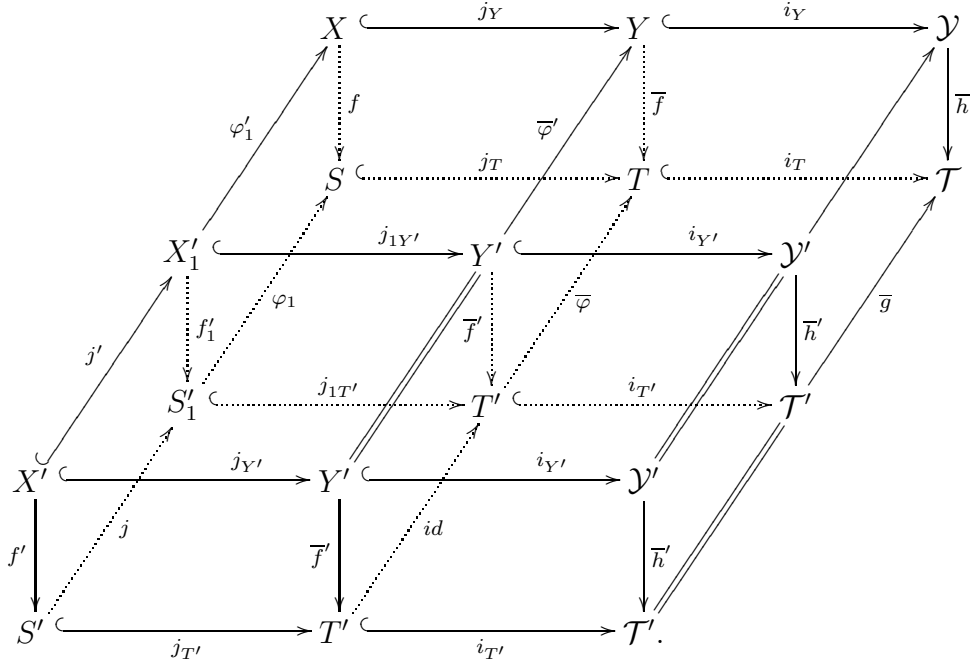
$$
\begin{array}{c}
 \shorthandoff{;:!?}
 \xymatrix@!0 @R=1cm @C=2cm{
&&X\ar@{.>}[dd]^{f}\  \ar @{^{(}->}[rr]^{j_{Y}}&&Y  \ar@{.>}[dd]^{\overline{f}}\  \ar@{^{(}->}[rr]^{i_{Y}}&&\mathcal{Y} \ar[dd]^{\overline{h}}\\
&&&&&&\\
&&S\  \ar@{^{(}.>}[rr]^{j_{T}}  && T\  \ar@{^{(}.>}[rr]^{i_{T}}&& \mathcal{T}\\
&X'_{1}\ar@{.>}[dd]^{f_{1}^{\prime}} \ar[uuur]^{\varphi'_{1}}\  \ar@{^{(}->}[rr]^(.7){j_{1Y'}}&&Y'\ar@{.>}[dd]_{\overline{f}'} \  \ar@{^{(}->}[rr]^(.7){i_{Y'}} \ar[uuur]^{\overline{\varphi}'}&& \mathcal{Y'} \ar[dd]^{\overline{h}'} \ar[uuur]\\
&&&&&&\\
& S'_{1}\  \ar@{^{(}.>}[rr]^{j_{1T'}} \ar@{.>}[uuur]_{\varphi_{1}} &&T' \  \ar@{^{(}.>}[rr] ^{i_{T'}} \ar@{.>}[uuur]_{\overline{\varphi}}&& \mathcal{T'} \ar[uuur]_{\overline{g}} &\\
X' \ar@{^{(}->}[rr]^(.7){j_{Y'}} \ar[dd]_{f'}\   \ar@{^{(}->}[uuur]^{j'}&&Y'\ar[dd]_{\overline{f}'}\  \ar@{^{(}->}[rr]^(.7){i_{Y'}} \ar@{=}[uuur]&& \mathcal{Y'}\ar[dd]^{\overline{h}'} \ar@{=}[uuur]&&\\
&&&&&&\\
S' \  \ar@{^{(}->}[rr]_{j_{T'}}  \ar@{.>}[uuur]_{j}&&T' \ar@{^{(}->}[rr]_{i_{T'}} \  \ar@{.>}[uuur]_{id}&&\mathcal{T'}. \ar@{=}[uuur]&&
}
\end{array}
$$
Pour $E \in Isoc^{\dag}((X,Y)/ \mathcal{W})$ on notera $(\varphi, \overline{\varphi})^{\ast}(E)$ son image inverse par le couple $(\varphi, \overline{\varphi})$ [B 3, (2.3.2)(iv)] pour pr\'eciser la d\'ependance en $\varphi$ et $\overline{\varphi}$: d'autres images inverses seront utilis\'ees, que le contexte pr\'ecisera (cf. d\'ef. (3.3.1.5)).\\

\noindent\textbf{Etape $\  \rondI$.} Avec les notations de (3.2) il existe un voisinage strict  $W$ de $]X[_{\mathcal{Y}}$ dans $]Y[_{\mathcal{Y}}$ tel que $E_{\mathcal{Y}}:=j_{Y}^{\dag}E_{W}$ soit une r\'ealisation de E, et d'apr\`es (2.2.2.2) on peut supposer que $W$ est de la forme $W=\overline{h}_{Y}^{-1}(V)$ pour un voisinage strict $V$ de $]S[_{\mathcal{T}}$ dans $]T[_{\mathcal{T}}$; de plus on a un diagramme commutatif
$$
\begin{array}{c}
\xymatrix{
]X[_{\mathcal{Y}}\ \ar@{^{(}->}[r] \ar[d]_{\overline{h}_{X}} & W\
\ar@{^{(}->}[r]^{\alpha_{W}} \ar[d]_{{h}_{V}} & \mathcal{Y}_{K}  \ar[d]^{\overline{h}_{K}} \\
]S[_{\mathcal{T}}\ \ar@{^{(}->} [r]&V\ \ar@{^{(}->} [r]_{\alpha_{V}} & \mathcal{T}_{K}
} 
\end{array}
$$
dans lequel les carr\'es sont cart\'esiens et o\`u les fl\`eches horizontales sont des immertions ouvertes, ${\overline{h}_{K}}$ est propre et ${\overline{h}_{X}}$ est propre et lisse [(2.2.3.2)(i)]: quitte \`a restreindre $V$ on peut supposer via [(2.2.3.2)(ii)] que ${h}_{V}$ est propre et lisse.\\
Or $R^{i+j}\overline{f}_{rig^{\ast}}((X,Y)/\mathcal{T};E)$ est l'aboutissement d'une suite spectrale de terme $E_{1}^{i,j}$ donn\'e par
$$E_{1}^{i,j}=R^{j}{\overline{h}_{K^{\ast}}}(j_{Y}^{\dag}E_{W}\otimes_{\mathcal{O}_{]Y[_{\mathcal{Y}}}}\Omega^{i}_{]Y[_{\mathcal{Y}}/\mathcal{T}_{K}})$$
avec filtration
\begin{eqnarray*}
Fil^{i}&:= &Fil^{i}(j_{Y}^{\dag}E_{W}\otimes\Omega^{\bullet}_{]Y[_{\mathcal{Y}}/\mathcal{T}_{K}})\\
	&=&j_{Y}^{\dag}E_{W}\otimes\Omega^{\geqslant i}_{]Y[_{\mathcal{Y}}/\mathcal{T}_{K}},
\end{eqnarray*}
et on a une suite d'isomorphismes
$$ \begin{array}{lll@{\ }l}
j_{Y}^{\dag}E_{W}\otimes_{\mathcal{O}_{]Y[_{\mathcal{Y}}}}\Omega^{i}_{]Y[_{\mathcal{Y}}/\mathcal{T}_{K}}&=&(\alpha_{W^{\ast}}j_{W}^{\dag}E)\otimes_{\mathcal{O}_{]Y[_{\mathcal{Y}}}}\Omega^{i}_{]Y[_{\mathcal{Y}}/\mathcal{T}_{K}}& \\
&\simeq&\alpha_{W^{\ast}}(j_{W}^{\dag}E_{W}\otimes_{\mathcal{O}_{W}}\alpha_{W}^{\ast}(\Omega^{i}_{]Y[_{\mathcal{Y}}/\mathcal{T}_{K}}))&\\
&\simeq&\alpha_{W^{\ast}}(j_{W}^{\dag}E_{W}\otimes_{\mathcal{O}_{W}}\Omega^{i}_{W/\mathcal{T}_{K}})&\mbox{car}\ \alpha_{W}\ \mbox{est \'etale}\\
&\simeq&\alpha_{W^{\ast}}(j_{W}^{\dag}E_{W}\otimes_{\mathcal{O}_{W}}\Omega^{i}_{W/V})&\mbox{car}\ \alpha_{V}\  \mbox{est \'etale}\\
&\simeq&\alpha_{W^{\ast}}(j_{W}^{\dag}\mathcal{O}_{W}\otimes_{\mathcal{O}_{W}}E_{W}\otimes_{\mathcal{O}_{W}}\Omega^{i}_{W/V})&[B\ 3,(2.1.3)(ii)]\\
&\simeq&\alpha_{W^{\ast}}j_{W}^{\dag}(E_{W}\otimes_{\mathcal{O}_{W}}\Omega^{i}_{W/V})&[\mbox{loc. cit.}]\\
&=&j_{Y}^{\dag}(E_{W}\otimes_{\mathcal{O}_{W}}\Omega^{i}_{W/V}),&\\
\end{array}$$
o\`u $\Omega^{i}_{W/V}$ est un $\mathcal{O}_{W}$-module coh\'erent et localement libre [(2.2.3.2(i)]. D'apr\`es (3.1.4) on en d\'eduit un isomorphisme
\begin{eqnarray*}
E_{1}^{i,j}&=&R^{j}{\overline{h}_{K^{\ast}}}({j}_{Y}^{\dag}E_{W}\otimes_{\mathcal{O}_{]Y[_{\mathcal{Y}}}}\Omega^{i}_{]Y[_{\mathcal{Y}}/\mathcal{T}_{K}})\\
&\overset{\sim}{\rightarrow}&j_{T}^{\dag}R^{j}{{h}_{V^{\ast}}}(E_{W}\otimes_{\mathcal{O}_{W}}\Omega^{i}_{W/V});
\end{eqnarray*}
or $R^{j}{{h}_{V^{\ast}}}(E_{W}\otimes_{\mathcal{O}_{W}}\Omega^{i}_{W/V})$ est un $\mathcal{O}_{V}$-module coh\'erent d'apr\`es le th\'eor\`eme (1.2), donc $E_{1}^{i,j}$ est un $j_{T}^{\dag}\mathcal{O}_{]T[_{\mathcal{T}}}$-module coh\'erent. Comme la filtration $Fil^{i}$ est finie, il en r\'esulte que l'aboutissement $R^{i+j}\overline{f}_{rig \ast}((X,Y)/\mathcal{T};E)$ est un $j_{T}^{\dag}\mathcal{O}_{]T[_{\mathcal{T}}}$-module coh\'erent. Remarquons que pour prouver cette coh\'erence on aurait pu appliquer le (1) de (3.3.2.2), puisque $E_{W}$ est un $\mathcal{O}_{W}$-module plat.\\
\\
\noindent\textbf{Etape \  \rondII}. On a vu \`a l'\'etape \  \rondI \ que 
\begin{eqnarray*}
R^{i}\overline{f}_{rig \ast}((X,Y)/\mathcal{T};E)&:=&R^{i}{\overline{h}_{K^{\ast}}}({j}_{Y}^{\dag}E_{W}\otimes_{\mathcal{O}_{]Y[_{\mathcal{Y}}}}\Omega^{\bullet}_{]Y[_{\mathcal{Y}}/\mathcal{T}_{K}})\\
&\simeq &R^{i}{\overline{h}_{{K}^{\ast}}}({j}_{Y}^{\dag}(E_{W}\otimes_{\mathcal{O}_{W}}\Omega^{\bullet}_{W/V})),
\end{eqnarray*}
et les $E_{W}\otimes_{\mathcal{O}_{W}}\Omega^{j}_{W/V}$ sont des $\mathcal{O}_{W}$-module coh\'erents et plats sur $\mathcal{O}_{V}$ car $h_{V}$ est lisse.\\
Puisque ou bien $\overline{g}$ est plat, ou bien $\overline{g}$ est lisse sur un voisinage de $S'$ dans $\mathcal{T}'$, on peut d'apr\`es (3.3.1.4) et quitte \`a restreindre $V$ supposer $g_{V}$ plat: le (2) de (3.3.2.2) nous fournit alors l'isomorphisme de changement de base
$$\overline{g}_{K}^{\ast}R^{i}\overline{f}_{rig \ast}((X,Y)/\mathcal{T};E)\ \overset{\sim}{\rightarrow}\ R^{i}\overline{f'}_{rig \ast}((X'_{1},Y')/\mathcal{T}';{(\varphi'_{1},\overline{\varphi}')}^{\ast}(E)).
\leqno{(3.4.4.1.1)}
$$\\
\\
\noindent\textbf{Etape \  \rondIII }.  En reprenant la construction de la connexion de Gau\ss -Manin d\'ecrite explicitement par Katz dans [K 3, 3.4 et 3.5] on prouve que cette connexion agit aussi bien sur le terme $E_{1}^{i,j}$ que sur l'aboutissement [loc. cit., theo. 3.5]: ainsi $R^{i}{f}_{rig \ast}(X/\mathcal{T};E)$ est muni d'une connexion $\nabla^{i}$ (de Gau\ss -Manin) int\'egrable. Pour construire cette connexion de Gau\ss -Manin $\nabla^{i}$ on peut aussi \'etablir des isomorphismes [B 3, (2.2.5.1)] v\'erifiant la condition de cocycles: cette construction co\"{\i}ncide avec la pr\'ec\'edente [B 1,V, 3.6.3, 3.6.4, 3.6.5], et cette deuxi\`eme construction va nous permettre d'\'etablir la surconvergence de $\nabla^{i}$.\\
On a un diagramme commutatif \`a carr\'es cart\'esiens
$$
\begin{array}{c}
\xymatrix{
{\mathcal{Y}}\ar@{->}[rr]^(.45){\Gamma_{\overline{h}}=(1_{\mathcal{Y}}\times\overline{h})} \ar[d]_{\overline{h}} & &{\mathcal{Y}}\times_{\mathcal{W}}\mathcal{T}\ar@{->}[rr]^(.55){p'_{1}} \ar[d]^{\overline{h}\times1_{\mathcal{T}}} & &\mathcal{Y} \ar[d]^{\overline{h}} \\
{\mathcal{T}}\ar@{->} [rr]_{\Delta_{\mathcal{T}}}& &\mathcal{T}\times_{\mathcal{W}}\mathcal{T}\ar@{->} [rr]_{p_{1\mathcal{T}}} & &\mathcal{T},
}
\end{array}
$$
o\`u $p_{1\mathcal{T}},\ p'_{1}$ sont les premi\`eres projections, $\Gamma_{\overline{h}}$ est le morphisme graphe de $\overline{h}$ et $\Delta_{\mathcal{T}}$ est le morphisme diagonal: $\Gamma_{\overline{h}}$ est une immersion ferm\'ee puisque $\mathcal{T}$ est s\'epar\'e sur $\mathcal{W}$; on note alors $]Y[_{\mathcal{Y}\times_{\mathcal{W}}\mathcal{T}}$ le tube de $Y$ dans ${\mathcal{Y}}\times_{\mathcal{W}}\mathcal{T}$ pour l'immersion ferm\'ee compos\'ee ${Y}\hooklongrightarrow\mathcal{Y}\overset{\Gamma_{\overline{h}}}{\hooklongrightarrow}{\mathcal{Y}}\times_{\mathcal{W}}\mathcal{T}$, et $ ]T[_{\mathcal{T}^{2}}$ celui de $T$ dans ${\mathcal{T}^{2}}$ pour l'immersion ferm\'ee $T\hooklongrightarrow\mathcal{T}\underset{\Delta_{\mathcal{T}}}{\hooklongrightarrow}{\mathcal{T}}\times_{\mathcal{W}}\mathcal{T}$.\\
D'apr\`es (2.1.2) on a alors un diagramme commutatif \`a carr\'e cart\'esien
$$
\begin{array}{c}
\xymatrix{
 &  & &  ]Y[_{\mathcal{Y}^{2}}\ar@/^1.5pc/[ddl]^{\overset{\sim}{h}} \ar@/_1.5pc/[dlll]_{p_{1 \mathcal{Y}}} \ar@{->}[dl]^{h''_{1}}\\
]Y[_{\mathcal{Y}}\ar@{->}[d]_{\overline{h}_{Y}} & &]Y[_{\mathcal{Y}\times\mathcal{T}} \ar@{->}[d]^{h'_{1}} \ar@{->}[ll]_{p'_{1}}&\\
]T[_{\mathcal{T}} & & ]T[_{\mathcal{T}^{2}} \ar@{->}[ll]^(.45){p_{1_{\mathcal{T}}}}&
	}
\end{array}
\leqno{(3.4.4.1.2)}
$$
o\`u $h'_{1}$, $h''_{1}$, $ \overset{\sim}{h}$ sont induits respectivement par $\overline{h}_{K}\times1_{\mathcal{T}_{K}}$, $1_{\mathcal{Y}_{K}}\times\overline{h}_{K}$ et $\overline{h}_{K}\times\overline{h}_{K}$.\\
On a de m\^eme un diagramme analogue $(3.4.4.1.2)' $ avec $p_{2}$ \`a la place de $p_{1}$ et $h'_{2},h''_{2}$...\\
D'apr\`es (3.2.2.1) et [B 5, (3.2.3) et (3.1.11)(i)] on a les isomorphismes canoniques
 \begin{eqnarray*}
R^{i}\overline{f}_{rig \ast}((X,Y)/\mathcal{T}^{2};E)&=&R^{i}{\overset{\sim}{h}_{\ast}}({p}_{1\mathcal{Y}}^{\ast}(E_{\mathcal{Y}})\otimes\Omega^{\bullet}_{]Y[_{\mathcal{Y}^{2}}/\mathcal{T}_{K}^{2}})\\
&\simeq &R^{i}{{h'}_{1\ast}}({p'}_{1}^{\ast}(E_{\mathcal{Y}})\otimes\Omega^{\bullet}_{]Y[_{\mathcal{Y}\times\mathcal{T}}/\mathcal{T}_{K}^{2}})\\
&\overset{\sim}{\leftarrow}&{p}_{1\mathcal{T}}^{\ast}R^{i}\overline{h}_{Y}^{\ast}(E_{\mathcal{Y}}\otimes\Omega^{\bullet}_{]Y[_{\mathcal{Y}}/\mathcal{T}_{K}})\qquad [\mbox{\'etape \  \rondII} ]\\
&=&{p}_{1\mathcal{T}}^{\ast}R^{i}\overline{f}_{rig \ast}((X,Y)/\mathcal{T};E);
\end{eqnarray*}
l'avant-dernier isomorphisme ci-dessus est r\'ealisable via l'\'etape \ \rondII \ car $\rho$ \'etant lisse sur un voisinage de $S$ dans $\mathcal{T}$, $p_{1\mathcal{T}}:\mathcal{T}\times_{\mathcal{W}}\mathcal{T}\rightarrow\mathcal{T}$ est lisse sur un voisinage de $S$ dans $\mathcal{T}\times_{\mathcal{W}}\mathcal{T}$.\\
Or l'isomorphisme [B 3, (2.2.5.1)]
$$
p_{2\mathcal{Y}}^{\ast}(E_{\mathcal{Y}})\overset{\sim}{\longrightarrow}p_{1\mathcal{Y}}^{\ast}(E_{\mathcal{Y}})
$$
assurant l'existence d'une connexion (surconvergente) sur $E_{\mathcal{Y}}$ fournit aussi les isomorphismes
$$
R^{i}{\overset{\sim}{h}_{\ast}}({p}_{2\mathcal{Y}}^{\ast}(E_{\mathcal{Y}})\otimes\Omega^{\bullet}_{]Y[_{\mathcal{Y}^{2}}/\mathcal{T}_{K}^{2}})\overset{\sim}{\longrightarrow}R^{i}{\overset{\sim}{h}_{\ast}}({p}_{1\mathcal{Y}}^{\ast}(E_{\mathcal{Y}})\otimes\Omega^{\bullet}_{]Y[_{\mathcal{Y}^{2}}/\mathcal{T}_{K}^{2}}),
$$
c'est-\`a-dire, par les m\^emes arguments que ci-dessus, des isomorphismes\\
$$
(3.4.4.1.3)\ {p}_{2\mathcal{T}}^{\ast}R^{i}\overline{f}_{rig\ast}((X,Y)/\mathcal{T};E)\overset{\sim}{\rightarrow}R^{i}\overline{f}_{rig\ast}((X,Y)/\mathcal{T}^{2};E)\overset{\sim}{\leftarrow}{p}_{1\mathcal{T}}^{\ast}R^{i}\overline{f}_{rig\ast}((X,Y)/\mathcal{T};E)\\
$$
satisfaisant aux conditions de [B 3, (2.2.5)]. De plus la connexion surconvergente sur $R^{i}\overline{f}_{rig\ast}((X,Y)/\mathcal{T};E)$ obtenue par l'isomorphisme (3.4.4.1.3) est bien celle de Gau\ss-Manin [B 1,V, 3.6.4].\\
 Ceci ach\`eve la preuve de l'\'etape \ \  \rondIII.  
  \\
  
\noindent\textbf{Etape \  \rondIV }. Puisque $R^{i}\overline{f}_{rig\ast}((X,Y)/\mathcal{T};E)$ est un isocristal surconvergent le long de $T\setminus S$, son image inverse par $(\varphi, \overline{\varphi})$ est, d'apr\`es l'\'etape\  \  \rondII , l'image inverse par $(j, id_{T'})$ de l'isocristal surconvergent (le long de $T'\backslash S^{\prime}_{1}$) \\
$$R^{i}\overline{f'}_{rig\ast}((X'_{1},Y')/\mathcal{T}';{(\varphi'_{1}, \overline{\varphi}')}^{\ast}(E))=j_{1_{T'}}^{\dag}(R^{i}h'_{{V'}^{\ast}}({g'}_{W}^{\ast}(E_{W})\otimes_{\mathcal{O}_{W'}}\Omega^{\bullet}_{W'/V'}))\\\mbox{(cf (3.3.1.3))}.$$
Le $\mathcal{O}_{V'}$-module
$$\mathcal{E}_{V'}^{i}:=R^{i}h'_{{V'}^{\ast}}({g}_{W}^{\prime\ast}(E_{W})\otimes_{\mathcal{O}_{W'}}\Omega^{\bullet}_{W'/V'})
$$
est coh\'erent, et par d\'efinition des images inverses [B 3, (2.3.2)(iv)] on a
$$
{(j, id_{T'})}^{\ast}(j_{1_{T'}}^{\dag}(\mathcal{E}_{V'}^{i}))=j_{{T'}}^{\dag}(\mathcal{E}_{V'}^{i});
$$
or d'apr\`es [(3.3.2.2)(1)] on a:
\begin{eqnarray*}
j_{{T'}}^{\dag}(\mathcal{E}_{V'}^{i})&\simeq&R^{i}\overline{h}'_{{K}^{\ast}}(j_{Y'}^{\dag}({g}_{W}^{\prime\ast}(E_{W}))\otimes_{\mathcal{O}_{W'}}\Omega^{\bullet}_{W'/V'})\\
&=:&R^{i}\overline{f}'_{rig^{\ast}}((X',Y')/\mathcal{T}';{(\varphi, \overline{\varphi})}^{\ast}(E)),
\end{eqnarray*}
d'o\`u l'isomorphisme de changement de base
$$
(3.4.4.1.4)\ {(\varphi, \overline{\varphi})}^{\ast}R^{i}{\overline{f}}_{rig\ast}((X,Y)/\mathcal{T};E)\simeq R^{i}\overline{f}'_{rig\ast}((X',Y')/\mathcal{T}';{(\varphi', \overline{\varphi'})}^{\ast}(E)),
$$
qui est celui de (3.4.4.1)(ii).\\
\\
\noindent\textbf{Etape  \  \rondV}. Pla\c{c}ons nous sous les hypoth\`eses (3.4.3) et prouvons (3.4.4.1) dans ce cas.\\
Consid\'erons les parall\'el\'epip\`edes commutatifs suivants dans lesquels les faces verticales sont cart\'esiennes

 $$
\shorthandoff{;:!?}
\begin{array}{l}
\xymatrix@!0 @R=1,4cm @C=1,4cm{
& X \ar@{.>}[dd]^(.3){f}  \ar@{^{(}->}[rr]^{j_{Y}} & & Y\ar@{.>}[dd]^(.3){\overline{f}}  \ar@{^{(}->}[rr]^{i_{Y}} & &\mathcal{Y} \ar[dd]^{\overline{h}} \\
X^{'}  \ar@{^{(}->}[rr]^(.8){j_{Y'}} \ar[dd]_{f^{'}} \ar[ur]^{\varphi'} & & Y' \ar@{^{(}->}[rr]^(.8){i''_{Y'}} \ar[dd]_(.7){\overline{f}'} \ar[ur]^{{\overline{\varphi}}'} & & \mathcal{Y}'' \ar[dd]_(.7){\overline{h}''} \ar[ur]^{p_{1\mathcal{Y}}} \\
& S\ar@{^{(}.>}[rr]^(.7){j_{T}}  & &T\ar@{^{(}.>}[rr]^(.8){i_{T}}   & &  \mathcal{T}\ar[rr]_{\rho}&& \mathcal{W}\\
S^{'}  \ar@{^{(}->}[rr]_{j_{T'}} \ar@{.>}[ur]^{\varphi}& & T' \ar@{^{(}->}[rr]_{i''_{T'}} \ar@{.>}[ur]_{\overline{\varphi}} & & **[r]\mathcal{T''}=\mathcal{T}\times_{\mathcal{W}}\mathcal{T}' \ar[ur]_{p_{1 \mathcal{T}''}} 
}
\end{array}
\leqno{(3.4.4.1.5)}
  $$
 $$
  \begin{array}{c}
\xymatrix{
& X' \ar@{.>}[dd]^(.3){f^{'}}  \ar@{^{(}->}[rr]^{j_{Y'}} & & Y'\ar@{.>}[dd]^(.3){\overline{f}'}  \ar@{^{(}->}[rr]^{i_{Y'}'''} & &\mathcal{Y}''' \ar[dd]^{\overline{h}\ '''} \ar[dl]  \\
X^{'}  \ar@{^{(}->}[rr]^(.8){j_{Y'}} \ar[dd]_{f^{'}} \ar@{=}[ur] & & Y' \ar@{^{(}->}[rr]^(.8){i_{Y'}} \ar[dd]_(.7){\overline{f}'} \ar@{=}[ur] & & \mathcal{Y}' \ar[dd] _(.7){\overline{h}'}\\
& S'\ar@{^{(}.>}[rr]^(.7){j_{T'}}  \ar@{.>}[dl]_{id}& &T'\ar@{^{(}.>}[rr]^(.8){i''_{T'}}  \ar@{.>} [dl] ^{id}& &  \mathcal{T}'' \ar[dl]^{p_{2 \mathcal{T}''}}\\
S^{'}  \ar@{^{(}->}[rr]_{j_{T'}} & & T' \ar@{^{(}->}[rr]_{i_{T'}} & &\mathcal{T}'  \ar[rr]_{\rho'}&& \mathcal{W}' 
}
\end{array}
\leqno{(3.4.4.1.6)}
  $$
o\`u $p_{i}$ est la projection sur le $i \ieme $  facteur et $i''_{T'}=(i_{T}\circ\overline{\varphi},i_{T'})$.\\
On forme aussi le carr\'e cart\'esien
$$
\begin{array}{c}
\xymatrix{
\overset{\sim}{\mathcal{Y}} \ar[rr]^{u_2} \ar[d]_{u_{1}}& &\mathcal{Y}'' \ar[d]^{\overline{h}''}\\
\mathcal{Y}''' \ar[rr]_{\overline{h}'''}& &\mathcal{T}'' \   .
}
\end{array}
\leqno{(3.4.4.1.7)}
$$
Comme $u_{1}$ et $u_{2}$ sont propres, et lisses sur un voisinage de $X'$ dans $\widetilde{\mathcal{Y}}$, on peut d'apr\`es [B 5, (3.1.2)] et [LS, 7.4.2] faire le calcul de 
$$R^{i}\overline{f}'_{rig\ast}((X',Y')/\mathcal{T}';{(\varphi', \overline{\varphi'})}^{\ast}(E))$$
\`a l'aide de la cohomologie de de Rham, aussi bien avec $\overline{h}'''$, $\overline{h}''' \circ u_{1} = \overline{h}'' \circ u_{2}$ ou $\overline{h}''$, qui sont tous les trois propres, et lisses sur un voisinage de $X'$ respectivement dans $\mathcal{Y}''', \widetilde{\mathcal{Y}}$ et $\mathcal{Y}''$. Or $p_{1 \mathcal{T}''}$ (resp $ p_{2\mathcal{T}''}$) \'etant lisse sur un voisinage de $S'$ dans $\mathcal{T}''$, on a d'apr\` es l'\'etape\  \  \rondIV  \ des isomorphismes de changement de base (le premier calcul\'e via $\overline{h}''$ et le second via $\overline{h}'''$)
$$ {(\varphi, \overline{\varphi}, p_{1\mathcal{T}''})}^{\ast}R^{i}{\overline{f}}_{rig\ast}((X,Y)/\mathcal{T};E)\simeq R^{i}\overline{f}'_{rig\ast}((X',Y')/\mathcal{T}'';{(\varphi', \overline{\varphi'})}^{\ast}(E)) $$
et
$$ {(id_{S'}, id_{T'}, p_{2\mathcal{T}''})}^{\ast}R^{i}\overline{f}'_{rig\ast}((X',Y')/\mathcal{T}' ;{(\varphi', \overline{\varphi'})}^{\ast}(E))\simeq R^{i}\overline{f}'_{rig\ast}((X',Y')/\mathcal{T}';{(\varphi', \overline{\varphi'})}^{\ast}(E))$$
et les seconds membres co\"{\i}ncident en faisant le calcul de la cohomologie de de Rham via $\overline{h}''' \circ u_{1} = \overline{h}'' \circ u_{2}$. Ceci ach\`eve la preuve de l'isomorphisme de changement de base.\\

\noindent\textbf{Etape \  \rondVI }. L'\'etape \ \rondIII \ reste inchang\'ee puisqu'on peut faire les changements de base par $p_{1\mathcal{T}}$ et $p_{2\mathcal{T}}$ car ils sont lisses sur un voisinage de $S$ dans $\mathcal{T}\times_{\mathcal{W}}\mathcal{T}$: en particulier $R^{i}{\overline{f}}_{rig\ast}((X,Y)/\mathcal{T};E)$ est un isocristal surconvergent le long de $T\setminus S$. Soit $\mathcal{W}_{0}$ la r\'eduction sur $k$ de $\mathcal{W}$. On note \\
$$
\begin{array}{c}
\xymatrix{
X'' \ar@{^{(}->}[r]^{j_{Y''}} \ar[d]_{f''}&Y'' \ar@{^{(}->}[r]^{i_{Y''}} \ar[d]_{\overline{f}''} & **[r] \mathcal{Y}''=\mathcal{Y}\times_{\mathcal{W}}\mathcal{T}' \ar[d]_{\overline{h}''}\\
**[l]S''=S\times_{\mathcal{W}_o}S' \ar@{^{(}->}[r]_(.3){j_{T''}=j_{T}\times j_{T'}}&T'':=T\times_{\mathcal{W}_o}T'\ar@{^{(}->}[r]_(.6){i_{T''}=i_{T}\times i_{T'}}&**[r] \mathcal{T}'':=\mathcal{T}\times_{\mathcal{W}}\mathcal{T}'
}
\end{array}
\leqno{(3.4.4.2.1)}
$$
l'image inverse de (2.2.1) par le diagramme commutatif
$$
\begin{array}{c}
\xymatrix{
S\ \ar@{^{(}->}[rr]^{j_{T}} &&T\ \ar@{^{(}->}[rr]^{i_{T}}&&\mathcal{T}\\
S''\ \ar@{^{(}->}[rr]_{j_{T''}} \ar[u]^{p_{1S}} &&T'' \ar@{^{(}->}[rr]_{i_{T''}} \ar[u]_{p_{1T}}&&\mathcal{T}''\ \ar[u]_{p_{1\mathcal{T}''}}
}
\end{array}
\leqno{(3.4.4.2.2)}
$$
o\`u les $p_{1}$ sont les projections sur le premier facteur: $\overline{h}''$ est propre et $\overline{h}''$ est lisse sur un voisinage de $X''$ dans $\mathcal{Y}''$.\\
De m\^eme 
$$
\begin{array}{c}
\xymatrix{
X'\ \ar@{^{(}->}[rr]^{j_{Y'}} \ar[d]_{f'}&&Y'\ \ar@{^{(}->}[rr]^{i_{Y''}\circ\psi_Y} \ar[d]^{\overline{f}'} &&\mathcal{Y}'' \ar[d]^{\overline{h}''}\\
S'\ \ar@{^{(}->}[rr]_{j_{T'}}&&T'\ \ar@{^{(}->}[rr]_{i_{T''}\circ\psi_T}&&\mathcal{T}''
}
\end{array}
\leqno{(3.4.4.2.3)}
$$
est l'image inverse de (3.4.4.2.1) par le diagramme commutatif
$$
\begin{array}{c}
\xymatrix{
S''\ \ar@{^{(}->}[rr]^{j_{T''}} &&T''\ \ar@{^{(}->}[rr]^{i_{T''}}&&\mathcal{T}''\\
S'\ \ar@{^{(}->}[rr]_{j_{T'}} \ar[u]^{\psi_{S}=(\varphi,1_{S'})
} &&T'\ \ar@{^{(}->}[rr]_{i_{T''}\circ\psi_T} \ar[u]_{\psi_T=(\overline{\varphi},1_{T'})}&&\mathcal{T}''\ \ar@{=}[u]_{id_{\mathcal{T}''}}
}
\end{array}
\leqno{(3.4.4.2.4)}
$$
o\`u $\psi_{S}, \psi_{T}$ sont des immersions ferm\'ees et $\psi_{Y}$ est l'immersion ferm\'ee d\'efinie par le carr\'e cart\'esien
$$
\xymatrix{
Y'\ \ar[d]_{\overline{f}'} \ar@{^{(}->}[rr]^{\psi_Y} && Y'' \ar[d]^{\overline{f}''}\\
T'\ \ar@{^{(}->}[rr]_{\psi_T} && T''\ .
}
$$
En appliquant (3.4.4.1)(ii) au changement de base par $p_{1\mathcal{T}}$ et $\psi_{T}$ on obtient un isomorphisme 
$$ \psi_{T}^{\ast}p_{1\mathcal{T}''}^{\ast}R^{i}{\overline{f}}_{rig\ast}((X,Y)/\mathcal{T};E)\simeq R^{i}\overline{f}'_{rig\ast}((X',Y')/\mathcal{T}'';{(\varphi', \overline{\varphi'})}^{\ast}(E)) .$$
Puisque le carr\'e suivant commute
$$
\xymatrix{
T\ \ar@{^{(}->}[rr]^{i_T} && \mathcal{T} \\
T'\ \ar@{^{(}->}[rr]_{i_{T''}\circ\psi_T} \ar[u]^{p_{1T}\circ\psi_T= \overline{\varphi}}&& \mathcal{T}'' \ar[u]_{p_{1\mathcal{T}''}} \ ,
}
$$
pour prouver (3.4.4.2) il suffit d'apr\`es [B 3, (2.3.2)(iv)] de prouver que la projection sur le second facteur $p_{2\mathcal{T}''}: \mathcal{T}'' \rightarrow \mathcal{T}'$ induit un isomorphisme 
$$ p_{2\mathcal{T}''}^{\ast}R^{i}\overline{f}'_{rig\ast}((X',Y')/\mathcal{T}';{(\varphi', \overline{\varphi'})}^{\ast}(E))\simeq R^{i}\overline{f}'_{rig\ast}((X',Y')/\mathcal{T}'';{(\varphi', \overline{\varphi'})}^{\ast}(E)) .$$
Consid\'erons alors le diagramme commutatif \`a carr\'e cart\'esien

$$
\begin{array}{c}
\xymatrix{
X'\ar@{^{(}->}[r]^{j_{Y'}}&Y'\ar@{^{(}->}[r]^{\psi_{Y}}&Y''\ar@{^{(}->}[r]^{i_{Y''}}&\mathcal{Y}'' \ar@/_1.5pc/[dddr]_{\overline{h}''} \ar[dr]_{i_{\mathcal{Y}''}} \ar@{=}[drrr]^{Id_{\mathcal{Y}''}}\\
&&&&\mathcal{Y}'''\ar[rr]_{p_{3\mathcal{Y}}} \ar[dd]^{\overline{h}'''}&&\mathcal{Y}''\ar[d]^{\overline{h}''}\\
&&&&&&\mathcal{T}'' \ar[d]^{p_{2\mathcal{T}''}}\\
&&&&\mathcal{T}''\ar[rr]_{p_{2\mathcal{T}''}}&&\mathcal{T}'\ .
}
\end{array}
$$
D'apr\`es les hypoth\`eses faites sur $\rho$, $p_{2\mathcal{T}''}$ est propre et $p_{2\mathcal{T}''}$ est lisse sur un voisinage de $S''$ dans $\mathcal{T}''$: ainsi $p_{2\mathcal{T}''}\circ \overline{h}''$ est propre et $p_{2\mathcal{T}''}\circ \overline{h}''$ est lisse sur un voisinage de $X'$ dans $\mathcal{Y}''$; par suite en appliquant (3.4.1.1) on en d\'eduit que 
$$R^{i}\overline{f}'_{rig\ast}((X',Y')/\mathcal{T}';{(\varphi', \overline{\varphi'})}^{\ast}(E)) \ \mbox{et} \ R^{i}\overline{f}'_{rig\ast}((X',Y')/\mathcal{T}'';{(\varphi', \overline{\varphi'})}^{\ast}(E))$$
sont des isocristaux surconvergents et puisque $  i_{\mathcal{Y}''} $ est propre, que l'on a des isomorphismes (o\`u $p_{1\mathcal{Y}}: \mathcal{Y}'' \rightarrow \mathcal{Y}$ est la premi\`ere projection)
$$
\begin{array}{rcl}
p_{2\mathcal{T}''}^{\ast}R^{i}\overline{f}'_{rig\ast}((X',Y')/\mathcal{T}';{(\varphi', \overline{\varphi'})}^{\ast}(E))& \simeq& R^{i}\overline{h}_{\ast}'''{p}_{3\mathcal{Y}}^{\ast}((p_{1\mathcal{Y}}^{\ast}E_{\mathcal{Y}})\otimes\Omega^{\bullet}_{]Y'[_{\mathcal{Y}''}/ \mathcal{T}'_{K}})\\
&\simeq & R^{i}\overline{h}'''_{\ast}({p}_{3\mathcal{Y}}^{\ast}p_{1\mathcal{Y}}^{\ast}(E_{\mathcal{Y}})\otimes\Omega^{\bullet}_{]Y'[_{\mathcal{Y}^{'''}}/ \mathcal{T}_{K}^{''}}) \ [(2.1.2)]\\
&\simeq & R^{i}\overline{h}''_{\ast}( {i}_{\mathcal{Y}^{''}}^{\ast}{p}_{3\mathcal{Y}}^{\ast}{p}_{1\mathcal{Y}}^{\ast}(E_{\mathcal{Y}})\otimes\Omega^{\bullet}_{]Y'[_{\mathcal{Y}^{''}}/ \mathcal{T}_{K}^{''}})\  [B\ 5, (3.2.2)] \\
&= & R^{i}\overline{f}'_{rig\ast}((X',Y')/\mathcal{T}^{''};{(\varphi', \overline{\varphi'})}^{\ast}(E))\  [(3.2.2)]. 
\end{array}
$$
D'o\`u (3.4.4.2), compte tenu des d\'efinitions (3.2).  $\square$  \\

\noindent\textbf{(3.4.5)}\\

(3.4.5.1) Supposons donn\'e un diagramme commutatif
$$
\xymatrix{
X\ \ar@{^{(}->}[rr]^{j_{\mathcal{Y}}} \ar[d]_f \ar @{}[drr] |{\square}&& \mathcal{Y}\ar[d]^{\overline{h}}&&\\
S\ \ar@{^{(}->}[rr]_{j_{\mathcal{T}}} &&\mathcal{T}\ar[rr]^{\rho}&&\mathcal{W},
}
$$
dans lequel le carr\'e est cart\'esien, $f$ est un morphisme propre de $k$-sch\'emas s\'epar\'es de type fini, $\overline{h}$ et $\rho$ sont des morphismes de $\mathcal{V}$-sch\'emas formels s\'epar\'es plats de type fini, $\overline{h}$ est propre, $\overline{h}$ (resp $\rho$) est lisse sur un voisinage de $X$ (resp $S$) dans $\mathcal{Y}$ (resp $\mathcal{T}$), $j_{\mathcal{Y}} \ \mbox{et} \ j_{\mathcal{Y}}$ sont des immersions. Soit $T$ l'adh\'erence sch\'ematique de $S$ dans $\mathcal{T}$ et 
$$S\underset{j_{T}}{\hooklongrightarrow} T\underset{i_{T}}{\hooklongrightarrow}{\mathcal{T}}$$
la factorisation de $j_{\mathcal{T}}$: $j_{T}$ est une immersion ouverte dominante et $i_{T}$ est une immersion ferm\'ee. On note $Y= \overline{h}^{-1}(T)$ et $\overline{f}:= \overline{h}_{|_{Y}}: Y \rightarrow T.$ Avec les notations pr\'ec\'edentes on est donc dans la situation (2.2.1) avec $\overline{f}^{-1}(S) = X.$\\

(3.4.5.2) Supposons de plus donn\'e un diagramme commutatif
$$
\begin{array}{c}
\xymatrix{
S\ \ar@{^{(}->}[rr]^{j_{\mathcal{T}}} &&\mathcal{T} \ar[rr]^{\rho}&&\mathcal{W}\\
S'\ \ar@{^{(}->}[rr]^{j_{\mathcal{T}'}} \ar[u]^{\varphi} &&\mathcal{T}' \ar[rr]^{\rho'} \ar[u]_{\overline{g}}&&\mathcal{W}'\ \ar[u]_{\theta}
}
\end{array}
$$
dans lequel $\varphi$ est un morphisme de $k$-sch\'emas s\'epar\'es de type fini, $\overline{g}, \ \rho'$ sont des morphismes s\'epar\'es de $\mathcal{V}$-sch\'emas formels s\'epar\'es plats de type fini, $\overline{g}$ et $\rho'$ sont lisses sur un voisinage de $S'$ dans $\mathcal{T}'$, $\theta$ est lisse et $j_{\mathcal{T}'}$ est une immersion.\\
On d\'efinit $f':X'\rightarrow S'$ par le carr\'e cart\'esien
$$
\xymatrix{
X' \ar[rr]^{\varphi'} \ar[d]_{f'} && X \ar[d]^{f}\\
S' \ar[rr]_{\varphi}&&S\ .
}
$$
On note $T'$ l'adh\'erence sch\'ematique de $S'$ dans $\mathcal{T'}$ et 
$$S'\underset{j_{T'}}{\hooklongrightarrow} T' \underset{i_{T'}}{\hooklongrightarrow}{\mathcal{T}'}$$
la factorisation de $j_{\mathcal{T}'}$ par l'immersion ouverte dominante $j_{T'}$ et l'immersion ferm\'ee $i_{T'}$: $\overline{g}$ induit un $k$-morphisme $\overline{\varphi}:T' \rightarrow T\ .$\\
Par image inverse de (3.4.5.1) par (3.4.5.2) on obtient un parall\'el\'epip\`ede commutatif tel que (3.3.1.1), dont on reprend les notations; on en d\'eduit:\\ 

\noindent\textbf{Corollaire (3.4.5.3).}
\textit { Sous les hypoth\`eses (3.4.5.1) et (3.4.5.2) les parties (3.4.4.1) et (3.4.4.3) du th\'eor\`eme (3.4.4) sont valides.}\\

(3.4.5.4) Supposons donn\'e cette fois un diagramme commutatif
$$
\xymatrix{
S\ \ar@{^{(}->}[rr]^{j_{\mathcal{T}}}&&\mathcal{T}\ar[rr]^{\rho}&&\mathcal{W}\\
S'\ \ar[u]^{\varphi} \ar@{^{(}->}[rr]_{j_{\mathcal{T}'}}&&\mathcal{T}'\ar[rr]_{\rho'}&&\mathcal{W}' \ar[u]_{\theta}
}
$$
dans lequel $\varphi$ est un morphisme de $k$-sch\'emas s\'epar\'es de type fini, $\rho, \rho' ,\theta$ sont des morphismes (s\'epar\'es) de $\mathcal{V}$-sch\'emas formels s\'epar\'es plats de type fini, $\theta$ est lisse, $\rho$ (resp $\rho'$) est lisse sur un voisinage de $S$ dans $\mathcal{T}$ (resp de $S'$ dans $ \mathcal{T}' $), $ j_{\mathcal{T}} $ et $ j_{\mathcal{T}' }$ sont des immersions. Et on suppose de plus que $ \rho $ est propre. D'apr\`es l'\'etape \ \rondVI \  \ de la d\'emonstration du th\'eor\`eme (3.4.4) on a montr\'e:\\

\noindent\textbf{Corollaire (3.4.5.5).}
\textit { Sous les hypoth\`eses (3.4.5.1) et (3.4.5.4) la partie (3.4.4.2) du th\'eor\`eme (3.4.4) s'applique.}\\

\noindent\textbf{3.4.6.} Sous les hypoth\`eses de (3.4.5.1) supposons que $\mathcal{W}= Spf\ \mathcal{V}$ et que $\rho:\  \mathcal{T} \rightarrow Spf\ \mathcal{V}$ est lisse sur un ouvert $\mathcal{S}$ de $\mathcal{T}$, avec $S$ contenu dans $\mathcal{S}$.\\

Soient $i_{s}: s= Spec\ k(s) \hookrightarrow S$ un point ferm\'e de $S$ et $f_{s}:X_{s} \rightarrow s$ la fibre de $f$ en $s$. On note $\mathcal{V}(s)= W(k(s))\otimes_{W} \mathcal{V}$, o\`u $W=W(k),\ W(k(s)) $ sont les anneaux de vecteurs de Witt \`a coefficients dans $k$ et $k(s)$ respectivement, et $K(s)$ le corps des fractions de $\mathcal{V}(s)$. Le morphisme $i_{s}$ d\'efinit des foncteurs images inverses [B 3, (2.3.6)]
$$
\xymatrix{
\  i_{s}^{\dag \ast}: \ Isoc^{\dag}((S,T)/ \mathcal{V})\ \ar[rr]&&\ Isoc^{\dag} (Spec\ k(s) / K(s))\ar [d]_{\simeq}\\
{}&&\ Isoc (Spec\ k(s) / K(s)),
}
\leqno (3.4.6.1)
$$
qui, pour $\rho$ propre, devient
$$
\xymatrix{
\  i_{s}^{\dag \ast}: \ Isoc^{\dag}(S /  K)\ \ar[rr]&&\ Isoc^{\dag} (Spec\ k(s) / K(s))\ar [d]_{\simeq}\\
{}&&\ Isoc (Spec\ k(s) / K(s)),
}
\leqno (3.4.6.2)
$$
et qui, pour $S=T$, devient
$$
\xymatrix{
\  \hat{i}_{s}^{\ast}: \ Isoc (S /  K)\ \longrightarrow Isoc (Spec\ k(s) / K(s)).
}
\leqno (3.4.6.3)
$$
Nous allons donner maintenant une r\'ealisation explicite de ces foncteurs: dans la notation $i_{s}^{\dag \ast}$, le $()^{\dag}$ n'est pas utile, sauf pour insister que l'on travaille avec la cat\'egorie surconvergente, et pas seulement la convergente.\\

Dans le carr\'e cart\'esien de $\mathcal{V}$-sch\'emas formels plats et s\'epar\'es 
$$
\xymatrix{
\mathcal{U}(s) \ar[rr]^{u'} \ar[d]_{v'} && \mathcal{S} \ar[d]^{v}\\
Spf\ \mathcal{V}(s)\ar[rr]_{u}&&Spf\ \mathcal{V}\ 
}
$$
les morphismes $u$ et $u'$ (resp $v$ et $v'$) sont finis \'etales (resp lisses et s\'epar\'es). Par lissit\'e de $\mathcal{S}$ sur $\mathcal{V}$ le morphisme compos\'e 

$$Spec\ k(s)\ \overset{i_{s}}{\hooklongrightarrow} S\ \hooklongrightarrow \mathcal{S}$$

\noindent se rel\`eve en un morphisme s\'epar\'e 
$$g_{s}\ :\ Spf\ \mathcal{V}(s) \longrightarrow \ \mathcal{S}$$
tel que $g_{s}\circ v'\ = \ u'$ donc $v'$ est fini et $v'$ est fini \'etale puisque $v'$ est lisse. Par la propri\'et\'e universelle du produit fibr\'e, $g_{s}$ se factorise en 
$$g_{s}\ :\ Spf\ \mathcal{V}(s)\ \overset{\tilde{g}_{s}}{\longrightarrow} \mathcal{U}(s)\ \overset{u'}{\hooklongrightarrow} \mathcal{S}$$
avec $v'\circ \tilde{g}_{s} = Id_{Spf\ \mathcal{V}(s)}$: comme $v'$ est \'etale, on en d\'eduit que $\tilde{g}_{s}$ est une immersion \`a la fois ouverte et ferm\'ee. Ainsi $g_{s}= u' \circ \tilde{g}_{s}$ est \'etale. Notons $\overline{g}_{s}$ le morphisme compos\'e
$$\overline{g}_{s}\ :\ Spf\ \mathcal{V}(s)\ \overset{g_{s}}{\longrightarrow} \mathcal{S}\ \hooklongrightarrow \mathcal{T};$$
$\overline{g}_{s}$ est \'etale et les foncteurs $ i_{s}^{\dag \ast}$ et $\hat{i}_{s}^{\ast}$ sont induits par $\overline{g}_{s}^{\ast}$ [B 3, (2.3.6)].\\

\noindent\textbf{Th\'eor\`eme (3.4.7).}
\textit{ On se place sous les hypoth\`eses et notations de (3.4.6).
\begin{enumerate}
\item [(3.4.7.1)] Soit $E \in Isoc^{\dag}((X,Y) / \mathcal{V})$. Alors, pour tout entier $i\geqslant 0$ et tout point ferm\'e $s$ de $S$, on a des isomorphismes
$$
\begin{array}{rcl}
i_{s}^{\dag \ast}R^{i}\overline{f}_{rig\ast}((X,Y)/\mathcal{T};\ E)& \simeq& R^{i}f_{s\ rig\ast}(X_{s}/\mathcal{V}(s); E_{X_{s}})\\
&\simeq & R^{i}f_{s\ conv\ast}(X_{s}/\mathcal{V}(s);\widehat{E_{X_{s}}})\\
&\simeq &H^{i}_{rig}( X_{s}/ K(s); E_{X_{s}})\\
&= & H^{i}_{conv}( X_{s}/ K(s); \widehat{E_{X_{s}}}) 
\end{array}
$$
o\`u l'isocristal convergent $\widehat{E_{X_{s}}}\in Isoc\ (X_{s}/K(s))= Isoc^{\dag}\ (X_{s}/K(s))$ co\"\i ncide avec l'isocristal surconvergent $E_{X_{s}}\in Isoc^{\dag}\ (X_{s}/K(s)).$
\item[(3.4.7.2)] Supposons $\rho\:\ \mathcal{T}\rightarrow\ Spf\ \mathcal{V}$ propre et soit $E \in Isoc^{\dag}\ (X/ K).$ Alors, pour tout entier $i\geqslant 0$ et tout point ferm\'e $s$ de $S$, on a des isomorphismes
$$
\begin{array}{rcl}
i_{s}^{\dag \ast}R^{i}f_{rig\ast}(X/\mathcal{T};\ E)& \simeq& R^{i}f_{s\ rig\ast}(X_{s}/\mathcal{V}(s);E_{X_{s}})\\
&\simeq & R^{i}f_{s\ conv\ast}(X_{s}/\mathcal{V}(s); \widehat{E_{X_{s}}})\\
&\simeq &H^{i}_{rig}( X_{s}/ K(s); E_{X_{s}})\\
&= & H^{i}_{conv}( X_{s}/ K(s); \widehat{E_{X_{s}}}) 
\end{array}
$$
o\`u l'isocristal convergent $\widehat{E_{X_{s}}}\in Isoc\ (X_{s}/K(s))= Isoc^{\dag}\ (X_{s}/K(s))$ co\"\i ncide avec l'isocristal surconvergent $E_{X_{s}}\in Isoc^{\dag}\ (X_{s}/K(s)).$
\item[(3.4.7.3)] Supposons que $S=T$ et soit $\mathcal{E}\in Isoc\ (X/K)$.  Alors, pour tout entier $i\geqslant 0$ et tout point ferm\'e $s$ de $S$, on a des isomorphismes
$$
\begin{array}{rcl}
\hat{i}_{s}^{ \ast}R^{i}f_{conv\ast}(X/\mathcal{T};\ \mathcal{E})& \simeq& R^{i}f_{s\ conv\ast}(X_{s}/\mathcal{V}(s); \mathcal{E}_{X_{s}})\\
&= & H^{i}_{conv}( X_{s}/ K(s); \mathcal{E}_{X_{s}}). 
\end{array}
$$
\end{enumerate}
}
\noindent\textit{D\'emonstration.} Puisque $\overline{g}_{s}$ est \'etale on peut donc appliquer le changement de base par ($i_{s},\ \overline{g}_{s}$) et les corollaires (3.4.5.3) et (3.4.5.5). $\square$\\

\noindent\textbf{ 3.4.8.}
Soient $S$ un $k$-sch\'ema lisse et s\'epar\'e et $f\ :  X\rightarrow S$ un $k$-morphisme projectif et lisse. Notons $S= \underset{\alpha}{\bigcup} \ S_{\alpha,0}$ une d\'ecomposition de $S$ en r\'eunion d'ouverts connexes affines $S_{\alpha,0}= Spec(A_{\alpha,0}), \ A_{\alpha}= \mathcal{V}[t_{1},...,t_{d_{\alpha}}]/J_{\alpha}$ une $\mathcal{V}$-alg\`ebre lisse relevant $A_{\alpha,0}$ dont on a fix\'e une pr\'esentation et $S_{\alpha}= Spec (A_{\alpha})$. On d\'esigne par $\overline{S}_{\alpha}$ l'adh\'erence sch\'ematique de $S_{\alpha}$ dans $\mathbb{P}_{\mathcal{V}}^{d_{\alpha}}$, par $\mathcal{S}_{\alpha}$ (resp $\overline{\mathcal{S}}_{\alpha}$) le compl\'et\'e formel $\mathfrak{m}$-adique de $S_{\alpha}$ (resp $\overline{S}_{\alpha}$) et par 
$$
f_{\alpha}:X_{\alpha,0}=X\times_{S}S_{\alpha,0}\longrightarrow S_{\alpha,0}
$$
la restriction de $f$. Quitte \`a d\'ecomposer $X_{\alpha}$ en somme disjointe de ses composantes connexes on peut supposer $X_{\alpha}$ connexe. D'apr\`es [I, (3.3.1)] il existe un carr\'e cart\'esien
$$
\begin{array}{c}
\xymatrix{
X_{\alpha} \ar@{^{(}->} [rr] \ar [d]_{h_{\alpha}} && \overline{X}_{\alpha}\ar[d]^{\overline{h}_{\alpha}}\\
S_{\alpha} \ar@{^{(}->} [rr] && \overline{S}_{\alpha} 
}
\end{array}
$$
dans lequel $\overline{h}_{\alpha}$ est projectif, $h_{\alpha}$ est un rel\`evement projectif de $f_{\alpha}$ et les fl\`eches horizontales sont des immersions ouvertes. Le compl\'et\'e formel de ce carr\'e est d'apr\`es [I, th\'eo (3.3)] un carr\'e cart\'esien de $\mathcal{V}$-sch\'emas formels 
$$
\begin{array}{c}
\xymatrix{
\mathcal{X}_{\alpha} \ar@{^{(}->} [rr] \ar [d]_{\hat{h}_{\alpha}} && \overline{\mathcal{X}}_{\alpha}\ar[d]^{\hat{\overline{h}}_{\alpha}}\\
\mathcal{S}_{\alpha} \ar@{^{(}->} [rr] && \overline{\mathcal{S}}_{\alpha} 
}
\end{array}
\leqno (3.4.8.1)
$$
dans lequel $\hat{\overline{h}}_{\alpha}$ est projectif, $\hat{h}_{\alpha}$ est un rel\`evement projectif de la restriction $f_{\alpha}$ de $f$ et les fl\`eches horizontales sont des immersions ouvertes.\\

\noindent\textbf{Th\'eor\`eme (3.4.8.2).}
\textit{Sous les hypoth\`eses (3.4.8) supposons que, pour tout $\alpha$, $X_{\alpha}$ est plat sur $\mathcal{V}$. Alors, pour tout entier $i \geqslant 0$, on a un diagramme commutatif de foncteurs naturels induits par $f$ et d\'efinis ci-apr\`es en (3.4.8.5)
$$
\xymatrix{
 Isoc^{\dag}(X/ K) \ar[rr]^{R^{i}f_{rig\ast}}\ar[d]&&Isoc^{\dag}(S/K)\ar [d]\\
Isoc (X/ K)\ar[rr]^{R^{i}f_{conv\ast}}&&Isoc(X/ K)
}
$$
o\`u les fl\`eches verticales sont les foncteurs d'oubli.\\
De plus ces foncteurs commutent \`a tout changement de base $S'\rightarrow S$ entre $k$-sch\'emas lisses et s\'epar\'es; en particulier ils commutent aux passages aux fibres en les points ferm\'es de $S$.
}\\

\noindent\textit{D\'emonstration.} Puisque, pour tout $\alpha$, $X_{\alpha}$ est plat sur $\mathcal{V}$, le th\'eor\`eme [I, (3.3)] prouve que $\hat{h}_{\alpha}$ est un rel\`evement projectif et lisse de la restriction $f_{\alpha}$ de $f$. Soient $E \in Isoc^{\dag}(X/K)$ et $E_{\alpha}$ sa restriction \`a $X_{\alpha}$: gr\^ace \`a l'existence du carr\'e cart\'esien (3.4.8.1) dans lequel $\hat{h}_{\alpha}$ est un rel\`evement projectif et lisse de  $f_{\alpha}$ on conclut \`a l'aide du th\'eor\`eme (3.4.4) que
$$
R^{i}f_{\alpha rig\ast}(X_{\alpha,0}/\overline{\mathcal{S}}_{\alpha}; E_{\alpha}) \in Isoc^{\dag}(S_{\alpha,0}/K).
$$
Montrons que celui-ci ne d\'epend que de $S_{\alpha,0}$ et non de $\overline{\mathcal{S}}_{\alpha}$.\\

Supposons donn\'es un $k$-sch\'ema propre $\overline{S}'_{\alpha,0}$, un $\mathcal{V}$-sch\'ema formel propre $\overline{\mathcal{S}}'_{\alpha}$, une immersion ouverte dominante $j'_{\alpha,0}: S_{\alpha,0}\hookrightarrow \overline{S}'_{\alpha,0} $ et une immersion ferm\'ee $i'_{\alpha,0}: \overline{S}'_{\alpha,0}\hookrightarrow \overline{\mathcal{S}}'_{\alpha,0} $ tels que le morphisme $\overline{\mathcal{S}}'_{\alpha}\longrightarrow Spf \mathcal{V}$ soit lisse sur un voisinage de $S_{\alpha,0}$ dans $ \overline{\mathcal{S}}'_{\alpha}$. Notons $\overline{T}_{\alpha,0}$ l'adh\'erence sch\'ematique de $S_{\alpha,0}$ plong\'e diagonalement dans $\overline{S}''_{\alpha}=\overline{S}_{\alpha,0}\times_{k} \overline{S}'_{\alpha,0}$ et $\overline{\mathcal{S}}''_{\alpha}=\overline{\mathcal{S}}_{\alpha}\hat{\times}_{\mathcal{V}} \overline{\mathcal{S}}'_{\alpha}$. On a un diagramme commutatif
$$
\begin{array}{c}
\xymatrix{
S_{\alpha,0}\ \ar@{^{(}->}[rr]^{j_{T}} \ar[d]_{id}&& \overline{T}_{\alpha,0}\ \ar@{^{(}->}[rr]^{i_{T}} \ar[d]^{p_{\alpha}} &&\overline{\mathcal{S}}''_{\alpha} \ar[d]^{\overline{p}_{1,\alpha}}\\
S_{\alpha,0}\ \ar@{^{(}->}[rr]_{j_{S}}&&\overline{S}_{\alpha,0}\ \ar@{^{(}->}[rr]_{i_{S}} &&\overline{\mathcal{S}}
}
\end{array}
\leqno{(3.4.8.3)}
$$
dans lequel les fl\`eches verticales sont induites par la premi\`ere projection, les $i$ (resp les $j$) sont des immersions ferm\'ees (resp ouvertes). Comme $p_{\alpha}$ est propre, le foncteur $(p_{\alpha},\overline{p}_{1,\alpha})^{\ast}$ est une \'equivalence de cat\'egorie [B 3, (2.3.5)] de la cat\'egorie des isocristaux sur $S_{\alpha,0}/K$ surconvergents le long de $Z_{\alpha,0}=\overline{S}_{\alpha,0}\setminus S_{\alpha,0}$ dans la cat\'egorie des isocristaux sur $S_{\alpha,0}/K$ surconvergents le long de $Z''_{\alpha,0}=\overline{T}_{\alpha,0}\setminus S_{\alpha,0}$. De plus $\overline{p}_{1,\alpha}$ est propre, et lisse sur un voisinage de $S_{\alpha,0}$ dans $\overline{\mathcal{S}}''_{\alpha}$, donc par le th\'eor\`eme (3.4.4) on a isomorphisme de changement de base
$$\overline{p}_{1\alpha}^{\ \ast}R^{i}f_{\alpha rig \ast}(X_{\alpha,0}/\overline{\mathcal{S}}_{\alpha};E_{\alpha}) \overset{\sim}	{\longrightarrow} 	R^{i}{f}_{\alpha rig \ast}(X_{\alpha,0}/\overline{\mathcal{S}}''_{\alpha};E_{\alpha}),$$
et de m\^eme
$$\overline{p}_{2\alpha}^{\ \ast}R^{i}f_{\alpha rig \ast}(X_{\alpha,0}/\overline{\mathcal{S}}'_{\alpha};E_{\alpha}) \overset{\sim}	{\longrightarrow} 	R^{i}{f}_{\alpha rig \ast}(X_{\alpha,0}/\overline{\mathcal{S}}''_{\alpha};E_{\alpha}),$$
d'o\`u un isomorphisme canonique
$$R^{i}f_{\alpha rig \ast}(X_{\alpha,0}/\overline{\mathcal{S}}_{\alpha};E_{\alpha}) \overset{\sim}	{\longrightarrow} 	R^{i}{f}_{\alpha rig \ast}(X_{\alpha,0}/\overline{\mathcal{S}}'_{\alpha};E_{\alpha}).
\leqno (3.4.8.4)
$$
Ainsi $R^{i}f_{\alpha rig \ast}(X_{\alpha,0}/\overline{\mathcal{S}}_{\alpha};E_{\alpha})$ ne d\'epend que de $S_{\alpha,0}$; de plus sur $S_{\alpha,0}\cap S_{\beta,0}$ ces isocristaux se recollent car on a des isomorphismes analogues \`a (3.4.8.4) et v\'erifiant la condition de cocycles. Ces donn\'ees qui se recollent fournissent un isocristal surconvergent sur $S/K$ [B 3, (2.3.2)] not\'e
$$
R^{i}f_{rig \ast}(E),
\leqno (3.4.8.5)
$$
et c'est le seul possible d'apr\`es le th\'eor\`eme de pleine fid\'elit\'e de Kedlaya  [Ked 3, Theo 5.2.1].\\
On raisonne de mani\`ere analogue pour $R^{i}f_{conv \ast}(E)$ et le carr\'e du th\'eor\`eme est clairement commutatif.\\
L'assertion sur les changements de base $S'\rightarrow S$ r\'esulte de (3.4.4.2). $\square$\\

\noindent\textbf{Corollaire (3.4.8.6).}
\textit{Sous les hypoth\`eses (3.4.8), supposons que $f$ d\'efinit $X$ comme une intersection compl\`ete relativement \`a $S$ dans un espace projectif sur $S$ [I, (3.3.6)]. Alors les conclusions du th\'eor\`eme (3.4.8.2) demeurent valides.} \\
 
\noindent \textit{D\'emonstration}. Avec les notations de [I, (3.3.6.2)], chaque $f_{\alpha, \beta} :\  X_{\alpha, \beta} \rightarrow S_{\alpha, \beta}$ se rel\`eve d'apr\`es [I, (3.3.7)] en un morphisme projectif et lisse au-dessus de $\mathcal{V}$ et donne lieu \`a un diagramme commutatif tel que (3.4.8.1) en rempla\c cant $\mathcal{S}_{\alpha}$ par $\mathcal{S}_{\alpha, \beta}$, et on conclut comme pour (3.4.8.2). $\square$\\

On a la variante suivante du th\'eor\`eme (3.4.8.2):\\

\noindent\textbf{Th\'eor\`eme (3.4.9).}
\textit{Supposons que le th\'eor\`eme (3.4.8.2) est vrai pour tout morphisme projectif et lisse $f :  X\rightarrow S$ avec $S$ un $k$-sch\'ema affine et lisse (donc sans faire l'hypoth\`ese $X_{\alpha}$ plat sur $\mathcal{V}$). Soient $S$ un $k$-sch\'ema lisse et s\'epar\'e et $f :  X\rightarrow S$ un $k$-morphisme projectif et lisse. Alors, pour tout entier $i \geqslant 0$, on a un diagramme commutatif de foncteurs naturels induits par $f$ et d\'efinis \`a la mani\`ere de(3.4.8.5)
$$
\xymatrix{
 Isoc^{\dag}(X/ K) \ar[rr]^{R^{i}f_{rig\ast}}\ar[d]&&Isoc^{\dag}(S/K)\ar [d]\\
Isoc (X/ K)\ar[rr]^{R^{i}f_{conv\ast}}&&Isoc(S/ K)
}
$$
o\`u les fl\`eches verticales sont les foncteurs d'oubli.\\
De plus ces foncteurs commutent \`a tout changement de base $S'\rightarrow S$ entre $k$-sch\'emas lisses et s\'epar\'es; en particulier ils commutent aux passages aux fibres en les points ferm\'es de $S$.
}\\

\noindent\textit{D\'emonstration.} Il suffit de d\'ecomposer $S$ en sch\'emas affines et lisses sur $k$ et d'op\'erer les recollements comme dans la preuve de (3.4.8.2).
 $\square$\\


\cleardoublepage


\vskip 10mm
\chapter*{III.  $F$-isocristaux convergents sur un sch\'ema lisse et images directes}
\markboth{\sc j.-y. etesse}{\sc III.  $F$-isocristaux convergents et images directes}
\noindent On conserve les notations du II (0).\\

\section*{1. $F$-isocristaux convergents sur un sch\'ema affine et lisse}

\subsection*{1.1. Notations}

Soient $X = \mbox{Spec}\  A_{0}$ un $k$-sch\'ema affine et lisse, $\mathcal{A}$ une $C(k)$-alg\`ebre lisse relevant $A_{0}$ [E$\ell$, th\'eo 6] et $A = \mathcal{A}\ \otimes_{C(k)} \mathcal{V}$. Fixons une pr\'esentation

$$\mathcal{A} = C(k) [t_{1}, ..., t_{n}]\ /\ (f_{1}, ... , f_{s}) ; $$
soient $P$ le compl\'et\'e formel de la fermeture projective de $\mathcal{X} = \mbox{Spec}\  A$ dans $\mathbb{P}^n_{\mathcal{V}}$, $Y$ sa r\'eduction sur $k$, $j : X\ \hookrightarrow\ Y. $ Alors $]Y[_{P}\  = P_{K}$ ; de plus $]X[_{P}$, qui est l'intersection de $\mathcal{X}^{an}_{K}$ avec la boule unit\'e $B(0,1^+) \subset \mathbb{A}^n_{K}$, est l'affino¬\"{\i}de 

$$]X[_{P}\  = \mbox{Spm}\ (K\{t_{1},...,t_{n}\}\ /\ (f_{1},...,f_{s})). $$ 
Notons $\hat{\mathcal{A}}$ (resp. $\hat{A}$) le s\'epar\'e compl\'et\'e $p$-adique de $\mathcal{A}$ (resp. $A$), $\mathcal{A}^{\dag} \subset \hat{\mathcal{A}}$ (resp. $A^{\dag} \subset \hat{A})$ le compl\'et\'e faible de $\mathcal{A}$ au-dessus de $(C(k), (p))$ (resp. de $A$ au-dessus de $(\mathcal{V},(\pi))$ [M-W, \S\ 1], et posons $A^{\dag}_{K} = A^{\dag}\ \otimes_{\mathcal{V}}\ K, \hat{A}_{K}\ =\ \hat{A}\ \otimes_{\mathcal{V}}\ K$. On a des isomorphismes\\

\qquad $\hat{\mathcal{A}}\ \simeq\ C(k) \{t_{1},...,t_{n} \} / (f_{1},...,f_{s}) $ , \\

\qquad $\hat{A}\ \simeq\ \mathcal{V} \{t_{1},...,t_{n} \} / (f_{1},...,f_{s})\ \simeq\  \hat{\mathcal{A}} \otimes_{C(k)} \mathcal{V}$ ,\\

\qquad $\mathcal{A}^{\dag}\ \simeq\ C(k) [t_{1},...,t_{n}]^{\dag} / (f_{1},...,f_{s})$ ,\\

\qquad $A^{\dag}\ \simeq\ \mathcal{V}[t_{1},...,t_{n}]^{\dag} / (f_{1},...,f_{s})\ \simeq\ \mathcal{A}^{\dag} \otimes_{C(k)} \mathcal{V}$ ,\\

\noindent et aussi [B 3, (2.1.2.4]\\

\qquad $\Gamma(P_{K}, j^{\dag} \mathcal{O}_{P_{K}})\ \simeq\ \Gamma(\mathcal{X}^{an}_{K}, j^{\dag} \mathcal{O}_{P_{K}})\ \simeq\ A^{\dag}_{K} $,\\

\qquad $\Gamma(] X [,  j^{\dag} \mathcal{O}_{P_{K}})\ \simeq\ \Gamma (]X[, \mathcal{O}_{P_{K}})\ =\  \Gamma(]X[, \mathcal{O}_{]X[})\ \simeq\ \hat{A}_{K}$.\\

On fixe un rel\`evement $F_{\mathcal{A}^{\dag}} : \mathcal{A}^{\dag}\ \rightarrow\ \mathcal{A}^{\dag}$ de l'\'el\'evation \`a la puissance $p^a$, $F_{A_{0}} : A_{0}\ \rightarrow\ A_{0}$, au-dessus de $\sigma$ [vdP, cor 2.4.3] : on peut choisir un tel rel\`evement $F_{\mathcal{A}^{\dag}}$ de mani\`ere compatible \`a une extension $k\ \hookrightarrow\ k'$ du corps de base [Et 5, I, 1.2]. Posons $F_{A^{\dag}} = F_{\mathcal{A}^{\dag}} \otimes_{C(k)} \mathcal{V}$, $F_{\hat{A}}\ = F_{A^{\dag}} \otimes_{A^\dag}  {\hat{A}}$ ; d'o\`u des morphismes $F_{A^{\dag}_{K}} : A^{\dag}_{K}\ \rightarrow\  A^{\dag}_{K}$,
$F_{\hat{A}_{K}} : \hat{A}_{K}\ \rightarrow\  \hat{A}_{K}$ au-dessus de $\sigma : K \rightarrow\ K$. De m\^eme pour $\mathcal{V}'$ comme ci-dessus et $A' := A\ \otimes_{\mathcal{V}}\ \mathcal{V}'$ il existe d'apr\`es [vdP, cor 2.4.3] un carr\'e commutatif

$$
\xymatrix{
A'^{\dag} \ar[r]^{F_{A'^{\dag}}} & A'^{\dag} \\
A^{\dag} \ar@{^{(}->}[u] \ar[r]_{F_{A^{\dag}}} & A^{\dag} \ar@{^{(}->}[u] 
}
$$

\noindent au-dessus du carr\'e commutatif

$$
\xymatrix{
\mathcal{V} \ar[r]^{\sigma}' & \mathcal{V}' \\
\mathcal{V} \ar@{^{(}->}[u] \ar[r]^{\sigma} & \mathcal{V}.  \ar@{^{(}->}[u] 
}
$$

On d\'esigne par $\mathbf{F^a \textrm{\bf-Mod}(A^{\dag}_{K})}$ (resp. $ \mathbf{F^a \textrm{\bf-Modloc}(A^{\dag}_{K})} $; 
resp. $\mathbf{F^a\textrm{\bf-Modlib}(A^{\dag}_{K})}$) la cat\'egorie des $A^{\dag}_{K}$-modules de type fini (resp. $A^{\dag}_{K}$-modules projectifs de type fini ; resp. $A^{\dag}_{K}$-modules libres de type fini) $M$ munis d'un morphisme

$$\phi_{M} : F^{\ast}_{A^{\dag}_{K}}(M) =: M^{\sigma} \longrightarrow\ M$$

\noindent appel\'e morphisme de Frobenius. Par analogie avec Wan [W 2, def 2.8] [W 3, def 2.1] nous dirons que $M$ est un $\mathbf{F^a\textrm{\bf-module surconvergent}}$ (resp. et $\textbf{projectif}$ ; resp. et $\textbf{libre}$) sur $A^{\dag}_{K}$. En consid\'erant les m\^emes d\'efinitions sur $\hat{A}_{K}$ au lieu de $A^{\dag}_{K}$ on dira que $\mathcal{M} \in \mathbf{F^a \textrm{\bf-Mod}(\hat{A}_{K})}$ (resp. 
$\mathcal{M} \in \mathbf{F^a\textrm{\bf-Modloc}(\hat{A}_{K})}$ ; resp. $\mathcal{M} \in \mathbf{F^a \textrm{\bf-Modlib}(\hat{A}_{K}))}$ muni de

$$\phi_{\mathcal{M}} : F^{\ast}_{\hat{A}_{K}}(\mathcal{M}) =: \mathcal{M}^{\sigma}\ \longrightarrow\ \mathcal{M}$$

\noindent est un $\mathbf{F^a\textrm{\bf-module convergent}}$ (resp. et \textbf{projectif} ; resp. et \textbf{libre}) sur $\hat{A}_{K}$.\\

Lorsque le Frobenius est un isomorphisme on dira que l'on a une \textbf{structure de Frobenius forte}. On a des notions analogues sur $A^{\dag}$ et $\hat{A}$, sans tensoriser par $K$ : on utilisera alors les notations $\mathbf{F^a\textrm{\bf-Mod}(A^{\dag})},... \mathbf{F^a\textrm{\bf-Modlib}(\hat{A})}$.\\

Soit $\Omega^1_{A^{\dag}}$ le module des $\mathcal{V}$-diff\'erentielles de $A^{\dag}$ au sens de Monsky-Washnitzer [M-W, theo 4.2]

$$\Omega^1_{A^{\dag}} :=  \Omega^1_{A^{\dag}/\mathcal{V}}\ /\ \displaystyle \mathop{\bigcap_{n}}\ \mathfrak{m}^n\ \Omega^1_{A^{\dag}/\mathcal{V}}\  ,$$

$$\Omega^1_{A^{\dag}_{K}} :=  \Omega^1_{A^{\dag}} \otimes_{\mathcal{V}} K\ ,\  \Omega^1_{\hat{A}} : = \widehat{\Omega^1_{\hat{A}/\mathcal{V}}}\ ,\ \Omega^1_{\hat{A}_{K}} := \Omega^1_{\hat{A}} \otimes_{\mathcal{V}} K.$$

\noindent Notons $\mathbf{F^a\textrm{\bf-Conn}^{\dag}(A^{\dag}_{K})}$ (resp. $\mathbf{F^a\textrm{\bf-Conn}(\hat{A}_{K}))}$ la cat\'egorie des $A^{\dag}_{K}$-modules (resp. $\hat{A}_{K}$-modules) projectifs de type fini $M$ (resp. $\mathcal{M})$ \`a connexion int\'egrable

$$\nabla : M \longrightarrow M \otimes_{A^{\dag}_{K}}  \Omega^1_{A^{\dag}_{K}}$$

$$( \mbox{resp.} \hat{\nabla} : \mathcal{M} \longrightarrow \mathcal{M} \otimes_{\hat{A}_{K}}  \Omega^1_{\hat{A}_{K}})$$

\noindent et munis d'un isomorphisme horizontal

$$\phi^{\dag} : (F^{\ast}_{A^{\dag}_{K}}(M), \ F^{\ast}_{A^{\dag}_{K}}(\nabla)) =: (M^{\sigma}, \nabla^{\sigma})\  \tilde{\longrightarrow}\ (M, \nabla) $$

$$(\mbox{resp.}\  \hat{\phi} : (F^{\ast}_{\hat{A}_{K}}(\mathcal{M}),\  F^{\ast}_{\hat{A}_{K}}(\hat{\nabla}) =: (\mathcal{M}^{\sigma}, \hat{\nabla}^{\sigma})\ \tilde{\longrightarrow}\ (\mathcal{M}, \hat{\nabla})    ).$$

\noindent On note $\mathbf{\textrm{\bf Conn}^{\dag}({A^{\dag}_{K}})}$ (resp. $\mathbf{\textrm{\bf Conn}\hat{ }(\hat{A}_{K})}$) la cat\'egorie des $A^{\dag}_{K}$-modules (resp. $\hat{A}_{K}$-modules) projectifs de type fini $M$ (resp. $\mathcal{M})$ \`a connexion int\'egrable dont la s\'erie de Taylor converge sur un voisinage strict du tube de la diagonale dans $\mathcal{X}^{\mbox{an}}_{K} \times \mathcal{X}^{\mbox{an}}_{K}$ (resp. dont la s\'erie de Taylor converge sur le tube de la diagonale dans $\mathcal{X}^{\mbox{an}}_{K} \times \mathcal{X}^{\mbox{an}}_{K}).$\\

\noindent On utilisera un exposant $(\ )^{\circ}$ pour sp\'ecifier les sous-cat\'egories des objets unit\'es (i.e. tels qu'en tout point g\'eom\'etrique les pentes du Frobenius sont nulles), $F^{a}$-Mod$(A^{\dag}_{K})^{\circ}$, $F^{a}$-Mod$(\hat{A}_{K})^{\circ}$, $F^{a}$-Conn$^{\dag}(A^{\dag}_{K})^{\circ}$, $F^{a}$-Conn$(\hat{A}_{K})^{\circ}$, ... ou la restriction de foncteurs aux objets unit\'es.\\

 On dispose de foncteurs naturels rendant commutatif le diagramme [Et 5, I, \S\ 5]

$$
\xymatrix{
F^{a}\mbox{-Isoc}^{\dag} (X/K) \ar[r]^{\Gamma^{\dag}} \ar[d]^{\mathcal{F}} & F^{a}\mbox{-Conn}^{\dag}(A^{\dag}_{K}) \ar[r] \ar[d]^{\mathcal{G}} & F^{a}\mbox{-Modloc}(A^{\dag}_{K}) \ar[d]^{\mathcal{H}}\\
F^{a}\mbox{-Isoc} (X/K) \ar[r]^{\hat{\Gamma}}& F^{a}\mbox{-Conn}{(\hat{A}}_{K}) \ar[r] & F^{a}\mbox{-Modloc} {(\hat{A}}_{K}), 
}
$$

\noindent o\`u $\Gamma^{\dag} := \Gamma(\mathcal{X}^{\mbox{an}}_{K}, -)$ est une \'equivalence de cat\'egories [B 3, cor (2.5.8)], $\hat{\Gamma} := \Gamma(]X[, -)$ est pleinement fid\`ele [O\ 2, 2.15, 2.23)] et [Et 5, I, (5.2.2)] et le foncteur $\mathcal{G}$ (resp. $\mathcal{H}$) envoie un ${A^{\dag}_{K}}$-module projectif de type fini $M$ sur son s\'epar\'e compl\'et\'e $p$-adique $\mathcal{M} = M \otimes_{A^{\dag}_{K}} \hat{A}_{K} : \mathcal{G}$ et $\mathcal{H}$ sont fid\`eles [Bour, A II, \S\ 5, $\mbox{n}^\circ\  3$ , prop 7] et [Bour, AC\ I, \S\ 3, $\mbox{n}^\circ\  5$, prop 9\ c)].\\

La restriction $\mathcal{F}^\circ$ de $\mathcal{F}$ \`a $F^a$-Isoc$^{\dag}(X / K)^\circ$ est un foncteur pleinement fid\`ele [Et 5, th\'eo 5]

\centerline{ $\mathcal{F}^{\circ} : F^a$-Isoc$^{\dag}  (X / K)^{\circ}\ \longrightarrow\ F^a$-Isoc$ (X / K)^{\circ}$ ;}

\noindent comme $\Gamma^{\dag \circ}$ est une \'equivalence de cat\'egories et $\hat{\Gamma}^{\circ}$ pleinement fid\`ele, le foncteur $\mathcal{G}^{\circ}$ est pleinement fid\`ele.\\

\noindent Nous allons montrer en 1.2 ci-apr\`es que $\hat{\Gamma}$ est en fait une \'equivalence de cat\'egories.\\

\subsection*{1.2. Des \'equivalences de cat\'egories}

Avec les notations de 1.1 nous allons montrer que la donn\'ee d'un $F^a$-isocristal convergent sur $X$ \'equivaut \`a celle d'un $\hat{A}_{K}$-module projectif de type fini $\mathcal{M}$, muni d'une connexion int\'egrable

$$\nabla : \mathcal{M}\  \longrightarrow\  \mathcal{M}\ \otimes_{\hat{A}_{K}}\ \Omega^1_{{\hat{A}_{K}}}$$

\noindent et d'un isomorphisme horizontal

$$\phi : (\mathcal{M}^{\sigma}, \nabla^{\sigma})\ \tilde{\longrightarrow}\ (\mathcal{M}, \nabla), $$

\noindent o\`u $(\mathcal{M}^{\sigma}, \nabla^{\sigma})$ provient de $(\mathcal{M}, \nabla)$ en \'etendant les scalaires par $F_{\hat{A}_{K}}$.

\vskip 3mm
\noindent \textbf{Proposition (1.2.1)}. \textit{Avec les notations ci-dessus, le foncteur $\Gamma(]X[, -)$ induit une \'equivalence entre  }

\begin{enumerate}
\item[(i)] \textit{La cat\'egorie des $\mathcal{O}_{]X[}$-modules coh\'erents (resp. et localement libres), et celles des $\hat{A}_{K}$-modules de type fini (resp. et projectifs) ; }
\item[(ii)] \textit{La cat\'egorie des  $\mathcal{O}_{]X[}$-modules coh\'erents \`a connexion int\'egrable (resp. des isocristaux convergents sur $X$), et celle des $\hat{A}_{K}$-modules projectifs de type fini munis d'une connexion int\'egrable (resp. et dont la s\'erie de Taylor converge sur $]X[_{\mathcal{X}^2}$).}

\end{enumerate}

\vskip 3mm
\noindent \textbf{Remarque}. Berthelot a fourni une description analogue pour les $j^{\dag} \mathcal{O}_{]X[-}$-modules coh\'erents [B 3, (2.5.2)].\\

\noindent \textit{D\'emonstration}. La d\'emonstration est semblable \`a celle de loc. cit. Remarquons simplement que $]X[$ \'etant affino¬\"{\i}de, la donn\'ee d'un $\mathcal{O}_{]X[}$-module coh\'erent $\mathcal{E}$ \'equivaut \`a celle du $\hat{A}_{K}$-module de type fini $\mathcal{M} = \Gamma(]X[, \mathcal{E})$ [B-G-R, 9.4.2]. De m\^eme l'assertion ``projectif" en (ii) r\'esulte du fait que $K$ est de caract\'eristique 0 [B 3, (2.2.3) (ii)] [P, 10.3.1]. $\square$

\vskip 3mm
\noindent \textbf{Th\'eor\`eme (1.2.2)}. \textit{Avec les notations de 1.1, soient $\mathcal{M}$ un $\hat{A}_{K}$-module de type fini, muni d'une connexion int\'egrable $\nabla$, et $(\mathcal{M}^{\sigma}, \nabla^{\sigma})$ le module \`a connexion int\'egrable d\'eduit de $(\mathcal{M}, \nabla)$ par l'extension des scalaires $F_{\hat{A}_{K}}$. On suppose qu'il existe un isomorphisme horizontal }

$$\phi : (\mathcal{M}^{\sigma}, \nabla^{\sigma})\ \tilde{\longrightarrow}\  (\mathcal{M}, \nabla).$$
\textit{Alors, si $(\mathcal{E, \nabla})$ est le $\mathcal{O}_{]X[}$-module correspondant \`a $(\mathcal{M}, \nabla)$ par l'\'equivalence de (1.2.1), la connexion $\nabla$ de $\mathcal{E}$ est convergente.}\\

\noindent \textit{D\'emonstration}. La preuve suit celle de [B 3, (2.5.7)]. Comme l'assertion est locale sur $X$ [B 3, (2.2.11)] on peut supposer $\Omega^1_{X}$ libre de base $d\overline{z}_{1},...,d\overline{z}_{m}$. Si $z_{1},...,z_{m}  \in A$ rel\`event $\overline{z}_{1},...\overline{z}_{m}$, $\Omega^1_{\hat{A}}$ est un $\hat{A}$-module libre de base $dz_{1},...,dz_{m}$. Soient $\partial_{1},...,\partial_{m}$ les d\'erivations correspondantes, et $e_{1},...,e_{r}$ une base de $\mathcal{M}$ ($\mathcal{M}$ est n\'ecessairement projectif d'apr\`es [P, lemme 10.3.1]). Pour tout $i$ et tout $k$, posons $\partial_{i}\ e_{k} = \displaystyle \mathop{\Sigma}_{j}\  b_{ijk}\ e_{j}$, et soit $B_{i} \in M_{r}(\hat{A}_{K})$ la matrice des $b_{i,j,k}$ ; la  connexion $\nabla$ est d\'etermin\'ee par les $B_{i}$. En utilisant l'isomorphisme horizontal $\phi$, on peut supposer, comme dans la preuve de [B 3, (2.5.7)] que les matrices $B_{i}$ sont \`a coefficients dans $\hat{A}$.\\

Notons $\eta_{0} =\  \mid \pi \mid$. Pour tout $e \in \mathcal{M} = \Gamma(]X[, \mathcal{E}),\  \underline{\partial}\ ^{\underline{k}}\  e$ est \`a coeffficients dans $\hat{A}$, donc $\parallel \underline{\partial}\ ^{\underline{k}}\  e \parallel\  \leqslant 1$, et comme $\eta_{0} < p^{-(1/p-1)}$ on en d\'eduit\\

\noindent (1.2.2.1) \qquad  $\parallel \frac{1}{\underline{k}!}\ \underline{\partial}\ ^{\underline{k}}\  e \parallel\ \eta_{0}^{\mid \underline{k} \mid}\ \longrightarrow 0$ quand $\mid \underline{k} \mid\ \longrightarrow + \infty .$\\

\noindent Posons $\zeta_{i} = 1 \otimes z_{i} - z_{i} \otimes 1\ ,\ \tau_{j} = 1 \otimes t_{j} - t_{j} \otimes 1.$\\

\noindent Pour $\eta < 1$, soient\\

$U = \{x \in \mathcal{X}^{an}_{K} \times \mathcal{X}^{an}_{K} /\  \forall j,~ \mid(1 \otimes t_{j}(x) \mid\  \leqslant 1, ~Ê\mid (t_{j} \otimes 1)(x) \mid\  \leqslant 1\} ,$\\

$W_{\eta} = \{x \in \mathcal{X}^{an}_{K} \times \mathcal{X}^{an}_{K} /\ \forall j,  \mid \tau_{j}(x) \mid\  \leqslant \eta \},$\\

$W'_{\eta} = \{x \in \mathcal{X}^{an}_{K} \times \mathcal{X}^{an}_{K} /\ \forall i,  \mid \zeta_{i}(x) \mid\  \leqslant 1 \},$\\

\noindent et $V_{\eta} = W_{\eta} \cap U$. Pour toute suite croissante $\underline{\eta}$ de limite 1, l'ouvert $V_{\underline{\eta}} :=  \displaystyle\mathop{\cup}_{{n}}\ V_{{\eta}_{n}}$  est \'egal \`a $] X [_{\mathcal{X}^2}$, d'apr\`es l'exemple de [B 3, (1.3.10)]. Nous allons construire une suite $\underline{\eta}$ telle qu'il existe sur l'ouvert $V_{\underline{\eta}} = ] X [_{\mathcal{X}^2}$ un isomorphisme $\epsilon : p^{\ast}_{2}\ \mathcal{E} \tilde{\longrightarrow} p^{\ast}_{1}\  \mathcal{E}$ induisant sur les voisinages infinit\'esimaux de la diagonale les isomorphismes $\epsilon_{n}$ d\'efinis par $\nabla$.\\

Puisque les $t_{j}$ engendrent $A$, il existe des relations $\zeta_{i} = \displaystyle \mathop{\Sigma}_{j}\ \beta_{ij}\  \tau_{j}$, avec $\beta_{ij} \in\ A \otimes_{\mathcal{V}} A$ ; d'o\`u $\parallel \zeta_{i} \parallel\  \leqslant \mbox{Sup}_{j} \parallel \tau_{j} \parallel$ , la norme \'etant la norme spectrale sur $U$. Par suite on a $V_{\eta_{0}} \subset W'_{\eta_{0}} \cap U$. On proc\`ede alors comme dans la d\'emonstration de [B 3, (2.2.13)] pour d\'efinir sur $V_{\eta_{0}}$ un isomorphisme $\varepsilon : p^{\ast}_{2}\ \mathcal{E}\  \tilde{\longrightarrow}\ 
p^{\ast}_{1}\ \mathcal{E}$ induisant sur les voisinages infinit\'esimaux de la diagonale les isomorphismes $\varepsilon_{n}$ d\'efinis par $\nabla$, en posant

$$\varepsilon(p^{\ast}_{2}(e)) = \displaystyle \mathop{\Sigma}_{\underline{k}}\ \frac{1}{\underline{k}!}\ \underline{\partial}^{\underline{k}}\ e \otimes \underline{\zeta}^{\underline{k}} ,$$

\noindent la s\'erie convergeant dans $\Gamma(V_{\eta_{0}}, p^{\ast}_{1} (\mathcal{E}))$ gr\^ace \`a (1.2.2.1).  On va utiliser ensuite l'action de Frobenius pour prolonger l'isomorphisme $\varepsilon$ de $V_{\eta_{0}}$ \`a $V_{\underline{\eta}}$, pour une suite $\underline{\eta}$ convenable : tout d'abord, on peut supposer d'apr\`es (1.2.1) qu'il existe un morphisme $F_{K} : \ ]X[ \ = \mbox{Spm} (\hat{A}_{K}) \longrightarrow\  ]X [$ tel que l'homomorphisme $\Gamma(]X[, \mathcal{O}_{]X[}) \longrightarrow\ \Gamma(]X[, \mathcal{O}_{]X[})$ induit par $F_{K}$ soit \'egal \`a $F_{\hat{A}_{K}} \otimes\  \sigma$. Le Frobenius $\phi$ de $\mathcal{E}$ est alors d\'efini par un isomorphisme $\phi : F^{\ast}_{K}\ \mathcal{E}\  \tilde{\longrightarrow}\   \mathcal{E}$ sur $]X[$, ce qui fournit des isomorphismes

$$\phi_{i} = p^{\ast}_{i}(\phi) : (F_{K} \times F_{K})^{\ast}\  (p^{\ast}_{i}\  \mathcal{E})\ \tilde{\longrightarrow}\ p^{\ast}_{i}\ \mathcal{E}\  \mbox{sur}\  V_{\eta_{0}} .$$

\noindent On d\'efinit une suite croissante $\underline{\eta}$ de limite 1 en posant

$$\eta_{n+1} = \eta^{1/p^{a}}_{n}\ .$$

\noindent Montrons que $(F_{K} \times F_{K}) (V_{\eta_{n+1}}) \subset V_{\eta_{n}}$. Posons

$$F_{\hat{A}}(t_{j}) = t^{p^{a}}_{j} + \pi\ a_{j}, $$

\noindent avec $a_{j} \in \hat{A}$. Pour $x \in U$ on a

$$\mid (t_{j} \otimes 1) ((F_{K} \times F_{K})(x)) \mid\  =\  \mid (F_{\hat{A}}(t_{j}) \otimes 1) (x) \mid $$

$$=\ \mid ((t^{p^ a}_{j} + \pi\ a_{j}) \otimes 1)(x) \mid\ \leqslant 1, $$

\noindent de m\^eme pour $1 \otimes t_{j}$, de sorte que  $(F_{K} \times F_{K})(x) \in U$. D'autre part, si $J = \mbox{Ker} (A \otimes_{\mathcal{V}} A \rightarrow A)$, on a dans $(A \otimes_{\mathcal{V}} A)^{\wedge}Ê$  la relation\\

\qquad  \qquad   ${(F_{\hat{A}} \times F_{\hat{A}}) (\tau_{j}) = 1 \otimes (t^{p^ a}_{j} +  \pi\  a_{j})\ -\ (t^{p^ a}_{j} +  \pi\  a_{j})\ \otimes\ 1} $

$\qquad \qquad \qquad \qquad \quad \quad ~ = \tau^{p^ a}_{j} + \pi\ \alpha_{j} ,$\\

\noindent avec $\alpha_{j} \in J(A \otimes_{\mathcal{V}} A)^{\wedge} $. Dans $(A \otimes_{\mathcal{V}} A)^{\wedge}$ on peut \'ecrire $\alpha_{j}$ sous la forme $\alpha_{j} = \displaystyle \mathop{\Sigma}_{i}\  \gamma_{ij}\  \tau_{i}$, avec $\| \gamma_{ij} \| \leqslant 1$. Alors, pour $x \in V_{\eta_{n+1}}$, on obtient

$$| \pi\ \alpha_{j}(x) | \leqslant | \pi |\  \eta_{n+1} \leqslant \eta_{n} $$

\noindent et

$$| \tau^{p^ a}_{j} (x) | \leqslant \eta^{p^a}_{n+1} = \eta_{n}.$$

\noindent Par cons\'equent on a bien

$$(F_{K} \times F_{K}) (V_{\eta_{n+1}}) \subset V_{\eta_{n}}. $$

\vskip 2mm
Supposons construit sur $V_{n} := \displaystyle \mathop{\cup}_{i \leqslant n} V_{\eta_{i}}$ un isomorphisme $\epsilon^{(n)} : p^{\ast}_{2}\  \mathcal{E} \tilde{\longrightarrow}\  p^{\ast}_{1}\ \mathcal{E}$ ; on conclut alors comme dans la preuve de [B 3, (2.5.7)] \`a l'existence d'un isomorphisme sur $V = \displaystyle \mathop{\cup}_{i } V_{\eta_{i}}  =\  ]X[_{\mathcal{X}^2}$, d'o\`u la convergence de $\nabla .~ \square$

\vskip 3mm
\noindent \textbf{Corollaire (1.2.3)}. \textit{Soient $X = \mbox{Spec}\  A_{0}$ un $k$-sch\'ema affine et lisse, $A$ une $\mathcal{V}$-alg\`ebre lisse relevant $A_{0}$, $\hat{A}$ le s\'epar\'e compl\'et\'e $p$-adique de $A$, $F_{\hat{A}} : \hat{A} \rightarrow \hat{A}$ un rel\`evement de la puissance $a^{\mbox{i\`eme}}$ de l'endomorphisme de Frobenius de $A_{0}$ au-dessus de $\sigma$. Alors la cat\'egorie des $F^{a}$-isocristaux convergents sur $X$ est \'equivalente \`a la cat\'egorie des $\hat{A}_{K}$-modules (n\'ecessairement projectifs) de type fini $\mathcal{M}$, munis d'une connexion int\'egrable $\nabla$ et d'un isomorphisme horizontal $\phi : (\mathcal{M}^{\sigma}, \nabla^{\sigma})\ \tilde{\longrightarrow}\ (\mathcal{M}, \nabla)$. }\\

\noindent \textit{D\'emonstration}. D'apr\`es la proposition (1.2.1) la cat\'egorie des isocristaux convergents sur $X$ est \'equivalente \`a celle des $\hat{A}_{K}$-modules projectifs de type fini, munis d'une connexion convergente $\nabla$. De plus cette \'equivalence est fonctorielle par rapport \`a $\hat{A}$ par des arguments analogues \`a [B 3, (2.5.6)]. On conclut par le th\'eor\`eme (1.2.2).\ $\square$\\

\section*{2. $F$-isocristaux convergents sur un sch\'ema lisse formellement relevable}

Nous allons g\'en\'eraliser au cas relevable l'\'equivalence de cat\'egories du corollaire (1.2.3) pr\'ec\'edent.\\

\subsection*{2.1. Espaces rigides associ\'es aux sch\'emas formels}

Dans ce \S\ 2 on se donne $f : X \rightarrow \mbox{Spec}\  k$ un $k$-sch\'ema lisse tel qu'il existe un $\mathcal{V}$-sch\'ema formel lisse $h : \mathcal{X} \rightarrow Spf\  \mathcal{V}$ relevant $f$. Par une construction de Raynaud on sait associer \`a $\mathcal{X}$ (resp. \`a $h$) un espace rigide analytique not\'e $\mathcal{X}_{K}$ [Bo-L¬\"u 1 et 2] [B 3] (resp. un morphisme $h_{K} : \mathcal{X}_{K}
\rightarrow \mbox{Spm}\ K)$. Si $h$ est propre alors $h_{K}$ est propre [L¬\"u]. L'espace $\mathcal{X}_{K}$ est muni d'une topologie de Grothendieck [B 3, (0.1.2), (0.2)] et d'un faisceau d'anneaux $\mathcal{O}_{\mathcal{X}_{K}}$ : nous dirons que $(\mathcal{X}_{K}, \mathcal{O}_{\mathcal{X}_{K}})$ est un $G$-espace annul\'e [B-G-R, 9.3.1].

\vskip 3mm
\noindent \textbf{Proposition (2.1.1)}. \textit{Sous les hypoth\`eses 2.1 il existe un $\mathcal{V}$-sch\'ema formel lisse $h' : \mathcal{X}' \rightarrow \mbox{Spf}\  \mathcal{V}$, un $\mathcal{V}$-isomorphisme $\mathcal{X'} \tilde{\longrightarrow} \mathcal{X}$ et des recouvrements par des ouverts lisses $\mathcal{X}' = \displaystyle \mathop{\cup}_{\alpha}\ \mbox{Spf}\ \hat{A}_{\alpha}$, $X = \displaystyle \mathop{\cup}_{\alpha}\ Spec\ A_{\alpha,0}$, o\`u les $A_{\alpha}$ sont des $\mathcal{V}$-alg\`ebres lisses et $A_{\alpha,0} := A_{\alpha} / \pi A_{\alpha}$.}\\

\noindent \textit{D\'emonstration}. Dire que $\mathcal{X}$ est un $\mathcal{V}$-sch\'ema formel lisse signifie qu'il existe un recouvrement $\mathcal{X} = \displaystyle \mathop{\cup}_{\alpha}\  \mbox{Spf}\  \mathcal{B}_{\alpha}$, o\`u les $\mathcal{B}_{\alpha}$ sont des $\mathcal{V}$-alg\`ebres plates s\'epar\'ees et compl\`etes pour la topologie $\pi$-adique et formellement lisses pour les topologies discr\`etes sur $\mathcal{V}$ et $\mathcal{B}_{\alpha} $ : les  $\mathcal{B}_{\alpha}$ sont donc des $\mathcal{V}$-alg\`ebres formellement lisses pour les topologies $\pi$-adiques sur $\mathcal{V}$ et $\mathcal{B}_{\alpha}$ [EGA O$_{IV}$, (19.3.1)]. Pour tout $\alpha$ notons $A_{\alpha,0} : = \mathcal{B}_{\alpha}\ /\ \pi\ \mathcal{B}_{\alpha}$, et $A_{\alpha}$ une $\mathcal{V}$-alg\`ebre lisse relevant $A_{\alpha,0}$ [E$\ell$, th\'eo 6] : d'apr\`es [Et 4, cor 1 du th\'eo 4] il existe, pour tout $\alpha$, un $\mathcal{V}$-isomorphisme

$$\hat{A}_{\alpha}\  \tilde{\longleftarrow}\ \mathcal{B}_{\alpha}, \mbox{o\`u}\  \hat{A}_{\alpha} := 
\displaystyle \mathop{\mbox{lim}}_{{\leftarrow} \atop{n}} A_{\alpha}\ /\  \pi^n\ A_{\alpha}.$$

\noindent D\'esignons par

$$P_{\alpha} = \mbox{Spf}\ \mathcal{B}_{\alpha}, \ P'_{\alpha} = \mbox{Spf}\ \hat{A}_{\alpha},\  P'_{\alpha \beta}\  := P'_{\alpha}\ \times_{P_{\alpha}} (P_{\alpha}\  \cap P_{\beta}),$$

\noindent $\psi_{\alpha}$ le $\mathcal{V}$-isomorphisme $\psi_{\alpha} : P'_{\alpha}\ \tilde{\rightarrow}\ P_{\alpha}$, et $\psi_{{\alpha}_{\beta}}$ le $\mathcal{V}$-isomorphisme 

$$P'_{\alpha \beta}\ \tilde{\rightarrow}\ P_{\alpha}\ \cap\ P_{\beta}$$

 \noindent d\'eduit de $\psi_{\alpha}$ par le changement de base $P_{\alpha}\ \cap\ P_{\beta}\ \hookrightarrow\ P_{\alpha}.$ Le $\mathcal{V}$-isomorphisme

$$\varphi_{{\alpha}_{\beta}} :  = \psi^{-1}_{\beta \alpha}\ \circ\  \psi_{{\alpha}_{\beta}} : P'_{\alpha \beta}\ \tilde{\rightarrow}\ P'_{\beta \alpha} $$

\noindent induit des $\mathcal{V}$-isomorphismes

$$\varphi_{{\alpha \beta \gamma}} : P'_{\alpha \beta}\ \cap\ P'_{{\alpha \gamma}}\ \tilde{\rightarrow}\  P'_{\beta \alpha}\ \cap\ P'_{\beta \gamma}$$

\noindent et on a les identit\'es

$$\varphi_{\alpha \beta}\ \circ\ \varphi_{\beta \alpha} = Id\ ,\ \varphi_{\alpha \alpha} = Id, $$

$$\varphi_{\alpha \beta \gamma} =  \varphi_{\gamma \beta \alpha} \circ\ \varphi_{\alpha \gamma \beta}\ ,$$

 \noindent gr\^ace au fait que les $P_{\alpha}$ se recollent pour former $\mathcal{X}$. Par cons\'equent on peut recoller les $P'_{\alpha}$ le long des $P'_{\alpha \beta}$ : le sch\'ema formel ainsi obtenu est le sch\'ema formel $\mathcal{X}'$ cherch\'e. $\square$

\vskip 3mm
\noindent \textbf{Corollaire (2.1.2)}. \textit{Sous les hypoth\`eses 2.1 le morphisme d'espaces rigides $h_{K} : \mathcal{X}_{K}\ \rightarrow\ \mbox{Spm}\ K$ est lisse. }\\

\noindent \textit{D\'emonstration}. Compte-tenu de l'isomorphisme $\mathcal{X}'\ \tilde{\rightarrow}\ \mathcal{X}$ de (2.1.1) il suffit d'appliquer le crit\`ere jacobien [B 3, (0.1.11)]. $\square$

\vskip 3mm
\noindent \textbf{Corollaire (2.1.3)}. \textit{Soient $S$ un $k$-sch\'ema lisse et $f : X\ \rightarrow\ S$ un $k$-morphisme lisse et supposons donn\'es un $\mathcal{V}$-sch\'ema formel lisse $\mathcal{S}$ et un $\mathcal{V}$-morphisme lisse $h : \mathcal{X}\ \rightarrow\ \mathcal{S}$ de sch\'emas formels relevant $f$. Alors }

\begin{enumerate}
\item[(i)] \textit{Il existe un $\mathcal{V}$-sch\'ema formel lisse $\mathcal{S}' = \displaystyle \mathop{\cup}_{\alpha}\ \mathcal{S}'_{\alpha}$\ , $\mathcal{S}'_{\alpha} = \mbox{Spf}\ \hat{A}_{\alpha}$ o\`u les $A_{\alpha}$ sont des $\mathcal{V}$-alg\`ebres lisses et un $\mathcal{V}$-isomorphisme $\mathcal{S}'\ \tilde{\rightarrow}\ \mathcal{S}$.}
\item[(ii)] \textit{Si l'on pose $\mathcal{X}' : = \mathcal{X} \times_{\mathcal{S}}\ \mathcal{S}'$, $\mathcal{X}'_{\alpha} : = \mathcal{X}' \times_{\mathcal{S}} \mathcal{S}'_{\alpha} $, il existe de plus un $\mathcal{S}'$-sch\'ema formel lisse $\mathcal{X}''$ et un $\mathcal{S}'$-isomorphisme $\mathcal{X}''\ \tilde{\rightarrow}\ \mathcal{X}'$ tel que}

$$\mathcal{X}'' = \displaystyle \mathop{\cup}_{\alpha} \mathcal{X}''_{\alpha}\ , \mbox{\textit{{o\`u}}}\  \mathcal{X}''_{\alpha} = \displaystyle \mathop{\cup}_{\beta}\ \mathcal{X}''_{\alpha, \beta} = \displaystyle \mathop{\cup}_{\beta}\  \mbox{Spf}\  \hat{B}_{\alpha, \beta}\ ,$$ 

\textit{les $B_{\alpha, \beta}$  \'etant des $\hat{A}_{\alpha}$-alg\`ebres lisses (resp. des $A^{\dag}_{\alpha}$-alg\`ebres lisses)}

\item[(iii)] \textit{Avec les notations du (ii) il existe aussi un $\mathcal{V}$-sch\'ema formel lisse $\mathcal{X}'''$ et un $\mathcal{V}$-isomorphisme $\mathcal{X}'''\ \tilde{\rightarrow}\ \mathcal{X}''$ tel que  }

$$\mathcal{X}'''\ = \displaystyle \mathop{\cup}_{\alpha}\ \mathcal{X}'''_{\alpha}\ ,\ \mbox{o\`u}\  \mathcal{X}'''_{\alpha} = \displaystyle \mathop{\cup}_{\beta}\ \mbox{Spf}\ \hat{C}_{\alpha, \beta}$$

\noindent \textit{et $ C_{\alpha,\beta}$ est une $\mathcal{V}$-alg\`ebre lisse munie d'un $\mathcal{V}$-isomorphisme}

$$\hat{B}_{\alpha, \beta} \simeq\ \hat{C}_{\alpha, \beta}.$$

\end{enumerate}

\noindent \textit{D\'emonstration}.
 La proposition (2.1.1.) fournit le (i).\\
 Pour (ii) et (iii) on utilise encore [Et 4, th\'eo 4 et son cor 1]: on peut recoller les $\mathcal{X}''_{\alpha,\beta}$ (resp. les $\mathcal{X}''_{\alpha}$, resp. les $\mathcal{X}'''_{\alpha})$ gr\^ace \`a l'existence globale de $\mathcal{X}'_{\alpha}$ (resp. de $\mathcal{X}'$, resp. de $\mathcal{X}'').\  \square$

\vskip 3mm
\noindent \textbf{Corollaire (2.1.4)}. \textit{Sous les hypoth\`eses (2.1.3) le morphisme d'espaces rigides $h_{K} : \mathcal{X}_{K} \rightarrow\mathcal{S}_{K}$ est lisse. Si de plus $h$ est propre alors $h_{K}$ est propre. }\\

\noindent \textit{D\'emonstration}. La lissit\'e de $h_{K}$ r\'esulte de (2.1.3) (ii) et du crit\`ere jacobien [B 3, (0.1.11]. \\
La d\'efinition d'un morphisme propre d'espaces rigides est donn\'ee dans [L¬\"u, 2.4] : la propret\'e de $h$ entra\^{\i}ne celle de $h_{K}$ [L¬\"u, theo 3.1].\ $\square$\\

\subsection*{2.2. F-isocristaux convergents et $\mathcal{O}_{\mathcal{X}_{K}}$-modules }
\vskip 2mm
Rappelons que les hypoth\`eses (2.1) sont satisfaites.\\

La donn\'ee de $\mathcal{E} \in \mbox{Isoc}(X/K)$ \'equivaut  \`a celle d'un $\mathcal{O}_{\mathcal{X}_{K}}$-module localement libre de type fini $\mathcal{E}_{\mathcal{X}} $ [B 3, (2.3.2), (2.2.3) (ii)] muni d'une connexion $\nabla$ relativement \`a $K$, int\'egrable et convergente. D'apr\`es (2.1.1) on a un $\mathcal{V}$-isomorphisme $\mathcal{X} \simeq \displaystyle \mathop{\cup}_{\alpha}\ \mbox{Spf}\ \hat{A}_{\alpha} = : \displaystyle \mathop{\cup}_{\alpha}\ \mathcal{X}_{\alpha}$ o\`u les $A
_{\alpha}$ sont des $\mathcal{V}$-alg\`ebres lisses ; si $F_{\alpha_{1}}, F_{\alpha_{2}}: \mathcal{X}_{\alpha} \rightarrow \mathcal{X}_{\alpha}$ sont deux rel\`evements de la puissance $a^{\mbox{i\`eme}}$ du Frobenius absolu de $X_{\alpha} = \mathcal{X}_{\alpha}\ \mbox{mod}\ \pi$, alors on a un isomorphisme canonique [B 3, (2.2.17)]

$$F^{\ast}_{\alpha_{1}}(\mathcal{E}_{\mathcal{X}_{\alpha}})\  \tilde{\longrightarrow}\ F^{\ast}_{\alpha_{2} }(\mathcal{E}_{\mathcal{X}_{\alpha}})\ , $$

\noindent o\`u $(\mathcal{E}_{\mathcal{X}_{\alpha}}, \nabla_{\alpha})$ est la restriction de $(\mathcal{E}_{\mathcal{X}}, \nabla)$ \`a $\mathcal{X}_{\alpha K} = \mbox{Spm}\ (\hat{A}_{\alpha K})$. \\

\noindent Supposons fix\'e pour chaque $\alpha$ un rel\`evement $F_{\alpha} : \mathcal{X}_{\alpha} \rightarrow \mathcal{X}_{\alpha}$, de la puissance $a^{\mbox{i\`eme}}$ du Frobenius absolu de $X_{\alpha} = \mathcal{X}_{\alpha}\  \mbox{mod}\  \pi$, et soit $\mathcal{E} \in F^a{\mbox{-}}\textrm{Isoc}(X/K)$. La structure de Frobenius sur $\mathcal{E}$ fournit pour tout $\alpha$ un isomorphisme [cor (1.2.3)]

$$
\phi_{\alpha} : (F^{\ast}_{\alpha}(\mathcal{E}_{\mathcal{X}_{\alpha}}), F^{\ast}_{\alpha}(\nabla_{\alpha}))  \displaystyle \mathop{\longrightarrow}^{\sim} (\mathcal{E}_{\mathcal{X}_{\alpha}}, \nabla_{\alpha}),
$$

\noindent avec compatibilit\'es \'evidentes quand $\alpha$ varie : d'apr\`es (1.2.3) la connaissance de $\mathcal{E}_{\mathcal{X}_{\alpha}}$ \'equivaut \`a celle de $\Gamma(\mathcal{X}_{\alpha K}, \mathcal{E}_{\mathcal{X}_{\alpha}})$ avec m\^emes donn\'ees.\\

Notons $\mathbf{F^a\textrm{\bf -Conn}(\mathcal{X}_{K})}$ la cat\'egorie des $\mathcal{O}_{\mathcal{X}_{K}}$-modules localement libres de type fini $\mathcal{M}$ munis d'une connexion $\nabla$ relativement \`a $K$, int\'egrable et d'une famille compatible d'isomorphismes

$$
\phi_{\alpha} : (F^{\ast}_{\alpha}(\mathcal{M}_{\alpha}), F^{\ast}_{\alpha}(\nabla_{\alpha}))  \displaystyle \mathop{\longrightarrow}^{\sim} (\mathcal{M}_{\alpha}, \nabla_{\alpha}) := (\mathcal{M}, \nabla)_{|\mathcal{X}_{\alpha K}}.
$$

\vskip 3mm
\noindent \textbf{Th\'eor\`eme (2.2.1)}. \textit{ Avec les notations de (2.2) la fl\`eche naturelle
$$
F^{a}{\mbox{-}}Isoc(X/K) \longrightarrow F^{a}{\mbox{-}}Conn(\mathcal{X}_{K})
$$
$$ \mathcal{E}  \longmapsto \mathcal{E}_{\mathcal{X}}$$
est une \'equivalence de cat\'egories.}\\

\noindent \textit{D\'emonstration}. Nous venons de montrer que si $\mathcal{E} \in F^{a}$-$Isoc(X/K)$ alors $\mathcal{E}_{\mathcal{X}} \in F^{a}\mbox{-}Conn(\mathcal{X}_{K})$.\\

R\'eciproquement soit $\mathcal{M} \in F^{a}\mbox{-}Conn(\mathcal{X}_{K})$ ; puisque $\mathcal{M}$ est localement libre de type fini et que $\mathcal{O}_{\mathcal{X}_{K}}$ est coh\'erent [B 3, (2.1.9)], alors $\mathcal{M}$ est un $\mathcal{O}_{\mathcal{X}_{K}}$-module coh\'erent.  D'apr\`es (1.2.3) l'existence des isomorphismes horizontaux

$$
\phi_{\alpha} : (F^{\ast}_{\alpha}(\mathcal{M}_{\alpha}), F^{\ast}_{\alpha}(\nabla_{\alpha})) \displaystyle \mathop {\longrightarrow}^{\sim} (\mathcal{M}_{\alpha}, \nabla_{\alpha})
$$
prouve que la connexion est convergente, car la donn\'ee de $(\mathcal{M}_{\alpha}, \nabla_{\alpha}, \phi_{\alpha})$ \'equivaut \`a celle d'un \'el\'ement de $F^{a}\mbox{-}Isoc(X_{\alpha}/K)$ : puisque ces donn\'ees se recollent la connexion $\nabla$ est convergente [B 3, (2.2.11)], d'o\`u un \'el\'ement $\mathcal{E}$ de Isoc$(X/K)$ [B 3, (2.3.2)] tel que $\mathcal{E}_{\mathcal{X}} = \mathcal{M}$. Enfin les $\phi_{\alpha}$ d\'efinissent un isomorphisme de Isoc$(X/K)$

$$
\phi_{\mathcal{E}} : F^{\ast}_{\sigma}\  \mathcal{E} \displaystyle \mathop {\longrightarrow}^{\sim} \mathcal{E}
$$
qui est le Frobenius de $\mathcal{E}$ [B 3, (2.3.7)] et $(\phi_{\mathcal{E}})_{\mathcal{X}} = \phi_{\mathcal{M}}. \ \square $\\

Rappelons que pour $\mathcal{E} \in Isoc(X/K)$ la cohomologie convergente de Ogus [O 3] et [B 5, (3.1.11) (i), (3.1.12) (ii)] est donn\'ee par

$$ H^i_{conv}(X/K, \mathcal{E}) = \mathbb{H}^i(\mathcal{X}_{K}, \mathcal{E}_{\mathcal{X}} \otimes \Omega^{\bullet}_{\mathcal{X}_{K}/K}). \leqno{(2.2.2)}$$

\section*{3. Images directes de $F$-isocristaux convergents}

Dans le cas d'une famille propre et lisse nous allons \'etendre ici le th\'eor\`eme 7.9 de Ogus [O 3] \'etablissant la finitude de la cohomologie convergente $H^i_{conv}(X/K, \mathcal{E})$ pour $X$ propre et lisse sur $k$ et $\mathcal{E}$ un $F$-isocristal convergent.\\

Sauf mention du contraire, on suppose dans ce paragraphe 3 que $k$ est un corps parfait de caract\'eristique $p>0, q = p^a, \mathcal{V}$ est un anneau de valuation discr\`ete complet, d'id\'eal maximal $\mathfrak{m}$ et corps r\'esiduel $k$. On suppose le corps des fractions $K$ de $\mathcal{V}$ de caract\'eristique 0,  on fixe une uniformisante $\pi$ et on note $e$ l'indice de ramification: on suppose, sauf mention expresse du contraire, que $e \leqslant p - 1$.\\

On rel\`eve la puissance $q$ sur $k$ en un automorphisme $\sigma$ de $\mathcal{V}$, tel que $\sigma(\pi) = \pi$, suivant la m\'ethode [Et 5, I 1.1] : on note encore $\sigma$ son extension \`a $K$ et ce $\sigma$ est un automorphisme de $K$ puisque $k$ est parfait [II,0].\\

\subsection*{3.1. Rel\`evement de Teichm¬\"uller}

Soit $X = \textrm{Spec}\ A_{0}$ un $k$-sch\'ema lisse. Pour $x \in\  \vert X \vert\  =  \{ \textrm{points ferm\'es de}\  X \}$, soit $i_{x} = Spec\ k(x) \hookrightarrow X$ l'immersion ferm\'ee canonique : $k(x) = A_{0}/\mathfrak{m}_{x}$ est une extension finie \'etale de $k$ de degr\'e deg $x = [k(x) : k]$. Notons $W = W(k)$ (resp. $W(x) = W(k(x))$ l'anneau des vecteurs de Witt \`a coefficients dans $k$ (resp. $k(x)$),

$$
\mathcal{V}(x) = W(x) \otimes_{W} \mathcal{V} \simeq W(x) [\pi]\ ,
$$

\noindent $K_{0} = \textrm{Frac}\ W$, $K_{0}(x) = \textrm{Frac}(W(x)), K(x) = \textrm{Frac}(\mathcal{V}(x))$, $\sigma_{x}$ la puissance $p^a$ sur $k(x)$, $\sigma_{W(x)} = W(\sigma_{x})$ le rel\`evement canonique de $\sigma_{x}$ \`a $W(x)$, $\sigma_{\mathcal{V}(x)} = \sigma_{W(x)} \otimes_{W} \mathcal{V}$ et $\sigma_{K(x)}$ (resp. $\sigma_{K_{0}(x)})$ son extension naturelle \`a $K(x)$ (resp. $K_{0}(x))$ d\'efinie par $\sigma_{K(x)}(u/v) = \sigma_{\mathcal{V}(x)}(u)/\sigma_{\mathcal{V}(x)}(v)$ (resp. $\sigma_{K_{0}(x)}(u/v) = \sigma_{W(x)}(u)/\sigma_{W(x)}v))$. Le morphisme $\sigma_{K(x)}$ co¬\"{\i}ncide, d'apr\`es [Et 5, I.1.1] et [B-M 2, (1.2.7) (ii)], avec le morphisme $\sigma' : K' \rightarrow K'$ (au-dessus de $\sigma : K \rightarrow K)$ de [Et 5, I.1.1] pour $k' = k(x)$.\\

Soit $A$ une $\mathcal{V}$-alg\`ebre lisse relevant $A_{0}$ et fixons une pr\'esentation $A = \mathcal{V}[t_{1},...,t_{n}]/(f_{1},...,f_{m})$. On d\'esigne par $\hat{A}$ ((resp. $A^{\dag}$) le s\'epar\'e compl\'et\'e (resp. le compl\'et\'e faible) $\mathfrak{m}$-adique de $A$. Soient $P$ le compl\'et\'e formel de la fermeture projective de $\mathcal{X} = Spec\ A$ dans $\mathbb{P}^n_{\mathcal{V}}$, $P_{1}
= Spf(\mathcal{V}(x))$, $X_{1} = Spec(k(x))$. D'apr\`es [Et 5, (1.2.1)] il existe un carr\'e commutatif

$$
\xymatrix{
A^{\dag} \otimes_{\mathcal{V}} \mathcal{V}(x) =: A^{\dag}(x)  \ar[rr]^{\qquad \quad F_{A^{\dag}(x)}} & & A^{\dag}(x) \\
A^{\dag} \ar[u] \ar[rr]^{F_{A^{\dag}}}   & & A^{\dag}  \ar[u]
}
$$

\noindent au-dessus du carr\'e commutatif

$$
\xymatrix{
\mathcal{V}(x)  \ar[r]^{\sigma_{\mathcal{V}(x)}}  & \mathcal{V}(x) &\\
\mathcal{V} \ar[u] \ar[r]^{\sigma_{\mathcal{V}}}  & \mathcal{V} \ar[u] &,
}
$$

\noindent o\`u $F_{A^{\dag}}$ est un rel\`evement \`a $A^{\dag}$ du Frobenius (puissance $q$) de $A_{0}$ ; d'o\`u un diagramme commutatif

$$
\xymatrix{
\mathcal{V}(x)  \ar[r] \ar[d]_{\sigma_{\mathcal{V}(x)}} & A^{\dag}(x) \ar[r] \ar[d]_{F_{A^{\dag}(x)}}  & \hat{A}(x) := \hat{A} \otimes_{\mathcal{V}} \mathcal{V}(x) \ar[d]^{F_{\hat{A}(x)}} &\\
\mathcal{V}(x) \ar[r]   & A^{\dag}(x) \ar[r] & \hat{A}(x) & .
}
$$

Par cons\'equent le morphisme $s : A_{0} \rightarrow k(x)$ se rel\`eve de mani\`ere unique d'apr\`es Katz [K 1] en un morphisme

$$
\tau(x) : \hat{A}(x) \rightarrow \mathcal{V}(x)
$$

\noindent tel que le diagramme

$$
\xymatrix{
\hat{A}(x)  \ar[r]^{\tau (x)} \ar[d]_{F_{\hat{A}(x)}} & \mathcal{V}(x) \ar[d]^{\sigma_{\mathcal{V}(x)}}\\
\hat{A}(x)  \ar[r]^{\tau (x)}  & \mathcal{V}(x)
}
$$

\noindent commute : $\tau (x)$ est appel\'e le \textbf{rel\`evement de Teichm¬\"uller de $s$ (ou de $x$)}. Les morphismes compos\'es

$$
\hat{\tau}(x) : \hat{A} \hookrightarrow \hat{A}(x) \displaystyle \mathop{\longrightarrow}^{\tau (x)} \mathcal{V}(x)\ ,
$$

$$
\tau^{\dag}(x) : A^{\dag} \hookrightarrow \hat{A} \displaystyle \mathop{\longrightarrow}^{\hat{\tau} (x)} \mathcal{V}(x)
$$

\noindent sont surjectifs car la r\'eduction de $\tau^{\dag}(x) \ \textrm{mod}\  \pi$ est le morphisme surjectif $s : A_{0} \rightarrow k(x)$ de d\'epart [M-W, theo 3.2] ; donc $\mathcal{V}(x)$ est un quotient de $\hat{A}$ et $\mathcal{V}(x) \simeq W(x) [\pi]$, qui est un anneau de valuation discr\`ete, est une extension finie \'etale de $\mathcal{V}$ de rang deg $x$. Le noyau du morphisme

$$
\xymatrix{
\hat{\tau}_{K}(x) := \hat{\tau}(x) \otimes_{\mathcal{V}} K : \hat{A}_{K} \ar@{->>}[r] & K(x) = \textrm{Frac} (\mathcal{V}(x))
}
$$

\noindent est ainsi un id\'eal maximal $\mathfrak{q}_{x}$ de $\hat{A}_{K}$ et le diagramme

$$
\begin{array}{c}
\xymatrix{
A^{\dag}_{K} \ar@{^{(}->}[r] \ar[d]_{F_{A^{\dag}_{K}}} & \hat{A}_{K} \ar[d]_{F_{\hat{A}_{K}}} \ar@{->>}[r]^{\hat{\tau}_{K }(x)} & K(x) \ar[d]^{\sigma_{K (x)}}\\
A^{\dag}_{K} \ar@{^{(}->}[r] & \hat{A}_{K} \ar@{->>}[r]^{\hat{\tau}_{K} (x)} & K(x)
}
\end{array}
\leqno{(3.1.1)}
$$

\noindent commute. On notera $\tau^{\dag}_{K}(x)$ la fl\`eche compos\'ee
$$
\tau^{\dag}_{K}(x) = \tau^{\dag}(x) \otimes_{\mathcal{V}} K :  A^{\dag}_{K} \displaystyle \mathop{\longrightarrow}^{\varphi}  \hat{A}_{K} \displaystyle -\hspace{-10pt}-\hspace{-10pt}-\hspace{-10pt}-\hspace{-10pt}-\hspace{-10pt} \mathop{\twoheadrightarrow}^{\hspace{-10pt}\hat{\tau}_{K }(x)}  K(x)\ ; \leqno{(3.1.2)}
$$\\
remarquons que $\tau^{\dag}_{K}(x)$ est aussi surjectif  car $\varphi$ induit une bijection entre les id\'eaux maximaux de $\hat{A}_{K}$ et ceux de $A^{\dag}_{K}$ [G-K 2, theo 1.7]. Par le morphisme de sp\'ecialisation [B 3, (0.2.2.1)]

$$
sp : Spm\  \hat{A}_{K} \rightarrow Spec\  A_{0}
$$

\noindent l'image de $\{ \mathfrak{q}_{x} \}$ n'est autre que $\{ \mathfrak{m}_{x} \}$ ; de plus $\hat{\tau}(x) : \hat{A} \twoheadrightarrow \mathcal{V}(x)$ est localis\'e en $\{ \mathfrak{p}_{x} \} \in Spec\ \hat{A}$, o\`u
$\mathfrak{p}_{x}$ est l'unique id\'eal maximal de $\hat{A}$ au-dessus de $\mathfrak{m}_{x}$. En d\'efinissant l'application (encore appel\'ee \textbf{rel\`evement de Teichm¬\"uller})\\

\noindent (3.1.3) $\qquad \qquad \qquad \qquad \hat{T}_{K} : Spm\ A_{0} \rightarrow Spm\ \hat{A}_{K}$\\

\noindent par $\hat{T}_{K}(x) = \{ \textrm{Ker}\ \hat{\tau}_{K}(x) \} = \{ \mathfrak{q}_{x} \}$, on vient de prouver que $\hat{T}_{K}$ est \textbf{une section du  morphisme de sp\'ecialisation}, consid\'er\'e comme une application

$$
sp : \vert Spm\ \hat{A}_{K} \vert\ \rightarrow \vert Spec\ A_{0} \vert .
$$
\vskip 3mm

\subsection*{3.2. Amplitude des $F$-isocristaux}

Soit $X$ un $k$-sch\'ema lisse. Berthelot a montr\'e [B 3, (2.4.2)] que pour tout $N \in F^a\mbox{-}Isoc(X/K)$ il existe un $F$-cristal non d\'eg\'en\'er\'e $M$ sur $X$ (i.e. un cristal $M$ muni d'un Frobenius $\phi : F^{\ast} M \rightarrow M$ et d'un morphisme $V : M \rightarrow F^{\ast} M$ tels que $\phi \circ V$ et $V \circ \phi$ soient la multiplication par $p^b$, pour un entier $b \geqslant 0$) de type fini sur $\mathcal{O}_{X/\mathcal{V}}$ et sans $p$-torsion, et un entier $r \geqslant 0$ tel que $N \simeq M^{an}(r)$ (o\`u $M^{an}(r)$ est le twist \`a la Tate de $M^{an}$). L'entier $b$ est appel\'e l'amplitude de $M$ (``width'' dans la terminologie de Ogus [O 4, 5.1.1]). Notons\\

\noindent (3.2.1) $\qquad \qquad \qquad \qquad  F^{a}\mbox{-}Isoc(X/K)_{plat}$ \\

\noindent la sous-cat\'egorie pleine de  $F^{a}\mbox{-}Isoc(X/K)$ form\'ee des $N$ tels que le $M$ ci-dessus soit plat sur $\mathcal{O}_{X/\mathcal{V}}$ (ce qui \'equivaut \`a $M$ localement libre puisqu'on est en situation noeth\'erienne).\\

\subsection*{3.3. Convergence des images directes}

\vskip 3mm
\noindent \textbf{Th\'eor\`eme (3.3.1)}. \textit{Supposons $k$ parfait et $e\leqslant p-1$. Soient $S$ un $k$-sch\'ema lisse et $f : X \rightarrow S$ un $k$-morphisme propre et lisse (resp. relevable en un morphisme propre et lisse $h : \mathcal{X} \rightarrow \mathcal{S}$ de $\mathcal{V}$-sch\'emas formels). Alors}\\
\begin{enumerate}
\item[(3.3.1.1)] \textit{Pour tout entier $i \geqslant 0$, $f$ induit un foncteur 
$$
R^{i} f_{conv^{\ast}} : F^{a}\mbox{-}Isoc(X/K)_{plat} \longrightarrow F^{a}\mbox{-}Isoc(S/K)
$$
$$
(resp. \ R^{i} f_{conv^{\ast}} : F^{a}\mbox{-}Isoc(X/K) \longrightarrow F^{a}\mbox{-}Isoc(S/K)).
$$}
\item[ (3.3.1.2)] \textit{Le foncteur pr\'ec\'edent est compatible aux changements de base entre $k$-sch\'emas lisses, c'est-\`a-dire : pour tout carr\'e cart\'esien}
$$
\xymatrix{
X' \ar[r]^{g'} \ar[d] _{f'} & X \ar[d]^{f} \\
S' \ar[r]_{g} & S
}
$$
\textit{avec $S'$ lisse sur $k$ et $\mathcal{E} \in F^{a}\mbox{-}Isoc(X/K)_{\mbox{plat}}$
(resp. $\mathcal{E} \in F^{a}\mbox{-}Isoc(X/K))$ on a un isomorphisme de changement de base
$$
g^{\ast} R^{i} f_{conv^{\ast}}(\mathcal{E)}\ \displaystyle \mathop{\longrightarrow}^{\sim}   R^{i} f'_{conv^{\ast}}(g'^{\ast}(\mathcal{E}))
$$
compatible aux connexions et aux Frobenius.}
\end{enumerate}
\newpage
\noindent\textit{D\'emonstration}.\\

\noindent\textit{Premi\`ere partie de la d\'emonstration: le cas propre et lisse, non n\'ecessairement relevable}.\\

Dans un premier temps nous ne supposerons pas que le corps $k$ est parfait, seulement qu'il est de caract\'eristique $p>0$ et que $e\leqslant p-1$: uniquement lorsqu'il faudra montrer que le Frobenius est un isomorphisme, nous supposerons que $k$ est parfait.\\

 Soient $\mathcal{E} \in F^{a}\mbox{-}Isoc(X/K)_{plat}$ et $E$ un $F$-cristal non d\'eg\'en\'er\'e tel que $\mathcal{E} \simeq E^{an}(r)$. On va associer \`a $E$, $f$ et l'entier $i \geqslant 0$,  un objet de $F^{a}\mbox{-}Isoc(S/K)$, i.e. pour tout ouvert $U$ (affine et lisse) de $S$ on va construire un \'el\'ement de $F^{a}\mbox{-}Isoc(U/K)$ avec donn\'ees de recollement [B 3, (2.3.2) (iii)].\\

 Si $\mathcal{T}$ est un $\mathcal{V}$-sch\'ema formel on pose $T_{n} = \mathcal{T}/ \pi^{n+1}\ \mathcal{T}$, et pour un $T_{0}-$sch\'ema de type fini $Y$ on d\'esigne par $(Y/\mathcal{T})_{\textrm{cris}}$ le topos cristallin [B-O, 7.17] [O 4, \S\  3.0]. Le morphisme $f$ induit un morphisme de topos 
 
 $$
 f_{\mbox{cris}} : (X/\mathcal{V})_{\mbox{cris}} \longrightarrow (S/\mathcal{V})_{\mbox{cris}}\ .
 $$
 
 \vskip 3mm
 Soit $U = \mbox{Spec} A_{0} \hookrightarrow S$ un ouvert affine et lisse et $A$ une $\mathcal{V}$-alg\`ebre lisse relevant $A_{0}$ : posons $\mathcal{T} = \mbox{Spf}\  \hat{A}$, $T_{n} = \mbox{Spec}\  (A/ \pi^{n+1}\ A)$. Alors $(U, T_{n}, \delta)$, o\`u $\delta$ sont les puissances divis\'ees canoniques sur 
$\pi\  \mathcal{O}_{T_{n}} (e \leqslant p-1$ : [B 1, I, \S\ 1]), est un ouvert de Cris $(S/ \mathcal{V})$. Le Frobenius $F_{U}$ de $U$ (\'el\'evation \`a la puissance $p^a$ sur $\mathcal{O}_{U})$ se rel\`eve en un endomorphisme $F_{\mathcal{T}}$ de $\mathcal{T}$ (resp. $F_{T_{n}} := F_{\mathcal{T}}\  \mbox{mod}\ \pi^{n+1}$, de $T_{n}$) et on consid\`ere les diagrammes commutatifs 
 
 $$
 \begin{array}{c}
\xymatrix{
X'_{U} \ar[r]^{\varphi} \ar[d]_{f'_{U}} & X_{U} \ar@{^{(}->}[r] \ar[d]^{f_{U}} & X \ar[d]^{f} & \\
U \ar[r]^{F_{U}} \ar[d]_{i_{T_{n}}} & U \ar@{^{(}->}[r] \ar[d]^{i_{T_{n}}} & S \ar[d] &\\
(T_{n}, \pi, \delta) \ar[r]_{F_{T_{n}}} & (T_{n}, \pi, \delta) \ar[r] & \mbox{Spec}\ \mathcal{V}_{n} \ ,
}
\end{array}
\leqno{(3.3.1.3)}
$$

 $$
 \begin{array}{c}
\xymatrix{
X'_{U} \ar[r]^{\varphi} \ar[d]_{f'_{U}} & X_{U} \ar@{^{(}->}[r] \ar[d]^{f_{U}} & X \ar[d]^{f} &\\
U \ar[r]^{F_{U}} \ar[d]_{i_{\mathcal{T}}} & U \ar@{^{(}->}[r] \ar[d]^{i_{\mathcal{T}}} & S \ar[d]  \\
(\mathcal{T}, \pi, \delta) \ar[r]_{F_{\mathcal{T}}} & (\mathcal{T}, \pi, \delta) \ar[r] & \mbox{Spf}\ \mathcal{V}\ ,
}
\end{array}
\leqno{(3.3.1.4)} 
$$

\noindent dont les carr\'es sup\'erieurs sont cart\'esiens. Soient $f_{X_{U}/T_{n}}$ [B 1, V, 3.5.2] et $f_{X_{U}/\mathcal{T}}$ [B-O, \S\ 7.21] les morphismes de topos

$$
f_{X_{U}/T_{n}} : (X_{U}/T_{n})_{\mbox{cris}} \displaystyle \mathop{\longrightarrow}^ {u_{X_{U}/T_{n}}} X_{U_{Zar}} \displaystyle \mathop{\longrightarrow}^ {i_{T_{n}}\circ f_{U}} T_{n_{Zar}}
$$
 
 $$
f_{X_{U}/\mathcal{T}} : (X_{U}/\mathcal{T})_{\mbox{cris}} \displaystyle \mathop{\longrightarrow}^ {u_{X_{U}/\mathcal{T}}} X_{U_{Zar}} \displaystyle \mathop{\longrightarrow}^ {i_{\mathcal{T}}\circ f_{U}} \mathcal{T}_{Zar}.
$$\\
Avec les notations de [B 1, V, 3.2.3] on a

$$ \mathbb{R} f_{\mbox{cris}^{\ast}} (E)_{(U,T_{n})} = \mathbb{R} f_{X_{U}/T_{n}^{\ast}} (\omega^{\ast}_{T_{n}}(E) ) \leqno{(3.3.1.5)}$$ 

\noindent et [B-O, 7.2.2.2] :

$$  \mathbb{R} f_{X_{U}/\mathcal{T}^{\ast}}  (\omega^{\ast}_{\mathcal{T}}(E))\ \simeq\ \mathbb{R} \displaystyle \mathop{\mbox{lim}}_{\leftarrow \atop{n}} \mathbb{R}  f_{X_{U}/T_{n}^{\ast}} (\omega^{\ast}_{T_{n}}(E)). \leqno{(3.3.1.6)}$$

Posons

$$
\displaystyle \widetilde{\mathop{R^i f_{\textrm{cris}}(E)_{(U,\mathcal{T})}}} := \displaystyle \mathop{\mbox{lim}}_{\leftarrow \atop{n}}\ R^{i}  f_{X_{U}/T_{n}^{\ast}} (\omega^{\ast}_{T_{n}}(E)) . \leqno{(3.3.1.7)}
$$

Le sch\'ema formel $\mathcal{T}$, identifi\'e \`a la limite inductive $\{T_{n}\}$ des $T_{n}$, est ce que Ogus appelle un ``\'epaississement fondamental de $U$ relativement \`a Spf $\mathcal{V}$'' [O 4, \S\ 3.0].\\

Puisque $E$ est localement libre de type fini et $f$ propre et lisse, le complexe $\mathbb{R}  f_{X_{U}/T_{n}^{\ast}} (\omega^{\ast}_{T_{n}}(E))$ est un complexe parfait de $\mathcal{O}_{{T}_{n}}$-modules [B 1, VII, 1.1.1] ; or, par le th\'eor\`eme de changement de base en cohomologie cristalline [B 1, V, 3.5.2] on a un isomorphisme

$$
\mathbb{R}f_{X_{U}/T_{n}^{\ast}}  (\omega^{\ast}_{T_{n}}(E)) \displaystyle \mathop{\otimes}^{\mathbb{L}}{_{\mathcal{O}_{{T}_{n}}}} \mathcal{O}_{T_{n-1}} \displaystyle \mathop{\longrightarrow}^{\sim}  \mathbb{R} f_{X_{U}/T_{n-1}^{\ast}}(\omega^{\ast}_{T_{n-1}}(E)) :
$$
par suite [B-O, def B-4] $\mathbb{R} f_{X_{U}/T_{\bullet}^{\ast}}(\omega^{\ast}_{T_{\bullet}}(E))$ est un objet ''consistant'' au sens de [loc. cit.], donc [B-O, prop B-7]

$$R^{i}  f_{X_{U}/\mathcal{T}^{\ast}} (\omega^{\ast}_{\mathcal{T}}(E)) =  H^{i}(\mathbb{R} f_{X_{U}/\mathcal{T}^{\ast}} (\omega^{\ast}_{\mathcal{T}}(E)))$$

$\qquad  \qquad  \qquad \qquad  \qquad  \qquad \qquad \simeq H^{i} (\mathbb{R} \displaystyle \mathop{\mbox{lim}}_{\leftarrow \atop{n}} \mathbb{R} f_{X_{U}/T_{n}^{\ast}} (\omega^{\ast}_{T_{n}}(E)))$ 

$\qquad  \qquad  \qquad \qquad  \qquad  \qquad \qquad  \simeq \displaystyle \mathop{\mbox{lim}}_{\leftarrow \atop{n}} R^{i} f_{X_{U}/T_{n}^{\ast}} (\omega^{\ast}_{T_{n}}(E))\ ,$

\noindent d'o\`u, via (3.3.1.7), et [B-O, prop B-7], un isomorphisme de $\mathcal{O}_{\mathcal{T}}$-modules coh\'erents

$$
\displaystyle \widetilde{\mathop{R^{i} f_{\textrm{cris}}(E)_{(U,\mathcal{T})}}} \simeq \ R^{i}  f_{X_{U}/\mathcal{T}^{\ast}} (\omega^{\ast}_{T}(E)) . \leqno{(3.3.1.8)}
$$

Soit $V_{\bullet} = (V_{\alpha})_{\alpha}\ \xrightarrow{u_{\bullet}=(u_{\alpha})_{\alpha}} X_{U}$ un hyperrecouvrement fini de $X_{U}$ par des ouverts affines,  $V_{\alpha} = Spec(B_{\alpha, 0})$ : soient $B_{\alpha}$ une $\hat{A}$-alg\`ebre lisse relevant $B_{\alpha, 0}$, $\mathcal{P}_{\alpha} := Spf\ (\hat{B}_{\alpha})$,\\

$h_{\alpha, n} : P_{\alpha, n} = Spec(B_{\alpha}/\pi^{n+1}\ B_{\alpha}) \longrightarrow T_{n},$\\

$
h_{\bullet, n} : P_{\bullet, n} \longrightarrow T_{n}, $\\

$h_{\alpha} : \mathcal{P}_{\alpha} \longrightarrow \mathcal{T}$ la limite inductive des $\{ h_{\alpha, n} \}_{n}$ et posons

$$
E_{P_{\alpha, n}} = u_{\alpha}^{\ast} (\omega^{\ast}_{T_{n}}(E)),\ E_{P_{\bullet, n}} = u_{\bullet}^{\ast} (\omega^{\ast}_{T_{n}}(E))\ ,
$$

$$
E_{\mathcal{P}_{\alpha}} = \displaystyle \mathop{\lim}_{\leftarrow \atop{n}}\ E_{P_{\alpha, n}} \quad , \quad E_{\mathcal{P}_{\bullet}}
= \{E_{\mathcal{P}_{\mathcal{\alpha}}} \}_{\alpha}\ ,
$$

$$
h_{\bullet} = \{ h_{\alpha} \}_{\alpha} : \mathcal{P}_{\bullet} = \{ \mathcal{P}_{\alpha} \}_{\alpha} \longrightarrow\ \mathcal{T}\ ,
$$

$$
C_{\alpha, n}^{\bullet} = E_{P_{\alpha, n}} \otimes \Omega_{P_{\alpha, n} /T_{n}
}^{\bullet}\ ,
$$

$$
 \Omega^{\bullet}_{P_{\bullet} /\mathcal{T}} = \{ \Omega^{\bullet}_{P_{\alpha} /T} \}_{\alpha}\ ,
 $$
 
 \noindent et remarquons que l'on a un isomorphisme
 
 $$
 \Omega^{\bullet}_{\mathcal{P}_{\alpha} /\mathcal{T}} \simeq \displaystyle \mathop{\lim}_{\leftarrow \atop{n}}\ \Omega_{P_{\alpha, n} /T_{n}}^{\bullet}.
 $$
 
 \noindent D'apr\`es (I$\ell$ 1, (0.3.2.6.2), (0.3.2.2), (0.3.2.4)] il existe un isomorphisme
 
 $$
 \mathbb{R} f_{X_{U}/T^{\ast}_{n}} (\omega^{\ast}_{T_{n}}(E)) \simeq \mathbb{R} h_{\bullet,n^{\ast}} (E_{P_{\bullet,n}} \otimes 
  \Omega^{\bullet}_{P_{\bullet, n} /T_{n}}) \ , 
  $$
  
  \vskip 2mm
  \noindent et avec les notations de [B-O, Appendix B] on a un diagramme commutatif

 $$
  \xymatrix{
  K(\mathbb{N},\mathcal{O}_{P_{\alpha, \bullet}}) \ar[rrr]^{Rh_{\alpha \bullet^{\ast}} =  \{ Rh_{\alpha, n^{\ast}} \}_{n}} \ar[dd]_{R \Gamma(\mathbb{N},-)} & &  &  K(\mathbb{N},\mathcal{O}_{T_{\bullet}})  \ar[dd]^{R \Gamma(\mathbb{N}, -)} &\\
\\
  K(\mathcal{O}_{\mathcal{P}_{\alpha}}) \ar[rrr]^{Rh_{\alpha \ast}} & & & K(\mathcal{O_{\mathcal{T}}}) & ,
 } 
$$

\noindent puisqu'il en est ainsi avant de passer aux foncteurs d\'eriv\'es [Et 2, d\'em. de III, 3.1.1]. Or le complexe $\{ C^{\bullet}_{\alpha, n} \}_{n} \in K (\mathbb{N}, \mathcal{O}_{P_{\alpha, \bullet}})$ a des fl\`eches de transition surjectives, donc v\'erifie la condition de Mittag-Leffler ; par suite

$$
\mathbb{R} \displaystyle \mathop{\lim}_{\longleftarrow \atop{n}} \mathbb{R} h_{\alpha, n^{\ast}} (C^{\bullet}_{\alpha, n}) = \mathbb{R}h_{\alpha^{\ast}} (E_{\mathcal{P}_{\alpha}} \otimes \Omega^{\bullet}_{P_{\alpha} /\mathcal{T}} )
$$

\noindent et donc

$$
\displaystyle \widetilde{\mathop{R^{i} f_{\textrm{cris}^{\ast}}(E)_{(U,\mathcal{T})}}} \simeq  R^{i}  f_{X_{U}/ \mathcal{T}^{\ast}} (\omega^{\ast}_{\mathcal{T}}(E)) \simeq H^{i}(\mathbb{R} h_{\bullet \ast} (E_{\mathcal{P}_{\bullet}} \otimes \Omega^{\bullet}_{\mathcal{P}_{\bullet} /\mathcal{T}} )). \leqno{(3.3.1.9)}
$$

Montrons \`a pr\'esent que  $\displaystyle \widetilde{\mathop{R^{i} f_{\textrm{cris}^{\ast}}(E)_{(U,\mathcal{T})}}}$ est muni d'un morphisme de Frobenius

$$
F^{\ast}_{\mathcal{T}} (\displaystyle \widetilde{\mathop{R^i f_{\textrm{cris}^{\ast}}(E)_{(U,\mathcal{T})}}}) \longrightarrow \displaystyle \widetilde{\mathop{R^i f_{\textrm{cris}^{\ast}}(E)_{(U,\mathcal{T})}}}\ . 
$$

\noindent Consid\'erons la factorisation usuelle du Frobenius $F_{X}$ (avec $F_{X}^{\ast}(x)=x^{q}$)

$$
\xymatrix{X \ar@/^1pc/[rrd]^{F_{X}} \ar@/_/[rdd]_{f}  \ar@{.>}[rd]^{F_{X/S}} \\
&  X' \ar[d]^{f'} \ar[r]^{\pi_{X/S}} & X \ar[d]^f\\
& S \ar[r] ^{F_{S}} & S
}
\leqno{(3.3.1.10)}
$$

\noindent o\`u le carr\'e est cart\'esien ; le morphisme de Frobenius de $E$, $\phi : F^{\ast}_{X/S}\  E' =F^{\ast}_{X}\  E \rightarrow E$, o\`u $E' := \pi^{\ast}_{X/S}(E)$, d\'efinit par image inverse 
$$
\omega^{\ast}_{T_{n}}(\phi) : F^{\ast}_{X_{U}/U}\  \varphi^{\ast}\  \omega^{\ast}_{T_{n}}\  E = \omega^{\ast}_{T_{n}}\ F^{\ast}_{X/S}\ E' = \omega^{\ast}_{T_{n}}\ F^{\ast}_{X}\ E \rightarrow\  \omega^{\ast}_{T_{n}}\  E .
$$

\noindent D'o\`u, en notant $\omega^{\ast}_{\mathcal{T}}(E) = \displaystyle \mathop{\lim}_{\leftarrow \atop_{n}}\ \omega^{\ast}_{T_{n}}(E)$ et $\omega^{\ast}_{\mathcal{T}}(\phi) =  \displaystyle \mathop{\lim}_{\leftarrow \atop_{n}}\  \omega^{\ast}_{T_{n}}(\phi)$ , une fl\`eche
$$\omega^{\ast}_{\mathcal{T}}(\phi) : F^{\ast}_{X_{U}/U}\ \varphi^{\ast}\  \omega^{\ast}_{\mathcal{T}}(E)=F^{\ast}_{X_{U}}\  \omega^{\ast}_{\mathcal{T}}(E) = \omega^{\ast}_{\mathcal{T}}F^{\ast}_{X}(E)\ \rightarrow\  \omega^{\ast}_{\mathcal{T}}(E)\ .
$$

\noindent La fl\`eche $\mathbb{R} f_{X_{U}/\mathcal{T}^{\ast}} (\omega^{\ast}_{\mathcal{T}}(\phi))$ est un morphisme
$$ \mathbb{R} f_{X_{U}/\mathcal{T}^{\ast}}(F^{\ast}_{X_{U}/U}\ \varphi^{\ast}\  \omega^{\ast}_{\mathcal{T}}(E)) = \mathbb{R} f_{X_{U}/\mathcal{T}^{\ast}}(F^{\ast}_{X_{U}}\  \omega^{\ast}_{\mathcal{T}}(E)) \ \rightarrow\  \mathbb{R} f_{X_{U}/\mathcal{T}^{\ast}}(\omega^{\ast}_{\mathcal{T}}(E));
$$

\noindent par composition avec le morphisme

$$
\mathbb{R}  f_{X'_{U}/\mathcal{T}^{\ast}} (\varphi^{\ast}  \omega^{\ast}_{\mathcal{T}}(E)) \longrightarrow\ \mathbb{R}  f_{X_{U}/\mathcal{T}^{\ast}} (F^{\ast}_{X_{U}/U}\ \varphi^{\ast}\omega^{\ast}_{\mathcal{T}}(E)) 
$$ 

\noindent provenant  du passage \`a la limite dans le morphisme de changement de base en cohomologie cristalline [B 1, V, (3.5.3)] et avec l'isomorphisme

$$
F^{\ast}_{\mathcal{T}} (\mathbb{R} f_{X_{U}/\mathcal{T}^{\ast}} (\omega^{\ast}_{\mathcal{T}}(E))) \displaystyle \mathop{\longrightarrow}^{\sim}\ \mathbb{R}  f_{X'_{U}/\mathcal{T}^{\ast}} (\varphi^{\ast}\  \omega^{\ast}_{\mathcal{T}}(E))
$$

\noindent provenant de la platitude de $F_{\mathcal{T}}$ et du passage \`a la limite dans l'isomorphisme de changement de base en cohomologie cristalline [B 1, V, (3.5.3)], on obtient un morphisme\\

\noindent (3.3.1.11) $\quad \displaystyle \tilde{\mathop{\phi}} : F^{\ast}_{\mathcal{T}} (\mathbb{R}  f_{X_{U}/\mathcal{T}^{\ast}} (\omega^{\ast}_{\mathcal{T}}(E))) \longrightarrow\ \mathbb{R}  f_{X_{U}/\mathcal{T}^{\ast}} (\omega^{\ast}_{\mathcal{T}}(E)).$\\

 Par passage \`a la cohomologie, $\phi^i = H^i(\displaystyle \tilde{\mathop{\phi}})$ est le morphisme de Frobenius recherch\'e

$$
\begin{array}{c}
\xymatrix{
\phi^i : F^{\ast}_{\mathcal{T}} (R^{i}  f_{X_{U}/\mathcal{T}^{\ast}} (\omega^{\ast}_{\mathcal{T}}(E))) \ar[r] \ar[d]^{\simeq} & R^{i}  f_{X_{U}/\mathcal{T}^{\ast}} (\omega^{\ast}_{\mathcal{T}}(E)) \ar[d]^{\simeq}\\
F^{\ast}_{\mathcal{T}} (\displaystyle \widetilde{\mathop{R^{i} f_{\textrm{cris}^{\ast}}(E)_{U,\mathcal{T}}}}) &  \displaystyle \widetilde{\mathop{R^{i} f_{\textrm{cris}^{\ast}}(E)_{(U,\mathcal{T})}}}.
}
\end{array}
\leqno{(3.3.1.12)}
$$

Consid\'erons \`a pr\'esent le point de vue rigide analytique. Notons $E^{an}$ le $F$-isocristal convergent, \'el\'ement de $F^{a}\mbox{-} \textrm{Isoc} (X/K)$, associ\'e par la construction de Berthelot [B 3, 2.4] au $F$-cristal non d\'eg\'en\'er\'e de type fini $E$.\\

Par d\'efinition des images directes sup\'erieures en cohomologie rigide [B 5, 3.2.3, 3.2.3.2, 3.1.12], [C-T, \S\ 10] et [II, 3.2], on a\\

\noindent (3.3.1.13) $\quad  \mathbb{R}f_{U \textrm{rig}^{\ast}}(X_{U}/ \mathcal{T}, E^{an}_{\mid X_{U}}) = \mathbb{R} h_{\bullet K \ast} (E_{\mathcal{P}_{\bullet K}} \otimes \Omega^{\bullet}_{\mathcal{P}_{\bullet K}/ \mathcal{T}_{K}})$\\

\noindent o\`u $h_{\bullet K} = \{ h_{\alpha} \otimes_{\mathcal{V}} K \}_{\alpha}$ \quad , \quad  $E_{\mathcal{P_{\bullet K}}} = \{ (E_{\mathcal{P_{\alpha}}})^{an}  \}_{\alpha} ,$

$$
\Omega^{\bullet}_{\mathcal{P}_{\bullet K}/ \mathcal{T}_{K}} = \{ ( \Omega^{\bullet}_{\mathcal{P}_{\alpha}/ \mathcal{T}})^{an} \}_{\alpha} .
$$

Vu les descriptions pr\'ec\'edentes, on a donc [B 5, (3.2.3.2)]

$$
\begin{array}{c}
\xymatrix{
R^{i} f_{U \textrm{rig}^{\ast}}(X_{U}/ \mathcal{T}, E^{an}_{\mid X_{U}}) = R^{i} f_{U \textrm{conv}^{\ast}}(X_{U}/ \mathcal{T}, E^{an}_{\mid X_{U}}) &\\
 \qquad \qquad \qquad \qquad = (\displaystyle \widetilde{\mathop{R^{i} f_{\textrm{cris}^{\ast}}(E)_{(U,\mathcal{T})}}})^{an} \ ,
}
\end{array}
\leqno{(3.3.1.14)}
$$

\noindent et ce $\mathcal{O}_{\mathcal{T}_{K}}$-module coh\'erent (d'apr\`es (3.3.1.8)) est muni de la connexion de Gau\ss-Manin via la description (3.3.1.13) (cf [B 1, V, 3.6]) : par la fonctorialit\'e de la construction de $\displaystyle \tilde{\mathop{\phi}}$ cette connexion est compatible au Frobenius $ \phi^i_{K}$ induit par $\phi^i$:
$$
\phi^i_{K} : F^{\ast}_{\mathcal{T}_{K}} R^{i} f_{U \textrm{rig}^{\ast}}(X_{U}/ \mathcal{T}, E^{an}_{\mid X_{U}}) \rightarrow R^{i} f_{U \textrm{rig}^{\ast}}(X_{U}/ \mathcal{T}, E^{an}_{\mid X_{U}}).
\leqno{(3.3.1.15)}
$$
 Nous allons montrer plus bas que $\phi^{i}_{K}$ est un isomorphisme lorsque $k$ est fini: puisque la source et le but de $\phi^{i}_{K}$ sont des $\mathcal{O}_{\mathcal{T}_{K}}$-modules coh\'erents, il suffit de montrer cet isomorphisme fibre \`a fibre aux points ferm\'es de $\mathcal{T}_{K}$ [B-G-R, 9.4.2, cor 7].\\

Soient 
$$
i_{s} : s = \textrm{Spec}\  k(s) \hookrightarrow S
$$

\noindent un point ferm\'e de $S$ et $f_{s} : X_{s} \rightarrow s$ la fibre de $f$ en $s$. On note $\mathcal{V}(s) = W(k(s)) \otimes_{W} \mathcal{V}$ et $K(s)$ le corps des fractions de $\mathcal{V}(s)$. Le morphisme $i_{s}$ d\'efinit un foncteur image inverse [B 3, (2.3.6), (2.3.7)]

$$
\xymatrix{
i^{\ast}_{s} : F^{a}\mbox{-}\textrm{Isoc}(S/K) \ar[r] & F^{a}\mbox{-}\textrm{Isoc}(\textrm{Spec}(k(s))/K(s)) \ar[d]^\simeq\\
&  F^{a}\mbox{-}\textrm{Isoc}^{\dag}(\textrm{Spec}(k(s))/K(s)) ,
}
$$

\noindent et un morphisme de topos [B 1, III, 2.2.3]

$$
i_{s\  \textrm{cris}} :(X_{s}/\mathcal{V}(s))_{\textrm{cris}} \rightarrow (X/\mathcal{V})_{\textrm{cris}}.
$$

\noindent Le point ferm\'e $s$ est contenu dans un ouvert affine et lisse $U = \textrm{Spec}\  A_{0}$ de $S$ : on consid\'erera $i_{s}$ comme un morphisme $i_{s} : \textrm{Spec}k(s)\hookrightarrow U$. On note 
$$
\hat{\tau}(s) : \mathcal{T}_{s} = \textrm{Spf}(\mathcal{V}(s)) \hookrightarrow \mathcal{T}
$$
\noindent le rel\`evement de Teim\"uller de $i_{s}$ [(3.1)]. On obtient alors un diagramme commutatif \`a carr\'es cart\'esiens
$$
\xymatrix{
X_{s} \ar@{^{(}->}[r]^{i_{X}} \ar[d]_{f_{s}} & X_{U} \ar[d]^{f_{U}}\\
\textrm{Spec}\ k(s) \ar@{^{(}->}[r]^{i_{s}} \ar[d]_{i_{\mathcal{T}_{s}}}  & U = \textrm{Spec}\ A_{0} \ar[d]^{i_{\mathcal{T}}}\\
\mathcal{T}_{s} = \textrm{Spf}(\mathcal{V}(s)) \ar@{^{(}->}[r]_{\hat{\tau}(s)} & \mathcal{T} = \textrm{Spf}(\hat{A})\  ;
}
$$
\noindent $\hat{\tau}(s)$ est un PD-morphisme pour les id\'eaux $\pi \mathcal{V}(s)$ et $\pi\ \hat{A}$ munis des puissances divis\'ees canoniques (car $e \leqslant p-1)$. \\

Par analogie avec (3.1) notons $\hat{\tau}_{K}(s)$
$$
\xymatrix{
\hat{\tau}_{K}(s): \ \mathcal{T}_{K(s)} = \textrm{Spm}(K(s)) \ar@{^{(}->}[r]& \mathcal{T}_{K} = \textrm{Spm}(\hat{A}_{K})\  
}
$$
le morphisme induit par $\hat{\tau}(s)$.\\

Comme $\hat{\tau}(s)$ d\'efinit un morphisme de topos $\mathcal{T}_{s \textrm{Zar}} \rightarrow \mathcal{T}_{ \textrm{Zar}}$, le foncteur $\hat{\tau}^{\ast}(s)$ est exact [SGA 4, T1, IV 3.1.2] ; ainsi le passage \`a la limite dans l'isomorphisme de changement de base de [B 1, V, (3.5.1)] fournit un isomorphisme
$$
\hat{\tau}^{\ast}(s) R^{i} f_{X_{U}/\mathcal{T}^{\ast}} (\omega^{\ast}_{\mathcal{T}}(E)) \displaystyle \mathop{\rightarrow}^{\sim} R^{i} f_{X_{s}/\mathcal{T}^{\ast}_{s}} (\omega^{\ast}_{\mathcal{T}_{s}}(E)).
$$
\noindent Comme $i^{\ast}_{s}$ est induit par $\hat{\tau}^{\ast}(s)$ [B 3, (2.3.6) p 72], en passant en rigide analytique on en d\'eduit l'isomorphisme
$$
i^{\ast}_{s}\  R^{i} f_{U \textrm{rig}^{\ast}}(X_{U}/ \mathcal{T}, E^{an}_{\mid X_{U}})=\hat{\tau}^{\ast}_{K}(s)\  R^{i} f_{U \textrm{rig}^{\ast}}(X_{U}/ \mathcal{T}, E^{an}_{\mid X_{U}}) \displaystyle \mathop{\rightarrow}^{\sim}R^{i} f_{s \textrm{rig}^{\ast}}(X_{s}/ \mathcal{T}_{s}, E^{an}_{\mid X_{s}}) ,
$$

\noindent o\`u, par d\'efinition [B 5] et [II, (3.2)], l'on a:
$$
R^{i} f_{s \textrm{rig}^{\ast}}(X_{s}/ \mathcal{T}_{s}, E^{an}_{\mid X_{s}}) 
 =\ H^{i}_{\textrm{rig}}(X_{s}/K(s), E^{an}_{\mid X_{s}});
$$
\noindent d'o\`u l'isomorphisme
$$
\hat{\tau}^{\ast}_{K}(s)\  R^{i} f_{U \textrm{rig}^{\ast}}(X_{U}/ \mathcal{T}, E^{an}_{\mid X_{U}}) \displaystyle \mathop{\rightarrow}^{\sim} H^{i}_{\textrm{rig}}(X_{s}/K(s), E^{an}_{\mid X_{s}}).
\leqno{(3.3.1.16)}
$$
\textit{Supposons dor\'enavant dans cette premi\`ere partie que $k$ est parfait}. Compte tenu de (3.1.1) et de (3.3.1.16) la fibre en $s$ (identifi\'e \`a $\hat{\tau}_{K}(s)$) de $\phi^i_{K}$ est donc un morphisme

$$
\xymatrix{
\phi^i_{K(s)} :&\hat{\tau}^{\ast}_{K}(s) F^{\ast}_{\mathcal{T}_{K}} R^{i} f_{U \textrm{rig}^{\ast}}(X_{U}/ \mathcal{T}, E^{an}_{\mid X_{U}}) \ar[r] \ar@{=}[d]& \hat{\tau}^{\ast}_{K}(s)R^{i} f_{U \textrm{rig}^{\ast}}(X_{U}/ \mathcal{T}, E^{an}_{\mid X_{U}})\ar@{=}[d]\\
& \sigma^{\ast}_{K(s)} \hat{\tau}^{\ast}_{K}(s)R^{i} f_{U \textrm{rig}^{\ast}}(X_{U}/ \mathcal{T}, E^{an}_{\mid X_{U}}) \ar[r] \ar[d]^{\simeq}& \hat{\tau}^{\ast}_{K}(s)R^{i} f_{U \textrm{rig}^{\ast}}(X_{U}/ \mathcal{T}, E^{an}_{\mid X_{U}})\ar[d]^{\simeq}\\
& \sigma^{\ast}_{K(s)}H^{i}_{\textrm{rig}}(X_{s}/K(s), E^{an}_{\mid X_{s}})\ar[r]&H^{i}_{\textrm{rig}}(X_{s}/K(s), E^{an}_{\mid X_{s}})\ .
}
\leqno{(3.3.1.17)}
$$
\noindent Or, on a la g\'en\'eralisation suivante de  [E-LS 1, 2.1]:\\
\vskip 2mm
\noindent \textbf{Proposition (3.3.1.18)}. \textit{
Soient $X$ un $k$-sch\'ema s\'epar\'e de type fini, $F_{X}$ l'it\'er\'e $a$-i\`eme du Frobenius absolu de $X$, $F_{X}=\pi_{X/k}\circ F_{X/k}$ sa factorisation (3.3.1.10)}
$$
\xymatrix{
X\ar[r]^{F_{X/k}}&X'\ar[r]^{\pi_{X/k}}&X
}
$$
 \textit{et $\sigma: K\rightarrow K'=K$ le rel\`evement choisi de la puissance $q$ de $k$ (cf[II, 0]).}
 \textit{Pour $E\in Isoc^{\dag}(X/K)$, on a:}
\begin{enumerate}

\item[(i)]\textit{ Pour tout entier $i\geqslant 0$, $F_{X/k}$ induit une injection $K$-lin\'eaire
$$F^{\ast}_{X/k}: H^{i}_{rig,c}(X'/K',\pi^{\ast}_{X/k}(E))\rightarrow H^{i}_{rig,c}(X/K,F^{\ast}_{X}(E)).$$
 Si de plus $X$ est lisse sur $k$, la m\^eme assertion vaut pour la cohomologie rigide sans supports compacts.}
\item[(ii)] \textit{Supposons de plus $k$ parfait. Alors, sous les hypoth\`eses du (i), le morphisme $F^{\ast}_{X/k}$ pr\'ec\'edent est un isomorphisme.}
\item[(iii)]\textit{Supposons $k$ parfait. Alors, pour tout entier $i\geqslant 0$, $F_{X}$ induit une bijection $\sigma$-lin\'eaire
$$
F_{X}^{\ast}: H^{i}_{rig,c}(X/K,E)\rightarrow H^{i}_{rig,c}(X/K,F^{\ast}_{X}(E))
$$
\noindent c'est-\`a dire un isomorphisme
$$
\sigma^{\ast}(H^{i}_{rig,c}(X/K,E))\displaystyle \mathop{\rightarrow}^{\sim} H^{i}_{rig,c}(X/K,F^{\ast}_{X}(E)).
$$
\noindent Si de plus $X$ est lisse sur $k$, la m\^eme assertion vaut pour la cohomologie rigide sans supports compacts.
}
\end{enumerate}

\noindent \textit{D\'emonstration de (3.3.1.18)}. La preuve suit celle de [E-LS 1, 2.1]. \\
\noindent \textit{Pour (i)}. On notera $H^{i}(X/K,E)$ la cohomologie rigide avec ou sans supports compacts.
En utilisant la suite exacte longue de localisation en cohomologie rigide \`a supports compacts (resp. la suite spectrale de localisation lorsque $X$ est lisse sur $k$) on se ram\`ene au cas o\`u $X$ est un sous-sch\'ema de $\mathbb{P}_{k}^{n}$ qui ne rencontre pas les hyperplans de coordonn\'ees; comme la cohomologie rigide commute aux extensions finies de $K$ on peut supposer que $K$ contient les racines $q$-i\`emes de l'unit\'e. On note $F_{\mathbb{P}}$ l'endomorphisme $\sigma$-lin\'eaire de $\mathbb{P}=\mathbb{P}_{\mathcal{V}}^{n}$ d\'efini par $F_{\mathbb{P}}^{^{\ast}}(t_{i})=t_{i}^{q}$ pour $i\in \llbracket0,n\rrbracket$; on a la factorisation usuelle de $F_{\mathbb{P}}$

$$
\xymatrix{\mathbb{P}_{\mathcal{V}}^{n} \ar@/^1pc/[rrd]^{F_{\mathbb{P}}} \ar@/_/[rdd]_{g}  \ar@{.>}[rd]^{F=F_{\mathbb{P}/\mathcal{V}}} \\
&  \mathbb{P}_{\mathcal{V}}^{n(\sigma)} \ar[d]^{g^{(\sigma)}} \ar[r]^{\pi_{\mathbb{P}/\mathcal{V}}} & \mathbb{P}_{\mathcal{V}}^{n} \ar[d]^g\\
& Spec(\mathcal{V}) \ar[r] ^{\sigma} & Spec(\mathcal{V})
}
$$

\noindent o\`u le carr\'e est cart\'esien. En dehors des hyperplans de coordonn\'ees, le morphisme
$$
F_{K}=(F_{\mathbb{P}/\mathcal{V}})_{K}: \mathbb{P}_{K}^{n}\rightarrow  \mathbb{P}_{K}^{n(\sigma)}
$$ 
est un rev\^etement \'etale galoisien de groupe $\mu_{q}^{n}\simeq (\mathbb{Z}/q\mathbb{Z})^{n}$. Or, si $V$ d\'esigne un voisinage strict suffisamment petit du tube de $X'$, alors $F_{K}$ induit un rev\^etement \'etale galoisien encore not\'e $F_{K}:W\rightarrow V$ de groupe $\mu_{q}^{n}\simeq (\mathbb{Z}/q\mathbb{Z})^{n}$. De plus, on peut supposer que $E'= \pi_{X/k}^{\ast}(E)$ provient d'un module \`a connexion $\mathcal{M}$ sur $V$: tout $F$-automorphisme $\psi$ de $\mathbb{P}_{\mathcal{V}}^{n}$  induit un automorphisme $\psi^{\ast}$ de $F_{K}^{\ast}\mathcal{M}\otimes \Omega^{\bullet}$ et l'endomorphisme $F_{\ast}(\sum \psi^{\ast})$ de $F_{K^{\ast}}F_{K}^{\ast}\mathcal{M}\otimes \Omega^{\bullet}_{V}$ se factorise de mani\`ere unique par le morphisme de complexes
$$
F_{K}^{\ast}: \mathcal{M}\otimes \Omega^{\bullet}_{V} \rightarrow F_{K^{\ast}}F_{K}^{\ast}\mathcal{M}\otimes \Omega^{\bullet}_{V}
$$
pour donner l'application trace
$$
Tr:F_{K^{\ast}}F_{K}^{\ast}\mathcal{M}\otimes \Omega^{\bullet}_{V} \rightarrow \mathcal{M}\otimes \Omega^{\bullet}_{V}
$$
(cf [E-LS 1, 2.1] et [Mi, V, lemma 1.12]). Cette application induit des homomorphismes\\

$tr: H^{i}(W, j_{W}^{\dag}F_{K}^{\ast}\mathcal{M}\otimes \Omega^{\bullet}) \rightarrow H^{i}(V, j_{V}^{\dag}\mathcal{M}\otimes \Omega^{\bullet})$\\

et $ tr: H^{i}_{]X[}(W,F_{K}^{\ast}\mathcal{M}\otimes \Omega^{\bullet}) \rightarrow H^{i}_{]X'[}(V, \mathcal{M}\otimes \Omega^{\bullet}).$\\

\noindent Puisque $F$ prolonge $F_{X/k}$, l'application

$$F^{\ast}_{X/k}: H^{i}(X'/K',\pi^{\ast}_{X/k}(E))\rightarrow H^{i}(X/K,F^{\ast}_{X}(E))$$

\noindent est induite par le morphisme de complexes $F_{K}^{\ast}$, et comme $Tr\circ F_{K}^{\ast}=q^{n}$ sur $\mathcal{M}\otimes \Omega^{\bullet}_{V}$, on en d\'eduit que $tr\circ F_{X/k}^{\ast}=q^{n}$ sur $H^{i}(X'/K',\pi^{\ast}_{X/k}(E))$. D'o\`u l'assertion du (i).\\

\noindent \textit{Pour (ii)}. Par d\'efinition on a $F_{K}\circ Tr= F_{\ast}(\sum \psi^{\ast})$ sur $F_{K^{\ast}}F_{K}^{\ast}\mathcal{M}\otimes \Omega^{\bullet}_{V}$. D'autre part, si l'on note $\psi_{0}$ la r\'eduction mod $\pi$ d'un $F$-automorphisme $\psi$ de $\mathbb{P}_{\mathcal{V}}^{n}$, alors, pour tout $i \in \llbracket 0,n\rrbracket$ et $\lambda_{i}\in k$, $\psi_{0}$ v\'erifie
$$
\psi_{0}^{\ast}(\lambda_{i}^{q}t_{i}^{q})=\psi_{0}^{\ast}F_{X/k}^{\ast}(\lambda_{i}\otimes t_{i})=F_{X/k}^{\ast}(\lambda_{i}\otimes t_{i})=F_{X/k}^{\ast}(1\otimes\lambda_{i}^{q} t_{i})=\lambda_{i}^{q} t_{i}^{q};
$$
puisque $k$ est parfait on en d\'eduit que $\psi_{0}$ est l'identit\'e de $X$. Ainsi on voit que
 $$
 F_{X/k}^{\ast}\circ tr= \sum Id^{\ast}=q^{n}
 $$
sur $H^{i}(X/K,F^{\ast}_{X}(E))$. On a donc montr\'e que $(1/q^{n})tr$ est un inverse pour $F_{X/k}^{\ast}$.\\

\noindent\textit{Pour (iii)}. Le morphisme du (iii) est le compos\'e
$$
\xymatrix{
\sigma^{\ast}(H^{i}(X/K, E))\ar[r]^{ u}&H^{i}(X'/K',\pi^{\ast}_{X/k}(E))\ar[r]^{F_{X/k}^{\ast}}&H^{i}(X/K,F^{\ast}_{X}(E))
}
$$
o\`u le morphisme $u$ est induit par le changement de base $F_{k}$ (puissance $q$ sur $k$) du carr\'e cart\'esien

$$
\xymatrix{
X'\ar[r]^{\pi_{X/k}}\ar[d]^{g'}&X\ar[d]^{g}&\\
Spec\ k\ar[r]_{F_{k}}&Spec\ k&.
}
$$
\noindent Puisque $k$ est parfait, $F_{k}$ et $\pi_{X/k}$ sont des isomorphismes, donc $u$ en est un aussi. Ceci ach\`eve la preuve de (3.3.1.18).  $\square$\\

\noindent\textbf{Remarque (3.3.1.19)}.
 Dans la preuve ci-dessus que $u$ est un isomorphisme le corps $k$ a \'et\'e supposÕe parfait: il est un autre cas o\`u $u$ est un isomorphisme. Supposons que $e\leqslant p-1$, $X$ propre et lisse sur $k$ et $E\in F^{a}\mbox{-}Isoc^{\dag}(X/K)_{plat}$: alors $u$ est un isomorphisme. En effet, d'apr\`es [B 3, (2.4.2)] il existe un $F$-cristal non d\'eg\'en\'er\'e $M$ sur $X$ et un entier $r\geqslant 0$ tels que $E \simeq M^{an}(r)$ et par l'hypoth\`ese faite sur $E$, $M$ est localement libre de type fini: le fait que $u$ soit un isomorphisme r\'esulte alors du th\'eor\`eme de changement de base en cohomologie cristalline [B 1, V, 3.5.7] et de l'isomorphisme

$$
H_{rig}^{i}(X/K, E) \simeq H_{cris}^{i}(X/\mathcal{V}, M)\otimes_{\mathcal{V}}K(r)
$$
( resp. de son analogue sur $X'$) puisque $X$ est propre et lisse sur $k$.  \\

Revenons \`a la preuve de (3.3.1): d'apr\`es (3.3.1.18) la derni\`ere fl\`eche horizontale du diagramme (3.3.1.17) est un isomorphisme; par cons\'equent on a prouv\'e que  $\phi^i_{K}$ est un isomorphisme.\\

Montrons \`a pr\'esent que ces constructions se recollent pour $U$ variable. Si $U_{1}$ et $U_{2}$ sont deux ouverts de $S$ et

$$
j_{1} : U_{3} = U_{1} \cap U_{2} \hookrightarrow U_{1}, j_{2} : U_{3} \hookrightarrow U_{2}
$$

\noindent d\'esignent les immersions ouvertes, le morphisme \'etale $j_{1}$ se rel\`eve de mani\`ere unique en un morphisme \'etale $j^n_{1} : T_{3,n} \rightarrow T_{1,n}$ o\`u $T_{i,n}$ correspond aux $T_{n}$ pr\'ec\'edents [EGA IV, (18.1.2)] : par passage \`a la limite on obtient un morphisme $j^{\infty}_{1} = \displaystyle \mathop{\lim}_{\rightarrow \atop{n}} j^n_{1} : \mathcal{T}_{3} = \displaystyle \mathop{\lim}_{\rightarrow \atop{n}} T_{3,n} \rightarrow \mathcal{T}_{1} = \displaystyle \mathop{\lim}_{\rightarrow \atop{n}} T_{1,n}$ s'ins\'erant dans un diagramme commutatif 

$$
\xymatrix{
(X_{U_{3}}/\mathcal{T}_{3})_{\textrm{cris}} \ar[rr]^{(j_{1X})_{\textrm{cris}}} \ar[d]_{u_{X_{U_{3}}/\mathcal{T}_{3}}} && (X_{U_{1}}/\mathcal{T}_{1})_{\textrm{cris}} \ar[d]^{u_{X_{U_{1}}/\mathcal{T}_{1}}}  &\\
X_{U_{3}} \ar@{^{(}->}[rr]^{j_{1X}} \ar[d]_{f_{U_{3}}}  & & X_{U_{1}} \ar[d]^{f_{U_{1}}}  &\\ 
U_{3}  \ar@{^{(}->}[rr]^{j_{1}} \ar[d]_{i_{\mathcal{T}_{3}}} & & U_{1} \ar[d]^{i_{\mathcal{T}_{1}}}  &\\
\mathcal{T}_{3} \ar@{^{(}->}[rr]_{j^{\infty}_{1}} && \mathcal{T}_{1} & .
}
$$

D'o\`u, par passage \`a la limite dans l'isomorphisme de changement de base en cohomologie cristalline, un isomorphisme

$$
\xymatrix{
j^{\infty \ast}_{1}\ R^{i} f_{X_{U_{1}}/\mathcal{T}^{\ast}_{1}} (\omega^{\ast}_{\mathcal{T}_{1}}(E)) & \overset{\sim}{\longrightarrow} & R^{i} f_{X_{U_{3}}/\mathcal{T}^{\ast}_{3}} (j^{\ast}_{1X \textrm{cris}}\ \omega^{\ast}_{\mathcal{T}_{1}}(E)) \ar@{=}[d]\\
&  &R^{i} f_{X_{U_{3}}/\mathcal{T}^{\ast}_{3}}  (\omega^{\ast}_{\mathcal{T}_{3}}(E)).
}
$$

\noindent En passant en rigide analytique, et par unicit\'e du rel\`evement $j_{1}^{\infty}$, les 

$$
R^{i} f_{U \textrm{rig}^{\ast}} (X_{U}/\mathcal{T}, E^{an}_{\mid X_{U}})
$$
\noindent se recollent pour $U$ variable et d\'efinissent donc un $F$-isocristal convergent sur $(S/K)$ gr\^ace au corollaire (1.2.3) et [B 3, (2.2.11) et 2.3] : on note celui-ci

$$
R^{i} f_{\textrm{conv}^{\ast}} (E^{an}) \in F^{a}\mbox{-}\textrm{Isoc} (S/K),
$$

\noindent et si $R^{i} f_{\textrm{cris}^{\ast}}(E)$ \'etait localement libre ce serait le $F$-isocristal convergent associ\'e, gr\^ace \`a (3.3.1.14. Pour notre $\mathcal{E} \in F^{a}\mbox{-}\textrm{Isoc} (X/K)_{\textrm{plat}}$ initial, donn\'e par $\mathcal{E} \simeq E^{an}(r)$, on pose donc

$$
R^{i} f_{\textrm{conv}^{\ast}}(\mathcal{E}) = R^{i} f_{\textrm{conv}^{\ast}} (E^{an}) (r),
$$

\noindent et sa formation commute aux changements de base $S' \rightarrow S$ puisqu'il en est ainsi pour   $R^{i} f_{X_{U}/\mathcal{T}^ \ast} (\omega^{\ast}_{\mathcal{T}}(E))$. Ceci ach\`eve la preuve de (3.3.1) dans le cas propre et lisse.\\

\noindent\textit{Deuxi\`eme partie de la d\'emonstration: cas o\`u $f$ est relevable en $h$}.\\

Dans le cas o\`u $f$ est relevable en $h : \mathcal{X} \rightarrow \mathcal{S}$ propre et lisse de (3.3.1.1) il n'est plus n\'ecessaire de supposer $E$ localement libre, $E$ sans $p$-torsion suffit.  

\noindent La fl\`eche $R^i f_{U\textrm{rig}^{\ast}}(X_{U}/\mathcal{T};(\omega^{\ast}_{\mathcal{T}_{K}}(\phi)))(r)$ est un morphisme
$$ R^i f_{U\textrm{rig}^{\ast}}(X_{U}/\mathcal{T};(F^{\ast}_{X_{U}}\  \omega^{\ast}_{\mathcal{T}}(E))^{an})(r) \ \rightarrow\  R^i f_{U\textrm{rig}^{\ast}}(X_{U}/\mathcal{T};(  \omega^{\ast}_{\mathcal{T}}(E))^{an})(r);
$$
 par composition avec le morphisme de changement de base [B 5, (3.1.11.3)]
$$
F^{\ast}_{\mathcal{T}_{K}}R^{i} f_{U \textrm{rig}^{\ast}}(X_{U}/\mathcal{T}, (\omega^{\ast}_{\mathcal{T}}(E))^{an})(r) \longrightarrow\ R^i f_{U\textrm{rig}^{\ast}}(X_{U}/\mathcal{T};(  \omega^{\ast}_{\mathcal{T}}(E))^{an})(r)
$$ 

\noindent on obtient le morphisme de Frobenius (qu'on prouve \^etre un isomorphisme comme en (3.3.1.17)) \\

\noindent (3.3.1.20)$\qquad \phi^i_{K} : F^{\ast}_{\mathcal{T}_{K}} R^{i} f_{U \textrm{rig}^{\ast}} (X_{U}/\mathcal{T}, \mathcal{E}_{\mid X_{U}}) \displaystyle \mathop{\rightarrow}^{\sim} R^{i} f_{U \textrm{rig}^{\ast}} (X_{U}/\mathcal{T}, \mathcal{E}_{\mid X_{U}}) $\\

\noindent valable pour $\mathcal{E} \in F^{a}\mbox{-}\textrm{Isoc}(X/K)$ et $f$  relevable en $h$ propre et lisse.\\

On conclut alors \`a la mani\`ere de [II, (3.4.8.2)] via [II, (3.4.4)]. $\square$\\

\noindent\textbf{Remarque (3.3.1.21)}. Nous avons prouv\'e ci-dessus que $R^{i} f_{U \textrm{rig}^{\ast}} (X_{U}/\mathcal{T}, E^{an}_{\mid X_{U}})$ est un $\mathcal{O}_{\mathcal{T}_{K}}$-module coh\'erent, qu'il est muni d'une connexion int\'egrable (de Gau\ss Manin) et  d'un morphisme de Frobenius sous la seule hypoth\`ese que le corps $k$ est de caract\'eristique $p>0$ et $e\leqslant p-1$: c'est seulement pour prouver que le Frobenius est un isomorphisme que nous avons \'et\'e amen\'es \`a supposer que le corps $k$ est parfait.\\

Le th\'eor\`eme suivant pr\'ecise (3.3.1) en l'\'etendant :\\

\noindent\textbf{Th\'eor\`eme (3.3.2)}. \textit{Supposons $k$ parfait. Soient $S$ un $k$-sch\'ema lisse et $f : X \rightarrow S$ un $k$-morphisme projectif et lisse satisfaisant aux hypoth\`eses de [II, (3.4.8.2)] ou [II, (3.4.8.6)] ou [II, (3.4.9)]. Alors, pour tout entier $i \geqslant 0$, $f$ induit un foncteur}
$$
R^{i} f_{conv\ast}\  : F^{a}\mbox{-}\textrm{Isoc}(X/K) \rightarrow F^{a}\mbox{-}\textrm{Isoc}(S/K)
$$
\textit{qui commute \`a tout changement de base $S' \rightarrow S$ entre $k$-sch\'emas lisses.} 

\vskip 3mm
\noindent \textit{D\'emonstration}. Compte tenu de [II, (3.4.8.2), (3.4.8.6) et (3.4.9)] il s'agit de v\'erifier que le Frobenius (qui est d\'efini par fonctorialit\'e [L.S, 8]) est un isomorphisme. D'apr\`es [B-G-R, (9.4.2/7)] il suffit de v\'erifier l'isomorphisme sur les points ferm\'es de $S$ et le r\'esultat provient alors de [(3.3.1.18)]. $\square$

\subsection*{3.4. Fibres des $F$-isocristaux convergents}

Sous les hypoth\`eses du th\'eor\`eme (3.3.1) avec $f$ propre et lisse, soient 

$$
i_{s} : s = \textrm{Spec}\  k(s) \hookrightarrow S
$$

\noindent un point ferm\'e de $S$ et $f_{s} : X_{s} \rightarrow s$ la fibre de $f$ en $s$. 
Pour $\mathcal{E} \in F^{a}\mbox{-}\textrm{Isoc}(X/K)_{\textrm{plat}}$ on pose\\

\noindent (3.4.1) $\qquad  R^{i} f_{\textrm{conv}^{\ast}}(\mathcal{E})_{s} := i^{\ast}_{s}\ R^{i} f_{\textrm{conv}^{\ast}} (\mathcal{E})).$\\

\noindent Comme $i^{\ast}_{s}$ est induit par $\hat{\tau}^{\ast}(s)$ [B 3, (2.3.6) p 72], on a d\'emontr\'e en (3.3.1.15) l'isomorphisme\\

\noindent (3.4.2) \qquad \qquad $i^{\ast}_{s}\  R^{i} f_{\textrm{conv}^{\ast}} (\mathcal{E})) \displaystyle \mathop{\rightarrow}^{\sim} R^{i} f_{s\ \textrm{conv}^{\ast}} (\mathcal{E}_{\mid X_{s}})= H^i_{\textrm{rig}}(X_{s}/K(s);\mathcal{E}_{\mid X_{s}}).$\\

Compte tenu de [II, (3.4.7) et (3.4.8.2)] on a prouv\'e:

\vskip 3mm
\noindent \textbf{Proposition (3.4.3)}. \textit{Sous les hypoth\`eses du th\'eor\`eme (3.3.1), avec $f$ propre et lisse (resp. sous les hypoth\`eses du th\'eor\`eme (3.3.2)), et pour $\mathcal{E} \in F^{a}\mbox{-}Isoc(X/K)_{\textrm{plat}}$ (resp. $\mathcal{E} \in F^a\mbox{-}\textrm{Isoc}(X/K)$) et tout point ferm\'e $s$ de $S$, la fibre en $s$ de $R^{i} f_{\textrm{conv}^{\ast}}(\mathcal{E})$ est un $K(s)$-espace vectoriel de dimension finie donn\'e par}

$$
\xymatrix{
R^{i}  f_{\textrm{conv}^{\ast}} (\mathcal{E})_{s} \simeq H^i_{\textrm{rig}}(X_{s}/K(s);\mathcal{E}_{\mid X_{s}}).
}
$$

Comme corollaire du th\'eor\`eme (3.3.1), on retrouve un cas particulier d'un th\'eor\`eme de Ogus sur la finitude de la cohomologie convergente [O 3, 0.7.9] :

\vskip 3mm
\noindent \textbf{Corollaire (3.4.4)}. \textit{Supposons $k$ parfait et $e \leqslant p-1$. Soient $X$ un $k$-sch\'ema propre et lisse, et $\mathcal{E} \in F^{a}\mbox{-}Isoc(X/K)_{\textrm{plat}}$. Alors, pour tout entier $i \geqslant 0$, les groupes de cohomologie convergente $H^i_{\textrm{conv}}(X/K; \mathcal{E})$ sont des $K$-espaces vectoriels de dimension finie.}

\vskip 3mm
\noindent \textit{D\'emonstration}. Il suffit d'appliquer (3.3.1) avec $S = \textrm{Spec}\ k$ et [B 5, (3.2.3)]. $\square$\\

\subsection*{3.5. Cas fini \'etale}

Dans le cas fini \'etale on n'a pas lieu de supposer $k$ parfait ni $e \leqslant p-1$ et de se restreindre, comme dans le th\'eor\`eme (3.3.1), au cas des $F$-cristaux localement libres:

\vskip 3mm
\noindent \textbf{Th\'eor\`eme (3.5.1)}. \textit{Soient $S$ un $k$-sch\'ema lisse et $f : X \rightarrow S$ un $k$-morphisme fini \'etale. Alors}\\

\noindent (3.5.1.1) \textit{Pour tout entier $i \geqslant 0$, on a des foncteurs }\\

\begin{itemize}
\item[(i)] $R^{i}  f_{\textrm{conv}^{\ast}} : \textrm{Isoc}(X/K) \rightarrow \textrm{Isoc}(S/K),$

\vskip 1mm
\item[(ii)] $R^{i}  f_{\textrm{conv}^{\ast}} : F^{a}\mbox{-}\textrm{Isoc}(X/K) \rightarrow F^{a}\mbox{-}\textrm{Isoc}(S/K),$

\vskip 1mm
\item[(iii)] \textit{Pour $\mathcal{E} \in \textrm{Isoc}(X/K)$ et $i \geqslant 1$ on a }
$$
R^{i} f_{\textrm{conv}^{\ast}} \mathcal{(E)} = 0.
$$
\end{itemize}
\noindent (3.5.1.2) \textit{Supposons de plus $f$ galoisien de groupe $G$. Pour $\mathcal{E} \in \textrm{Isoc}(X/K)$ on a des isomorphismes canoniques}\\
\begin{itemize}
\item[(i)] $\mathcal{E} \displaystyle \mathop{\longrightarrow}^{\sim}\ (f_{\textrm{conv}^{\ast}}\ f^{\ast}(\mathcal{E}))^G,$
\vskip 1mm
\item[(ii)] $H^i_{\textrm{conv}}(S/K, \mathcal{E}) \displaystyle \mathop{\longrightarrow}^{\sim}\ H^i_{\textrm{conv}}(X/K, f^{\ast}(\mathcal{E}))^G,$
\vskip 1mm
\item[(iii)] \textit{Si  $\mathcal{E} \in F^{a}\mbox{-}\textrm{Isoc}(S/K)$, ces isomorphismes sont compatibles aux Frobenius.}
\end{itemize}

\vskip 3mm
\noindent \textit{D\'emonstration.\\
Pour (3.5.1.1)}. Le (i) est l\`a pour m\'emoire, car prouv\'e en [II, 3.4.8]. On a vu dans la d\'emonstration de [loc. cit.] que la d\'efinition de $R^{i}  f_{\textrm{conv}^{\ast}}\mathcal{(E)}$ est locale sur $S$ : on peut donc supposer $S = \textrm{Spec}\ A_{0}$ affine et lisse sur $k$.\\

Posons $\mathcal{S} = \textrm{Spf}\ \hat{A}$ o\`u $A$ est une $\mathcal{V}$-alg\`ebre lisse relevant $A_{0}$, et relevons $f : X \rightarrow S$ en un morphisme fini \'etale de $\mathcal{V}$-sch\'emas formels $h : \mathcal{X} \rightarrow \mathcal{S}$ [EGA IV, (18.3.2) ou (18.3.4)], et soit $F_{\mathcal{S}} : \mathcal{S} \rightarrow \mathcal{S}$ un rel\`evement du Frobenius de $S$. Puisque $f$ est \'etale, dans la d\'ecomposition classique du Frobenius $F_{X}$ de $X$

$$
\xymatrix{X \ar@/^1pc/[rrd]^{F_{X}} \ar@/_/[rdd]_{f}  \ar@{.>}[rd]^{F_{X/S}} \\
&  X^{(q)} \ar[d]^{f^{(q)}} \ar[r]^{\pi_{X/S}} & X \ar[d]^f\\
& S \ar[r] _{F_{S}} & S
}
$$

\noindent le morphisme $F_{X/S}$ est un isomorphisme et se rel\`eve de mani\`ere unique en un isomorphisme $F_{\mathcal{X}/\mathcal{S}}$ s'ins\'erant dans le diagramme commutatif \`a carr\'e cart\'esien

$$
\xymatrix{\mathcal{X}  \ar@/_/[rdd]_{h}  \ar[rd]^{F_{\mathcal{X}/\mathcal{S}}}  &\\
&  \mathcal{X}^{(q)} \ar[d]^{h^{(q)}} \ar[r]^{\pi_{\mathcal{X}/\mathcal{S}}} & \mathcal{X} \ar[d]^h &\\
& \mathcal{S} \ar[r] _{F_{\mathcal{S}}} & \mathcal{S} & .
}
$$\\
On pose $F_{\mathcal{X}} = \pi_{\mathcal{X}/\mathcal{S}} \circ F_{\mathcal{X}/\mathcal{S}}.$\\

Pour $\mathcal{E} \in \textrm{Isoc}(X/K)$, soit $\mathcal{E}_{\mathcal{X}}$ une r\'ealisation de $\mathcal{E}$ sur $\mathcal{X}_{K}$ ; par d\'efinition on a

$$
\xymatrix{
R^{i} f_{\textrm{conv}^{\ast}}(X/ \mathcal{S};\mathcal{E}) = H^{i}(\mathbb{R} h_{K^{\ast}}
(\mathcal{E}_{\mathcal{X}} \otimes \Omega^{\bullet}_{\mathcal{X}_{K}/\mathcal{S}_{K}}))\\
= R^{i} h_{K^{\ast}} (\mathcal{E}_{\mathcal{X}}),
}
$$

\noindent car $\mathcal{X}$ est \'etale sur $\mathcal{S}$. D'o\`u le (iii) par le th\'eor\`eme B de Kiehl car $h$ est affine. \\

Soit $\mathcal{E} \in F^{a}\mbox{-}\textrm{Isoc}(X/K)$ et $\phi : F^{\ast}_{\mathcal{X}}
(\mathcal{E}_{\mathcal{X}}) \displaystyle \mathop{\longrightarrow}^{\sim} \mathcal{E}_{\mathcal{X}}$ le Frobenius.\\
D'apr\`es [III, (1.2.3)] il suffit de construire un isomorphisme $\phi^i$ (de Frobenius) sur $R^{i} f_{\textrm{conv}^{\ast}}(\mathcal{E})$, compatible aux connexions. Comme $F_{\mathcal{S}}$ est plat, le morphisme de changement de base 
$$
F^{\ast}_{\mathcal{S}_{K}} R^{i} f_{\textrm{conv}^{\ast}}(X/\mathcal{S}, \mathcal{E}) \rightarrow R^{i} f_{\textrm{conv}^{\ast}}^{(q)}(X^{(q)}/\mathcal{S},\mathcal{E}) \simeq R^{i} h_{K^{\ast}}^{(q)}(\pi^{\ast}_{\mathcal{X}_{K}/\mathcal{S}_{K}}(\mathcal{E}_{\mathcal{X}}))
$$
\noindent est un isomorphisme [II, (3.4.4)] ; par composition avec les isomorphismes

$$
R^{i} h_{K^{\ast}}^{(q)}(\pi^{\ast}_{\mathcal{X}_{K}/\mathcal{S}_{K}}(\mathcal{E}_{\mathcal{X}})) \simeq R^{i} h_{K^{\ast}}(F^{\ast}_{\mathcal{X}_{K}}(\mathcal{E}_{\mathcal{X}}))
$$

 \noindent (puisque $F_{\mathcal{X}/\mathcal{S}}$ est un isomorphisme)

\noindent et

$$
R^{i}  h_{K^{\ast}}(\phi) : R^{i}  h_{K^{\ast}}(F^{\ast}_{\mathcal{X}_{K}}(\mathcal{E}_{\mathcal{X}})) \simeq R^{i}  h_{K^{\ast}}(\mathcal{E}_{\mathcal{X}}) ,
$$

\noindent on obtient le Frobenius $\phi^i$ cherch\'e

$$
\phi^i : F^{\ast}_{\mathcal{S}_{K}}  R^{i} f_{\textrm{conv}^{\ast}}(X/\mathcal{S}, \mathcal{E}) \simeq R^{i} f_{\textrm{conv}^{\ast}}(X/\mathcal{S}, \mathcal{E}).
$$

En reprenant la preuve de [II, (3.4.4)] on v\'erifie que $\phi^i$ est compatible aux connexions d'o\`u le (ii).\\

\textit{Pour (3.5.1.2)}. Soit $\mathcal{E}_{\mathcal{S}}$ une r\'ealisation de $\mathcal{E}$ sur $\mathcal{S}_{K}$. Par d\'efinition on a

$$
f_{\textrm{conv}^{\ast}}(X/\mathcal{S}, f^{\ast}(\mathcal{E})) = h_{K^{\ast}}
(h^{\ast}_{K}(\mathcal{E}_{\mathcal{S}})).
$$

\noindent Comme $h_{K}$ est fini \'etale galoisien de groupe G [II, (2.3.1)], la fl\`eche canonique
$$
\mathcal{E}_{\mathcal{S}} \rightarrow (h_{K^{\ast}} (h^{\ast}_{K}(\mathcal{E}_{\mathcal{S}})))^G
$$

\noindent est un isomorphisme, d'o\`u (i).\\

L'isomorphisme du (ii) est alors une cons\'equence classique du (i) [Et 2, III, 3.1.1].\\

La fonctorialit\'e des constructions pr\'ec\'edentes prouve le (iii).  $\square$


\cleardoublepage

\vskip 10mm
\chapter*{IV.  Images directes de $F$-isocristaux surconvergents}
\markboth{\sc j.-y. etesse}{\sc IV.  Images directes de $F$-isocristaux surconvergents}

\section*{1. Frobenius}

\textbf{1.1.} On fixe dans ce paragraphe 1 un entier $a \in \mathbb{N}^{\ast}$ ; on pose $q = p^a$. Pour tout $k$-sch\'ema $S$ on notera $F_{S}$ le Frobenius de $S$ induit par la puissance $q$ sur le faisceau $\mathcal{O}_{S}$.\\

On fixe un rel\`evement $\sigma : \mathcal{V} \rightarrow \mathcal{V}$ de la puissance $q$ sur $k$ \`a la mani\`ere de [Et 5, I, 1.1].\\

Si $S$ est lisse sur $k$ et $e \leqslant p-1$, on notera  $F^{a}\mbox{-}\textrm{Isoc}^{\dag}(S/K)_{\textrm{plat}}$ la sous-cat\'egorie pleine de $F^{a}\mbox{-}\textrm{Isoc}^{\dag}(S/K)$ form\'e des objets dont l'image par le foncteur d'oubli

$$
F^{a}\mbox{-}\textrm{Isoc}^{\dag}(S/K) \longrightarrow F^{a}\mbox{-}\textrm{Isoc}(S/K) 
$$

\noindent est dans  $F^{a}\mbox{-}\textrm{Isoc}(S/K)_{\textrm{plat}}$, cf [III, 3.1].\\

\textbf{1.2.} Soit $S$ un $k$-sch\'ema affine et lisse. En utilisant les notations du [I, th\'eo (3.4) (3)] il existe une $\mathcal{V}$-alg\`ebre lisse $A$ telle que $Spec\ A$ rel\`eve $S$ et un $\mathcal{V}$-morphisme fini $\psi$ relevant le Frobenius $F_{S}$, s'ins\'erant dans un diagramme commutatif \`a carr\'es cart\'esiens

$$
\begin{array}{c}
\xymatrix{
Spec\ B_{T} \ar[r] \ar[d]_{\psi_{T}}& Spec\ B\ \ar@{^{(}->}[r]^(.62){j_{Z'}} \ar[d]^{\psi} & Z' \ar[d]^{\overline{\psi}} \\
Spec\ A_{T} \ar[r]  & Spec\ A\ \ar@{^{(}->}[r]^(.62){j_{Z}}& Z
}
\end{array}
\leqno{(1.2.1)}
$$

\noindent o\`u les $j$ sont des immersions ouvertes, $\overline{\psi}$ est fini,$\psi_{T}$ est fini et plat, et $Z$ est un $\mathcal{V}$-sch\'ema propre, normal.\\

Soit $\hat{A}$ le s\'epar\'e compl\'et\'e $\mathfrak{m}$-adique de $A$ : c'est aussi le s\'epar\'e compl\'et\'e de $A_{T}$, et on a un isomorphisme [I, th\'eo (3.4) (2) (i)]

$$
B_{T} \otimes_{A_{T}} \hat{A} \simeq \hat{B} \simeq \hat{A} 
$$

\noindent tel que dans le diagramme commutatif \`a carr\'es cart\'esiens\\

$$
\begin{array}{c}
\xymatrix{
Spec\ \hat{A} \ar[r]^{\sim} \ar[rd]_{\varphi} \ar@/^2pc/[rr]^{\rho_{B}} & Spec\ \hat{B}_{T} \ar[r] \ar[d]^{\hat{\psi}_{T}} & Spec\ B\ \ar@{^{(}->}[r]^(.62){j_{Z'}} \ar[d]^{\psi} & Z' \ar[d]^{\overline{\psi}} \\
& Spec\ \hat{A}_{T} = Spec \hat{A} \ar[r]_(.62){\rho_{A}}  & Spec\ A\ \ar@{^{(}->}[r]_(.62){j_{Z}}& Z
}
\end{array}
\leqno{(1.2.2)}
$$

\noindent $\varphi$ est un rel\`evement de $F_{S}$.\\

Le morphisme diagonal $Spec\ \hat{A} \longrightarrow Spec\ \hat{A}\  \times_{\mathcal{V}}\  Spec\ \hat{A}$ est une immersion ferm\'ee, donc $Spec\ \hat{A}$ est isomorphe \`a son image sch\'ematique $\hat{\Delta}$ par ce morphisme. Consid\'erons l'image sch\'ematique de $\hat{\Delta}$ par le morphisme compos\'e

$$
Spec\ \hat{A}\ \times_{\mathcal{V}}\ Spec\ \hat{A} \displaystyle \mathop{\twoheadrightarrow}_{\rho=\rho_{_{A}}\times \rho_{B}}\ Spec\ A \times_{\mathcal{V}}\  Spec\ B \displaystyle \mathop{\hooklongrightarrow}_{j_{\mathcal{Z}} \times j_{\mathcal{Z}'} = j} \mathcal{Z} \times_{\mathcal{V}} \mathcal{Z}'\ ;
$$

\noindent notons $\Delta$ (resp. $\mathcal{Z}''$) l'image sch\'ematique de $\hat{\Delta}$  (resp. de $\Delta$) par $\rho_{A} \times \rho_{B}$ (resp. par $j_{Z} \times j_{Z'}$). L'immersion ouverte $j$ induit une immersion ouverte $j_{Z}'' : \Delta \hookrightarrow Z''$ [EGA I, (5.4.4)]. \\

Montrons que $\rho(\hat{\Delta}) = \Delta$. Quitte \`a d\'ecomposer la $\mathcal{V}$-alg\`ebre lisse (donc normale) $A$ en somme de ses composantes connexes, on peut supposer $A$ int\`egre, donc int\'egralement clos : ainsi $\hat{A}$ est int\'egralement clos [I, prop (1.6) (4) (iv)]. Soit $I$ l'id\'eal de $\hat{A} \otimes_{\mathcal{V}} \hat{A}$ d\'efinissant $\hat{\Delta} = Spec(\hat{A} \otimes_{\mathcal{V}} \hat{A}/I)$ : comme $\hat{A}$ est int\`egre, $I$ est un id\'eal premier. L'image de $A$ par $\rho$ est donc l'ensemble des id\'eaux premiers de $A \otimes_{\mathcal{V}} B$ contenant $\rho (I)$ : c'est donc $Spec(A \otimes_{\mathcal{V}} B/ \tilde{\rho}^{-1}(I))$ o\`u $\tilde{\rho} : A \otimes_{\mathcal{V}} B \rightarrow \hat{A} \otimes_{\mathcal{V}} \hat{A}$ induit $\rho$ ; comme c'est d\'ej\`a un sous-sch\'ema ferm\'e de $Spec\ A \times_{\mathcal{V}} Spec\ B$, il est \'egal \`a $\Delta$.\\

En remarquant que $\rho_{n} = \rho\ \textrm{mod}\  \mathfrak{m}^{n+1}$ est l'identit\'e, $\rho_{n}$ induit un isomorphisme

$$
(Spec\ \hat{A}\  \textrm{mod}\ \mathfrak{m}^{n+1}) \displaystyle \mathop{\longrightarrow}^{\sim} (\hat{\Delta} \ \textrm{mod}\ \mathfrak{m}^{n+1})  \displaystyle \mathop{\longrightarrow}^{\sim}_{\rho_{n}}
(\Delta\  \textrm{mod}\   \mathfrak{m}^{n+1}).
$$

Notons alors $\mathcal{S}, \mathcal{Z}, \mathcal{Z'}, \mathcal{Z''}$ les sch\'emas formels associ\'es respectivement \`a $Spec\ \hat{A}, \ Z, \ Z', \ Z''$. Ce qui pr\'ec\`ede fournit un diagramme commutatif \`a carr\'e cart\'esien\\

$$
\begin{array}{c}
\xymatrix{
& \mathcal{Z}\\
& \mathcal{Z}''  \ar@{^{(}->}[r] \ar[d] \ar[d]^{v_{\mathcal{Z}'}}  \ar[u]_{v_{\mathcal{Z}}} & \mathcal{Z} \times_{\mathcal{V}} \mathcal{Z}' \ar[dl]^{\textrm{proj}}  \ar[ul]_{\textrm{proj}}\\
\mathcal{S}   \ar@{^{(}->}[ur]_{j_{\mathcal{Z}''}}    \ar@{^{(}->}[uur]^{j_{\mathcal{Z}}}  \ar@{^{(}->}[r]_{j_{\mathcal{Z}'}}    \ar[d]_{F_{\mathcal{S}}:=\hat{\varphi}} &  \mathcal{Z}' \ar[d]^{\hat{\overline{\psi}}=:F_{\mathcal{Z}'}} &\\
\mathcal{S}  \ar@{^{(}->}[r]^{j_{\mathcal{Z}}}   & \mathcal{Z} & 
}
\end{array}
\leqno{(1.2.3)}  
$$

\noindent o\`u $\mathcal{Z}$ est propre sur $\mathcal{V}$ et normal, [I, th\'eo (3.4) (3) (ii)], $\hat{\overline{\psi}}$ est fini ; $F_{\mathcal{S}}$ est un rel\`evement fini et plat du Frobenius $F_{S}$ ; $j_{\mathcal{Z}}, j_{\mathcal{Z'}}, j_{\mathcal{Z''}}$ sont des immersions ouvertes induites respectivement par $j_{Z}, j_{Z'}, j_{Z''}$ ; $i$ est l'immersion ferm\'ee induite par l'immersion ferm\'ee $\mathcal{Z}'' \hookrightarrow \mathcal{Z}\times_{\mathcal{V}}\mathcal{Z}'  ; \  v_{\mathcal{Z}}, v_{\mathcal{Z}'}$ sont des morphismes propres par composition de morphismes propres.\\

Avec les notations de [II, (2.3.1) (2)] soit $V_{\lambda} = Spm\ A_{\lambda}$: il existe $ \lambda_{0} > 1 $ tel que, pour $1< \lambda \leqslant \lambda_{0}$, $V_{\lambda}$ est lisse sur $K$ ; notons $W'_{\lambda} = F^{-1}_{\mathcal{Z}'_{K}}(V_{\lambda})$ et $W''_{\lambda} = v^{-1}_{\mathcal{Z}'_{K}}(W'_{\lambda})$. \\

\noindent Puisque $(V_{\lambda})_{\lambda}$ d\'ecrit un syst\`eme fondamental de voisinages stricts de $\mathcal{S}_{K}$ dans $\mathcal{Z}_{K}$, alors $(W'_{\lambda})_{\lambda}$ d\'ecrit un syst\`eme fondamental de voisinages stricts de $\mathcal{S}_{K}$ dans $\mathcal{Z}'_{K}$ [II, prop (2.1.2)]. Comme $v_{\mathcal{Z}'}$ est \'etale au voisinage de $\mathcal{S}$, il existe  $\lambda_{1} > 1$ tel que pour tout $\lambda$, $1 < \lambda \leqslant \lambda_{1} \leqslant\lambda_{0}$, on ait un isomorphisme $W''_{\lambda} \displaystyle \mathop{\rightarrow}^{\sim} W'_{\lambda}$ induit par $v_{\mathcal{Z}'}$ [B 3, (1.3.5)]. De m\^eme $v_{\mathcal{Z}}$ qui est \'etale au voisinage de $\mathcal{S}$, induit un isomorphisme entre un syst\`eme fondamental de voisinages stricts de $\mathcal{S}_{K}$ dans $\mathcal{Z}''_{K}$ et un syst\`eme fondamental de voisinages stricts de $\mathcal{S}_{K}$ dans $\mathcal{Z}_{K}$ : par composition il existe $\mu$, $1 < \mu \leqslant \lambda \leqslant \lambda_{1}$, et un morphisme fini $F_{\lambda \mu}$ rendant cart\'esien le carr\'e

$$
\begin{array}{c}
\xymatrix{
\mathcal{S}_{K} \ar@{^{(}->}[r] \ar[d]_{F_{\mathcal{S}_{K}}} & V_{\mu} \ar[d]^{F_{\lambda \mu}}\\
\mathcal{S}_{K} \ar@{^{(}->}[r] & V_{\lambda}
}
\end{array}
\leqno{(1.2.4)}
$$

\noindent o\`u les fl\`eches horizontales sont des immersions ouvertes et $V_{\lambda}$ est lisse sur $K$.

\vskip 6mm
\section*{2. Cas relevable}

\noindent \textbf{Th\'eor\`eme (2.1)}. \textit{Soient $S$ un $k$-sch\'ema lisse et s\'epar\'e et $f : X \rightarrow S$ un $k$-morphisme propre et lisse. On suppose qu'il existe un carr\'e cart\'esien de $\mathcal{V}$- sch\'emas formels}

$$
\begin{array}{c}
\xymatrix{
\mathcal{X} \ar@{^{(}->}[r]^{j_{\overline{\mathcal{X}}}} \ar[d]_{h} & \overline{\mathcal{X}} \ar[d]^{\overline{h}} &\\
\mathcal{S} \ar@{^{(}->}[r]_{j_{\overline{\mathcal{S}}}} & \overline{\mathcal{S}} & ,
}
\end{array}
\leqno{(2.1.1)}
$$

\noindent \textit{de r\'eduction mod $\mathfrak{m}$ \'egale \`a}

$$
\begin{array}{c}
\xymatrix{
X \ar@{^{(}->}[r]^{j_{\overline{X}}} \ar[d]_{f} & \overline{X} \ar[d]^{\overline{f}} &\\
S \ar@{^{(}->}[r]_{j_{\overline{S}}} & \overline{S} & ,
}
\end{array}
\leqno{(2.1.2)}
$$

\noindent \textit{o\`u $\overline{\mathcal{S}}$ est propre sur $\mathcal{V}$, $\overline{h}$ est propre, $h$ est propre et lisse et les $j$ sont des immersions ouvertes.}\\

\textit{Soit $E \in F^{a}\mbox{-}\textrm{Isoc}^{\dag}(X/K)$ et $\hat{E} = \mathcal{E}$ son image dans $F^{a}\mbox{-}\textrm{Isoc}(X/K)$. Alors, pour tout entier $i \geqslant 0$}

\begin{itemize}
\item[(1)] $E_{i} := R^{i} f_{\textrm{rig}^{\ast}}(X/ \overline{\mathcal{S}}, E) \in F^{a}\mbox{-}\textrm{Isoc}^{\dag}(S/K)$
\item[(2)] \textit{Soient $\hat{E}_{i} = j^{\ast}_{\overline{S}}(E_{i}), \mathcal{E}_{i} = R^{i} f_{\textrm{conv}^{\ast}} (X/ \mathcal{S}, \mathcal{E})$ et $\phi_{E_{i}} : F_{S}^{\ast}\ E_{i} \rightarrow E_{i}$, $\phi_{\mathcal{E}_{i}} : F_{S}^{\ast}\ \mathcal{E}_{i} \rightarrow \mathcal{E}_{i}$ les isomorphismes de Frobenius. Le diagramme commutatif d'isomorphismes ci-dessous d\'efinit $\phi_{\hat{E}_{i}}$ et permet les identifications canoniques} 
\end{itemize}

$$
j^{\ast}_{\overline{S}} (\phi_{E_{i}}) = \phi_{\hat{E}_{i}} = \phi_{\mathcal{E}_{i}}
$$

$$ 
\xymatrix{
j^{\ast}_{\overline{S}}\ F^{\ast}_{S}\ R^{i} f_{\textrm{rig} \ast}(X/ \overline{\mathcal{S}}, E) \ar[r]^{\sim}_(.53){j^{\ast}_{\overline{S}} (\phi_{E_{i}})} \ar[d]^{\simeq} & j^{\ast}_{\overline{S}}\ R^{i} f_{\textrm{rig}\ast}(X/ \overline{\mathcal{S}}, E) \ar@{=}[d]\\
F^{\ast}_{S}\ j^{\ast}_{\overline{S}}\  R^{i} f_{\textrm{rig} \ast}(X/ \overline{\mathcal{S}}, E) \ar[r]^{\sim}_{\phi_{\hat{E}_{i}}}   \ar[d]^{\simeq} & j^{\ast}_{\overline{S}}\  R^{i} f_{\textrm{rig}\ast}(X/ \overline{\mathcal{S}}, E) \ar[d]^{\simeq}\\
F^{\ast}_{S}\  R^{i} f_{\textrm{conv}^{\ast}}(X/ \mathcal{S}, \mathcal{E}) \ar[r]^{\sim}_{\phi_{\mathcal{E}_{i}}}  & R^{i}  f_{\textrm{conv}^{\ast}}(X/ \mathcal{S}, \mathcal{E}) .
}
$$

\vskip 3mm
\noindent \textit{D\'emonstration}. L'image inverse $F^{\ast}_{\sigma}\ E_{i}$ par Frobenius s'obtient [B 3, (2.3.7)] en appliquant le foncteur de changement de base
$$
\sigma^{\ast} : \textrm{Isoc}^{\dag}(S/K) \longrightarrow \textrm{Isoc}^{\dag}(S^{(q)}/K)
$$
\noindent puis le foncteur image inverse par le Frobenius $F_{S/k} : S \rightarrow S^{(q)}$.\\

\noindent Comme $\sigma$ est fix\'e on notera $F^{\ast}_{\sigma}\ E = F^{\ast}_{S} \ E$. Il nous reste donc \`a d\'efinir l'isomorphisme  de Frobenius $\phi_{E_{i}}$ de $E_{i}$.\\

Quitte \`a d\'ecomposer $S$ en somme de ses composantes connexes il suffit de d\'efinir $\phi_{E_{i}}$ sur chacune de ces composantes connexes. Soit $S_{\alpha}$ un ouvert affine d'une composante connexe $S_{0}$ de $S$ : comme le foncteur
$$
F^{a}\mbox{-} \textrm{Isoc}^{\dag}(S_{0}/K) \longrightarrow F^{a}\mbox{-} \textrm{Isoc}^{\dag}(S_{\alpha}/K)
$$
\noindent est pleinement fid\`ele [Et 5, th\'eo 4], il suffit de d\'efinir $\phi_{E_{i}}$ sur $S_{\alpha}$.\\

Soit $j_{s_{\alpha}} : S_{\alpha} = Spec\ A_{0} \hookrightarrow S$ l'immersion ouverte et $A$ une $\mathcal{V}$-alg\`ebre lisse relevant $A_{0}$. D'apr\`es (1.2.3) on a un diagramme commutatif de $\mathcal{V}$-sch\'emas formels de type fini, \`a carr\'e cart\'esien\\

$$
\begin{array}{c}
\xymatrix{
& & & &  \overline{\mathcal{S}}_{\alpha}\\
& & & \overline{\mathcal{S}}''_{\alpha}   \ar[ur]_{v_{\overline{\mathcal{S}}_{\alpha}}}  \ar[d]^{v_{\overline{\mathcal{S}}'_{\alpha}}} & &\\
S_{\alpha}  \ar@{^{(}->}[rrr]_{j_{\overline{\mathcal{S}}'_{\alpha}}}    \ar@{^{(}->}[urrr]_{j_{\overline{\mathcal{S}}''_{\alpha}}}      \ar@{^{(}->}[uurrrr]^{j_{\overline{\mathcal{S}}_{\alpha}}} \ar[d]_{F_{\alpha}} &  & &  \overline{\mathcal{S}}'_{\alpha} \ar[d]^{\overline{F}_{\alpha}}    &  \\
\mathcal{S}_{\alpha}  \ar@{^{(}->}[rrr]^{j_{\overline{\mathcal{S}}_{\alpha}}}  & & & \overline{\mathcal{S}}_{\alpha} & 
}
\end{array}
\leqno{(2.1.3)}  
$$

\noindent o\`u $\mathcal{S}_{\alpha} = Spf \hat{A}$, $\overline{\mathcal{S}}_{\alpha}$ est propre sur $\mathcal{V}$, $v_{\overline{\mathcal{S}}'_{\alpha}}$ et $v_{\overline{\mathcal{S}}_{\alpha}}$ sont propres, $F_{\alpha}$ est un rel\`evement fini et plat du Frobenius $F_{S_{\alpha}}$ de $S_{\alpha}$, $\overline{F}_{\alpha}$ est fini et les $j$ sont des immersions ouvertes. Notons $j_{\overline{\mathcal{T}}_{\alpha}} : \mathcal{T}_{\alpha} := \mathcal{S}_{\alpha} \times_{\mathcal{V}} \mathcal{S} \longrightarrow \overline{\mathcal{T}}_{\alpha} : = \overline{\mathcal{S}}_{\alpha} \times_{\mathcal{V}} \overline{\mathcal{S}}$ l'immersion ouverte et $\overline{u}_{\alpha} : \overline{\mathcal{T}}_{\alpha} \longrightarrow \overline{\mathcal{S}}_{\alpha}, \overline{v}_{\alpha} : \overline{\mathcal{T}}_{\alpha} \longrightarrow \overline{\mathcal{S}}$ les projections ; 
soient $\overline{T}_{\alpha}$ (resp. $T_{\alpha}$) la r\'eduction de $\overline{\mathcal{T}_{\alpha}}$ (resp. $\mathcal{T}_{\alpha}$) mod $\mathfrak{m}$ et $\tilde{S}_{\alpha}$ l'image sch\'ematique de $S_{\alpha}$ plong\'e diagonalement dans $\overline{T}_{\alpha}$ :

$$
\xymatrix{
S_{\alpha} \ar@{^{(}->}[r]_{j_{\tilde{S}_{\alpha}}} & \tilde{S}_{\alpha}  \ar@{^{(}->}[r]^{i_{\overline{T}_\alpha}}  \ar@/_1pc/[rr]_{i_{\tilde{S}_{\alpha}}} & \overline{T}_{\alpha} \ar@{^{(}->}[r]^{i_{\overline{\mathcal{T}}_\alpha}} & \overline{\mathcal{T}}_{\alpha}
}
$$

\noindent $j_{\tilde{S}_{\alpha}}$ est une immersion ouverte et les $i$ des immersions ferm\'ees. On a alors un diagramme commutatif \`a carr\'es verticaux cart\'esiens\\

$$
\begin{array}{c}
 \shorthandoff{;:!?}
 \xymatrix@!0 @R=1cm @C=2cm{
&&X\ar@{.>}[dd]^{f}\  \ar @{^{(}->}[rr]^{j_{\overline{X}}}&&\overline{X}  \ar@{.>}[dd]^{\overline{f}}\  \ar@{^{(}->}[rr]^{i_{\overline{X}}}&&\overline{\mathcal{X} }\ar[dd]^{\overline{h}}\\
&&&&&&\\
&&S\  \ar@{^{(}.>}[rr]^{j_{\overline{S}}}  && \overline{S}\  \ar@{^{(}.>}[rr]^{i_{\overline{S}}}&& \overline{\mathcal{S}}\\
&X_{\alpha}\ar@{.>}[dd]^{f_{\alpha}} \ar@{^{(}->}[uuur]^{j_{X_{\alpha}}} \ar@{^{(}->}[rr]  &&\overline{X}\ar@{.>}[dd]_{\overline{f}} \  \ar@{^{(}->}[rr]^(.7){i_{\overline{X}}} \ar@{=}[uuur]&& \overline{\mathcal{X}} \ar[dd]^{\overline{h}} \ar@{=}[uuur]\\
&&&&&&\\
& S_{\alpha}\  \ar@{^{(}.>}[rr]^{\overline{j}_{\alpha}} \ar@{.>}[uuur]_{j_{S_{\alpha}}} &&\overline{S} \  \ar@{^{(}.>}[rr] ^{i_{\overline{S}}} \ar@{.>}[uuur]_{id}&& \overline{\mathcal{S}} \ar@{=}[uuur]_{id}&\\
X_{\alpha} \ar@{^{(}->}[r] \ar[dd]_{f_{\alpha}}\   \ar@{=}[uuur]^{id}& \tilde{X}_{\alpha}\ar@{^{(}->}[r] \ar[dd]&\overline{Y}_{\alpha}\ar[dd]\  \ar@{^{(}->}[rr] \ar[uuur]&& \overline{\mathcal{Y}}_{\alpha}\ar[dd]^{\overline{h}_{\alpha}} \ar[uuur]&&\\
&&&&&&\\
S_{\alpha} \  \ar@{^{(}->}[r]_{j_{\tilde{S}_{\alpha}}}  \ar@{.>}[uuur]_{id}& \tilde{S}_{\alpha}\ar@{^{(}->}[r]&\overline{T}_{\alpha} \ar@{^{(}->}[rr] \  \ar@{.>}[uuur]&&\ \overline{\mathcal{T}}_{\alpha}\ . \ar[uuur]_{\overline{v}_{\alpha}}&&
}
\end{array}
\leqno{(2.1.4)}
$$

Ainsi on a une suite d'isomorphismes

$$
j^{\ast}_{S_{\alpha}}\ R^i f_{\textrm{rig}\ast}(X/\overline{\mathcal{S}}, E) \tilde{\longrightarrow} R^i f_{\alpha \textrm{rig}\ast} (X_{\alpha}/ \overline{\mathcal{S}}, j^{\ast}_{X_{\alpha}} E)\ [\textrm{II, th\'eo(3.4.4)}]
$$

\noindent (2.1.5) $\qquad \simeq \overline{v}^{\ast}_{\alpha_{K}} R^i f_{\alpha \textrm{rig}\ast}(X_{\alpha}/ \overline{\mathcal{S}}, j^{\ast}_{X_{\alpha}} E)\ [\textrm{B 3, (2.3.6), (2.3.2) (iv)}]$\\

$\qquad \qquad \tilde{\rightarrow} R^i f_{\alpha \textrm{rig}\ast}(X_{\alpha}/ \overline{\mathcal{T}_{\alpha}}, j^{\ast}_{X_{\alpha}} E)\ [\textrm{II, th\'eo (3.4.4)}].$\\

\noindent On est donc ramen\'e \`a construire un isomorphisme de Frobenius sur

$$
R^i f_{\alpha \textrm{rig}\ast}(X_{\alpha}/ \overline{\mathcal{T}_{\alpha}}, j^{\ast}_{X_{\alpha}} E) \in \textrm{Isoc}^{\dag}(S_{\alpha}/K).
$$

\noindent A partir du carr\'e commutatif

$$
\xymatrix{
\mathcal{T}_{\alpha} \ar@{^{(}->}[rr]^(.4){j_{\overline{\mathcal{T}}'_{\alpha}}} \ar[d]_{F_{\mathcal{T}_{\alpha}} := F_{\alpha} \times 1} && \overline{\mathcal{T}}'_{\alpha} : = \overline{\mathcal{S}}'_{\alpha} \times_{\mathcal{V}} \overline{\mathcal{S}} \ar[d]^{\overline{F}_{\alpha} \times 1 =: F_{\overline{\mathcal{T}}'_{\alpha}}} &\\
\mathcal{T}_{\alpha} \ar@{^{(}->}[rr]_(.4){j_{\overline{\mathcal{T}}_{\alpha}}}  &&   \overline{\mathcal{T}}_{\alpha} = \overline{\mathcal{S}}_{\alpha} \times_{\mathcal{V}} \overline{\mathcal{S}} & ,
}
$$

\noindent d\'eduit de (2.1.3), et du carr\'e cart\'esien

$$
\xymatrix{
\mathcal{Y}_{\alpha} \ar@{^{(}->}[r] \ar[d]_{h_{\alpha}} & \overline{\mathcal{Y}_{\alpha}} \ar[d]^{\overline{h}_{\alpha}} &\\
\mathcal{T}_{\alpha} \ar@{^{(}->}[r]_{j_{\overline{\mathcal{T}}_{\alpha}}} & \overline{\mathcal{T}}_{\alpha} & ,
}
$$

\noindent o\`u les fl\`eches $j$ sont des immersions ouvertes, on forme les diagrammes commutatifs \`a carr\'es cart\'esiens\\

$$
\begin{array}{c}
\xymatrix{
& X_{\alpha} \ar@{.>}[dd]^(.3){f_{\alpha}} |\hole \ar[rr] & & Y_{\alpha} \ar@{.>}[dd] |\hole \ar@{^{(}->}[rr] & & \mathcal{Y}_{\alpha} \ar@{.>}[dd]^(.3){h_{\alpha}}  \ar@{^{(}->}[rr] && \overline{\mathcal{Y}}_{\alpha}\ar[dd]^{\overline{h}_{\alpha}} \\
X^{(q/S_{\alpha})}_{\alpha} \ar[rr] \ar[dd]_{f^{(q)}_{\alpha}} \ar[ur] & & Y'_{\alpha}\ar@{^{(}->}[rr] \ar[dd] \ar[ur]  & & \mathcal{Y}'_{\alpha} \ar[dd]  \ar[ur]  \ar@{^{(}->}[rr]  && \overline{\mathcal{Y}}'_{\alpha} \ar[ur] \ar[dd] \\
& S_{\alpha} \ar@{.>}[rr]^(.7){i_{\alpha}} |\hole & & T_{\alpha}\ar@{^{(}.>}[rr]^(.4){i_{\mathcal{T}_{\alpha}}}  |\hole & &  \mathcal{T}_{\alpha}  \ar@{^{(}.>}[rr]^(.4){j_{\overline{\mathcal{T}}_{\alpha}}}  |\hole && \overline{\mathcal{T}_{\alpha}}\\
S_{\alpha} \ar[rr]_{i_{\alpha}} \ar@{.>}[ur]^{F_{S_{\alpha}}} & & T_{\alpha} \ar@{^{(}->}[rr]_{i_{\mathcal{T}_{\alpha}}} \ar@{.>}[ur]_{F_{S_{\alpha}} \times 1} & & \mathcal{T}_{\alpha} \ar@{^{(}->}[rr]_{j_{\overline{\mathcal{T}}'_{\alpha}}} \ar@{.>}[ur]_{F_{\mathcal{T}_{\alpha}}} && \overline{\mathcal{T}}'_{\alpha} \ar[ur]_{F_{\overline{\mathcal{T}}'_{\alpha}}}
}
\end{array}
\leqno{(2.1.6)}
$$

\noindent (o\`u $i_{\alpha}$ et $i_{\mathcal{T}_{\alpha}}$ sont des immersions ferm\'ees), et \\

$$
\begin{array}{c}
\xymatrix{
& X_{\alpha} \ar@{.>}[dd]^(.3){f_{\alpha}} |\hole \ar@{^{(}->}[rr]^{j_{\tilde{X}_{\alpha}}} & & \tilde{X}_{\alpha} \ar@{.>}[dd]^(.3){\tilde{f}_{\alpha}} |\hole \ar[rr]& & \overline{\mathcal{Y}}_{\alpha} \ar[dd]^{\overline{h}_{\alpha}} \\
X^{(q/S_{\alpha})}_{\alpha} \ar@{^{(}->}[rr] \ar[dd]_{f^{(q)}_{\alpha}} \ar[ur]^{\pi_{X_{\alpha}}} & & \tilde{X}'_{\alpha} \ar[rr] \ar[dd]_(.3){\tilde{f}'_{\alpha}} \ar[ur]  & & \overline{\mathcal{Y}}'_{\alpha} \ar[dd]^(.3){\overline{h}'_{\alpha}} \ar[ur] \\
& S_{\alpha} \ar@{^{(}.>}[rr]^(.7){j_{\tilde{S}_{\alpha}}} |\hole & &\tilde{S}'_{\alpha}\ar@{^{(}.>}[rr]^(.4){i_{\tilde{S}_{\alpha}}}  |\hole & &  \overline{\mathcal{T}}'_{\alpha}  \\
S_{\alpha} \ar@{^{(}->}[rr]_{j_{\tilde{S}'_{\alpha}}} \ar@{.>}[ur]^{F_{S_{\alpha}}} & & \tilde{S}'_{\alpha} \ar@{^{(}->}[rr]_{i_{\tilde{S}'_{\alpha}}} \ar@{.>}[ur]_{\tilde{F}_{\alpha}} & & \  \overline{\mathcal{T}}'_{\alpha}\ . \ar[ur]_{F_{\overline{\mathcal{T}}'_{\alpha}}} 
}
\end{array}
\leqno{(2.1.7)}
$$

D'apr\`es [II, th\'eo (3.4.4)] on a un isomorphisme\\

\noindent (2.1.8)  $\qquad  F^{\ast}_{S_{\alpha}} R^i f_{\alpha\  \textrm{rig}^{\ast}} (X_{\alpha}/\overline{\mathcal{T}_{\alpha}}, j^{\ast}_{X_{\alpha}} E) \tilde{\rightarrow} R^i f^{(q)}_{\alpha\ \textrm{rig}^{\ast}} (X^{(q/S_{\alpha})}_{\alpha} / \overline{\mathcal{T}}'_{\alpha}, \pi^{\ast}_{X_{\alpha}}\ j^{\ast}_{X_{\alpha}} E).$ \\

Notons $v_{\overline{\mathcal{T}}_{\alpha}} := v_{\overline{\mathcal{S}}_{\alpha}} \times 1_{\overline{\mathcal{S}}} : \overline{\mathcal{T}}''_{\alpha} = \overline{\mathcal{S}}''_{\alpha} \times_{\mathcal{V}}    \overline{\mathcal{S}} \rightarrow  \overline{\mathcal{T}}_{\alpha} = \overline{\mathcal{S}}_{\alpha} \times_{\mathcal{V}} \overline{\mathcal{S}}$ et $v_{\overline{\mathcal{T}}'_{\alpha}} := v_{\overline{\mathcal{S}}'_{\alpha}} \times 1_{\overline{\mathcal{S}}} : \overline{\mathcal{T}}''_{\alpha} = \overline{\mathcal{S}}''_{\alpha} \times_{\mathcal{V}}  \overline{\mathcal{S}} \rightarrow \overline{\mathcal{T}}'_{\alpha} = \overline{\mathcal{S}}'_{\alpha} \times_{\mathcal{V}} \overline{\mathcal{S}}$ ; de (2.1.3) on d\'eduit un diagramme commutatif\\

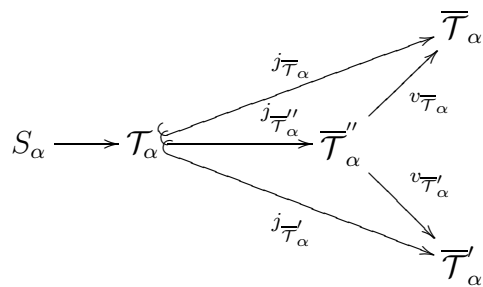
$$
\xymatrix{
& & & & \overline{\mathcal{T}}_{\alpha}\\
S_{\alpha} \ar[r] &  \mathcal{T}_{\alpha} \ar@{^{(}->}[rr]^(.7){j_{\overline{\mathcal{T}}^{''}_{\alpha}}} \ar@{^{(}->}[rrrd]_{j_{\overline{\mathcal{T}}^{'}_{\alpha}}}   \ar@{^{(}->}[urrr]^{j_{\overline{\mathcal{T}}_{\alpha}}}  && \overline{\mathcal{T}}^{''}_{\alpha} \ar[rd]^{v_{\overline{\mathcal{T}}'_{\alpha}}}  \ar[ur]_{v_{\overline{\mathcal{T}}_{\alpha}}} \\
&&& &  \overline{\mathcal{T}}'_{\alpha}
}
$$

\noindent o\`u les $j$ sont des immersions ouvertes et $S_{\alpha} \rightarrow \mathcal{T}_{\alpha}$ une immersion.\\

On a un diagramme commutatif dont les carr\'es verticaux sont cart\'esiens, de m\^eme que le carr\'e horizontal  $\ \rondI$  en bas \`a droite\\

$$
\begin{array}{c}
\xymatrix{
& X^{(q/S_{\alpha})} \ar@{.>}[dd]^(.3){f^{(q)}_{\alpha}} |\hole \ar@{^{(}->}[rr]& & \tilde{X}'_{\alpha} \ar@{.>}[dd]^(.3){\tilde{f}'_{\alpha}} |\hole \ar[rr]& & \overline{\mathcal{Y}}'_{\alpha} \ar[dd]^{\overline{h}'_{\alpha}} \\
X^{(q/S_{\alpha})} \ar@{^{(}->}[rr] \ar[dd]_{f^{(q)}_{\alpha}} \ar@{=} [ur]^{\textrm{id}} & & \tilde{X}''_{\alpha} \ar[rr] \ar[dd]_(.3){\tilde{f}''_{\alpha}} \ar[ur]  & & \overline{\mathcal{Y}}''_{\alpha} \ar[dd]^(.3){\overline{h}''_{\alpha}} \ar[ur]& \\
& S_{\alpha} \ar@{^{(}.>}[rr]^(.7){j_{\tilde{S}'_{\alpha}}} |\hole & &\tilde{S}'_{\alpha}\ar@{^{(}.>}[rr]^(.4){i_{\tilde{S}'_{\alpha}}}  |\hole & &  \overline{\mathcal{T}}'_{\alpha}  \\
S_{\alpha} \ar@{^{(}->}[rr]_{j_{\tilde{S}''_{\alpha}}} \ar@{==}[ur]^{\textrm{id}} & & \tilde{S}'' _{\alpha} \ar@{^{(}->}[rr]_{i_{\tilde{S}''_{\alpha}}} \ar@{.>}[ur] \ar@{}[urrr] |{\rondI} & & \  \overline{\mathcal{T}}''_{\alpha} \ .\ar[ur]_{v_{\overline{\mathcal{T}}'_{\alpha} }} 
}
\end{array}
\leqno{(2.1.9)}
$$

\noindent D'apr\`es [II, (3.4.4)], $v^{\ast}_{\overline{\mathcal{T}'}_{\alpha} K}$ induit un isomorphisme\\

\noindent (2.1.10)$ \qquad R^i f^{(q)}_{\alpha\  \textrm{rig}^{\ast}}(X^{(q/S_{\alpha})}_{\alpha} / \overline{\mathcal{T}}'_{\alpha}, \pi^{\ast}_{X_{\alpha}}\  j^{\ast}_{X_{\alpha}} E) \tilde{\rightarrow}  R^i f^{(q)}_{\alpha\  \textrm{rig}^{\ast}} (X^{(q/S_{\alpha})} / \overline{\mathcal{T}}''_{\alpha},  \pi^{\ast}_{X_{\alpha}}\  j^{\ast}_{X_{\alpha}} E) : $\\

\noindent en effet on peut appliquer [loc.cit.] car le morphisme propre $v_{\overline{\mathcal{T'}}_{\alpha}}$, \'etant \'etale au voisinage de $S_{\alpha}$, il induit [B 3, (1.3.5)] un isomorphisme entre un voisinage strict de $]S_{\alpha}[_{\overline{\mathcal{T}}''_{\alpha}}$ dans $\overline{\mathcal{T}}''_{\alpha K}$ et un voisinage strict de $]S_{\alpha}[_{\overline{\mathcal{T}}'_{\alpha}}$ dans $\overline{\mathcal{T}}'_{\alpha K}$.\\
Notons $\overline{T}''_{\alpha}$ la r\'eduction de $\overline{\mathcal{T}}''_{\alpha}\  \textrm{mod}\  \mathfrak{m}, v_{\overline{T}_{\alpha}}$ la r\'eduction de $v_{\overline{\mathcal{T}}_{\alpha}} \  \textrm{mod}\  \mathfrak{m}$ et $\tilde{S}'''_{\alpha }$ l'adh\'erence sch\'ematique de $S_{\alpha}$ dans ${\overline{T''}_{\alpha}}$ ; on a un diagramme commutatif o\`u les carr\'es $\  \rondI$ et $\ \rondII$  sont cart\'esiens\\

$$
\begin{array}{c}
\xymatrix{
& & \tilde{S}''_{\alpha} \ar@{^{(}->}[rrdd]^{i''} & & &&\\
& \tilde{S}'''_{\alpha} \ar@{^{(}->}[rd]^{i'''} \ar@{^{(}->}[ur]^{i'} & &&& &\\
& & v^{-1}_{\overline{T}_{\alpha}}(\tilde{S}_{\alpha}) \ar@{^{(}->}[rr] \ar[d]^{v_{\tilde{S}_{\alpha}}} \ar @{}[drr]|{\rondI}&& \overline{T}''_{\alpha} \ar@{^{(}->}[rr] \ar[d]^{v_{\overline{T}_{\alpha}}} \ar@{}[drr]|{\rondII}&& \overline{\mathcal{T}}''_{\alpha} \ar[d]^{v_{\overline{\mathcal{T}}_{\alpha}}} \\
S_{\alpha} \ar@{^{(}->}[rr]_{j_{\tilde{S}_{\alpha}}} \ar@{^{(}->}[uur]^{j_{\tilde{S}'''_{\alpha}}} & & \tilde{S}_{\alpha} \ar@{^{(}->}[rr]_{i_{\overline{T}_{\alpha}}}  && \overline{T}_{\alpha} \ar@{^{(}->}[rr]  && \overline{\mathcal{T}}_{\alpha}
}
\end{array}
\leqno{(2.1.11)}
$$

\noindent o\`u les $j$ (resp. les $i$) sont des immersions ouvertes (resp. ferm\'ees). Posons $v_{\alpha} = v_{\tilde{S}_{\alpha}} \circ i'''$.\\

\noindent Soit $\tilde{f}'''_{\alpha} : \tilde{X}'''_{\alpha} \rightarrow \tilde{S}'''_{\alpha}$ l'image inverse de $\tilde{f}''_{\alpha} : \tilde{X}''_{\alpha} \rightarrow \tilde{S}''_{\alpha}$ par $i' : \tilde{S}'''_{\alpha} \hookrightarrow \tilde{S}''_{\alpha}$. On a un diagramme commutatif \\

$$
\begin{array}{c}
\xymatrix{
X^{(q/S_{\alpha})} \ar@{^{(}->}[r] \ar[d]_{f^{(q)}_{\alpha}} & \tilde{X}'''_{\alpha} \ar@{^{(}->}[r] \ar[d]^{\tilde{f}'''_{\alpha}}  & \overline{\mathcal{Y}''}_{\alpha} \ar[d]^{\overline{h}''_{\alpha}} \\
S_{\alpha} \ar@{^{(}->}[r]^{j_{\tilde{S}'''_{\alpha}}} \ar@{=}[d]  & \tilde{S}'''_{\alpha} \ar[d]^{i_{\mathcal{T}}} \ar@{^{(}->}[r]^{i_{\tilde{S}'''_{\alpha}}} & \overline{\mathcal{T}''_{\alpha}} \ar[d]^{v_{\overline{\mathcal{T}''_{\alpha}}}}\\
S_{\alpha} \ar@{^{(}->}[r]^{j_{\tilde{S}_{\alpha}}} & \tilde{S}_{\alpha} \ar@{^{(}->}[r]^{i_{\tilde{S}_{\alpha}}} & \overline{\mathcal{T}_{\alpha}}
}
\end{array}
\leqno{(2.1.12)}
$$

\noindent o\`u les $j$ (resp. les $i$) sont des immersions ouvertes (resp. ferm\'ees). Comme $v_{\overline{\mathcal{T}}_{\alpha}}$ est \'etale au voisinage de $S_{\alpha}$ et que $v_{\alpha}$ est propre on d\'eduit de [B 3, th\'eo (1.3.5)] que $v_{\overline{\mathcal{T}}_{\alpha K}}$ induit un isomorphisme entre un voisinage strict de $]S_{\alpha}[_{\overline{\mathcal{T}}''_{\alpha}}$ dans $]\tilde{S}'''_{\alpha}[_{\overline{\mathcal{T}}''_{\alpha}}$ et un voisinage strict de $]S_{\alpha}[_{\overline{\mathcal{T}}_{\alpha}}$ dans  $]\tilde{S}_{\alpha}[_{\overline{\mathcal{T}}_{\alpha}}$. Par suite [II, th\'eo (3.4.4)] $v_{\overline{\mathcal{T}}_{\alpha K}}$ induit un isomorphisme\\

\noindent (2.1.13) $ R^i f^{(q)}_{\alpha\  \textrm{rig}^{\ast}}(X^{(q/S_{\alpha})}_{\alpha} / \overline{\mathcal{T}''}_{\alpha}, \pi^{\ast}_{X_{\alpha}}\  j^{\ast}_{X_{\alpha}} E) \tilde{\rightarrow}  R^i f^{(q)}_{\alpha\  \textrm{rig}^{\ast}} (X^{(q/S_{\alpha})} / \overline{\mathcal{T}}_{\alpha},  \pi^{\ast}_{X_{\alpha}}\  j^{\ast}_{X_{\alpha}} E).$\\

Par composition des isomorphismes (2.1.8), (2.1.10) et (2.1.13) on obtient un isomorphisme induit par $\pi_{X_{\alpha}} : X^{(q/S_{\alpha})}_{\alpha} \rightarrow X_{\alpha}$

$$ 
\hspace{-1cm}\begin{array}{c}
\xymatrix{
F^{\ast}_{S_{\alpha}}\, R^i\, f_{\alpha \textrm{rig}^{\ast}}(X_{\alpha}/ \overline{\mathcal{T}_{\alpha}}; j^{\ast}_{X_{\alpha}}\, E) \ar[r]^{\sim} \ar[d]^{\simeq} &  R^i\, f^{(q)}_{\alpha \textrm{rig}^{\ast}}(X^{(q/S_{\alpha}}/ \overline{\mathcal{T}}_{\alpha}; \pi^{\ast}_{X_{\alpha}} j^{\ast}_{X_{\alpha}} E)\\
F^{\ast}_{\overline{\mathcal{T}}'_{\alpha K}} R^i f_{\alpha \textrm{rig}^{\ast}}(X_{\alpha}/\overline{\mathcal{T}_{\alpha}}; j^{\ast}_{X_{\alpha}} E) . & 
}
\end{array}
\leqno{(2.1.14)}
$$

Si l'on prend l'image inverse de cet isomorphisme par $j_{\overline{\mathcal{T}}'_{\alpha}} : \mathcal{T}_{\alpha} \rightarrow \overline{\mathcal{T}}'_{\alpha}$, on voit ais\'ement en suivant les diagrammes pr\'ec\'edents que l'on obtient l'isomorphisme en cohomologie convergente, induit par $\pi_{X_{\alpha}}$ :\\

\noindent (2.1.15) $\quad F^{\ast}_{S_\alpha} R^i f_{\alpha\  \textrm{conv}^{\ast}} (X_{\alpha} / \mathcal{T}_{\alpha}, j^{\ast}_{X_{\alpha}} \mathcal{E})\  \tilde{\rightarrow}\  R^{i} f^{(q)}_{\alpha\  \textrm{conv}^{\ast}} (X^{(q/S_{\alpha}}) / \mathcal{T}_{\alpha},\pi^{\ast}_{X_{\alpha}}\  j^{\ast}_{X_{\alpha}} \mathcal{E})$\\

\noindent o\`u l'on a, comme pour (2.1.5), un isomorphisme\\

\noindent (2.1.16) $\qquad \qquad j^{\ast}_{S_{\alpha}} R^i f_{\textrm{conv}^{\ast}} (X / \mathcal{S}, \mathcal{E})\  \tilde{\rightarrow}\  R^i f_{\alpha\ \textrm{conv}^{\ast}} (X_{\alpha} / \mathcal{T}_{\alpha}, j^{\ast}_{X_{\alpha}} \mathcal{E}).$\\

Si $F_{\tilde{S}_{\alpha}}$ d\'esigne le Frobenius (puissance $q$) de $\tilde{S}_{\alpha}$, on notera $\tilde{X}_{\alpha}^{(q/\tilde{S}_{\alpha})}$ le produit fibr\'e  d\'efini par le diagramme \`a carr\'e cart\'esien

$$
\xymatrix{
\tilde{X}_{\alpha} \ar[r]^(.4){F_{\tilde{X}_{\alpha}/\tilde{S}_{\alpha}}} \ar[rd]_{\tilde{f}_{\alpha}} \ar@/^4pc/[rr]^{F_{\tilde{X}_{\alpha}}} & \tilde{X}_{\alpha}^{(q/\tilde{S}_{\alpha})} \ar[r]^{\pi_{\tilde{X}_{\alpha}}} \ar[d]^{\tilde{f}^{(q)}_{\alpha}} & \tilde{X}_{\alpha}  \ar[d]^{\tilde{f}_{\alpha}} \\
& \tilde{S}_{\alpha} \ar[r]_(.48){F_{\tilde{S}_{\alpha}}}  \ar@{^{(}->}[d]_{i_{\tilde{S}_{\alpha}}} & \tilde{S}_{\alpha}\\
& \overline{\mathcal{T}_{\alpha}} &.
}
$$

\noindent Le calcul de $R^i f^{(q)}_{\alpha\  \textrm{rig}^{\ast}}(X^{(q/S_{\alpha})} / \overline{\mathcal{T}}_{\alpha}, \pi^{\ast}_{X_{\alpha}}\  j^{\ast}_{X_{\alpha}} E)$ \'etant ind\'ependant de la compactification de $X^{(q/S_{\alpha}})$ choisie [B 5, (3.1.11), (3.1.12), (3.2.3)], nous choisirons dor\'enavant  $\tilde{X}_{\alpha}^{(q/ \tilde{S}_{\alpha})}$ [C-T, \S\ 10] comme compactification de $X^{(q/s_{\alpha})}$ au lieu de $\tilde{X}'''_{\alpha}$. Consid\'erons le diagramme commutatif

$$
\xymatrix{
X_{\alpha} \ar[rr]^{F_{X_{\alpha}/S_{\alpha}}} \ar@{^{(}->}[d]_{j_{\tilde{X}_{\alpha}}} && X_{\alpha}^{(q/S_{\alpha})} \ar[d]  \\
\tilde{X}_{\alpha} \ar[rr]^{F_{\tilde{X}_{\alpha}/\tilde{S}_{\alpha}}} \ar[d]_{i_{\tilde{S}_{\alpha}} \circ \tilde{f}_{\alpha}} && \tilde{X}^{(q/\tilde{S}_{\alpha})}_{\alpha}    \ar[d]^{i_{\tilde{S}_{\alpha}} \circ \tilde{f}^{(q)}_{\alpha}}   \\
\overline{\mathcal{T}}_{\alpha} \ar@{=}[rr] && \overline{\mathcal{T}}_{\alpha}   ;  \\
} 
$$

\noindent $F_{X_{\alpha}/S_{\alpha}}$ d\'efinit par fonctorialit\'e [B 5, (3.1.11) (ii), (3.1.12) (i)] ou [C-T, (10.5.2)]  un morphisme\\

$$
\begin{array}{c}
\xymatrix{
F^{\ast}_{X_{\alpha}/S_{\alpha}} :  R^i f^{(q)}_{\alpha\  \textrm{rig}^{\ast}} (X^{(q/S_{\alpha})} / \overline{\mathcal{T}}_{\alpha}; \pi^{\ast}_{X_{\alpha}}\  j^{\ast}_{X_{\alpha}} E) \ar[r] & R^i f_{\alpha\  \textrm{rig}^{\ast}}(X_{\alpha} / \overline{\mathcal{T}}_{\alpha}; F^{\ast}_{X_{\alpha}/S_{\alpha}}  \pi^{\ast}_{X_{\alpha}}  j^{\ast}_{X_{\alpha}} E) \ar[d]^{\simeq}\\
  R^i f_{\alpha\  \textrm{rig}^{\ast}} (X_{\alpha}/ \overline{\mathcal{T}}_{\alpha};  j^{\ast}_{X_{\alpha}} F^{\ast}_{X} E)\ar[r]^{ \simeq} &R^i f_{\alpha\  \textrm{rig}^{\ast}} (X_{\alpha}/  \overline{\mathcal{T}}_{\alpha}; F^{\ast}_{X_{\alpha}}\  j^{\ast}_{X_{\alpha}} E) .
}
\end{array}
\leqno{(2.1.17)}
$$

Par image inverse par $j_{\overline{\mathcal{T}}_{\alpha}} : \mathcal{T}_{\alpha} \hookrightarrow \overline{\mathcal{T}}_{\alpha}$ il fournit le morphisme en cohomologie convergente, induit par $F_{X_{\alpha}/S_{\alpha}}$\\

\noindent (2.1.18) $F^{\ast}_{X_{\alpha}/S_{\alpha}} : R^i f^{(q)}_{\alpha\  \textrm{conv}^{\ast}} (X^{(q/S_{\alpha})} / \mathcal{T}_{\alpha}, \pi^{\ast}_{X_{\alpha}}\  j^{\ast}_{X_{\alpha}} \mathcal{E}) \rightarrow R^i f_{\alpha\  \textrm{conv}^{\ast}}(X_{\alpha} / \mathcal{T}_{\alpha}, F^{\ast}_{X_{\alpha}}\  j^{\ast}_{X_{\alpha}} \mathcal{E}).$ \\

Enfin, l'isomorphisme de Frobenius de $E, \phi_{E} : F^{\ast}_{X} \tilde{\rightarrow} E$, fournit par fonctorialit\'e un isomorphisme\\

\noindent (2.1.19) $\qquad  R^i f_{\alpha\  \textrm{rig}^{\ast}} (X_{\alpha} / \overline{\mathcal{T}}_{\alpha}, j^{\ast}_{X_{\alpha}}  F^{\ast}_{X} E)\  \tilde{\rightarrow}\  R^{i} f_{\alpha\  \textrm{rig}^{\ast}} (X_{\alpha} / \overline{\mathcal{T}}_{\alpha},    j^{\ast}_{X_{\alpha}} \mathcal{E}),$\\

\noindent dont l'image inverse par $j_{\overline{\mathcal{T}}_{\alpha}} : \mathcal{T}_{\alpha} \rightarrow \overline{\mathcal{T}}_{\alpha}$ est l'isomorphisme induit en cohomologie convergente par $\phi_{\mathcal{E}} : F^{\ast}_{X}\  \mathcal{E}\  \tilde{\rightarrow}\  \mathcal{E}$,\\

\noindent (2.1.20) $\qquad R^i f_{\alpha\  \textrm{conv}^{\ast}} (X_{\alpha} / \mathcal{T}_{\alpha}, j^{\ast}_{X_{\alpha}}  F^{\ast}_{X} \mathcal{E})\  \tilde{\rightarrow}\  R^{i} f_{\alpha\  \textrm{conv}^{\ast}} (X_{\alpha} / \mathcal{T}_{\alpha},    j^{\ast}_{X_{\alpha}} \mathcal{E}).$\\

En composant les morphismes (2.1.14), (2.1.17) et (2.1.19) [resp. (2.1.15), (2.1.18) et (2.1.20)] on obtient le morphisme de Frobenius souhait\'e\\

\noindent (2.1.21) $\qquad  \phi_{E_{\alpha_{i}}} : F^{\ast}_{S_{\alpha}} R^i f_{\alpha\  \textrm{rig}^{\ast}} (X_{\alpha} / \overline{\mathcal{T}}_{\alpha},   j^{\ast}_{X_{\alpha}} E)\  \rightarrow\  R^i f_{\alpha\  \textrm{rig}^{\ast}}(X_{\alpha} / \overline{\mathcal{T}}_{\alpha}, \   j^{\ast}_{X_{\alpha}} E)$\\

\noindent [resp.\\

\noindent (2.1.22) $\qquad  \phi_{\mathcal{E}_{\alpha_{i}}} : F^{\ast}_{S_{\alpha}} R^i f_{\alpha\  \textrm{conv}^{\ast}} (X_{\alpha} / \mathcal{T}_{\alpha},   j^{\ast}_{X_{\alpha}} \mathcal{E})\  \rightarrow\  R^i f_{\alpha\  \textrm{conv}^{\ast}}(X_{\alpha} / \mathcal{T}_{\alpha}, \   j^{\ast}_{X_{\alpha}} \mathcal{E})$\\

\noindent qui est l'image inverse de $\phi_{E_{\alpha_{i}}}$ par $j_{\overline{\mathcal{T}}_{\alpha}}].$\\

D'apr\`es le [III, th\'eo (3.3.1)] $\phi_{\mathcal{E}_{\alpha_{i}}}$ est un isomorphisme : nous allons en d\'eduire que $\phi_{E_{\alpha_{i}}}$ est un isomorphisme, ce qui ach\`evera la preuve du th\'eor\`eme.\\

On a un diagramme commutatif\\

$$
\begin{array}{c}
\xymatrix{
& \tilde{S}_{\alpha} \ar[d]^{i_{\overline{T}_{\alpha}}} \ar@{^{(}->}[rd]^{i_{\tilde{S}_{\alpha}}}  \ar@/_1pc/[dd]_{v} &\\
S_{\alpha} \ar@{^{(}->}[ru]^{j_{\tilde{S}_{\alpha}}}    \ar@{^{(}->}[rd]_{j_{\overline{S}_{\alpha}}} & \overline{T}_{\alpha} \ar[d]^{\overline{u}'_{\alpha}} \ar@{^{(}->}[r]_{i_{\overline{T}_{\alpha}}}   & \overline{\mathcal{T}}_{\alpha} \ar[d]^{\overline{u}_{\alpha}}\\
& \overline{S}_{\alpha} \ar@{^{(}->}[r]_{i_{\overline{S}_{\alpha}}} & \overline{\mathcal{S}}_{\alpha}
}
\end{array}
\leqno{(2.1.23)}
$$

\noindent o\`u $j_{\overline{S}_{\alpha}}$ (resp. $\overline{u}'_{\alpha}$) est la r\'eduction mod $\mathfrak{m}$ de $j_{\overline{\mathcal{S}}_{\alpha}} : \mathcal{S}_{\alpha} \hookrightarrow \overline{\mathcal{S}}_{\alpha}$ (resp. de $\overline{u_{\alpha}} : \overline{\mathcal{T}_{\alpha}} \rightarrow \overline{\mathcal{S}_{\alpha}})$ : les $j$ (resp. les $i$) sont des immersions ouvertes (resp. ferm\'ees).\\

De m\^eme le triangle commutatif\\

$$
\xymatrix{
& \mathcal{T}_{\alpha} = \mathcal{S}_{\alpha} \times_{\mathcal{V}} \mathcal{S} \ar[dd]^{u_{\alpha}} \\
\mathcal{S}_{\alpha} \ar@{^{(}->} [ur]^{\Delta} \ar@{=}[rd]_{\textrm{id}}\\
& \mathcal{S}_{\alpha}
}
$$

\noindent o\`u $\Delta$ est le morphisme diagonal et $u_{\alpha}$ la projection ($u_{\alpha}$ est lisse), fournit un triangle commutatif\\

$$
\begin{array}{c}
\xymatrix{
&   \mathcal{T}_{\alpha} \ar[d]^{u_{\alpha}} \\
S_{\alpha} \ar@{^{(}->}[ur]^{i'} \ar@{^{(}->}[r]_{i} & \mathcal{S}_{\alpha}\\
}
\end{array}
\leqno{(2.1.24)}
$$

\noindent o\`u $i$ et $i'$ sont des immersions ferm\'ees.\\
D'apr\`es [B 3, (2.3.1)] le foncteur $u^{\ast}_{\alpha_{K}}$ induit une auto-\'equivalence de la cat\'egorie $F^a\mbox{-}\textrm{Isoc}(S_{\alpha}/K)$ : en composant un foncteur quasi-inverse canonique \`a $u^{\ast}_{K}$ (cf. [B 5, (3.1.10)]) avec l'\'equivalence de cat\'egories du [III, cor (1.2.3)] on constate que la donn\'ee de l'isomorphisme $\phi_{\mathcal{E}_{\alpha_{i}}}$ correspond \`a la donn\'ee d'un isomorphisme 

$$
\phi_{\mathcal{M}} : \mathcal{M}^{\sigma} \tilde{\longrightarrow}\  \mathcal{M}
$$

\noindent de $\hat{A}_{K}$-modules projectifs de type fini commutant aux connexions.\\

De m\^eme, d'apr\`es [B 3, (2.3.5)] le morphisme propre $v$ de (2.1.23) induit l'auto-\'equivalence $v^{\ast} = \overline{u}^{\ast}_{\alpha_{k}}$ de Isoc$^{\dag}(S_{\alpha}/K)$ : en composant un foncteur quasi-inverse \`a $v^{\ast}$ avec l'\'equivalence de cat\'egories de [B 3, (2;5;2) (ii)], on constate que la donn\'ee du morphisme $\phi_{E_{\alpha_{i}}}$ correspond \`a la donn\'ee d'un morphisme

$$
\phi_M : M^{\sigma} \longrightarrow M
$$

\noindent de $A^{\dag}_{K}$-modules projectifs de type fini commutant aux  connexions.\\

Ce qui pr\'ec\`ede peut \^etre formalis\'e par un diagramme commutatif de foncteurs entre cat\'egories\\

$$
\begin{array}{c}
\xymatrix{
\textrm{Conn}^{\dag} (A^{\dag}_{K}) \ar[d]^{\mathcal{G}} && \textrm{Isoc}^{\dag} (S_{\alpha}/K) \ar[ll]_{\Gamma(\overline{\mathcal{S}}_{\alpha_{K},-})} ^{\simeq} \ar[rr]^{\overline{u}^{\ast}_{\alpha_K}}_{\simeq} \ar[d]^{j^{\ast}_{\overline{\mathcal{S}}_{\alpha K}}} && \textrm{Isoc}^{\dag}(S_{\alpha}/K) \ar[d]^{j^{\ast}_{\overline{\mathcal{T}}_{\alpha K}}} \\
\textrm{Conn}^{\wedge} (\hat{A}_{K}) && \textrm{Isoc}(S_{\alpha}/K) \ar[ll]_{\Gamma(\mathcal{S}_{\alpha_{K},-})} ^{\simeq}   \ar[rr]^{u^{\ast}_{\alpha_K}}_{\simeq}  && \textrm{Isoc}(S_{\alpha}/K)
}
\end{array}
\leqno{(2.1.25)}
$$

\noindent o\`u les fl\`eches verticales sont les foncteurs d'oubli et les fl\`eches horizontales des \'equivalences de cat\'egories : pour les notations et r\'esultats cf. [III, (1.1) et prop. (1.2.1)]. Ainsi on a un carr\'e commutatif 

$$
\xymatrix{
\mathcal{G}(\phi_{M}) : \mathcal{G}(M)^{\sigma} \ar[r] \ar[d]^{\simeq} & \mathcal{G}(M) \ar[d]^{\simeq}\\
\phi_{\mathcal{M}} : \mathcal{M}^{\sigma} \ar[r]^{\sim} & \mathcal{M}\ ,
 }
$$

\noindent  donc $\mathcal{G}(\phi_{M})$ est un isomorphisme : par fid\`ele platitude de $\hat{A}_{K}$ sur $A^{\dag}_{K}$ on en d\'eduit que $\phi_{M}$ est un isomorphisme. Par suite $\phi_{E_{\alpha_{i}}}$ est un isomorphisme. $\square$

 \noindent \textbf{Remarque (2.2)}. 
\begin{itemize}
\item[(i)] En fait on a prouv\'e, plus pr\'ecis\'ement, que le morphisme (2.1.17) induit par $F_{X_{\alpha}/S_{\alpha}}$ est un isomorphisme.
\item[(ii)] Pour construire le morphisme de Frobenius $\phi_{E_{i}}$ nous n'avons pas suppos\'e l'existence d'un rel\`evement \`a $\overline{\mathcal{S}}$ du Frobenius de $\overline{S}$, contrairement \`a [C-T, 12.2]. 
\end{itemize}

\vskip 3mm
\section*{3. Cas propre et lisse}

\noindent \textbf{Th\'eor\`eme (3.1)}. \textit{Soient $S$ un $k$-sch\'ema lisse et s\'epar\'e et $f : X \rightarrow S$ un $k$-morphisme projectif et lisse satisfaisant aux propri\'et\'es de [II, (3.4.8.2)] ou [II, (3.4.8.6)] ou [II, (3.4.9)]. Alors} \\

\noindent (3.1.1) \textit{Pour tout entier $i \geqslant 0$, on a un diagramme commutatif de foncteurs naturels induits par $f$ et d\'efinis en [II, (3.4.8.5)]}
$$
\xymatrix{
 F^{a}\mbox{-}Isoc^{\dag}(X/ K) \ar[rr]^{R^{i}f_{rig\ast}}\ar[d]&&F^{a}\mbox{-}Isoc^{\dag}(S/K)\ar [d]\\
F^{a}\mbox{-}Isoc (X/ K)\ar[rr]^{R^{i}f_{conv\ast}}&&F^{a}\mbox{-}Isoc(S/ K)
}
$$
\textit{o\`u les fl\`eches verticales sont les foncteurs d'oubli.}\\

\noindent (3.1.2) \textit{Le foncteur $R^i f_{\textrm{rig}^{\ast}}$ pr\'ec\'edent est compatible aux changements de base entre $k$-sch\'emas lisses et s\'epar\'es (en particulier il commute aux passages aux fibres en les points ferm\'es de $S$), c'est-\`a-dire : pour tout carr\'e cart\'esien}

$$
\xymatrix{
X' \ar[r]^{g'} \ar[d]_{f'} & X \ar[d]^f\\
S' \ar[r]_{g} & S
}
$$

\noindent \textit{o\`u $S'$ est un $k$-sch\'ema lisse et s\'epar\'e et $E  \in F^a\mbox{-}\textrm{Isoc}^{\dag}(X/K)$ on a un isomorphisme de changement de base}

$$g^{\ast} R^i f_{\textrm{rig}^{\ast}}(E) \tilde{\longrightarrow} R^i f'_{\textrm{rig}^{\ast}}(g'^{\ast} (E)) $$

\noindent \textit{compatible aux connexions et aux Frobenius.}

\vskip 3mm
\noindent \textit{D\'emonstration}. \\
\noindent \textit{Pour (3.1.1)}. Vu la d\'efinition locale sur $S$ de $R^i f_{\textrm{rig}^{\ast}}$ [cf [II, (3.4.8.5)]) on peut supposer $S$ affine lisse et connexe et se ramener au cas de [II, (3.4.8.2)]: alors $f$ est relevable comme dans le th\'eor\`eme (2.1) ci-dessus qu'il suffit d'appliquer, d'o\`u la conclusion.\\

\noindent \textit{Pour (3.1.2)}. Soient $V = Spec\ A_{0} \displaystyle \mathop{\hookrightarrow}^{j_{U_{0}}}S$ un ouvert affine de $S$ et $Y = Spec\ B_{0} \displaystyle \mathop{\hookrightarrow}_{j'_{Y_{0}}}S'$ un ouvert affine de $V' = S' \times_{S} V$ ; on note $\psi_{VY} : Y \rightarrow V$, le morphisme induit par $g$. Soient $A = \mathcal{V}[t_{1},...,t_{n}]/(f_{1},...f_{r})$ une $\mathcal{V}$-alg\`ebre lisse relevant $A_{0}$, $U = Spec\ A$ et $\overline{U}$ la fermeture projective de $U$ dans $\mathbb{P}^n_{\mathcal{V}}$, $\mathcal{S}$ et $\overline{\mathcal{S}}$ leurs compl\'et\'es formels respectifs, $j_{\overline{\mathcal{S}}}$ l'immersion ouverte $\mathcal{S} \hookrightarrow \overline{\mathcal{S}}$, et posons $X_{V} = X \times_{S} V$, et $\overline{V} = \overline{U}\  \textrm{mod} . \pi$. Si $V$ (resp. $Y$) parcourt un recouvrement ouvert affine de $S$ (resp. de $V'$) alors $Y$ parcourt un recouvrement ouvert affine de $S'$ ; or la donn\'ee de $R^i f_{\textrm{rig}^{\ast}}(E)$ (resp. de $g^{\ast} R^i f_{\textrm{rig}^{\ast}}(E))$ \'equivaut \`a la donn\'ee des $j^{\ast}_{V} R^i f_{\textrm{rig}^{\ast}}(E)$ (resp. des $j_{Y}^{\prime\ast} g^{\ast} R^i f_{\textrm{rig}^{\ast}}(E))$, donc la donn\'ee de $g^{\ast} R^i f_{\textrm{rig}^{\ast}}(E)$ \'equivaut \`a celle des

$$
\psi^{\ast}_{VY}\  j^{\ast}_{V}\  R^i f_{\textrm{rig}^{\ast}}(E) = \psi^{\ast}_{VY}\  R^i f_{\textrm{rig}^{\ast}}(X_{V} / \overline{\mathcal{S}}, E_{X_{V}}).
$$
Puisque $\psi  := \psi_{VY}$ est de type fini on peut choisir une pr\'esentation $B_{0} = A_{0}[x_{1},...,x_{d}]/(g_{1},...,g_{s})$ de $B_{0}$ sur $A_{0}$ : notons $\mathcal{Y}$ (resp. $\overline{\mathcal{Y}}$) le compl\'et\'e formel de $\mathbb{A}^d_{U}$ (resp. de $\mathbb{P}^d_{\overline{U}})$, $\overline{Y}$ l'adh\'erence sch\'ematique de $Y$ dans $\mathbb{P}^d_{\overline{U}}$, $\overline{\psi} : \overline{Y} \rightarrow \overline{V}$ le morphisme canonique et $j_{\overline{\mathcal{Y}}} : \mathcal{Y} \hookrightarrow \overline{\mathcal{Y}}$, $j_{\overline{Y}} : Y \hookrightarrow \overline{Y}$ les immersions ouvertes ; $\theta : \mathcal{Y} \rightarrow \mathcal{S}$, $\overline{\theta} : \overline{\mathcal{Y}} \rightarrow \overline{\mathcal{S}}$ d\'esignerons les projections canoniques. D'o\`u un diagramme commutatif

$$
\xymatrix{
& \overline{Y} \ar@{^{(}->}[rr] \ar@{.>}[dd]^(.4){\overline{\psi}} |\hole && \overline{\mathcal{Y}} \ar[dd]^{\overline{\theta}} \\
Y \ar@{^{(}->}[rr] \ar[dd]_{\psi} \ar@{^{(}->}[ur] && \mathcal{Y} \ar[dd] \ar@{^{(}->}[ur]_{j_{\overline{\mathcal{Y}}}}\\
&\overline{V} \ar@{.>}[rr]^(.40){\theta}  && \overline{\mathcal{S}}  \\
V \ar@{^{(}->}[rr] \ar@{^{(}.>}[ur] &&\  \mathcal{S}\  .\ar@{^{(}->}[ur]_{j_{\overline{\mathcal{S}}}} \\
}
$$
\\

Ainsi [B 3, (2.3.2) (iv)]

$$
\psi^{\ast} R^i f_{\textrm{rig}^{\ast}}(X_{V}/ \overline{\mathcal{S}}, E_{X_{V}}) = \overline{\theta}^{\ast} R^i f_{\textrm{rig}^{\ast}} (X_{V}/ \overline{\mathcal{S}}, E_{X_{V}})
$$

\noindent et \quad $j^{\ast}_{\overline{\mathcal{Y}}}\  \overline{\theta}^{\ast} R^i f_{\textrm{rig}^{\ast}}(X_{V} / \overline{\mathcal{S}}, E_{X_{V}})$\\

$\quad =  \theta^{\ast} j^{\ast}_{\overline{\mathcal{S}}}\ R^i f_{\textrm{rig}^{\ast}}(X_{V} / \overline{\mathcal{S}}, E_{X_{V}})$\\

$\quad = \theta^{\ast} R^i f_{\textrm{conv}^{\ast}}(X_{V} / S, \hat{E}_{X_{V}})\ [\mbox{II}, (3.4.4)]$\\

$\quad = \psi^{\ast} R^i f_{\textrm{conv}^{\ast}}(X_{V} / S, \hat{E}_{X_{V}}) = \psi^{\ast} (R^i f_{\textrm{conv}^{\ast}} (\hat{E})_{\vert_{V}})$\\

$\qquad \qquad \qquad \qquad \qquad \qquad \quad = (g^{\ast} (R^i f_{\textrm{conv}^{\ast}} (\hat{E}))_{\vert_{Y}}$\\

\noindent o\`u $\hat{E}$ est l'isocristal convergent associ\'e \`a $E$.\\

De m\^eme on a\\

$j^{\ast}_{\overline{\mathcal{Y}}}\ j'^{\ast}_{Y}\  R^i f'_{\textrm{rig}^{\ast}} (g'^{\ast} (E)) \simeq R^i f'_{\textrm{conv}^{\ast}} (X'_{Y}/\mathcal{Y}, g'^{\ast}(\hat{E})_{X'_{Y}})$\\

$\qquad \qquad \qquad \qquad \qquad \simeq (R^i f'_{\textrm{conv}^{\ast}} (g'^{\ast}(\hat{E})))_{\vert_{Y}}$\ . \\

Or le th\'eor\`eme [III, (3.2.1)] fournit un isomorphisme de changement de base

$$
g^{\ast} R^i f_{\textrm{conv}^{\ast}}(\hat{E})\  \tilde{\rightarrow}\ R^i f'_{\textrm{conv}^{\ast}} (g'^{\ast}(\hat{E}))\  ; 
$$

\noindent donc par le th\'eor\`eme de pleine fid\'elit\'e pour les $F$-isocristaux de Kedlaya [Ked 2, theo 1.1] on en d\'eduit un isomorphisme

$$
g^{\ast} R^i f_{\textrm{rig}^{\ast}}(E)\  \tilde{\rightarrow}\  R^i f'_{\textrm{rig}^{\ast}} g'^{\ast}(E)\ 
$$

\noindent compatible aux connexions et aux Frobenius. $\square$

\vskip 3mm
\noindent \textbf{Th\'eor\`eme (3.2)}. \textit{Soient $S = \displaystyle \mathop{\coprod}_{\alpha = 1}^n S_{\alpha}$ un $k$-sch\'ema lisse et s\'epar\'e, d\'ecompos\'e en la somme de ses composantes connexes, et $f : X \rightarrow S$ un $k$-morphisme propre et lisse v\'erifiant la propri\'et\'e $\mathcal{P}$ suivante :}\\

 \noindent $\mathcal{P}$
 $\left\lbrace
 \begin{array}{l}
 Il\  existe\  un\  ouvert\  dense\  U \subset X quasi\mbox{-}projectif sur\  S\ 
tel\  que\\ 
pour\  tout\  \alpha \in \mathbb{[[}1,n \mathbb{]]}\   on\  ait\  S_{\alpha} \backslash f(X \backslash U) \neq \phi . 
\end{array}
 \right .$\\

 \noindent \textit{Alors}\\
\begin{itemize} 
\item[{(3.2.1)}] \textit{La propri\'et\'e ($\mathcal{P})$ \'equivaut \`a dire que $f$ est g\'en\'eriquement projective, i.e. qu'il existe un ouvert dense $V \subset S$ tel que l'application $f_{V} : X_{V} = X_{X_{S}} V \rightarrow V$ induite par $f$ soit projective et lisse.\\
}

\item[{(3.2.2)}]\textit{Supposons $k$ parfait, $e\leqslant p-1$ et que le $f_{V}$ de (3.2.1) satisfait aux hypoth\`eses de [II, (3.4.8.2)] ou [II,(3.4.8.6)] ou [II, (3.4.9)]}. $ Si\  E \in F^a\mbox{-}\textrm{Isoc}^{\dag}(X/K)_{\textrm{plat}} \ a\ pour\ image\ \mathcal{E} \in F^a\mbox{-}\textrm{Isoc}(X/K), alors :$

	\begin{itemize}
	\item[{(i)}] $\mathcal{E}^i = R^i   f_{\textrm{conv}^{\ast}}(\mathcal{E}) \in F^a\mbox{-}\textrm{Isoc}(S/K)$,

	\item[{(ii)}] \textit{Il existe $E^i \in F^a\mbox{-}\textrm{Isoc}^{\dag}(S/K)$, unique \`a isomorphisme pr\`es, tel que $\mathcal{E}^i$ soit l'image de $E^i$ par le foncteur d'oubli}

	$$F^a\mbox{-}\textrm{Isoc}^{\dag}(S/K) \longrightarrow F^a\mbox{-}\textrm{Isoc}(S/K) $$
	$$E^i \longmapsto \widehat{E^i} = \mathcal{E}^i . $$
	\end{itemize}
\end{itemize}

\vskip 3mm
\noindent \textit{D\'emonstration}. \\
\noindent \textit{Prouvons (3.2.1)}. Si $f$ est g\'en\'eriquement projective et lisse on prouve facilement qu'elle v\'erifie $(\mathcal{P})$.\\

R\'eciproquement supposons que $f$ v\'erifie $(\mathcal{P})$.\\

Par le lemme de Chow pr\'ecis de Gruson-Raynaud [R-G, I, cor 5.7.14] il existe un \'eclatement $U$-admissible $g : X' \rightarrow X$, avec $X'$ quasi-projectif sur $S$ : en particulier $g$ induit un isomorphisme

$$
g_{U} : U' = g^{-1}(U)\   \tilde{\longrightarrow}\  U\ .
$$

\noindent De plus, comme $f$ et $g$ sont propres, le morphisme compos\'e $f \circ g : X' \rightarrow X$ est projectif [EGA II, (5.5.3) (ii)]. L'image du ferm\'e $Z := X \backslash U$ par le morphisme propre $f$ est un ferm\'e $f(Z)$ de $S$ et l'ouvert $V = S \backslash f(Z) = \displaystyle \mathop{\coprod}_{\alpha} (S_{\alpha} \backslash f(X-U))$ est non vide par hypoth\`ese : comme $S_{\alpha}$ est connexe et int\`egre, l'ouvert non vide $V_{\alpha} := S_{\alpha} \backslash f(X  \backslash U)$ de $S_{\alpha}$ est dense, donc l'immersion ouverte $j : V \hookrightarrow S$ est dominante. D'autre part l'ouvert $X_{V} = X \times_{S} V$ de $X$ ne rencontre pas $f^{-1}(f(Z))$, donc $X_{V}$ est un ouvert de $U$ : par suite l'isomorphisme $g_{U}$ induit un isomorphisme 

$$
g_{V} : X'_{U} = g^{-1}(X_{V}) \tilde{\longrightarrow} X_{V} .
$$

\noindent Notons $f_{V} : X_{V} \rightarrow V$ la restriction de $f$ ; le morphisme compos\'e $f_{V} \circ g_{V}$, restriction du morphisme projectif $f \circ g$, est lui aussi projectif : ainsi $f_{V}$ est projectif, d'o\`u (3.2.1).\\

\noindent \textit{Prouvons le (3.2.2)}.  Le (i) est mis pour m\'emoire, car prouv\'e en [III, (3.3.1)]. Pour le (ii) consid\'erons le carr\'e cart\'esien

$$
\xymatrix{
X_{V} \ar@{^{(}->}[r]^{j'} \ar[d]_{f_{V}} & X \ar[d]^{f}\\
V \ar@{^{(}->}[r]^{j} & S   ; 
}
$$

\noindent on a un isomorphisme de changement de base en cohomologie convergente [III, (3.3.1)]\\

\noindent (3.2.2.1) $\qquad \qquad  j^{\ast} R^i   f_{\textrm{conv}^{\ast}}(\mathcal{E}) \tilde{\longrightarrow} R^i   f_{V{\textrm{conv}}^{\ast}}(j'^{\ast} (\mathcal{E})) =: \mathcal{E}^i_{V}, $\\

\noindent o\`u $\mathcal{E}$ d\'esigne l'image de $E$ par le foncteur d'oubli $F^a\mbox{-}\textrm{Isoc}^{\dag}(X/K) \rightarrow F^a\mbox{-}\textrm{Isoc}(X/K)$, $E \mapsto \hat{E} = \mathcal{E}$.\\
Pour la suite de la d\'emonstration on peut supposer $S$ connexe, int\`egre : quitte \`a restreindre $V$ on peut supposer $V$ affine, lisse et connexe, $V = Spec\ A_{0}$. On utilise alors les notations introduites dans la d\'emonstration du th\'eor\`eme (3.1) :  $j_{\overline{\mathcal{S}}} : \mathcal{S} = Spf\ \hat{A} \hookrightarrow \overline{\mathcal{S}}$.
De plus $f_{V}$ se rel\`eve en un morphisme projectif et lisse  $h : \mathcal{X} \rightarrow \mathcal{S}$ s'ins\'erant dans un carr\'e cart\'esien de $\mathcal{V}$-sch\'emas formels
 
$$
\xymatrix{
\mathcal{X} \ar@{^{(}->}[r] \ar[d]_{h} & X \ar[d]^{\overline{h}}\\
\mathcal{S} \ar@{^{(}->}[r]_{j_{\overline{\mathcal{S}}}} & \overline{\mathcal{S}}  
}
$$
\noindent o\`u $\overline{h}$ est projectif [I, th\'eo (3;3)].
En notant $E^i_{V} = R^i  f_{V{\textrm{rig}}^{\ast}}(X_{V}/\overline{\mathcal{S}}, j'^{\ast}(E))$, le th\'eor\`eme (3.1) prouve que $E^i_{V} \in F^a\mbox{-}\textrm{Isoc}^{\dag}(V/K)$.
On peut appliquer le [II, (3.4.4)] qui fournit un isomorphisme\\

\noindent (3.2.2.2) $\qquad \qquad \qquad \qquad \qquad \widehat{E^i_{V}}\  \tilde{\longrightarrow}\  \mathcal{E}^i_{V}$\\

\noindent compatible aux Frobenius. Par le th\'eor\`eme 2 de [Et 5], les isomorphismes (3.2.2.1) et (3.2.2.2) assurent l'existence de $E^i \in F^a\mbox{-}\textrm{Isoc}^{\dag}(S/K)$ tel que

$$
\mathcal{E}^i = \widehat{E^i} \qquad \textrm{et}   \qquad E^i_{V} = j^{\ast}(E^i).
$$

\noindent L'unicit\'e de $E^i$ \`a isomorphisme pr\`es provient de la pleine fid\'elit\'e du foncteur d'oubli $F^a\mbox{-}\textrm{Isoc}^{\dag}(S/K) \rightarrow F^a\mbox{-}\textrm{Isoc}(S/K)$ \'etabli par Kedlaya [Ked 2, theo 1.1]. D'o\`u le th\'eor\`eme. $\square$\\

\vskip 3mm
\section*{4. Cas fini \'etale}

\noindent \textbf{Th\'eor\`eme (4.1)}. \textit{Soient $S$ un $k$-sch\'ema lisse et s\'epar\'e et $f : X \rightarrow S$ un $k$-morphisme fini \'etale. Alors, pour tout entier $i \geqslant 0$, $f$ induit des foncteurs canoniques }

$$
R^i  f_{\textrm{rig}\ast} : \textrm{Isoc}^{\dag}(X/K) \longrightarrow \textrm{Isoc}^{\dag}(S/K)
$$

$$
R^i  f_{\textrm{rig}\ast} : F^a\mbox{-}\textrm{Isoc}^{\dag}(X/K) \longrightarrow F^a\mbox{-}\textrm{Isoc}^{\dag}(S/K)
$$

\vskip 3 mm
\noindent \textit{et $\qquad \qquad  \qquad \qquad R^i  f_{\textrm{rig} \ast}(E)  = 0$ pour tout $i\geqslant 1$.} 

\vskip 3mm
\noindent \textit{D\'emonstration}. Soient $S_{0} = Spec\ A_{0} \hookrightarrow S$ un ouvert affine et $A_{1}$, $A_{2}$ deux $\mathcal{V}$-alg\`ebres lisses relevant $A_{0}$. On pose $S_{1} = Spec\ A_{1}$, $S_{2} = Spec\ A_{2}$ ; par la m\'ethode du [I, th\'eo (3.1)] on a des compactifications $\overline{S}_{1} := P_{1}$, $\overline{S}_{2} := P_{2}$ de $S_{1}$ et $S_{2}$ et on note $\overline{S}_{0}$ l'adh\'erence sch\'ematique de $S_{0}$ plong\'e diagonalement dans $\overline{S}_{1} \times_{\mathcal{V}} \overline{S}_{2}$. En d\'esignant par $f_{0}$ la restriction de $f$ \`a $X_{0} = f^{-1}(S_{0})$ et par $\mathcal{S}_{1}$,  $\mathcal{S}_{2}$, $\overline{\mathcal{S}_{1}}$, $\overline{\mathcal{S}_{2}}$ les compl\'et\'es formels de $S_{1}$, $S_{2}$, $\overline{S}_{1}$, $\overline{S}_{2}$  respectivement, le th\'eor\`eme (3.1) du I fournit des carr\'es cart\'esiens, $i = 1,2$,

$$
\xymatrix{
\mathcal{X}_{i} \ar[r] \ar[d]_{h_{i}} & \overline{\mathcal{X}_{i}} \ar[d]^{\overline{h}_{i}}\\
\mathcal{S}_{i} \ar@{^{(}->}[r]  & \overline{\mathcal{S}_{i}}
}
$$

\noindent o\`u $\overline{h}_{i}$ est fini, $h_{i}$ est fini \'etale et r\'el\`eve $f_{0}$ ; d'o\`u deux cubes commutatifs ($i = 1,2$)\\

$$
\xymatrix{
&\mathcal{X}_{i} \ar@{^{(}->}[rr] \ar@{.>}[dd]^(.4){h_{i}} |\hole && \overline{\mathcal{X}}_{i} \ar[dd]^{\overline{h}_{i}} \\
 \mathcal{X}_{1} \times_{\mathcal{V}}  \mathcal{X}_{2} \ar@{^{(}->}[rr] \ar[dd]_{h_{1} \times h_{2}} \ar[ur]^{u_{\mathcal{X}_i }} && \overline{\mathcal{X}}_{1} \times_{\mathcal{V}}  \overline{\mathcal{X}}_{2} \ar[dd] \ar[ur]_{u_{\overline{\mathcal{X}}_i }} \\
&\mathcal{S}_{i} \ar@{.>}[rr]^{\overline{h}_{1} \times \overline{h}_{2}}  && \overline{\mathcal{S}_{i}}  \\
 \mathcal{S}_{1} \times_{\mathcal{V}}  \mathcal{S}_{2} \ar[rr] \ar@{.>}[ur]_{u_{\mathcal{S}_i}} && \overline{\mathcal{S}}_{1} \times_{\mathcal{V}}  \overline{\mathcal{S}}_{2} \ar[ur]_{u_{\overline{\mathcal{S}}_i}}\ . \\
}
$$

\noindent Par le th\'eor\`eme [II, (3.4.1)] on sait que pour $E \in \textrm{Isoc}^{\dag}(X/K)$ et $E_{0}$ sa restriction \`a $X_{0}$, alors $R^i f_{0 \textrm{rig}^{\ast}} (X_{0} / \overline{\mathcal{S}_{1}}, E_{0})$ et 
$R^i f_{0 \textrm{rig}^{\ast}}(X_{0} / \overline{\mathcal{S}_{2}}, E_{0})$ sont \'el\'ements de $\textrm{Isoc}^{\dag}(S_{0}/K)$ : de plus ils sont nuls pour $i \geqslant 1$ car $\overline{h}_{1}$ et $\overline{h}_{2}$ sont finis. De plus le th\'eor\`eme [II, (3.4.4)] fournit des isomorphismes de changement de base

$$
u^{\ast}_{{\overline{\mathcal{S}}}_{i}} : f_{0 \textrm{rig} \ast} (X_{0} / \overline{\mathcal{S}_{i}}, E_{0}) \tilde{\longrightarrow} f_{0 \textrm{rig} \ast}(X_{0} / \overline{\mathcal{S}_{1}} \times_{\mathcal{V}} \overline{\mathcal{S}_{2}}, u^{\ast}_{\overline{\mathcal{X}}_{i}} E_{0}) ;
$$

\noindent d'o\`u un isomorphisme

$$
f_{0 \textrm{rig} \ast}(X_{0} / \overline{\mathcal{S}_{1}}, E_{0}) \tilde{\longrightarrow} f_{0 \textrm{rig}\ast}(X_{0} / \overline{\mathcal{S}_{2}}, E_{0}) \ ,
$$

\noindent et cet isomorphisme v\'erifie la condition de cocycles pour trois r\'el\`evements $S_{1}, S_{2}, S_{3}$ de $S_{0}$.\\

Par suite $f$ induit un foncteur

$$
f_{\textrm{rig}^{\ast}} : \textrm{Isoc}^{\dag}(X/K) \longrightarrow \textrm{Isoc}^{\dag}(S/K)
$$

\noindent puisque les constructions se recollent sur les ouverts de $S$. On pouvait aussi conclure en appliquant [II, (3.4.8)].\\

La construction du Frobenius \'etant locale, on peut, pour montrer que $f_{\textrm{rig}^{\ast}}$ induit un foncteur 

$$
f_{\textrm{rig}^{\ast}} :  F^a\mbox{-}\textrm{Isoc}^{\dag}(X/K) \longrightarrow F^a\mbox{-}\textrm{Isoc}^{\dag}(S/K),
$$

\noindent supposer que $S$ est affine et lisse. La construction du th\'eor\`eme (2.1) s'applique ; le morphisme (2.1.17) est alors un isomorphisme car $F_{X/S}$ est un isomorphisme puisque $f$ est \'etale : l\`a on n'a pas besoin d'utiliser les r\'esultats de Ogus via le cas convergent (o\`u l'on avait suppos\'e $e \leqslant p-1$). On en d\'eduit directement que $\phi_{E_{i}}$ est un isomorphisme. D'o\`u le th\'eor\`eme. $\square$

\vskip 3mm

\noindent \textbf{Remarque (4.1.1)}. Tsuzuki a abord\'e dans [Tsu 1, theo (2.6.3)] la construction de $f_{\textrm{rig}^{\ast}}(X_{0}/\overline{\mathcal{S}},-)$ dans le cas fini \'etale, mais il n'\'etudie pas l'ind\'ependance par rapport \`a $\overline{\mathcal{S}}$ et ne prouve pas l'existence d'un $\mathcal{V}$- morphisme fini relevant le $f_{0}$ ci-dessus.

\vskip 3mm
\noindent \textbf{Th\'eor\`eme (4.2)}. \textit{Soient $S$ un $k$-sch\'ema s\'epar\'e de type fini, $E \in \textrm{Isoc}^{\dag}(S/K)$ et $f : X \rightarrow S$ un $k$-morphisme fini \'etale galoisien de groupe $G$}.\\

\noindent (4.2.1) \textit{Si $S$ est lisse sur $k$, alors, pour tout entier $i \geqslant 0$, on a des isomorphismes canoniques}\\

$\qquad \qquad H^i_{\textrm{rig}}(S/K, E)\  \tilde{\longrightarrow}\  (H^i_{\textrm{rig}}(S/K, f_{\textrm{rig}\ast} f^{\ast} E))^G$\\

$\qquad \qquad \qquad \qquad \qquad \tilde{\longrightarrow}\  (H^i_{\textrm{rig}}(X/K, f^{\ast} E))^G.$\\

\noindent (4.2.2) \textit{Si $k$ est parfait, ou si $S$ est affine et lisse sur $k$, alors, pour tout entier $i \geqslant 0$, on a des isomorphismes canoniques}\\

$\qquad \qquad H^i_{\textrm{rig},c}(S/K, E)\  \tilde{\longrightarrow}\  (H^i_{\textrm{rig},c}(S/K, f_{\textrm{rig}\ast} f^{\ast} E))^G$\\

$\qquad \qquad \qquad \qquad \qquad \tilde{\longrightarrow}\  (H^i_{\textrm{rig},c}(X/K, f^{\ast} E))^G.$\\

\noindent (4.2.3) \textit{Si $E \in F^a\mbox{-}\textrm{Isoc}^{\dag}(S/K)$ alors les isomorphismes de (4.2.1) et (4.2.2) sont compatibles \`a l'action du Frobenius.}

\vskip 3mm
\noindent \textit{D\'emonstration}. Par additivit\'e de la cohomologie rigide, avec ou sans supports, on peut supposer, pour le (1) et le (2), que $S$ est connexe.\\

\noindent \textit{Pour le (4.2.1)}, la suite spectrale de $\check{\mbox{C}}$ech en cohomologie rigide nous ram\`ene \`a $S$ affine et lisse sur $k$, $S = Spec\ A_{0}$. On choisit une $\mathcal{V}$-alg\`ebre lisse $A$ relevant $A_{0}$ et on reprend les notations utilis\'ees dans la preuve de (3.2.2): il existe un carr\'e cart\'esien de $\mathcal{V}$-sch\'emas formels

$$
\xymatrix{
\mathcal{X} \ar@{^{(}->}[r] \ar[d]_{h} & \overline{\mathcal{X}} \ar[d]^{\overline{h}} \\
\mathcal{S} \ar@{^{(}->}[r] & \overline{\mathcal{S}}
}
$$

\noindent et un syst\`eme fondamental $(V_{\lambda})_{\lambda} = (Spm\ A_{\lambda})_{\lambda}$ de voisinages stricts de $]S[_{\overline{\mathcal{S}}}$ dans $\overline{\mathcal{S}_{K}}$ et $\lambda_{0} > 1$ tel que pour $1 < \lambda \leqslant \lambda_{0}$ on ait un diagramme \`a carr\'es cart\'esiens

$$ \begin{array}{c}
\xymatrix{
\mathcal{X}_{K} \ar@{^{(}->}[r] \ar[d] _{h_{K}} & W_{\lambda} \ar@{^{(}->}[r] \ar[d]_{h_{\lambda}} & \overline{\mathcal{X}}_{K} \ar[d]^{\overline{h}_{K}}\\
\mathcal{S}_{K} \ar@{^{(}->}[r] & V_{\lambda} \ar@{^{(}->}[r] & \overline{\mathcal{S}}_{K}
}
\end{array}
\leqno{(4.2.1.1)}
$$

\noindent avec $\overline{h}_{K}$ fini, $h_{K}$ et $h_{\lambda}$ finis \'etales galoisiens de groupe $G$ [II, (2.3.1)(2)].\\

\noindent (4.2.1.2) Soit $E_{\lambda}$ un $\mathcal{O}_{V_{\lambda}}$-module localement libre de type fini. Pour $1 < \mu \leqslant \lambda$ on note

$$
\alpha_{\lambda_{\mu}} : V_{\mu} \hookrightarrow V_{\lambda} \quad , \quad \alpha_{\lambda} : V_{\lambda} \hookrightarrow \overline{\mathcal{S}}_{K}
$$

$$
\alpha'_{\lambda_{\mu}} : W_{\mu} \hookrightarrow W_{\lambda} \quad , \quad \alpha'_{\lambda} : W_{\lambda} \hookrightarrow \overline{\mathcal{X}}_{K},
$$

\noindent les immersions ouvertes et on pose

$$
j^{\dag}_{\lambda}\  E_{\lambda} = \displaystyle \mathop{\lim}_{\rightarrow} \alpha_{\lambda\mu^{\ast}}\  \alpha^{\ast}_{\lambda \mu}(E_{\lambda}),
$$

$$
j'^{\dag}_{\lambda}\  h^{\ast}_{\lambda}\  E_{\lambda} = \displaystyle \mathop{\lim}_{\rightarrow} \alpha'_{\lambda\mu^{\ast}}\  \alpha'^{\ast}_{\lambda \mu}\ h^{\ast}_{\lambda}\ (E_{\lambda}),
$$

$$
j^{\dag}\  E_{\lambda} =  \alpha_{\lambda^{\ast}}\  j'^{\dag}_{\lambda}\ E_{\lambda},
$$

$$
j'^{\dag}_{\lambda}\  h^{\ast}_{\lambda} (E_{\lambda}) = \alpha'_{\lambda^{\ast}}\  j'^{\dag}_{\lambda}\  h^{\ast}_{\lambda}\ (E_{\lambda}).
$$

\vskip 3mm
\noindent \textbf{Lemme (4.2.1.3)}. \textit{Avec les notations pr\'ec\'edentes on a des isomorphismes canoniques}

\begin{itemize}
\item[(i)] $(h_{\lambda^{\ast}}\  h^{\ast}_{\lambda}(E_{\lambda}))^G\  \tilde{\longrightarrow}\  E_{\lambda}.$
\item[(ii)] $(h_{\lambda \ast}\ h^{\ast}_{\lambda}\ j^{\dag}_{\lambda}\ E_{\lambda})^G \  \tilde{\longrightarrow}\  j^{\dag}_{\lambda}\ E_{\lambda}.$
\item[(iii)] $(\overline{h}_{K \ast}\ \overline{h}^{\ast}_{K}\ j^{\dag}\ E_{\lambda})^G\  \tilde{\longrightarrow}\  j^{\dag}\ E_{\lambda}.$
\end{itemize}

\vskip 3mm
\noindent \textit{D\'emonstration du lemme (4.2.1.3)}. 

\begin{itemize}
\item[(i)] Comme $E_{\lambda}$ est localement libre de type fini on a un isomorphisme
$$
h_{\lambda \ast}\ h^{\ast}_{\lambda}(E_{\lambda})\ \tilde{\longrightarrow}\  h_{\lambda \ast}\ h^{\ast}_{\lambda}(\mathcal{O}_{V_{\lambda}}) \otimes_{\mathcal{O}_{V_{\lambda}}} E_{\lambda},
$$
et l'action de $G$ sur le membre de gauche se fait par l'interm\'ediaire de $h_{\lambda \ast}\ h^{\ast}_{\lambda}(\mathcal{O}_{V_{\lambda}})$ puisque $G$ agit trivialement sur $E_{\lambda}$ : on est ramen\'e au cas $E = \mathcal{O}_{V_{\lambda}}$ qui a \'et\'e prouv\'e dans la proposition [II, (2.3.1)].
\item[(ii)] On a des isomorphismes\\
$$
\begin{array}{r@{\ \, }c@{\ \, }l}
h_{\lambda \ast}\ h^{\ast}_{\lambda}\ j^{\dag}_{\lambda}\ E_{\lambda} & \simeq & h_{\lambda \ast}\ j'^{\dag}_{\lambda}\ h^{\ast}_{\lambda}\ E_{\lambda}\  [B 3, (2.1.4.7)]\\
&\simeq & j^{\dag}_{\lambda}\ h_{\lambda \ast}\ h^{\ast}_{\lambda}\ E_{\lambda}\ [\mbox{II}, (3.1.4.1)]\\
&\simeq & h_{\lambda \ast}\ h^{\ast}_{\lambda}\ E_{\lambda}\ \otimes_{\mathcal{O}_{V_{\lambda}}} j^{\dag}_{\lambda}\ \mathcal{O}_{V_{\lambda}}\ [B 3, (2.1.3) (ii)].
\end{array}
$$
L'action de $G$ sur $h_{\lambda \ast}\ h^{\ast}_{\lambda}\ j^{\dag}_{\lambda}\ E_{\lambda}$ se fait par l'interm\'ediaire de $h_{\lambda \ast}\ h^{\ast}_{\lambda}\  E_{\lambda}$ puisque $G$ agit trivialement sur $j^{\dag}_{\lambda}(\mathcal{O}_{V_{\lambda}})$ : le (ii) r\'esulte alors du (i).
\item[(iii)] La preuve est semblable \`a celle du (ii) en utilisant cette fois [B 3, (2.1.4.8)] et [II (3.1.4.2)]. D'o\`u le lemme. $\square$
\end{itemize}

\vskip 3mm
Soit $E \in \textrm{Isoc}^{\dag}(S/K)$: on choisit le $V_{\lambda}$ comme ci-dessus de sorte qu'il existe un $\mathcal{O}_{V_{\lambda}}$-module localement libre et coh\'erent $E_{\lambda}$ tel que $j^{\dag} E_{\lambda}$ soit une r\'ealisation de $E$.\\

La cohomologie rigide $H^{\ast}_{\textrm{rig}}(S/K;E)$ est, pour $1 < \mu \leqslant \lambda$, la cohomologie des complexes

$$
\begin{array}{c}
\xymatrix{
 R \Gamma(V_{\mu};\  j^{\dag}_{\mu}  E_{\mu} \otimes_{\mathcal{O}_{V_{\mu}}} \Omega^{\bullet}_{V_{\mu}/K})\   \tilde{\leftarrow}\  R \Gamma (V_{\lambda};\  j^{\dag}_{\lambda}  E_{\lambda} \otimes_{\mathcal{O}_{V_{\lambda}}} \Omega^{\bullet}_{V_{\lambda}/K})\\
\qquad  \qquad \qquad \tilde{\rightarrow}\  M \otimes_{A^{\dag}_{K}} \Omega^{\bullet}_{A^{\dag}_{K}}\ , \\
}
\end{array} \leqno{(4.2.1.4)}
$$

\noindent o\`u la premi\`ere fl\`eche (resp. la deuxi\`eme) est un isomorphisme (resp. un quasi-isomorphisme), o\`u $M := \Gamma(V_{\lambda}; j^{\dag}_{\lambda}\ E_{\lambda})$ est un $A^{\dag}_{K}$-module projectif de type fini \`a connexion int\'egrable [B 3, (2.5.2)], o\`u $\Omega^1_{V_{\lambda}/K}$ est localement libre de type fini sur le faisceau coh\'erent d'anneaux $\mathcal{O}_{V_{\lambda}}$ [II (2.3.1) (2)] et o\`u $\Omega^1_{A^{\dag}_{K}}$ est un $A^{\dag}_{K}$-module projectif de type fini [Et 5, 1.3].\\

De m\^eme la cohomologie rigide

$$
H^{\ast}_{\textrm{rig}}(S/K ; f_{\textrm{rig}\ast}\  f^{\ast}\ E)\quad\quad  (\textrm{resp.} H^{\ast}_{\textrm{rig}}(X/K ;  f^{\ast}\ E))
$$

\noindent est la cohomologie des complexes\\

\noindent (4.2.1.5)  $ \qquad \qquad   R\Gamma(V_{\lambda} ; h_{\lambda^{\ast}}\ h^{\ast}_{\lambda} (j^{\dag}_{\lambda}  E_{\lambda}) \otimes_{\mathcal{O}_{V_{\lambda}}} \Omega^{\bullet}_{V_{\lambda}/K})$\\

\noindent [resp. des complexes\\

\noindent (4.2.1.6)  $ \qquad \qquad   R\Gamma(W_{\lambda} ;  h^{\ast}_{\lambda} (j^{\dag}_{\lambda}  E_{\lambda}) \otimes_{\mathcal{O}_{W_{\lambda}}} \Omega^{\bullet}_{W_{\lambda}/K})\ ].$\\

\noindent Or la formule de projection, jointe au fait que $h_{\lambda}$ est \'etale, fournit des isomorphismes

$$
\begin{array}{c}
\xymatrix{
h_{\lambda^{\ast}} h^{\ast}_{\lambda} (j^{\dag}_{\lambda} E_{\lambda}) \otimes_{\mathcal{O}_{V_{\lambda}}} \Omega^{\bullet}_{V_{\lambda}/K} \simeq h_{\lambda^{\ast}} (h^{\ast}_{\lambda}\  j^{\dag}_{\lambda} E_{\lambda} \otimes_{\mathcal{O}_{W_{\lambda}}} h^{\ast}_{\lambda} (\Omega^{\bullet}_{V_{\lambda}/K})) \\
\qquad \qquad \qquad \qquad  \qquad \simeq\  h_{\lambda^{\ast}}(h^{\ast}_{\lambda}\  j^{\dag}_{\lambda}  E_{\lambda} \otimes_{\mathcal{O}_{W_{\lambda}}} \Omega^{\bullet}_{W_{\lambda}/K})\ ;\\
}
\end{array} \leqno{(4.2.1.7)}
$$

\noindent donc les complexes (4.2.1.5) et (4.2.1.6) sont quasi-isomorphes puisque $h_{\lambda}$ est fini.\\

Compte tenu du lemme (4.2.1.3) les isomorphismes

$$
H^i_{\textrm{rig}}(S/K;E)\  \tilde{\rightarrow}\  H^i_{\textrm{rig}}(S/K; f_{\textrm{rig}^{\ast}} f^{\ast} E)^G
$$

$\qquad \qquad \qquad \qquad  \qquad \qquad \quad \tilde{\rightarrow}\  H^i_{\textrm{rig}}(X/K; f^{\ast}E)^G$\\

\noindent s'\'etablissent comme [Et 2, (3.1.1)].\\

\noindent \textit{Pour le (4.2.2)}, comme la cohomologie rigide ne d\'epend que du sch\'ema r\'eduit sous-jacent, on supposera $S$ r\'eduit : si $k$ est parfait il existe alors un ouvert dense $U \hookrightarrow S$ avec $U$ affine et lisse sur $k$ et $Z := S \setminus U$ de dimension strictement plus petite que celle de $S$ [Et 3, d\'em. du th\'eo 3]. De plus $H^j(G, H^i _{\textrm{rig},c}(S/K, f_{\textrm{rig}\ast}\ f^\ast E)) = 0$ pour $j \geqslant 1$ [S\ 2, VIII, \S\ 2, cor 1 de prop 4] ; par fonctorialit\'e en $E$ de la cohomologie rigide \`a supports on en d\'eduit un morphisme de suites exactes longues

$$
\xymatrix{
\ar[r] & H^i_{\textrm{rig},c}(U, E_{\vert U}) \ar[d] \ar[r] & H^i_{\textrm{rig},c}(S,E) \ar[d] \ar[r] & H^i_{\textrm{rig},c} (Z, E_{\vert Z}) \ar[d] \ar[r] & \\
 \ar[r] & (H^i_{\textrm{rig},c} (U, f_{U\textrm{rig}\ast}f_{U}^{\ast}E_{\vert U}))^G \ar[r] & (H^i _{\textrm{rig},c} (S, f_{\textrm{rig}\ast}f^{\ast}E))^G \ar[r]  & (H^i_{\textrm{rig},c} (Z, f_{Z\textrm{rig}\ast}f_{Z}^{\ast}E_{\vert Z}))^G   \ar[r]  & 
}
$$
\\

\noindent Par r\'ecurrence sur la dimension on est ramen\'e \`a montrer l'isomorphisme du (4.2.2) pour $S$ affine et lisse sur $k$.\\

Reprenons les notations utilis\'ees pour la d\'emonstration du (4.2.1) et consid\'erons le diagramme commutatif \`a carr\'es cart\'esiens

$$
\begin{array}{c}
\xymatrix{
\mathcal{X}_{K} \ar@{^{(}->}[r] \ar[d]_{h_{K}} & W_{\lambda}  \ar[d]_{h_{\lambda}} & W_{\lambda} \setminus \mathcal{X}_{K} \ar[d]^{h'_{\lambda}} \ar@{_{(}->}[l]_{i'_{\lambda}} \\
\mathcal{S}_{K} \ar@{^{(}->}[r] & V_{\lambda}  & V_{\lambda} \setminus \mathcal{S}_{K} . \ar@{_{(}->}[l]_{i_{\lambda}} 
}
\end{array}
\leqno{(4.2.2.1)}
$$

\noindent La cohomologie \`a supports $H^{\ast}_{\textrm{rig},c}(S/K;E)$

$$
[\textrm{resp.} H^{\ast}_{\textrm{rig},c}(S/K; f_{\textrm{rig} \ast}\ f^{\ast}\ E), \textrm{resp.} H^{\ast}_{\textrm{rig},c}(X/K; f^{\ast}\ E)]
$$\\
\noindent est la cohomologie du complexe\\

\noindent (4.2.2.2) $\qquad  R\Gamma(V_{\lambda} ; j^{\dag}_{\lambda}  E_{\lambda} \otimes_{\mathcal{O}_{V_{\lambda}}} \Omega^{\bullet}_{V_{\lambda}/K} \longrightarrow i_{\lambda \ast} i^{\ast}_{\lambda}(j^{\dag}_{\lambda} E_{\lambda} \otimes_{\mathcal{O}_{V_{\lambda}}} \Omega^{\bullet}_{V_{\lambda}/K}))$\\

\noindent [resp.\\

\noindent (4.2.2.3) $R\Gamma(V_{\lambda} ; h_{\lambda \ast} h^{\ast}_{\lambda}(j^{\dag}_{\lambda}  E_{\lambda}) \otimes_{\mathcal{O}_{V_{\lambda}}} \Omega^{\bullet}_{V_{\lambda}/K} \rightarrow i_{\lambda \ast} i^{\ast}_{\lambda}(h_{\lambda \ast} h^{\ast}_{\lambda }(j^{\dag}_{\lambda} E_{\lambda})\otimes_{\mathcal{O}_{V_{\lambda}}} \Omega^{\bullet}_{V_{\lambda}/K}));$\\

\noindent resp.\\

\noindent (4.2.2.4) $R\Gamma(W_{\lambda} ;  h^{\ast}_{\lambda}(j^{\dag}_{\lambda}  E_{\lambda} \otimes_{\mathcal{O}_{W_{\lambda}}} \Omega^{\bullet}_{W_{\lambda}/K} \rightarrow i'_{\lambda \ast} i'^{\ast}_{\lambda}(h^{\ast}_{\lambda }(j^{\dag}_{\lambda} E_{\lambda})\otimes_{\mathcal{O}_{W_{\lambda}}} \Omega^{\bullet}_{W_{\lambda}/K}))].$\\

Le th\'eor\`eme de changement de base pour un morphisme propre [II, th\'eor\`eme (3.3.2)] fournit des isomorphismes

$$
 i_{\lambda \ast} i^{\ast}_{\lambda}(h_{\lambda \ast} h^{\ast}_{\lambda }(j^{\dag}_{\lambda} E_{\lambda})\otimes_{\mathcal{O}_{V_{\lambda}}} \Omega^{\bullet}_{V_{\lambda}}) \simeq  i_{\lambda \ast} i^{\ast}_{\lambda} h_{\lambda \ast} (h^{\ast}_{\lambda }(j^{\dag}_{\lambda} E_{\lambda})\otimes_{\mathcal{O}_{W_{\lambda}}} \Omega^{\bullet}_{W_{\lambda}})
$$
$$
\simeq  i_{\lambda \ast} h'_{\lambda \ast} i'^{\ast}_{\lambda} (h^{\ast}_{\lambda} j^{\dag}_{\lambda} E_{\lambda} \otimes \Omega^{\bullet}_{W_{\lambda}}) \simeq h_{\lambda \ast} i'_{\lambda \ast} i'^{\ast}_{\lambda} (h^{\ast}_{\lambda} j^{\dag}_{\lambda} E_{\lambda} \otimes_{\mathcal{O}_{W_{\lambda}}} \Omega^{\bullet}_{W_{\lambda}})\ ;
$$

\noindent donc via (4.2.1.7) les complexes (4.2.2.3) et (4.2.2.4) sont quasi-isomorphes puisque $h_{\lambda}$ est fini. L'isomorphisme (4.2.2) du th\'eor\`eme (4.2) en r\'esulte, compte tenu de (4.2.1.3) (ii).\\

\noindent \textit{Pour le (4.2.3)}, on peut supposer $S$ connexe affine et lisse sur $\mathcal{V}$ comme ci-dessus, dont on reprend les notations ainsi que celles de [II, (2.3.1) (2)]. On fixe un rel\`evement $F_{A^{\dag}} : A^{\dag} \rightarrow A^{\dag}$
$$
\xymatrix{
\textrm{[resp.} F_{B^{\dag}} : B^{\dag} \ar[rr]_{1_{B^{\dag}}\otimes F_{A^{\dag}}}&& B^{\dag} \otimes A^{\dag}\ar[rr]^{\sim}_{F_{B^{\dag}/A^{\dag}}}&& B^{\dag}]\\
 }
 $$
 
 \noindent du Frobenius de $A_{0}$ [resp. de $B_{0}]$ comme dans [Et 5, (1.2)] : par extension des scalaires on en d\'eduit $F_{\hat{A}_{K}} : \hat{A}_{K} \rightarrow \hat{A}_{K}$ et $F_{\hat{B}_{K}} : \hat{B}_{K} \rightarrow \hat{B}_{K}$. On a vu en (1.2.4) qu'on dispose de carr\'es cart\'esiens o\`u $F_{\lambda \mu}$ et $F'_{\lambda \mu}$ sont finis :

$$
\begin{array}{c}
\xymatrix{
\mathcal{S}_{K} = Spm\  (\hat{A}_{K}) \ar@{^{(}->}[r] \ar[d]_{F_{\mathcal{S}_{K}}=Sp\ F_{\hat{A}_{K}}} 
& V_{\mu} = Spm(A_{\mu})  \ar[d]^{F_{\lambda \mu}}\\
\mathcal{S}_{K} = Spm\  (\hat{A}_{K}) \ar@{^{(}->}[r] & V_{\lambda} = Spm(A_{\lambda})
}
\end{array}
\leqno{(4.2.3.1)}
$$

\noindent et

$$
\begin{array}{c}
\xymatrix{
\mathcal{X}_{K} = Spm\  (\hat{B}_{K}) \ar@{^{(}->}[r] \ar[d]_{F_{\mathcal{X}_{K}}=Sp\ F_{\hat{B}_{K}}} 
& W_{\mu} = Spm(B_{\mu})  \ar[d]^{F'_{\lambda \mu}}\\
\mathcal{X}_{K} = Spm\  (\hat{B}_{K}) \ar@{^{(}->}[r] & V_{\mu} = Spm(B_{\lambda})\ ;
}
\end{array}
\leqno{(4.2.3.2)}
$$

\noindent plus pr\'ecis\'ement, \'etant donn\'e $\lambda$ on trouve $\mu$ de la fa\c{c}on suivante : en fixant des g\'en\'erateurs $\{ x_{i} \}$ de $B_{\lambda} $ sur $A_{\lambda}$ comme dans la preuve de [II, (2.3.1) (2)], les \'el\'ements $F_{\hat{B}_{K}}(x_{i})$ sont entiers sur $A_{\lambda} \subset A^{\dag}_{K}$, donc a fortiori sur $B^{\dag}_{K} = \displaystyle \mathop{\lim}_{\rightarrow \atop{n}}\  B_{\lambda'}$ : il existe donc $\mu$, $1 < \mu \leqslant \lambda$ tel que pour tout $i$ on ait $F_{\hat{B}_{K}}(x_{i}) \in B_{\mu}$. Comme dans la preuve de [II, (2.3.1) (2)] on peut aussi supposer que pour tout $i$ et tout $g \in G$ on a $g_{\hat{B}_{K}}(x_{i}) \in B_{\mu}$. Ainsi $F'_{\lambda \mu} : W_{\mu} \rightarrow W_{\lambda}$ (resp. $g_{\lambda} : W_{\lambda} \rightarrow W_{\lambda}$ est induit par $F_{B^{\dag}} : B^{\dag} \rightarrow B^{\dag}$ (resp. $g_{B^{\dag}} : B^{\dag} \rightarrow B^{\dag})$ ; pour prouver le lemme suivant:
\vskip 3mm

\noindent \textbf{Lemme (4.2.3.3)}. 

$$
g_{\lambda} \circ F'_{\lambda \mu} = F'_{\lambda \mu} \circ g_{\mu}.
$$
\noindent il suffit de prouver le 

\vskip 3mm

\noindent \textbf{Lemme (4.2.3.4)}. 

$$
g_{B^{\dag}} \circ F_{B^{\dag}} = F_{B^{\dag}} \circ g_{B^{\dag}}.
$$

\noindent Or $g \in G$ induit un morphisme $g_{X} : X \rightarrow X$ tel que $g_{X} \circ F_{X} = F_{X}
 \circ g_{X}$, puisque $g(x^q) = g(x)^q$ pour toute section $x$ de $O_{X}$ ; d'o\`u un diagramme commutatif\\
 
 $$
 \begin{array}{c}
 \xymatrix{
X  \ar@{}[drr] |{\rondI}  && X^{(q/S)}  \ar[ll]_{\pi_{X}} \ar@{}[drr] |{\rondII} && X  \ar[ll]_{F_{X/S}}^{\sim}  \ar@/_2pc/[llll]_{F_{X}} && \\
X \ar[u]^{g_{X}} && X^{(q/S)}   \ar[ll]^{\pi_{X}} \ar[u]^{g^{(q)}_{X}} && X \ar[ll]^{F_{X/S}}_{\sim}  \ar[u]_{g_{X}} &&\  .
}
\end{array}
\leqno{(4.2.3.5)}
$$

 \noindent Le carr\'e commutatif    $\ \rondI$ se rel\`eve en le carr\'e commutatif\\
 
 $$
\begin{array}{c}
\xymatrix{
B^{\dag} \ar[rr]^(.4){1_{B^{\dag}}  \otimes F_{A^{\dag}}} \ar[d]_{g_{B^{\dag}}} && B^{\dag} \otimes A^{\dag} \ar[d]^{g_{B^{\dag}} \otimes 1_{A^{\dag}}}\\
B^{\dag} \ar[rr]_(.4){1_{B^{\dag}}  \otimes F_{A^{\dag}}} && B^{\dag} \otimes A^{\dag} \ .
}
\end{array}
\leqno{(4.2.3.6)}
$$ 

\noindent Par l'\'equivalence de cat\'egories $B^{\dag} \longmapsto B_{0}$ de la cat\'egorie des 
$A^{\dag}$-alg\`ebres finies \'etales dans la cat\'egorie des $A_{0}$-alg\`ebres finies \'etales [Et 4, th\'eo 7], on rel\`eve le carr\'e commutatif     $\ \rondII$ en le carr\'e commutatif \\

 $$
\begin{array}{c}
\xymatrix{
B^{\dag} \otimes A^{\dag} \ar[d]_{g_{B^{\dag}} \otimes 1_{A^{\dag}}}\ar[rr]_{\sim}^(.6){F_{B^{\dag}/A^{\dag}}} && B^{\dag} \ar[d]^{g_{B^{\dag}}}&\\
B^{\dag} \otimes A^{\dag}  \ar[rr]^{\sim}_(.6){F_{B^{\dag}/A^{\dag}}} && B^{\dag} & .
}
\end{array}
\leqno{(4.2.3.7)}
$$ 

\noindent Par composition on a prouv\'e (4.2.3.4), donc (4.2.3.3).\\

Compte tenu de la commutation (4.2.3.3) et de la d\'efinition de la cohomologie rigide (resp. de la cohomologie rigide \`a supports compacts) donn\'ee en (4.2.1.4), (4.2.1.5), (4.2.1.6) [resp. en (4.2.2.2), (4.2.2.3), (4.2.2.4)] les isomorphismes (4.2.1) et (4.2.2) du th\'eor\`eme (4.2) sont compatibles \`a l'action du Frobenius. $\square$\\

\noindent Dans la preuve du th\'eor\`eme (4.2) on a montr\'e au passage :

\vskip 3mm
\noindent \textbf{Lemme (4.3)}. \textit{Si $S$ est un $k$-sch\'ema s\'epar\'e de type fini, $f : X \longrightarrow S$ est fini \'etale (non n\'ecessairement galoisien) et $E \in \textrm{Isoc}^{\dag}(X/K)$ on a des isomorphismes canoniques }

\begin{enumerate}
\item[(1)] $H^i_{\textrm{rig}}(X/K;E) \tilde{\longrightarrow} H^i_{\textrm{rig}}(S/K;f_{\textrm{rig}\ast}\ E)$.
\item[(2)] $H^i_{\textrm{rig},c}(X/K;E) \tilde{\longrightarrow} H^i_{\textrm{rig},c}(S/K;f_{\textrm{rig}\ast}\ E)$.
\item[(3)] \textit{Si de plus $E \in F^a\mbox{-}\textrm{Isoc}^{\dag}(X/K)$ les isomorphismes du (1) et (2) commutent \`a l'action de Frobenius}.
\end{enumerate}

\vskip 3mm
\noindent \textbf{Remarques (4.4)}.
\begin{itemize}
\item[(i)] Les r\'esultats du lemme (4.3) sont donn\'es par Tsuzuki dans [Tsu 1, cor (2.6.5) et (2.6.6)], sans pr\'ecisions de d\'emonstration, notamment pour le (2) du lemme : nous y avons utilis\'e le th\'eor\`eme de changement de base pour un morphisme propre [II, (3.3.2)], qui n'est pas mentionn\'e  par Tsuzuki.
\item[(ii)] Le (4.2.2) du th\'eor\`eme (4.2) est une \'etape essentielle pour \'etablir la finiture de la cohomologie rigide \`a supports compacts \`a coefficients dans un $F$-isocristal surconvergent unit\'e \`a partir de la finitude de la cohomologie cristalline via la suite exacte longue de localisation en cohomologie rigide, la preuve de cet isomorphisme crucial n'appara\^it pas dans la d\'emonstration du th\'eor\`eme 6.1.2 de [Tsu 1].
\end{itemize}

\vskip10mm
\section*{5. Cas plongeable}

\textbf{5.1.} On suppose donn\'e un diagramme commutatif\\

$$
\xymatrix{
X \ar@{^{(}->}[r]^{j_{\mathcal{Y}}} \ar[d]_{f} & \mathcal{Y} \ar[d]^{\overline{h}} &\\
S \ar@{^{(}->}[r]_{j_{\mathcal{T}}} 		  & \mathcal{T} \ar[r]_{\rho} & Spf\  \mathcal{V}
}
$$

\noindent dans lequel $f$ est un morphisme de $k$-sch\'emas s\'epar\'es de type fini, $\overline{h}$ et $\rho$ sont des morphismes propres de $\mathcal{V}$-sch\'emas formels, $\overline{h}$ (resp. $\rho$) est lisse sur un voisinage de $X$ dans $\mathcal{Y}$ (resp. un voisinage de $S$ dans $\mathcal{T}$), $j_{\mathcal{Y}}$ et $j_{\mathcal{T}}$ sont des immersions. D\'esignons par $T$ (resp. $Y$) l'adh\'erence sch\'ematique de $S$ dans $\mathcal{T}$ (resp. de $X$ dans $\mathcal{Y}$), $\overline{f} : Y \rightarrow T$ le morphisme induit par $\overline{h}$, $i_{Y} : Y \hookrightarrow \mathcal{Y}$ l'immersion ferm\'ee, $X_{1} := \overline{f}^{-1}(S)$ et $f_{1} : X_{1} \rightarrow S$ le morphisme induit par $\overline{f}$.\\

On note $F_{S}$ (resp. $F_{X}$) le Frobenius de $S$ (resp. de $X$) (\'el\'evation \`a la puissance $q = p^a$ sur le faisceau structural) ; d'o\`u le diagramme commutatif

$$
\xymatrix{X \ar@/^1pc/[rrd]^{F_{X}} \ar@/_1,5pc/[rdd]_{f} \ar[rd]_{F_{X/S}} \\
&  X^{(q)} \ar[r]^{\pi_{X/S}} \ar[d]^{f^{(q)}} & X \ar[d]^{f}&\\
& S \ar[r]_{F_{S}} & S &.
}
$$

\vskip 3mm
\noindent \textbf{Th\'eor\`eme (5.2)}. \\
\textit{(5.2.1) Sous les hypoth\`eses (5.1) supposons que $\overline{h}^{-1}(S)= \overline{f}^{-1}(S) = X$; alors, pour tout entier $i \geqslant 0$, le morphisme $f$ induit un foncteur}

$$
R^i f_{\textrm{rig} \ast}(X/ \mathcal{T}, -) : F^a\mbox{-}\textrm{Isoc}^{\dag}(X/K) \longrightarrow F^a\mbox{-}\textrm{Isoc}^{\dag}(S/K).
$$

\noindent \textit{(5.2.2) Supposons donn\'es des morphismes}

$$
\xymatrix{
S' \ar@{^{(}->}[r]^{j'} & T' \ar@{^{(}->}[r]^{\rho'} & Spf \mathcal{V}
}
$$

\noindent \textit{o\`u $\rho'$ est un morphisme propre de $\mathcal{V}$-sch\'emas formels, $S'$ est un $k$-sch\'ema s\'epar\'e de type fini, $j'$ est une immersion et $\rho'$ est lisse sur un voisinage de $S'$ dans $\mathcal{T}'$. Alors le foncteur de (5.2.1) commute \`a tout changement de base s\'epar\'e de type fini $S' \rightarrow S$ : en particulier ce foncteur commute aux passages aux fibres en les points ferm\'es de $S$}.\\

\noindent \textit{D\'emonstration}.\\
\textit{Pour (5.2.1)}, soit $(E, \phi) \in F^a\mbox{-}\textrm{Isoc}^{\dag}(X/K)$; pour tout entier $i \geqslant 0$, on a d'apr\`es [II, (3.4.4)] un isomorphisme

$$
F^{\ast}_{S} R^i f_{\textrm{rig} \ast} (X/ \mathcal{T}, E)  \displaystyle \mathop{\longrightarrow}^{\sim} R^if^{(q)}_{\textrm{rig} \ast}(X^{(q)}/\mathcal{T}, \pi^{\ast}_{X/S}(E)).
$$

\noindent L'identit\'e de $S$ induit un morphisme

$$
\xymatrix{
\theta^i : R^i f^{(q)}_{\textrm{rig} \ast}(X^{(q)}/ \mathcal{T}, \pi^{\ast}_{X/S}(E)) \ar[r] & R^i f_{\textrm{rig} \ast}(X/ \mathcal{T}, F^{\ast}_{X/S} \pi^{\ast}_{X/S}(E)) \ar@{=}[d]\\
&  R^i f_{\textrm{rig} \ast}(X/ \mathcal{T}, F^{\ast}_{X} (E)),
}
$$

\noindent et le Frobenius $\phi$ de $E$ induit un isomorphisme

$$
R^i f_{\textrm{rig} \ast}(X/ \mathcal{T}, F^{\ast}_{X} E) \displaystyle \mathop{\longrightarrow}^{\sim} R^i f_{\textrm{rig} \ast}(X/ \mathcal{T},  E).
$$

\noindent Par composition de ces trois morphismes on obtient le Frobenius de $R^i f_{\textrm{rig} \ast}(X/ \mathcal{T},  E)$\\

\noindent (5.2.3) $\qquad \qquad \qquad \phi^{i}: \  F^{\ast}_{S} R^i f_{\textrm{rig} \ast}(X/ \mathcal{T},  E) \longrightarrow R^i f_{\textrm{rig} \ast}(X/ \mathcal{T},  E) $\\

\noindent et il s'agit de prouver que $\phi^i$ est un isomorphisme : pour \c{c}a il suffit de prouver que c'est le cas pour $\theta^i$. On sait d\'ej\`a que $\theta^i$ est un morphisme d'isocristaux surconvergents  : d'apr\`es [B 3, (2.1.11) et (2.2.7)] il suffit de montrer que $\theta^i$ induit un isomorphisme dans la cat\'egorie convergente Isoc$(S/K)$ ; d'apr\`es [B-G-R, 9.4.3/3 et 9.4.2/7] il suffit de le v\'erifier apr\`es passage aux fibres de $\theta^i$ en les points ferm\'es $s$ de $S$. Pour un tel point $s$ notons $\mathcal{V}(s) = W(k(s)) \otimes_{W} \mathcal{V}$ et $K(s)$ le corps des fractions de $\mathcal{V}(s)$. D'apr\`es [II, (3.4.4)] on a un diagramme commutatif

$$
\xymatrix{
R^i f^{(q)}_{\textrm{rig} \ast}(X^{(q)}/\mathcal{T}, \pi^{\ast}_{X/S}(E))_{s} \ar[r]^{\theta_{s}} \ar[d]^{\simeq} & R^i f_{\textrm{rig} \ast}(X/\mathcal{T}, F^{\ast}_{X}(E))_{s} \ar[d]_{\simeq}\\
R^i f_{s\ \textrm{rig} \ast}^{(q)} (X^{(q)}_{s}/\mathcal{V}(s), E_{X_{s}^{(q)}}) \ar@{.>}[r] \ar@{=}[d] & R^i f_{s\ \textrm{rig} \ast} (X_{s}/\mathcal{V}(s),F^{\ast}_{X_{s}} (E_{X_{s}})) \ar@{=}[d]\\
H^i_{\textrm{rig}} (X^{(q)}_{s}/K(s), E_{X^{(q)}_{s}}) \ar@{.>}[r]  \ar@{=}[d] & H^i_{\textrm{rig}}(X_{s}/K(s), F^{\ast}_{X_{s}}(E_{X_{s}})) \ar@{=}[d] \\
H^i_{\textrm{rig}, c} (X^{(q)}_{s}/K(s), E_{X^{(q)}_{s}}) \ar@{.>}[r]   & H^i_{\textrm{rig}, c}(X_{s}/K(s), F^{\ast}_{X_{s}}(E_{X_{s}})) 
}
$$

\noindent o\`u les fl\`eches verticales sont des isomorphismes ; or la fl\`eche horizontale inf\'erieure est un isomorphisme par [E-LS 1, prop 2.1, o\`u il faut supposer $X$ lisse sur $\mathbb{F}_{q}$ dans le cas de la cohomologie sans support]. D'o\`u (5.2.1).\\
L'assertion (5.2.2) r\'esulte de [II, (3.4.4)]. $\square$


\cleardoublepage


\vskip 10mm
\chapter*{V.  Cohomologie syntomique}
\markboth{\sc j.-y. etesse}{\sc V.  Cohomologie syntomique}

\section*{1. Site syntomique}

\textbf{1.1.} Si $X$ est un sch\'ema quelconque, le gros site syntomique de $X$ est d\'efini comme suit [SGA 3, IV, 6.3] : la cat\'egorie sous-jacente est celle des sch\'emas sur $X$ et la topologie est engendr\'ee par les familles finies surjectives de morphismes syntomiques (i.e. plats et localement intersection compl\`ete). On rappelle que les morphismes syntomiques sont ouverts, demeurent syntomiques par changement de base, et sont localement relevables le long d'une immersion ferm\'ee. Le gros site (resp. petit site) syntomique de $X$ sera not\'e $SYNT(X)$ (resp. $synt(X)$) et le topos correspondant $X_{SYNT}$ (resp. $X_{synt}$) : lorsqu'on ne voudra pas distinguer entre les deux situations on notera $\mathcal{T}(X)$ (resp. $X_{\mathcal{T}}$) l'un ou l'autre de ces deux sites (resp. topos).\\

Soit $j : U \hookrightarrow X$ une immersion ouverte ; $j$ d\'efinit un couple de foncteurs adjoints $(j_{\ast}, j^{-1})$

$$
U_{\mathcal{T}} \displaystyle \mathop{\rightarrow}^{j^{-1} \atop{\longleftarrow}}_{j _{\ast}} X_{\mathcal{T}}.
$$

\noindent Dans la suite $X_{\mathcal{T}}$ sera annel\'e par un faisceau d'anneaux $\mathcal{A}$ et $U_{\mathcal{T}}$ sera annel\'e par $j^{-1} \mathcal{A}$, not\'e $\mathcal{A}_{\vert U}$ : on note $_{\mathcal{A}}X_{\mathcal{T}}$ (resp. $_{\mathcal{A}\vert U}U_{\mathcal{T}}$) la cat\'egorie des faisceaux des $\mathcal{A}$-modules \`a gauche sur $X_{\mathcal{T}}$ (resp. des $\mathcal{A}_{\vert U}$-modules \`a gauche sur $U_{\mathcal{T}}$). Le foncteur

$$
j^{\ast} : \ _{\mathcal{A}}X_{\mathcal{T}} \longrightarrow\  _{\mathcal{A}\vert U}U_{\mathcal{T}}
$$

\noindent admet un adjoint \`a gauche $j_{!}$ [Mi, II, Rk 3.18] et [SGA 4, IV, $\S$ 14] d\'efini par\\

\noindent
(1.1.1)
$\left\lbrace
\begin{array}{l}
j_{!} (\mathcal{F})(X') = \mathcal{F}(X')\  \textrm{si}\ X' \longrightarrow X\  \textrm{se factorise par}\  U\\
\textrm{et}\  j_{!} (\mathcal{F}) (X') = 0\ \textrm{sinon\ ;}
\end{array}
\right.
$\\

\noindent $j_{!}$ est exact [loc. cit].\\

\noindent Remarquons que $j^{\ast}$ est exact puisqu'il admet un adjoint \`a droite et un adjoint \`a gauche.\\

De m\^eme si $i : Z \hookrightarrow X$ est une immersion ferm\'ee, $i$ d\'efinit un couple de foncteurs adjoints $(i_{\ast}, i^{-1})$ :

$$
Z_{\mathcal{T}} \displaystyle \mathop{\rightarrow}^{i^{-1} \atop{\longleftarrow}}_{i _{\ast}} X_{\mathcal{T}}\ ;
$$

\noindent $Z_{\mathcal{T}}$ sera annel\'e par $i^{-1} \mathcal{A}$, not\'e $\mathcal{A}_{\vert Z}$.\\

\noindent Le foncteur
$$
i^{\ast} :\  _{\mathcal{A}}X_{\textrm{SYNT}} \longrightarrow\ _{\mathcal{A} \vert Z}Z_{\textrm{SYNT}}
$$

\noindent admet un adjoint \`a gauche $i_{!}$ [Mi, II, Rk 3.18] et $i_{!}$ est exact [loc. cit.] : en particulier 
$$i^{\ast} :\  _{\mathcal{A}}X_{\textrm{SYNT}} \longrightarrow\ _{\mathcal{A} \vert Z}Z_{\textrm{SYNT}}$$ est exact.\\
Le foncteur $$i^{\ast} :\  _{\mathcal{A}}X_{\textrm{synt}} \longrightarrow\ _{\mathcal{A} \vert Z}Z_{\textrm{synt}}$$ est lui aussi exact gr\^ace \`a [Mi, II, 2.6 et 3.0 p 68] et [EGA $0_{I}$, 1.4.12] car les morphismes syntomiques demeurent syntomiques par changement de base.\\
De plus le foncteur 
$$i_{\ast} :\  _{\mathcal{A} \vert Z}Z_{\mathcal{T}} \longrightarrow\ _{\mathcal{A}}X_{\mathcal{T}}$$
  est exact, car tout morphisme syntomique se rel\`eve, localement le long d'une immersion ferm\'ee, en un morphisme syntomique.\\

\textbf{1.2.} Soient $i : Z \hookrightarrow X$ une immersion ferm\'ee, et $j : U \hookrightarrow X$ l'immersion ouverte du compl\'ementaire de $Z$.
Pour un $X$-sch\'ema $X'$ on note $U' = X' \times_{X} U$, $Z' = X' \times_{X} Z$; tout recouvrement syntomique $W \twoheadrightarrow Z'$ se rel\`eve localement en $\tilde{W} \rightarrow X'$ syntomique : pour all\'eger l'\'ecriture on supposera le rel\`evement global. Par suite, si $\tilde{U} \twoheadrightarrow U'$ est surjectif syntomique et $\tilde{W}$ tel que ci-dessus, alors $(\tilde{U}, \tilde{W})$ est un recouvrement syntomique de $X'$.

\vskip 3mm
\noindent \textbf{Lemme (1.2.1)}. 
\textit{Avec les notations de (1.2) et pour $\mathcal{F} \in {_{\mathcal{A}}X_{\mathcal{T}}}$ le carr\'e suivant, o\`u les fl\`eches sont les fl\`eches canoniques}

$$
\xymatrix{
\mathcal{F} \ar[r] \ar[d] & j_{\ast} j^{\ast}(\mathcal{F)} \ar[d] \\
i_{\ast} i^{\ast}(\mathcal{F}) \ar[r] & i_{\ast} j_{\ast} j^{\ast}(\mathcal{F})
}
$$
\noindent \textit{est cart\'esien.}

\vskip 3mm
\noindent \textit{D\'emonstration}. On notera $\mathcal{G}$ le produit fibr\'e.\\
Pour un $X$-sch\'ema $X'$ on consid\`ere un recouvrement $(\tilde{U}, \tilde{W})$ de $X'$ du type pr\'ec\'edent : comme $\tilde{U} \rightarrow X'$ et $\tilde{W} \rightarrow X'$ sont deux morphismes syntomiques, on est ramen\'e \`a \'etablir l'isomorphisme $\mathcal{F}\  \tilde{\rightarrow}\  \mathcal{G}$ au-dessus d'un $X'$-sch\'ema syntomique $\tilde{W}$ ; la d\'emonstration se fait donc sur le petit site de $X'$. On pose $W = \tilde{W} \times_{X'}Z'$.\\

Le faisceau $j_{\ast} j^{\ast}(\mathcal{F})$ est le faisceau associ\'e au pr\'efaisceau

$$
\tilde{W} \longmapsto \mathcal{F}(V) \qquad , \quad \textrm{avec}\  V := \tilde{W} \times_{X'} U' ,
$$

\noindent et $i_{\ast} i^{\ast}(\mathcal{F})$ est le faisceau associ\'e au pr\'efaisceau 

$$
\tilde{W} \longmapsto \displaystyle \mathop{\lim}_{\rightarrow \atop{W'}} \mathcal{F}(\mathcal{W}')\ ,
$$
\noindent la limite \'etant prise sur les diagrammes commutatifs

$$\begin{array}{c}
\xymatrix{
\tilde{W} \ar[d] & W' \ar[l] & W \ar[l] \ar[d]\\
X' & &Z' \ar[ll]
}
\end{array}
\leqno{(1.2.2)}
$$

\noindent avec $W' \rightarrow \tilde{W}$ syntomique.\\
De m\^eme $i_{\ast} i^{\ast} j_{\ast} j^{\ast}(\mathcal{F})$ est le faisceau associ\'e au pr\'efaisceau

$$
\tilde{W} \longmapsto \displaystyle \mathop{\lim}_{\rightarrow \atop{W'}} \mathcal{F}(\mathcal{W}' \times_{\tilde{W}} V) =  \displaystyle \mathop{\lim}_{\rightarrow \atop{W'}} \mathcal{F} (\mathcal{W}' \times_{X'} U'),
 $$
 
 \noindent avec $W'$ comme en (1.2.2). On remarque alors que $(V, W')$ est un recouvrement syntomique de $\tilde{W}$ : en effet le $\tilde{W}$-morphisme $W \rightarrow W'$ fournit une section du morphisme syntomique $W' \times_{\tilde{W}} W \rightarrow W$ ; ce dernier est donc surjectif et on conclut comme pour le recouvrement $(\tilde{U}, \tilde{W})$ de $X'$.\\
 Ainsi on a bien un isomorphisme $\mathcal{F}\  \tilde{\rightarrow}\  \mathcal{G}$ au-dessus de $\tilde{\mathcal{W}}$, d'o\`u le lemme. $\square$\\
 
Notons $\mathbb{T}(_{\mathcal{A}}X_{\mathcal{T}})$ la cat\'egorie des triplets $(\mathcal{F}_{1}, \mathcal{F}_{2}, \alpha)$ o\`u $\mathcal{F}_{1} \in _{\mathcal{A} \vert Z}Z_{\mathcal{T}}$, $\mathcal{F}_{2} \in\  _{\mathcal{A} \vert U}U_{\mathcal{T}}$ et $\alpha$ est un morphisme $\alpha : \mathcal{F}_{1} \rightarrow i^{\ast} j_{\ast}  \mathcal{F}_{2}$ ; les morphismes entre deux tels triplets sont d\'efinis de la mani\`ere naturelle, analogue \`a [Mi, II, \S\ 3].
 
 \vskip 3mm
\noindent \textbf{Th\'eor\`eme(1.3)}. 
\textit{Soient $i : Z \hookrightarrow X$ une immersion ferm\'ee de sch\'emas et $j : U \hookrightarrow X$ l'immersion ouverte du compl\'ementaire de $Z$. Le foncteur
$$\mathcal{F} \longmapsto (i^{\ast} \mathcal{F}, j^{\ast} \mathcal{F}, \alpha)$$ o\`u $\alpha$
est le morphisme canonique $\alpha : i^{\ast} \mathcal{F} \rightarrow i^{\ast} j_{\ast} j^{\ast} \mathcal{F}$, induit une \'equivalence de cat\'egories entre $_\mathcal{A}X_{\mathcal{T}}$ et $\mathbb{T}(_\mathcal{A}X_{\mathcal{T}})$.
 }

\vskip 3mm
\noindent \textit{D\'emonstration}. La d\'emonstration est analogue \`a celle de Fontaine-Messing [F-M, 4.4]. Le th\'eor\`eme r\'esulte du lemme (1.2.1) par la m\^eme m\'ethode que pour le site \'etale [Mi, II, theo 3.10]. $\square$\\

En identifiant $_\mathcal{A}X_{\mathcal{T}}$ et $\mathbb{T}(_\mathcal{A}X_{\mathcal{T}})$ via le th\'eor\`eme (1.3) on d\'efinit six foncteurs

$$
\begin{array}{c}
\xymatrix{
& \ar[l]_{i^{\ast}} \quad & \quad  \ar[l]_{j_{!}} &\\
_{\mathcal{A} \vert Z}Z_{\mathcal{T}} \ar[r]^{i_{\ast}} & _{\mathcal{A}}X_{\mathcal{T}} \ar[r]^{j^{\ast}} & _{\mathcal{A}}U_{\mathcal{T}}\\
& \ar[l]_{i^{!}} \quad & \quad  \ar[l]_{j_{\ast}} &
}
\end{array}
\leqno{(1.4)}
$$

\noindent dont la description est la suivante :\\
$i^{\ast} : \mathcal{F}_{1} \leftarrow (\mathcal{F}_{1}, \mathcal{F}_{2}, \alpha : \mathcal{F}_{1} \rightarrow i^{\ast} j_{\ast} \mathcal{F}_{2}),\  j_{!} : (0, \mathcal{F}_{2}, 0) \leftarrow \mathcal{F}_{2}$\\
$i_{\ast} : \mathcal{F}_{1} \mapsto (\mathcal{F}_{1}, 0, 0) \quad \qquad \qquad \qquad ,\ j^{\ast} : (\mathcal{F}_{1}, \mathcal{F}_{2}, \alpha : \mathcal{F}_{1} \rightarrow i^{\ast} j_{\ast} \mathcal{F}_{2}) \mapsto \mathcal{F}_{2}$\\
$i^! : \mbox{Ker} \alpha \leftarrow (\mathcal{F}_{1}, \mathcal{F}_{2}, \alpha : \mathcal{F}_{1} \rightarrow i^{\ast} j_{\ast} \mathcal{F}_{2})$,\\
 $j_{\ast} : (i^{\ast} j_{\ast} \mathcal{F}_{2}, \mathcal{F}_{2}, id : i^{\ast} j_{\ast} \mathcal{F}_{2} \rightarrow i^{\ast} j_{\ast} \mathcal{F}_{2}) \leftarrow \mathcal{F}_{2}. $

\vskip 3mm
\noindent \textbf{Th\'eor\`eme(1.5)}. 
\textit{Avec les notations pr\'ec\'edentes on a :
\begin{enumerate}
\item[(1)] Chaque foncteur est adjoint \`a gauche de celui \'ecrit la ligne au-dessous ; en particulier on a un isomorphisme de transitivit\'e $(j_{1} j_{2})_{!} = j_{1!} j_{2!}$.
\item[(2)] Les foncteurs $i^{\ast}, i_{\ast}, j^{\ast}, j_{!}$ sont exacts ; les foncteurs $j_{\ast}$, $i^{!}$ sont exacts \`a gauche.
\item[(3)] Les compos\'es $i^{\ast} j_{!}, i^{!} j_{!}, i^{!} j_{\ast}, j^{\ast} i_{\ast}$ sont nuls.
\item[(4)] Les foncteurs $i_{\ast}, j_{\ast}$ et $j_{!}$ sont pleinement fid\`eles.
\item[(5)] Les foncteurs $j_{\ast}, j^{\ast}, i^{!}, i_{\ast}$ envoient les injectifs sur les injectifs.
\item[(6)] Pour tout $\mathcal{F} \in _{\mathcal{A}}X_{\mathcal{T}}$  (resp. $\mathcal{F} \in _{\mathcal{A} \vert U}U_{\mathcal{T}})$ on a des suites exactes courtes
\begin{itemize}
\item[(6.1)] $0 \longrightarrow j_{!} j^{\ast} \mathcal{F} \longrightarrow \mathcal{F} \longrightarrow i_{\ast} i^{\ast} \mathcal{F} \longrightarrow 0$
\item[(6.2)] $0 \longrightarrow i_{\ast} i^{!} \mathcal{F} \longrightarrow \mathcal{F} \longrightarrow j_{\ast}j^{\ast} \mathcal{F}$
\item[[resp. (6.3)] $0 \longrightarrow j_{!} \mathcal{F} \longrightarrow j_{\ast} \mathcal{F} \longrightarrow i_{\ast} i^{\ast} j_{\ast} \mathcal{F} \longrightarrow 0].$\\
De plus le couple de foncteurs $(i^{\ast}, j^{\ast}) $ est conservatif.
\end{itemize}
\end{enumerate}
}

\vskip 3mm
\noindent \textit{D\'emonstration}.\par

Le \textit{(1)} et le \textit{(2)} ont d\'ej\`a \'et\'e vus, et l'isomorphisme de transitivit\'e $(j_{1} j_{2})_{!} = j_{1!} j_{2!}$ r\'esulte de la formule $(j_{1} j_{2})^{\ast} = j_{2}^{\ast} j_{1}^{\ast}$.\par

Le \textit{(3)} et le \textit{(4)} r\'esultent des descriptions (1.4).\par

Le \textit{(5)} r\'esulte de \textit{(1)} et \textit{(2)} et du fait qu'un foncteur avec un adjoint \`a gauche exact pr\'eserve les injectifs [Mi, III, 1.2].\par

La suite \textit{(6.3))} provient de \textit{(6.1)} en remarquant que $j^{\ast} j_{\ast} \mathcal{F} = \mathcal{F}.$\par

Compte tenu des identifications (1.4) la suite \textit{(6.1)} s'\'ecrit

$$
0 \longrightarrow (0, j^{\ast} \mathcal{F}, 0) \longrightarrow (i^{\ast} \mathcal{F}, j^{\ast} \mathcal{F}, \alpha) \longrightarrow (i^{\ast} \mathcal{F}, 0, 0) \longrightarrow 0.
$$

\noindent Pour d\'emontrer qu'elle est exacte il nous suffit donc de montrer qu'une suite de $_{\mathcal{A}}X_{\mathcal{T}}$\\

\noindent \textit{(6.4)} $\qquad \qquad  0 \longrightarrow \mathcal{F}' \longrightarrow \mathcal{F} \longrightarrow \mathcal{F}'' \longrightarrow 0$\\

\noindent est exacte si et seulement si les suites\\

\noindent \textit{(6.5)} $\qquad \qquad  0 \longrightarrow i^{\ast}(\mathcal{F}') \longrightarrow i^{\ast}(\mathcal{F}) \longrightarrow i^{\ast}(\mathcal{F}'') \longrightarrow 0$\\

\noindent \textit{(6.6)} $\qquad \qquad  0 \longrightarrow j^{\ast}(\mathcal{F}') \longrightarrow j^{\ast}(\mathcal{F}) \longrightarrow j^{\ast}(\mathcal{F}'') \longrightarrow 0$\\

\noindent de $_{\mathcal{A} \vert Z}Z_{\mathcal{T}}$ et $_{\mathcal{A} \vert U}U_{\mathcal{T}}$ respectivement, sont exactes, i.e. que le couple de foncteurs $(i^{\ast}, j^{\ast})$ est conservatif.\\

L'exactitude de \textit{(6.4)} entra\^{\i}ne celle de \textit{(6.5)} et \textit{(6.6)} puisque $i^{\ast}$ et $j^{\ast}$ sont exacts.\\

R\'eciproquement, l'exactitude de \textit{(6.5)} et \textit{(6.6)} entra\^{\i}ne celle des suites\\

\noindent \textit{(6.7)} $\qquad \qquad  0 \longrightarrow i_{\ast} i^{\ast}(\mathcal{F}') \longrightarrow i_{\ast} i^{\ast}(\mathcal{F}) \displaystyle \mathop{\longrightarrow} ^{\varphi_{Z}} i_{\ast} i^{\ast}(\mathcal{F}'') \longrightarrow 0,$\\

\noindent \textit{(6.8)} $\qquad \qquad  0 \longrightarrow j_{\ast} j^{\ast}(\mathcal{F}') \longrightarrow j_{\ast} j^{\ast}(\mathcal{F}) \displaystyle \mathop{\longrightarrow} ^{\varphi_{U}} j_{\ast} j^{\ast}(\mathcal{F}'') ,$\\

\noindent \textit{(6.9)} $\qquad \qquad  0 \longrightarrow i_{\ast} i^{\ast} j_{\ast} j^{\ast}(\mathcal{F}') \longrightarrow i_{\ast} i^{\ast} j_{\ast} j^{\ast}(\mathcal{F}) \longrightarrow i_{\ast} i^{\ast} j_{\ast} j^{\ast}(\mathcal{F}''),$\\

\noindent d'o\`u l'exactitude de\\

\noindent \textit{(6.10)} $\qquad \qquad \qquad 0 \longrightarrow \mathcal{F}' \longrightarrow \mathcal{F} \longrightarrow \mathcal{F}''$\\

\noindent gr\^ace au lemme (1.2.1).\\

Donc pour tout $\mathcal{F} \in _{\mathcal{A}}X_{\mathcal{T}}$ on a l'exactitude de la suite

$$
0 \longrightarrow j_{!}j^{\ast} \mathcal{F} \longrightarrow \mathcal{F} \longrightarrow i_{\ast} i^{\ast} \mathcal{F}.
$$

\noindent Montrons la surjectivit\'e de $\rho : \mathcal{F} \rightarrow i_{\ast}i^{\ast} \mathcal{F}$. Comme pour le lemme (1.2.1), dont on utilise les notations, on est ramen\'e au petit site de $X'$. Soient $\tilde{W}$ un $X'$-sch\'ema syntomique, $V := \tilde{W} \times_{X'}U'$ et $s \in i_{\ast} i^{\ast}(\mathcal{F})(\tilde{W}) = \displaystyle \mathop{\lim}_{\rightarrow \atop{W'}} \mathcal{F}(W')$, o\`u $W'$ est comme dans (1.2.2) : il existe un $W'$ et $s_{W'} \in \mathcal{F}(W')$ d'image $s$. On a vu dans la preuve du lemme (1.2.1) que $(V, W')$ est un recouvrement syntomique de $\tilde{W}$ : notons $s'$ l'image de $s$ par l'application restriction

$$
i_{\ast} i^{\ast} \mathcal{F}(\tilde{W}) \longrightarrow i_{\ast} i^{\ast} \mathcal{F}(W')
$$

$$
s \mapsto s'\ ;
$$

\noindent alors $s_{W'} \in \mathcal{F}(W')$ a pour image $s'$ par $\rho$.\\

\noindent Comme $i_{\ast} i^{\ast}(\mathcal{F})(V) = 0$, tout $s_{V} \in \mathcal{F}(V)$ est un rel\`evement de $s$ dans $i_{\ast} i^{\ast}(\mathcal{F})(V) = 0$. D'o\`u la surjectivit\'e de $\rho$, et l'exactitude de \textit{(6.1)}.\\

On a alors un diagramme commutatif \`a lignes exactes

$$
\xymatrix{
0 \ar[r]  & j_{!} j^{\ast} \mathcal{F} \ar[r] \ar[d]_{u} & \mathcal{F} \ar[r] \ar[d]_{v} & i_{\ast} i^{\ast} \mathcal{F} \ar[r] \ar[d]_{w} & 0 \\
0 \ar[r] & j_{!} j^{\ast} \mathcal{F}'' \ar[r] & \mathcal{F}'' \ar[r] & i_{\ast} i^{\ast} \mathcal{F}'' \ar[r] & 0.
}
$$

\noindent L'exactitude de \textit{(6.5)} et \textit{(6.6)} et celle des foncteurs $j_{!}$ et $i_{\ast}$ prouve que $u$ et $w$ sont surjectifs : d'o\`u la surjectivit\'e de $v$ ; jointe \`a l'exactitude de \textit{(6.10)} ceci prouve que le couple $(i^{\ast}, j^{\ast})$ est conservatif.\\

\noindent L'exactitude de \textit{(6.2)} en r\'esulte via la description de $i^{!}$ fournie en (1.4). $\square$

\newpage
\section*{2. Cohomologie syntomique \`a supports compacts}

Soit $X$ un sch\'ema s\'epar\'e de type fini sur un corps $k$ ; on sait par Nagata que $X$ est ouvert dans un $k$-sch\'ema propre $\overline{X}$ ; on note $j : X \hookrightarrow \overline{X}$ l'immersion ouverte.\\

On se place sous les notations de (1.1) : $\overline{X}_{\mathcal{T}}$ est annel\'e par un faisceau d'anneaux $\mathcal{A}$ et $\mathcal{F}$ est un faisceau de $\mathcal{A}$-modules.\\

Si $\mathcal{F}$ est \'element de $_{\mathcal{A}}\overline{X}_{\textrm{synt}}$, on notera encore $\mathcal{F}$ son image inverse dans $_{\mathcal{A}}\overline{X}_{\textrm{SYNT}}$ par le morphisme de topos $\overline{X}_{\textrm{SYNT}} \rightarrow \overline{X}_{\textrm{synt}}$ et pour tout entier $i \geqslant 0$, on a [E-LS 2, (1.2)]

$$
H^i(\overline{X}_{\textrm{synt}}, \mathcal{F}) = H^i(\overline{X}_{\textrm{SYNT}}, \mathcal{F}),
$$
\noindent et de m\^eme pour la topologie \'etale.

\vskip 3mm
\noindent \textbf{Proposition - D\'efinition (2.1)}. \textit{Sous les hypoth\`eses pr\'ec\'edentes, si $\mathcal{A}$ est de torsion et $\mathcal{F}$ un \'el\'ement de $_{\mathcal{A}_{\vert X}}X_{\mathcal{T}}$, le complexe $R \Gamma(\overline{X}_{\mathcal{T}}, j_{!}\  \mathcal{F})$ est ind\'ependant de la compactification $\overline{X}$ de $X$, et sera not\'e}

$$
R \Gamma_{\textrm{synt},c}(X, \mathcal{F})\ ;
$$

\noindent \textit{ses groupes de cohomologie seront not\'es}

$$
H^i_{\textrm{synt},c}(X, \mathcal{F})
$$

\noindent \textit{et appel\'es groupes de cohomologie syntomique \`a supports compacts.}

\vskip 3mm
\noindent \textit{D\'emonstration}. Si 
$$\xymatrix{
X \ar@{^{(}->}[rr]^{j'}  \ar[rd] & & \overline{X}' \ar[ld] \\
 &  Spec\ k &
 }$$
est une autre compactification de $X$, on note $\overline{X}''$ l'image sch\'ematique de $X$ plong\'e diagonalement dans $\overline{X}\times_{k} \overline{X}'$, $j'' : X \hookrightarrow \overline{X}''$ l'immersion ouverte et $g : \overline{X}'' \rightarrow \overline{X}$, $g' : \overline{X}'' \rightarrow \overline{X}'$ les deux projections (propres).\\

Il s'agit de montrer que\\

\noindent (2.2) $\qquad \qquad \qquad R \Gamma(\overline{X}_{\mathcal{T}}, j_{!}\  \mathcal{F}) = R \Gamma(\overline{X}''_{\mathcal{T}}, j''_{!}\  \mathcal{F}).$\\

Ceci va r\'esulter de la proposition plus g\'en\'erale suivante.

\vskip 3mm
\noindent \textbf{Proposition (2.3)}. 
\textit{Supposons donn\'e un carr\'e cart\'esien de $k$-sch\'emas}

$$
\xymatrix{
X' \ar@{^{(}->}[r]^{j'} \ar[d]_f & \overline{X}' \ar[d]^{\overline{f}}&\\
X \ar@{^{(}->}[r]^{j} & \overline{X}&,
}
$$

\noindent \textit{o\`u $\overline{X}$ est propre sur $k$, $\overline{f}$ est propre, $j$, $j'$ sont des immersions ouvertes. Si $\mathcal{F}$ est un faisceau ab\'elien de torsion sur $X'_{\mathcal{T}}$, alors on a des isomorphismes}\\

$\qquad \qquad R \Gamma(\overline{X}_{\mathcal{T}}, j_{!}\  Rf_{\ast}\ \mathcal{F}) \simeq R \Gamma(\overline{X}_{\mathcal{T}}, R \overline{f}_{\ast}\ j'_{!}\  \mathcal{F})$\\

$\qquad \qquad \qquad \qquad \qquad  \qquad \simeq R \Gamma(\overline{X}'_{\mathcal{T}},  j'_{!}\  \mathcal{F}).$\\

En effet (2.1) r\'esulte de (2.3) via le lemme suivant :

\vskip 3mm
\noindent \textbf{Lemme (2.4)}. 
\textit{Si} 
$$
\xymatrix{
& \overline{X}'' \ar[dd]^{g}\\
X \ar@{^{(}->}[ur]^{j''} \ar@{^{(}->}[rd]_{j}\\
& \overline{X}
}
$$
\textit{est un triangle commutatif de sch\'emas, avec $j$, $j''$ des immersions ouvertes dominantes, alors $X$ est le produit fibr\'e $X \times_{\overline{X}} \overline{X}''$.}

\vskip 3mm
\noindent \textit{D\'emonstration de (2.4)}. Soit $U''$ le produit fibr\'e

$$
\xymatrix{
X \ar@/^/[rrd]^{j''} \ar@/_/[rdd]_{\textrm{id}} \ar@{->}[rd]^{\varphi}\\
& U'' \ar[d]^{f} \ar@{^{(}->}[r]_{\overline{j}} & \overline{X}''  \ar[d]^{g}&\\
& X \ar@{^{(}->}[r]_{j} & \overline{X}&.
}
$$

\noindent Puisque $f \circ \varphi = \textrm{id}$, $\varphi$ est une immersion ferm\'ee [EGA I, (4.3.6) (iv)] ; or $\varphi$ est \'etale, car $\tilde{j} \circ \varphi = j''$ et $\tilde{j}$, $j''$ sont \'etales [EGA IV, (17.3.5)]. Ainsi $\varphi$ est une immersion ouverte [EGA IV, (17.9.1); EGA I, 4.2]. De plus $\varphi$ est dominante car $j''$ l'est ; donc $\varphi$ est surjective, car $\varphi$ est finie. Une immersion ouverte surjective est un isomorphisme. $\square$

\vskip 3mm
\noindent \textit{D\'emonstration de (2.3)}. Faisons la d\'emonstration dans le cas des gros topos syntomiques : pour les petits topos cela r\'esulte de l'\'egalit\'e 

$$
H^i(\overline{X}_{\textrm{SYNT}}, \mathcal{G}) = H^i(\overline{X}_{\textrm{synt}}, \mathcal{G})
$$

\noindent valable pour tout faisceau ab\'elien $\mathcal{G}$ sur synt$(X)$, et tout entier $i \geqslant 0$ [E-LS 2, 1.2] et [Mi, II, prop 3.1]. On d\'esigne par un "ET" en indice les gros topos \'etales [E-LS 2, \S\ 1].\\

On a un cube commutatif de morphismes de gros topos

$$
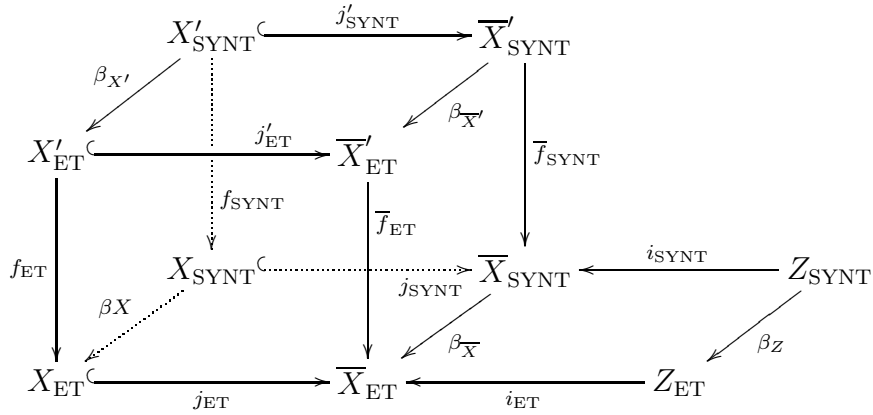
\xymatrix{
&X'_{\textrm{SYNT}} \ar@{^{(}->}[rr]^{j'_{\textrm{SYNT}}} \ar@{.>}[dd]^(.7){f_{\textrm{SYNT}}} |\hole \ar[dl]_{\beta_{X'}}  && \overline{X}'_{\textrm{SYNT}} \ar[dd]^{\overline{f}_{\textrm{SYNT}}}  \ar[dl]^{\beta_{\overline{X}'}} & & \\
 X'_{\textrm{ET}} \ar@{^{(}->}[rr]^(.7){j'_{\textrm{ET}}} \ar[dd]_{f_{\textrm{ET}}} && \overline{X}'_{\textrm{ET}}  \ar[dd]^(.3){\overline{f}_{\textrm{ET}}}  & &\\
& X_{\textrm{SYNT}} \ar@{^{(}.>}[rr]_(.7){j_{\textrm{SYNT}}} \ar@{.>}[dl]_{\beta{X}} && \overline{X}_{\textrm{SYNT}}   \ar[dl]^{\beta_{\overline{X}}} & & Z_{\textrm{SYNT}} \ar[dl]^{\beta_{Z}}  \ar[ll]_{i_{\textrm{SYNT}}} \\
X_{\textrm{ET}} \ar@{^{(}->}[rr]_{{j_{\textrm{ET}}}}  && \overline{X}_{\textrm{ET}} &&  Z_{\textrm{ET}} \ar[ll]^{i_{\textrm{ET}}}\\
}
$$

\noindent o\`u $Z = \overline{X} \setminus X$, et un isomorphisme

$$
R \Gamma(\overline{X}_{\textrm{SYNT}},\  j_{!}\ R f_{\textrm{SYNT}^{\ast}}\  \mathcal{F}) \simeq R \Gamma(\overline{X}_{\textrm{ET}}, R \beta_{\overline{X}^{\ast}}\  j_{\textrm{SYNT !}}\  R f_{\textrm{SYNT}^{\ast}}\  \mathcal{F}).
$$
\noindent Supposons \'etablie la proposition suivante :

\vskip 3mm
\noindent \textbf{Proposition (2.5)}. 
\textit{Avec les notations pr\'ec\'edentes, annelons $\overline{X}_{\textrm{SYNT}}$ par $\mathcal{A}$ et $\overline{X}_{\textrm{ET}}$ par $\mathcal{B} = \beta_{\overline{X}^{\ast}} \mathcal{A}$. Alors pour tout $\mathcal{H} \in {_{\mathcal{A} \vert X}X_{\textrm{SYNT}}}$ on a un isomorphisme}
$$
j_{\textrm{ET}!}\  R \beta_{X^{\ast}}(\mathcal{H}) \displaystyle \mathop{\rightarrow}^{\sim} R \beta_{\overline{X}^{\ast}}\  j_{\textrm{SYNT !}} (\mathcal{H}).
$$

Alors on a des isomorphismes \\

$R \Gamma(\overline{X}_{\textrm{SYNT}},\  j_{!}\ R f_{\textrm{SYNT}^{\ast}}\  \mathcal{F}) \simeq R \Gamma(\overline{X}_{\textrm{ET}}, j_{\textrm{ET}!}\ R \beta_{X^{\ast}}\  R f_{\textrm{SYNT}^{\ast}}\  \mathcal{F})$\\

$\quad \simeq R \Gamma(\overline{X}_{\textrm{ET}}, j_{\textrm{ET}!}\  R f_{\textrm{ET}^{\ast}}\ R \beta_{X'^{\ast}}\   \mathcal{F})$\\

$\quad \simeq R \Gamma(\overline{X}_{\textrm{ET}}, R \overline{f}_{\textrm{ET}^{\ast}}\ j'_{\textrm{ET}!}\   R \beta_{X'^{\ast}}\   \mathcal{F})$ [SGA 4, XVII, lemme 5.1.6] car $\mathcal{F}$ de torsion.\\

$\quad \simeq R \Gamma(\overline{X}_{\textrm{ET}}, R \overline{f}_{\textrm{ET}^{\ast}}\ R \beta_{\overline{X}'^{\ast}}\   j'_{\textrm{SYNT}!}\  \mathcal{F})$ [Prop (2.5)]\\

$\quad \simeq R \Gamma(\overline{X}_{\textrm{ET}},  R \beta_{\overline{X}^{\ast}}\   R \overline{f}_{\textrm{SYNT}^{\ast}}\ j'_{\textrm{SYNT}!}\  \mathcal{F})$ \\

$\quad \simeq R \Gamma(\overline{X}_{\textrm{SYNT}},  R \overline{f}_{\textrm{SYNT}^{\ast}}\ j'_{\textrm{SYNT}!}\  \mathcal{F})\ ;$ \\

\noindent d'o\`u la proposition (2.3).\\

\noindent\textit{Etablissons la proposition (2.5)}.\\
Puisque $j_{\textrm{ET}!}$ et $j_{\textrm{SYNT}!}$ sont exacts, il suffit de montrer l'isomorphisme \hfill\break $j_{\textrm{ET}!}\  \beta_{X^{\ast}} \displaystyle \mathop{\rightarrow}^{\sim} \beta_{\overline{X}^{\ast}} j_{\textrm{SYNT}!}$.\\
Par le lemme du serpent appliqu\'e au morphisme de suites exactes \\

$$
\begin{array}{c}
\xymatrix{
0 \ar[r] & j_{\textrm{ET}!} \beta_{X^{\ast}} \mathcal{H} \ar[r] \ar@{^{(}->}[ddd] & j_{\textrm{ET}^{\ast}} \beta_{X^{\ast}} \mathcal{H} \ar[r] \ar[ddd]^\simeq]& i_{\textrm{ET}^{\ast}} i^{\ast}_{\textrm{ET}} j_{\textrm{ET}^{\ast}} \beta_{X^{\ast}} \mathcal{H} \ar[r]  \ar[d]^\simeq & 0\\
& & &  i_{\textrm{ET}^{\ast}} i^{\ast}_{\textrm{ET}} \beta_{\overline{X}^{\ast}} j_{\textrm{SYNT}^{\ast}} \mathcal{H} \ar[d] &\\
& & & i_{\textrm{ET}^{\ast}} \beta_{Z^{\ast}} i_{\textrm{SYNT}}^{\ast} j_{\textrm{SYNT}^{\ast}} \mathcal{H} \ar[d]^\simeq &\\
0 \ar[r] & \beta_{\overline{X}^{\ast}} j_{\textrm{SYNT}!} \mathcal{H} \ar[r] & \beta_{\overline{X}^{\ast}} j_{\textrm{SYNT}^{\ast}} \mathcal{H} \ar[r]  & \beta_{\overline{X}^{\ast}} i_{\textrm{SYNT}^{\ast}}  i_{\textrm{SYNT}}^{\ast}    j_{\textrm{SYNT}^{\ast}} \mathcal{H}\ ,  & 
}
\end{array}
\leqno{(2.5.1)}
$$

\noindent il nous suffit de montrer que, pour tout faisceau $\mathcal{G}$, on a un isomorphisme

\noindent (2.5.2) $\qquad \qquad i_{\textrm{ET}^{\ast}}\  i_{\textrm{ET}}^{\ast}\  \beta_{\overline{X}^{\ast}}(\mathcal{G}) \displaystyle \mathop{\rightarrow}^{\varphi \atop{\sim}} i_{\textrm{ET}^{\ast}}\  \beta_{Z^{\ast}}\  i^{\ast}_{\textrm{SYNT}}(\mathcal{G}).$\\

\noindent Or $i_{\textrm{ET}^{\ast}}\ i_{\textrm{ET}}^{\ast}\ Ê\beta_{\overline{X}^{\ast}}(\mathcal{G})$ est le faisceau associ\'e au pr\'efaisceau\\

\vskip 4mm

$\qquad \qquad X' \mapsto i_{\textrm{ET}}^{\ast}\  \beta_{\overline{X}^{\ast}}(\mathcal{G}) (Z \times_{\overline{X}} \overline{X}') = \displaystyle \mathop{\lim}_{\rightarrow \atop{X''}} \beta_{\overline{X}^{\ast}}(\mathcal{G}) (X'')$

$\qquad \qquad  \qquad  \qquad \qquad \qquad \qquad \qquad  \quad = \displaystyle \mathop{\lim}_{\rightarrow \atop{X''}} \mathcal{G} (X'') \simeq \mathcal{G}(Z \times_{\overline{X}} \overline{X}')$\\

\noindent o\`u $\overline{X}'$ est un $\overline{X}$-sch\'ema et la limite inductive est prise sur les diagrammes commutatifs

$$
\xymatrix{
X''  \ar[d] & Z \times_{\overline{X}} \overline{X}' \ar[l]\ \ar[d]\\
\overline{X} & Z \ar[l]
}
$$

\noindent o\`u $X''$ est un $\overline{X}$-sch\'ema ; comme $i_{\textrm{ET}^{\ast}}\ \beta_{Z^{\ast}}\  i^{\ast}_{\textrm{SYNT}}(\mathcal{G})$ a la m\^eme description, il en r\'esulte que $\varphi$ est un isomorphisme. $\square$

\vskip 3mm
\noindent \textbf{Th\'eor\`eme (2.6)}. 
\textit{Soient $i_{1} : Z \hookrightarrow X$ une immersion ferm\'ee entre deux $k$-sch\'emas s\'epar\'es de type fini, $j_{1} : \cup \hookrightarrow X$ l'immersion ouverte du compl\'ementaire et $\mathcal{F} \in {_{\mathcal{A}}X_{\mathcal{T}}}$, o\`u $\mathcal{A}$ est un faisceau d'anneaux de torsion. Alors on a une suite exacte longue de cohomologie syntomique \`a supports compacts }\\

$ \rightarrow H^i_{\textrm{synt},c}(U, \mathcal{F}_{\vert U}) \rightarrow H^i_{\textrm{synt},c}(U, \mathcal{F}) \rightarrow H^i_{\textrm{synt},c}(Z, \mathcal{F}_{\vert Z}) \rightarrow H^{i+1}_{\textrm{synt},c}(U, \mathcal{F}_{\vert U}) \rightarrow $

\vskip 3mm
\noindent \textit{D\'emonstration}. Choisissons une compactification $\overline{X}$ de $X$ au-dessus de $k$, $\overline{j} : X \hookrightarrow \overline{X}$ l'immersion ouverte dominante et soit $\overline{Z}$ l'adh\'erence sch\'ematique de $Z$ dans $\overline{X}$. On a alors un diagramme commutatif \`a carr\'es cart\'esiens

$$
\begin{array}{c}
\xymatrix{
Z \ar@{^{(}->}[r]^{\overline{j}'} \ar@{^{(}->}[d]_{i_{1}} & \overline{Z} \ar@{^{(}->}[d]^{i}\\
X \ar@{^{(}->}[r]_{\overline{j}} & \overline{X}\\
U \ar@{=}[r] \ar@{^{(}->}[u]^{j_{1}} & U \ar@{^{(}->}[u]_{j}
}
\end{array}
\leqno{(2.6.1)}
$$

\noindent o\`u $i$ est une immersion ferm\'ee et $j$, $\overline{j}'$ des immersions ouvertes.\\
En appliquant le foncteur exact $\overline{j}_{!}$ \`a la suite exacte

$$
0 \longrightarrow j_{1 !}\  j^{\ast}_{1}\  \mathcal{F} \longrightarrow \mathcal{F} \longrightarrow\  i_{1^{\ast}}\  i^{\ast}_{1}\  \mathcal{F} \longrightarrow 0,
$$

\noindent on obtient la suite exacte\\

\noindent (2.6.2) $\qquad \qquad 0 \longrightarrow j_{!}\  j^{\ast}_{1}\   \mathcal{F} \longrightarrow \overline{j}_{!}\ \mathcal{F} \longrightarrow  \overline{j}_{!}\  i_{1^{\ast}}\  i^{\ast}_{1}\  \mathcal{F} \longrightarrow 0$\\

\noindent car $\overline{j}_{!}\ j_{1 !} = j_{!}$ [Th\'eo (1.5) (1)]. Or l'exactitude des foncteurs $i_{1 \ast}$ et $i_{\ast}$, jointe \`a la proposition (2.3), donne un isomorphisme

$$
R \Gamma(\overline{X}_{\mathcal{T}}, \overline{j}_{!}\ i_{1 \ast}\ i^{\ast}_{1}\ \mathcal{F} \displaystyle \mathop{\longrightarrow}^{\sim} R \Gamma (\overline{Z}_{\mathcal{T}}, \overline{j'}_{!}\ i^{\ast}_{1}\ \mathcal{F})\ ;
$$

\noindent ainsi, par application du foncteur $R \Gamma(\overline{X}_{\mathcal{T}}, -)$ \`a la suite exacte (2.6.2) on obtient, via (2.1), un triangle distingu\'e\\

\noindent (2.6.3) $\qquad R \Gamma_{\textrm{synt},c}(U, \mathcal{F}_{\vert U}) \longrightarrow R \Gamma_{\textrm{synt},c}(X, \mathcal{F}) \longrightarrow R \Gamma_{\textrm{synt},c}(Z, \mathcal{F}_{\vert Z}),$\\

\noindent qui fournit \`a son tour la suite exacte longue du th\'eor\`eme. $\square$\\

\vskip 3mm
\section*{3. Comparaison avec la cohomologie \'etale et la cohomologie rigide}

\textbf{3.0.} On suppose dans ce $\S$\  3 que le corps $k$ contient $\mathbb{F}_{q},\ q = p^a$. On d\'esigne par $C(k)$ un anneau de Cohen de $k$ de caract\'eristique 0 et par $K_{0}$ le corps des fractions de $C(k)$.\\

Comme en III (3.3), $\mathcal{V}$ est un anneau de valuation discr\`ete complet, d'uniformisante $\pi$ et $\sigma : \mathcal{V \rightarrow \mathcal{V}}$ un rel\`evement de la puissance $q$ de $k$ tel que $\sigma(\pi) = \pi$ construit via [Et 5, 1.1].\\

On note $e$ l'indice de ramification de $\mathcal{V}$, $K = \textrm{Frac}(\mathcal{V})$, $\mathcal{V}^{\sigma} = \textrm{Ker} \{1 - \sigma : \mathcal{V} \rightarrow \mathcal{V} \}, \mathcal{V}_{n} = \mathcal{V} / \pi^{n+1} \mathcal{V}, \mathcal{V}^{\sigma}_{n} = \mathcal{V}^{\sigma} / \pi^{n+1} \mathcal{V}^{\sigma}$ et $K^{\sigma} = \textrm{Frac} (\mathcal{V}^{\sigma})$.\\

Si $X$ est un sch\'ema on dit qu'un $\mathcal{V}^{\sigma}_{n}$-module $\mathcal{F}$ sur \'et$(X)$ est localement trivial (on dit aussi constant-tordu constructible ou encore localement constant constructible) s'il est localement isomorphe \`a une somme directe finie de copies de $\mathcal{V}^{\sigma}_{n}$ : c'est alors la m\^eme chose de dire qu'il est localement trivial sur SYNT$(X)$ [E-LS 2, 5.1].\\

Un $\mathcal{V}^{\sigma}$-faisceau lisse $\mathcal{F}$ (localement libre de rang fini) sur $X$ est un syst\`eme projectif $\mathcal{F} = (\mathcal{F}_{n})_{n \in \mathbb{N}}$ o\`u, pour tout $n$, $\mathcal{F}_{n}$ est un $\mathcal{V}^{\sigma}_{n}$-module localement trivial sur \'et$(X)$ et pour $n' \geqslant n$, $\mathcal{F}_{n} = \mathcal{F}_{n'} \otimes \mathcal{V}^{\sigma}_{n}$. Les $K^{\sigma}$-faisceaux lisses sont les $\mathcal{V}^{\sigma}$-faisceaux lisses \`a isog\'enie pr\`es.\\

Pour un $\mathcal{V}^{\sigma}$-faisceau lisse $\mathcal{F}$ (resp. un $K^{\sigma}$-faisceau lisse   $\mathcal{F_{\mathbb{Q}}}$) sur $X$ on pose, pour tout entier $i \geqslant 0$

$$
H^i_{\textrm{synt},c}(X, \mathcal{F}) := \displaystyle \mathop{\lim}_{\leftarrow \atop_{n}} H^i_{\textrm{synt},c}(X, \mathcal{F}_{n})
$$\\

[resp. $ \qquad \qquad  H^i_{\textrm{synt},c}(X, \mathcal{F}_{\mathbb{Q}}) = (\displaystyle \mathop{\lim}_{\leftarrow \atop_{n}} H^i_{\textrm{synt},c}(X, \mathcal{F}_{n})) \otimes_{\mathcal{V}^{\sigma}} K^{\sigma}]\ .
$

\vskip 3mm
\noindent \textbf{Proposition (3.1)}. 
\textit{Soient $X$ un sch\'ema et $\mathcal{F}$ un faisceau ab\'elien sur \'et$(X)$ ; on note encore $\mathcal{F}$ son image inverse par le morphisme de topos $X_{\textrm{synt}} \rightarrow X_{\textrm{\'et}}$. Alors, pour tout entier $i \geqslant 0$, on a un isomorphisme}

$$
H^i_{\textrm{\'et},c} (X, \mathcal{F}) \displaystyle \mathop{\longrightarrow}^{\sim} H^i_{\textrm{synt},c}(X, \mathcal{F}).
$$

\vskip 3mm
\noindent \textit{D\'emonstration}. R\'esulte de [E-LS 2, 1.3]. \ $\square$

\vskip 3mm
\noindent \textbf{Th\'eor\`eme (3.2)}. 
\textit{Supposons $k$ s\'eparablement clos. Soient $X$ un $k$-sch\'ema s\'epar\'e de type fini, $\mathcal{F}_{\mathbb{Q}}$ un $K^{\sigma}$-faisceau lisse sur $X$ et $f : Y \rightarrow X$ sur $k$-morphisme fini \'etale galoisien de groupe $G$. Alors, pour tout entier $i \geqslant 0$, on a des isomorphismes de $K^{\sigma}$-espaces vectoriels de dimension finie}

$$
H^i_{\textrm{\'et},c}(X, \mathcal{F}_{\mathbb{Q}}) \simeq H^i_{\textrm{\'et},c}(X, f_{\ast} f^{\ast} \mathcal{F}_{\mathbb{Q}})^G \simeq H^i_{\textrm{\'et},c}(Y, f^{\ast} \mathcal{F}_{\mathbb{Q}})^G
$$

$$
\simeq H^i_{\textrm{synt},c}(X, \mathcal{F}_{\mathbb{Q}}) \simeq H^i_{\textrm{synt},c} (X, f_{\ast} f^{\ast} \mathcal{F}_{\mathbb{Q}})^G \simeq H^i_{\textrm{synt},c}(Y, f^{\ast} \mathcal{F}_{\mathbb{Q}})^G.
$$

\vskip 3mm
\noindent \textit{D\'emonstration}. Compte tenu de (3.1) il suffit de montrer l'assertion pour la cohomologie \'etale. Puisque $f$ est fini, $f_{\ast}$ est exact, et on est ramen\'e \`a montrer l'isomorphisme\\

\noindent (3.2.1) $\qquad  H^i_{\textrm{\'et},c}(X, \mathcal{F}) \otimes_{\mathcal{V}^{\sigma}} K^{\sigma} \simeq [H^i_{\textrm{\'et},c}(X, f_{\ast} f^{\ast} \mathcal{F}) \otimes_{\mathcal{V}^{\sigma}} K^{\sigma}]^G $ \\

\noindent pour un $\mathcal{V}^{\sigma}$-faisceau lisse $\mathcal{F}$.\\
Comme $\mathcal{F}$ est localement libre on \'etablit le lemme suivant comme [Et 2, III (3.1.2)].

\vskip 3mm
\noindent \textbf{Lemme (3.2.2)}. 
\textit{Sous les hypoth\`eses pr\'ec\'edentes, on a un isomorphisme}

$$
\mathcal{F} \displaystyle \mathop{\longrightarrow}^{\sim} (f_{\ast} f^{\ast}(\mathcal{F}))^G.
$$    

Soient $\overline{X}$ une compactification de $X$ au-dessus de $k$ et $j : X \hookrightarrow \overline{X}$ l'immersion ouverte correspondante. Par exactitude du foncteur $j_{!}$ on d\'eduit de (3.2.2) des isomorphismes\\

\noindent (3.2.3) $\qquad \qquad  j_{!}   \mathcal{F}   \displaystyle \mathop{\longrightarrow}^{\sim} j_{!} (f_{\ast} f^{\ast} \mathcal{F})^G \simeq (j_{!} f_{\ast} f^{\ast} \mathcal{F})^G.$ \\

\noindent Le corps $k$ \'etant s\'eparablement clos, les groupes $H^i_{\textrm{\'et},c}(X, \mathcal{F}_{n})$ et   $H^i_{\textrm{\'et},c}(X, f_{\ast} f^{\ast}  \mathcal{F}_{n})$ sont des $\mathcal{V}^{\sigma}_{n}$-modules de type fini [SGA 4, XVII, 5.3.8] ; par suite les groupes $H^i_{\textrm{\'et},c}(X, \mathcal{F})$ et    $H^i_{\textrm{\'et},c}(X,  f_{\ast} f^{\ast}  \mathcal{F})$ sont des $\mathcal{V}^{\sigma}$-modules de type fini, engendr\'es par tout sous-ensemble qui les engendre mod.$\pi$. \\

On ach\`eve la d\'emonstration de (3.2) comme [Et 2, III, 3.1.1]. $\square$  

\vskip 3mm
\textbf{3.3.} Soient $X$ un $k$-sch\'ema et $\mathcal{H}$ un $\mathcal{V}^{\sigma}_{n}$-module localement trivial sur SYNT$(X)$. On consid\`ere le morphisme de topos annel\'es [E-LS 2, 5.3]

$$
u = u^{(m)}_{X/\mathcal{V}_{n}-\textrm{SYNT}} : ((X/\mathcal{V}_{n})^{(m)}_{\textrm{CRIS-SYNT}}, \mathcal{O}^{(m)}_{X/\mathcal{V}_{n}}) \longrightarrow (X_{\textrm{SYNT}}, \mathcal{V}^{\sigma}_{n})
$$

\noindent et on note

$$
T^{(m)}(\mathcal{H}) := u^{\ast}(\mathcal{H}) = \mathcal{H} \otimes_{\mathcal{V}^{\sigma}_{n}} \mathcal{O}^{(m)}_{X/\mathcal{V}_{n}},
$$

\noindent et
$$T^{(m)}(\mathcal{H})^{\textrm{cris}} := u_{\ast} (T^{(m)}(\mathcal{H})) = u_{\ast} u^{\ast}(\mathcal{H})\ [E\mbox{-}LS \2, 1.11]$$

\noindent (3.3.1) $\qquad \qquad \qquad = \mathcal{H}  \otimes_{\mathcal{V}^{\sigma}_{n}} u_{\ast}({\mathcal{O}}^{(m)}_{X/\mathcal{V}_{n}}) = \mathcal{H}  \otimes_{\mathcal{V}^{\sigma}_{n}} \mathcal{O}^{m-\textrm{cris}}_{n, X},$\\

\noindent o\`u l'on a pos\'e $\mathcal{O}^{m-\textrm{cris}}_{n, X} := u_{\ast}(\mathcal{O}^{(m)}_{X/\mathcal{V}_{n}}).$\\

Soient $i : Z \hookrightarrow X$  une immersion  ferm\'ee et $j : U \hookrightarrow X$ l'immersion ouverte du compl\'ementaire : $j$ et $i$ d\'efinissent respectivement des morphismes de topos [cf (1.1)]

$$
j : U_{\textrm{SYNT}} \longrightarrow X_{\textrm{SYNT}}\ ,
$$

$$
i : Z_{\textrm{SYNT}} \longrightarrow X_{\textrm{SYNT}}\ ,
$$

\noindent et m\^eme des morphismes de topos annel\'es\\

\noindent (3.3.2) $\qquad \qquad \qquad j_{\mathcal{V}} : (U_{\textrm{SYNT}}, \mathcal{V}^{\sigma}_{n,U}) \longrightarrow (X_{\textrm{SYNT}}, \mathcal{V}^{\sigma}_{n,X})$\\

\noindent (3.3.3) $\qquad \qquad \qquad i_{\mathcal{V}} : (Z_{\textrm{SYNT}}, \mathcal{V}^{\sigma}_{n,Z}) \longrightarrow (X_{\textrm{SYNT}}, \mathcal{V}^{\sigma}_{n,X})$\\

\noindent o\`u $\mathcal{V}^{\sigma}_{n,U} = j^{-1}(\mathcal{V}^{\sigma}_{n,X}) \quad, \quad \mathcal{V}^{\sigma}_{n,Z} = i^{-1}(\mathcal{V}^{\sigma}_{n,X}),$\\

\noindent $\mathcal{V}^{\sigma}_{n,X}$ \'etant le faisceau d'anneaux $\mathcal{V}^{\sigma}_{n}$ sur $X_{\textrm{SYNT}}.$\:

On d\'eduit alors de (1.4) l'existence de six foncteurs

$$
\begin{array}{c}
\xymatrix{
& \ar[l]_{i^{\ast}_{\mathcal{V}}} \qquad & \qquad  \ar[l]_{j_{\mathcal{V}!}} &\\
\ _{\mathcal{V}^{\sigma}_{n,Z} }Z_{\textrm{SYNT}} \ar[r]^{i_{\mathcal{V}^{\ast}}} & \ _{\mathcal{V}^{\sigma}_{n,X} }X_{\textrm{SYNT}} \ar[r]^{j^{\ast}_{\mathcal{V}}} & \  _{\mathcal{V}^{\sigma}_{n,U} }U_{\textrm{SYNT}}\\
 &\  \ar[l]_{i^{!}_{\mathcal{V}}} \  &  \  \ar[l]_{j_{\mathcal{V}^{\ast}}}  \ &
}
\end{array}
\leqno{(3.3.4)}
$$

\vskip 3mm
\noindent \textbf{Lemme (3.3.5)}. 
\textit{Sous les hypoth\`eses pr\'ec\'edentes on a des isomorphismes de faisceaux sur les gros sites syntomiques :}\\

\noindent (3.3.5.1) $\qquad \qquad j^{-1} \mathcal{O}^{m-\textrm{cris}}_{n,X} \displaystyle \mathop{\longrightarrow}^{\sim} \mathcal{O}^{m-\textrm{cris}}_{n,U}$\\

\noindent (3.3.5.2) $\qquad \qquad i^{-1} \mathcal{O}^{m-\textrm{cris}}_{n,X} \displaystyle \mathop{\longrightarrow}^{\sim} \mathcal{O}^{m-\textrm{cris}}_{n,Z}.$\\

\vskip 3mm
\noindent \textit{D\'emonstration}. Comme $j^{-1} \mathcal{O}^{m-\textrm{cris}}_{n,X}$ est le faisceau associ\'e au pr\'efaisceau qui \`a tout $U$-sch\'ema $U'$ associe\\
$$
\mathcal{O}^{m-\textrm{cris}}_{n,X}(U') = \Gamma((U'/\mathcal{V}_{n})^{(m)}_{\textrm{CRIS-SYNT}}, \mathcal{O}^{(m)}_{U'/\mathcal{V}n})\  [\textrm{E\mbox{-}LS\2, 1.10}] 
$$
$\qquad \qquad \qquad \qquad \qquad  = \mathcal{O}^{m-\textrm{cris}}_{n,U}(U'),$\\

\noindent on a bien (3.3.5.1). De m\^eme pour (3.3.5.2). $\square$\\

Gr\^ace au lemme (3.3.5) les morphismes de topos $j$ et $i$ ci-dessus induisent des morphismes de topos annel\'es\\

\noindent (3.3.6) $\qquad \qquad j_{\mathcal{O}} : (U_{\textrm{SYNT}}, \mathcal{O}^{m-\textrm{cris}}_{n,U}) \longrightarrow (X_{\textrm{SYNT}}, \mathcal{O}^{m-\textrm{cris}}_{n,X})$\\

\noindent (3.3.7) $\qquad \qquad i_{\mathcal{O}} : (Z_{\textrm{SYNT}}, \mathcal{O}^{m-\textrm{cris}}_{n,Z}) \longrightarrow (X_{\textrm{SYNT}}, \mathcal{O}^{m-\textrm{cris}}_{n,X})$\\

\noindent et six foncteurs

$$
\begin{array}{c}
\xymatrix{
& \ar[l]_{i^{\ast}_{\mathcal{O}}} \qquad & \qquad  \ar[l]_{j_{\mathcal{O}!}} &  {} &\\
\  _{\mathcal{O}^{m-\textrm{cris}}_{n,Z} } Z_{\textrm{SYNT}} \ar[r]^{i_{\mathcal{O}^{\ast}}} & \  _{\mathcal{O}^{m-\textrm{cris}}_{n,X} } X_{\textrm{SYNT}} \ar[r]^{j^{\ast}_{\mathcal{O}}} & \  _{\mathcal{O}^{m-\textrm{cris}}_{n,U} } U_{\textrm{SYNT}}&\\
 &\  \ar[l]_{i^{!}_{\mathcal{O}}} \ &\   \ar[l]_{j_{\mathcal{O}^{\ast}}}\ & .
}
\end{array}
\leqno{(3.3.8)}
$$
    
En utilisant la description (1.4), ou le th\'eor\`eme (1.5) \textit{(6.1)} on en d\'eduit la proposition suivante :    

\vskip 3mm
\noindent \textbf{Proposition (3.3.9)}. 
\textit{Sous les hypoth\`eses et notations de (3.3) on a :}
\begin{itemize}
\item[(1)] \textit{Si $\mathcal{G}$ est un $\mathcal{O}^{m-\textrm{cris}}_{n,U}$-module, alors on a un isomorphisme canonique}\\
$$j_{\mathcal{V}!}(\mathcal{G}) \displaystyle \mathop{\rightarrow}^{\sim} j_{\mathcal{O}!}(\mathcal{G}).$$
\item[(2)] \textit{Si $\mathcal{G}$ est un $\mathcal{V}^{\sigma}_{n,U}$-module, alors on a des isomorphismes canoniques}
$$j_{\mathcal{V}!}(\mathcal{G} \otimes_{\mathcal{V}^{\sigma}_{n,U}} \mathcal{O}^{m-\textrm{cris}}_{n,U}) \simeq j_{\mathcal{O}!} (\mathcal{G} \otimes_{\mathcal{V}^{\sigma}_{n,U}} \mathcal{O}^{m-\textrm{cris}}_{n,U})$$
$\qquad \qquad \simeq j_{\mathcal{V}!}(\mathcal{G}) \otimes_{\mathcal{V}^{\sigma}_{n,X}} j_{\mathcal{V}!} (\mathcal{O}^{m-\textrm{cris}}_{n,U}) \simeq j_{\mathcal{V}!} (\mathcal{G}) \otimes_{\mathcal{V}^{\sigma}_{n,X}} j_{\mathcal{O}!} (\mathcal{O}^{m-\textrm{cris}}_{n,U})$\\

$\qquad \qquad \simeq  j_{\mathcal{V}!}(\mathcal{G}) \otimes_{\mathcal{V}^{\sigma}_{n,X}} \mathcal{O}^{m-\textrm{cris}}_{n,X}.$

\end{itemize} 

\vskip 3mm
Les formules du (2) sont \`a comparer \`a celles de [SGA 4, IV, prop 12.11 (b)]. \\

Supposons \`a pr\'esent que $\mathcal{F}$ est un $\mathcal{V}^{\sigma}_{n}$-module localement trivial sur SYNT$(U) : \mathcal{F}$ \'etant localement trivial, le morphisme $F^{\ast} : \mathcal{F} \rightarrow \mathcal{F}^{(q)} = F^{-1}_{X}(\mathcal{F})$ est un isomorphisme ; on pose $F = F^{\ast-1} : \mathcal{F}^{(q)} \displaystyle \mathop{\rightarrow}^{\sim} \mathcal{F}$ et $\phi _{U}= T^{(m)}(F): T^{(m)}(\mathcal{F})^{(q)}=T^{(m)}(\mathcal{F}^{(q)})\displaystyle \mathop{\rightarrow}^{\sim} T^{(m)}(\mathcal{F}) $ qui munit $T^{(m)}(\mathcal{F})$ d'une structure de $F\mbox{-}m$-cristal localement trivial [E-LS 2, 5.3]. \\
\noindent

On note encore $\phi_{U} : T^{(m)}(\mathcal{F})^{\textrm{cris}} \rightarrow T^{(m)}(\mathcal{F})^{\textrm{cris}}$ l'homomorphisme obtenu en composant $\phi_{U}^{\textrm{cris}}$ avec  $F^{\ast} :  T^{(m)}(\mathcal{F})^{\textrm{cris}}   \rightarrow T^{(m)}(\mathcal{F})^{\textrm{cris}(q)}$ [E-LS 2, 5.2]. Alors la suite exacte de [E-LS 2, th\'eo 5.5] s'interpr\`ete, via (3.3.1), comme une suite exacte sur SYNT$(U)$ :  \\

\noindent (3.3.10) $\qquad \qquad 0 \longrightarrow \mathcal{F} \longrightarrow u_{\ast}\  u^{\ast}(\mathcal{F}) \displaystyle \mathop{\longrightarrow} _{1-\phi_{U}} u_{\ast}\ u^{\ast}(\mathcal{F}) \longrightarrow 0$ \\

\noindent ou encore\\

\noindent (3.3.11) $\qquad \quad 0 \longrightarrow \mathcal{F} \longrightarrow \mathcal{F} \otimes_{\mathcal{V}^{\sigma}_{n}}\  \mathcal{O}^{m-\textrm{cris}}_{n,U} \displaystyle \mathop{\longrightarrow} _{1-\phi_{U}}  \mathcal{F} \otimes_{\mathcal{V}^{\sigma}_{n}}\  \mathcal{O}^{m-\textrm{cris}}_{n,U} \longrightarrow 0$\ .\\

\noindent En lui appliquant le foncteur exact $j_{\mathcal{V}!}$ [th\'eo 1.5], on obtient encore une suite exacte, ce qui, compte tenu de (3.3.9), \'etablit le th\'eor\`eme suivant :

\vskip 3mm
\noindent \textbf{Th\'eor\`eme (3.3.12)}. 
\textit{Sous les notations de (3.3), si $\mathcal{F}$ est un $\mathcal{V}^{\sigma}_{n}$-module localement trivial sur SYNT$(U)$, alors on a des suites exactes de $\mathcal{V}^{\sigma}_{n}$-modules sur SYNT$(X)$ :}

$$
\xymatrix{
0 \ar[r] & j_{\mathcal{V}!}\ \mathcal{F} \ar[r] \ar@{=}[d]& j_{\mathcal{O}!}(\mathcal{F} \otimes_{\mathcal{V}^{\sigma}_{n}} \mathcal{O}^{m-\textrm{cris}}_{n,U}) \ar[r]_{1-\phi} \ar[d]^{\simeq} & j_{\mathcal{O}!} (\mathcal{F} \otimes_{\mathcal{V}^{\sigma}_{n}} \mathcal{O}^{m-\textrm{cris}}_{n,U})\ar[r] &  0\\
0 \ar[r] & j_{\mathcal{V}!}\ \mathcal{F} \ar[r] & j_{\mathcal{V}!}(\mathcal{F}) \otimes_{\mathcal{V}^{\sigma}_{n}} \mathcal{O}^{m-\textrm{cris}}_{n,X} \ar[r]_{1-\phi} & j_{\mathcal{V}!} (\mathcal{F}) \otimes_{\mathcal{V}^{\sigma}_{n}} \mathcal{O}^{m-\textrm{cris}}_{n,X} \ar[r] &  0
}
$$

\noindent \textit{o\`u $\phi = j_{\mathcal{V}!} (\phi_{U})$.}\\

Les suites exactes du th\'eor\`eme (3.3.12) vont nous permettre en passant \`a la cohomologie dans le th\'eor\`eme suivant de relier cohomologie \'etale et cohomologie rigide.\\

\noindent \textbf{Th\'eor\`eme (3.3.13)}. 
\textit{On suppose le corps $k$ parfait. Soient $X$ un $k$-sch\'ema s\'epar\'e de type fini, $\mathcal{F}_{\mathbb{Q}}$ un $K^{\sigma}$-faisceau lisse sur $X$ et $E_{K} \in F^a\mbox{-}\textrm{Isoc}(X/K)^{\circ}$ le $F$-isocristal convergent associ\'e \`a $\mathcal{F}_{\mathbb{Q}}$ [E-LS 2 ; 5.6) et on suppose que $E_{K}$ provient de $E^{\dag}_{K} \in F^a\mbox{-}\textrm{Isoc}^{\dag}(X/K)^{\circ}$ par le foncteur d'oubli $F^a\mbox{-}\textrm{Isoc}^{\dag}(X/K)^{\circ} \rightarrow F^a\mbox{-}\textrm{Isoc}(X/K)^{\circ}$. On note $\mathcal{F} = (\mathcal{F}_{n})_{n \in \mathbb{N}}$ un $\mathcal{V}^{\sigma}$-faisceau lisse associ\'e \`a $\mathcal{F}_{\mathbb{Q}}$ et $E^{m-\textrm{cris}}_{n} = T^{(m)}(\mathcal{F}_{n})^{\textrm{cris}}$ [cf (3.3.1)].\\
Alors on a :}
\begin{itemize}
\item[(1)] \textit{Il existe un isomorphisme canonique}
$$
R \Gamma_{\textrm{rig},c}(X/K, E^{\dag}_{K}) = R \displaystyle \mathop{\lim}_{\leftarrow \atop{m}} [(R  \displaystyle \mathop{\lim}_{\leftarrow \atop{n}} R\Gamma_{\textrm{synt},c}(X, E_{n}^{m-\textrm{cris}})) \otimes \mathbb{Q}].
$$
\item[(2)] \textit{Si de plus $k$ est alg\'ebriquement clos, il existe, pour tout entier $i \geqslant 0$, une suite exacte courte}
$$
O \rightarrow H^i_{\textrm{\'et},c}(X, \mathcal{F}_{\mathbb{Q}}) \rightarrow H^i_{\textrm{rig},c}(X/K,E^{\dag}_{K}) \displaystyle \mathop{\longrightarrow}_{1-\phi} H^i_{\textrm{rig},c}(X/K, E^{\dag}_{K}) \rightarrow 0.$$
\end{itemize}

\vskip 3mm
\noindent \textit{D\'emonstration}.\\
\textit{Prouvons le (1)}. Comme la cohomologie rigide \`a supports compacts ne d\'epend que du sch\'ema r\'eduit sous-jacent \`a $X$, on peut supposer $X$ r\'eduit. On va faire une r\'ecurrence sur la dimension de $X$.\\

Si dim $X = 0$, alors $X = \displaystyle \mathop{\bigcup}_{\textrm{finie}} \textrm{Spec} A_{i}$ o\`u $A_{i}$ est artinien [Eis., cor 9.1] car $X$ est de type fini sur $k$, et $A_{i}$ est un produit fini $\displaystyle \mathop \Pi_{j} A_{i,m_{j}}$ d'anneaux artiniens locaux r\'eduits ($X$ est r\'eduit) : ainsi $A_{i,m_{j}}$ est un corps [Bour, A VIII, $\S$\ 6, \no4, prop 9] $k_{ij}$ extension finie du corps parfait $k$ [Eis, cor 2.15]. En particulier $X$ est fini \'etale sur $k$, et alors l'assertion du th\'eor\`eme est prouv\'ee dans [E-LS 2, prop 3.11].\\

Si dim $X \geqslant 1$, il existe, puisque $k$ est parfait et $X$ r\'eduit, un ouvert non vide $U \hookrightarrow X$ qui est lisse sur $k$, de ferm\'e compl\'ementaire $Z$ tel que dim $Z < \textrm{dim}\ X$. Comme les deux foncteurs $R \Gamma_{\textrm{rig},c}(X/K, -)$ et $R \displaystyle \mathop{\lim}_{\leftarrow \atop{m}} \ [(R \displaystyle \mathop{\lim}_{\leftarrow \atop{n}} R\Gamma_{\textrm{synt},c}(X, (-)_{n}^{m-\textrm{cris}})) \otimes \mathbb{Q})]$ donnent lieu \`a des triangles distingu\'es faisant intervenir $X$, $U$ et $Z$ on est ramen\'e \`a prouver le th\'eor\`eme pour $U$. Donc on peut supposer $X$ lisse connexe et aussi $K' = K$ avec $e \leqslant p-1$, avec $\mathcal{F}$ un $\mathcal{V}^{\sigma}$-faisceau lisse sur $X$, associ\'e \`a un $F$-cristal unit\'e $E$ sur $X/\mathcal{V}$, d'isocristal convergent unit\'e $E_{K}$ par [B 3, (2.4.2)] suppos\'e provenir de $E^{\dag}_{K} \in F^a\mbox{-}\textrm{Isoc}(X/K)^{\circ}$.\\

Notons $\overline{X}$ une compactification de $X$ sur $k$. D'apr\`es le th\'eor\`eme de monodromie finie ``g\'en\'erique'' de Tsuzuki [Tsu 2, theo 3.1] il existe un $k$-sch\'ema projectif et lisse $\overline{X}'$, un $k$-morphisme propre surjectif $\overline{w} : \overline{X'} \rightarrow \overline{X}$ g\'en\'eriquement \'etale, tel qu'en posant $X' = \overline{w}^{-1}(X)$, $j' : X' \hookrightarrow \overline{X'}$ l'immersion ouverte, il existe un unique $N^{\dag} \in F^a\mbox{-}\textrm{Isoc}^{\dag}(\overline{X'}/K)^{\circ}$ avec $\overline{w}^{\ast}(E^{\dag}_{K}) \simeq (j')^{\dag}(N^{\dag})$.\\

\noindent Notons $U \hookrightarrow X$ un ouvert dense tel que la restriction $w : U' = U \times_{\overline{X}} \overline{X'} \rightarrow U$ de $\overline{w}$ soit finie \'etale : quitte \`a r\'etr\'ecir $U$ on peut supposer $U$ affine et int\'egralement clos, de m\^eme pour $U'$. Puisque $U$ et $U'$ sont connexes il existe un morphisme fini \'etale $s : U'' \rightarrow U'$ tel que le compos\'e $f : w \circ s : U'' \rightarrow U$ soit fini \'etale galoisien de groupe not\'e $G$ [Mi, I, Rk 5.4]. D\'esignons par $\overline{X''}$ la fermeture int\'egrale de $\overline{X'}$ dans $U''$ : on obtient un diagramme commutatif \`a carr\'es cart\'esiens

$$
\xymatrix{
U'' \ar@{^{(}->}[rr]^{j''} \ar@/_2pc/[dd]_{f} \ar[d]_{s}  &  & \overline{X''} \ar[d]^{\overline{s}}\\
U' \ar@{^{(}->}[r]  \ar[d]_{w} & X' \ar@{^{(}->}[r] \ar[d] & \overline{X'} \ar[d]^{\overline{w}}\\
U \ar@{^{(}->}[r]  \ar@/_1pc/[rr]_{j_{U}} & X \ar@{^{(}->}[r] ^{j} & \overline{X}
}
$$

\noindent o\`u $\overline{s}$ est un morphisme fini et les fl\`eches horizontales sont des immersions ouvertes : en particulier $\overline{X''}$ est une compactification de $U''$.
On note $\overline{Z}$ l'adh\'erence sch\'ematique de $Z$ dans $\overline{X}$ et $j_{Z} : Z \hookrightarrow \overline{Z}$ l'immersion ouverte dominante.\\

Puisque $\overline{X'}$ est propre et lisse sur $k$, il existe, d'apr\`es (3.3.13) un $F$-cristal unit\'e $M$ sur $\overline{X'}$ tel que $N^{\dag} \simeq M^{an}$. Posons $\overline{\mathcal{M}} = \overline{s}^{\ast}_{\textrm{CRIS}}(M)$, $\mathcal{M} = j''^{\ast}_{\textrm{CRIS}}(\overline{\mathcal{M}})$ ; d'apr\`es[B 3, (2.42)] $\mathcal{M}$ est isog\`ene \`a $f^{\ast}_{\textrm{CRIS}}(E_{\vert U})$, i.e. il existe un entier $r \geqslant 0$ et des morphismes

$$
\alpha : \mathcal{M} \longrightarrow f^{\ast}_{\textrm{CRIS}}(E_{\vert U}),
$$

$$
\beta : f^{\ast}_{\textrm{CRIS}}(E_{\vert U}) \longrightarrow \mathcal{M},
$$

\noindent tels que $\alpha \circ \beta = p^r$ et $\beta \circ \alpha = p^r$. On notera $\overline{\mathcal{M}}_{K}$ le $F$-isocristal (sur)convergent sur $\overline{X''}$ associ\'e \`a $\overline{\mathcal{M}}$ et $\mathcal{M}^{\dag}_{K} = j''^{\dag}(\overline{\mathcal{M}}_{K})$. D'apr\`es [B-M 1, cor du th\'eo 6] le $F$-cristal unit\'e $\overline{\mathcal{M}}$ est le cristal de Dieudonn\'e d'un groupe $p$-divisible \'etale, dont le $\mathcal{V}^{\sigma}$-faisceau lisse associ\'e sera not\'e $\overline{\mathcal{G}} = (\overline{\mathcal{G}}_{n})_{n \in \mathbb{N}}$ ; on pose $\mathcal{G} = j''^{\ast}(\overline{\mathcal{G}})$.

L'isog\'enie $\alpha$ (resp $\beta$) fournit une isog\'enie $\alpha_{\mathcal{G}} : \mathcal{G} \rightarrow f^{\ast}(\mathcal{F}_{\vert U})$ (resp $\beta_{\mathcal{G}} : f^{\ast}(\mathcal{F}_{\vert U}) \rightarrow \mathcal{G})$ telle que $\alpha_{\mathcal{G}} \circ \beta_{\mathcal{G}} = p^r$ et $\beta_{\mathcal{G}} \circ \alpha_{\mathcal{G}} = p^r$.

L'isog\'enie $\alpha$ (resp $\beta$) fournit, par la construction de Berthelot [B 3, (2.4.2)], un isomorphisme sur les $F$-isocristaux convergents associ\'es

$$
\alpha_{K} : \mathcal{M}_{K} \displaystyle \mathop{\longrightarrow}^{\sim} f^{\ast}_{\textrm{rig}}(E_{K \vert U})
$$

$$
(\textrm{resp}\  \beta_{K} : f^{\ast}_{\textrm{rig}}(E_{K \vert U}) \displaystyle \mathop{\longrightarrow}^{\sim} \mathcal{M}_{K}).
$$
\noindent D'apr\`es [Et 5, th\'eo 5] l'isomorphisme $\alpha_{K}$ (resp. $\beta_{K}$) se rel\`eve de mani\`ere unique en un isomorphisme
$$
\alpha_{K}^{\dag} : \mathcal{M}_{K}^{\dag} \displaystyle \mathop{\longrightarrow}^{\sim} f^{\ast}_{\textrm{rig}}(E_{K \vert U}^{\dag})
$$

$$
(\textrm{resp}\  \beta_{K}^{\dag} : f^{\ast}_{\textrm{rig}}(E_{K \vert U}^{\dag}) \displaystyle \mathop{\longrightarrow}^{\sim} \mathcal{M}_{K}^{\dag});
$$

\noindent de m\^eme l'action de $G$ sur $f^{\ast}_{\textrm{rig}}(E_{K \vert U})$ se rel\`eve de mani\`ere unique \`a $f^{\ast}_{\textrm{rig}}(E_{K \vert U}^{\dag})$.\\

D'autre part, par le th\'eor\`eme (3.3.12), on a un morphisme de suites exactes sur SYNT$(\overline{X''})$

$$
\xymatrix{
 (S_{1})\ 0 \ar[r]  & j_{!} f^{\ast}(\mathcal{F}_{n \vert U}) \ar[r] \ar[d]_{{j''_{!}}(\beta_{\mathcal{G}})_{n}} & j_{!} f^{\ast}_{\textrm{CRIS}}(E_{\vert U})^{m-\textrm{cris}}_{n} \ar[r]^{1-\phi} \ar[d]^{{j''_{!}}(\beta)^m_{n}}  &  j_{!} f^{\ast}_{\textrm{CRIS}}(E_{\vert U})^{m-\textrm{cris}}_{n} \ar[r] \ar[d]^{{j''_{!}}(\beta)^m_{n}} & 0 \\
  (S_{2})\ 0 \ar[r]  & j''_{!}  \mathcal{G}_{n} \ar[r]& j''_{!} \mathcal{M}^{m-\textrm{cris}}_{n} \ar[r]^{1-\phi} &  j''_{!} \mathcal{M}^{m-\textrm{cris}}_{n} \ar[r]  & 0
}
$$

Le groupe $G$ agit de mani\`ere \'equivariante sur la suite exacte $(S_{1})$, donc sur le triangle distingu\'e $\mathcal{C}''(S_{1})$ obtenu en lui appliquant le foncteur

$$
\mathcal{C}'' := R \displaystyle \mathop{\lim}_{\leftarrow \atop{m}} \{(R  \displaystyle \mathop{\lim}_{\leftarrow \atop{n}} R \Gamma(\overline{X''}_{\textrm{SYNT},} -)) \otimes \mathbb{Q} \}.
$$\\

Comme $\mathcal{C}''(j''_{!}(\beta)^m_{n}) =: \tilde{\beta}$ est un isomorphisme, et de m\^eme en rempla\c{c}ant $\beta$ par $\alpha$, ou $\alpha_{\mathcal{G}}$, $\beta_{\mathcal{G}}$, le triangle distingu\'e $\mathcal{C}''(S_{2})$ obtenu en appliquant $\mathcal{C}''$ \`a $(S_{2})$ est aussi $G$-\'equivariant par transport de structure par ces isomorphismes : compte tenu de [E-LS 2, (3.11)] ce morphisme de triangles s'identifie \`a\\

$$
\xymatrix{
(3.3.13.1)\ R \Gamma_{\textrm{\'et},c}(U'',f^{\ast}(\mathcal{F}_{\vert U})) \otimes \mathbb{Q} \ar[r]  \ar[d] _{R \Gamma_{\textrm{\'et},c} (\beta_{\mathcal{G}}) \otimes \mathbb{Q}}^{\simeq} & R \Gamma_{\textrm{rig},c}(U'',f^{\ast}_{\textrm{rig}}(E^{\dag}_{K_{\vert U}})) \ar[r]^{1-\phi} \ar[d]_{R \Gamma_{\textrm{rig},c} (\beta^{\dag}_{K})}^{\simeq} & R \Gamma_{\textrm{rig},c}(U'',f^{\ast}_{\textrm{rig}}(E^{\dag}_{K_{\vert U}})) \ar[d]_{R \Gamma_{\textrm{rig},c} (\beta^{\dag}_{K})}^{\simeq}\\
(3.3.13.2)\qquad R \Gamma_{\textrm{\'et},c}(U'',\mathcal{G}) \otimes \mathbb{Q} \ar[r]  & R \Gamma_{\textrm{rig},c}(U'',\mathcal{M}^{\dag}_{K})  \ar[r]_{1-\phi} & R \Gamma_{\textrm{rig},c}(U'',\mathcal{M}^{\dag}_{K}) .
}
$$

En prenant les points fixes sous $G$ dans l'isomorphisme

$$
R \Gamma_{\textrm{rig},c}(U'', f^{\ast}_{\textrm{rig}}(E^{\dag}_{K \vert U})) \displaystyle \mathop{\longrightarrow}^{\sim} \mathcal{C}''(j''_{!} f^{\ast}_{\textrm{CRIS}}(E_{\vert U})^{m-\textrm{cris}}_{n})
$$

\noindent et en prouvant \`a la mani\`ere du th\'eor\`eme (3.2) que les points fixes sous $G$ du membre de droite s'identifient \`a
$$
\mathcal{C}(j_{U!}\ E_{\vert U_{n}}^{m-\textrm{cris}}) :=  R \displaystyle \mathop{\lim}_{\leftarrow \atop{m}} \{ (R  \displaystyle \mathop{\lim}_{\leftarrow \atop{n}} R \Gamma(\overline{X}_{\textrm{SYNT}}, j_{U!}\  E_{\vert U_{n}}^{m-\textrm{cris}})) \otimes \mathbb{Q} \} 
$$

\noindent on a prouv\'e le (1) du th\'eor\`eme (3.3.13) pour $U$ gr\^ace \`a [IV, th\'eo (4.2)].\\

On a donc deux triangles distingu\'es reli\'es par des fl\`eches qui sont des isomorphismes

$$
\xymatrix{
R \Gamma_{\textrm{rig},c}(U, E^{\dag}_{K_{\vert U}}) \ar[r] \ar[d]^{\simeq} &  R \Gamma_{\textrm{rig},c}(X, E^{\dag}_{K }) \ar[r] &  R \Gamma_{\textrm{rig},c}(Z, E^{\dag}_{K_{\vert Z}}) \ar[d]^{\simeq} \\
\mathcal{C} (j_{U!}\ E_{\vert U_{n}}^{m-\textrm{cris}}) \ar[r] & \mathcal{C} (j_{!}\ E_{n}^{m-\textrm{cris}}) \ar[r] & \mathcal{C} (j_{Z!}\ E_{\vert Z_{n}}^{m-\textrm{cris}})\   ;
}
$$

\noindent d'apr\`es les axiomes des cat\'egories triangul\'ees [H 1, I, $\S1$] on peut compl\'eter par un isomorphisme au milieu. \\
Ceci ach\`eve la preuve du (1) du th\'eor\`eme.\\

\noindent \textit{Prouvons \`a pr\'esent le (2)}. Comme $k$ est alg\'ebriquement clos les points fixes sous $G$ de $R \Gamma_{\textrm{\'et},c}(U'', f^{\ast}(\mathcal{F}_{\vert U})) \otimes \mathbb{Q}$ sont \'egaux \`a $R \Gamma_{\textrm{\'et},c}(U, \mathcal{F}_{\vert U}) \otimes \mathbb{Q}$. Par suite les points fixes sous $G$ du triangle distingu\'e (3.3.13.1) $G$-\'equivariant fournissent une suite exacte longue de cohomologie\\

\noindent (3.3.13.3) $ \rightarrow H^i_{\textrm{\'et},c}(U, \mathcal{F}_{\vert U}) \otimes \mathbb{Q} \rightarrow H^i_{\textrm{rig},c}(U/K, E^{\dag}_{K_{\vert U}}) \displaystyle \mathop{\rightarrow}_{1-\phi} H^i_{\textrm{rig},c}(U/K, E^{\dag}_{K_{\vert U}}) \rightarrow .$\\

\noindent Les groupes de cohomologie rigide $H^i_{\textrm{rig},c}(U/K, E^{\dag}_{K_{\vert U}})$ \'etant de dimension finie sur $K$ [Tsu 1, theo 6.1.2], et $k$ alg\'ebriquement clos, ces suites exactes se scindent en suites exactes courtes [I$\ell$ 1; II, lemme 5.6]. Notons $Z \hookrightarrow X$ l'immersion du ferm\'e compl\'ementaire \`a $U$ (rappelons que $X$ est suppos\'e r\'eduit). Si dim $Z = 0$, $Z$ est fini \'etale sur $k$ et on a une suite analogue \`a (3.3.13.3) pour $Z$ : par r\'ecurrence sur la dimension de $X$ on peut donc supposer l'existence de suites exactes courtes telles que (3.3.13.3) pour $U$ et pour $Z$. En particulier on a trois triangles distingu\'es horizontaux reli\'es par des fl\`eches induites par (3.3.13.3) appliqu\'e \`a $U$ et $Z$ :

$$
\begin{array}{c}
\xymatrix{
R \Gamma_{\textrm{\'et},c}(U, \mathcal{F}_{U \mathbb{Q}}) \ar[r] \ar[d] & R \Gamma_{\textrm{\'et},c}(X, \mathcal{F}_{ \mathbb{Q}}) \ar[r]  & R \Gamma_{\textrm{\'et},c}(Z, \mathcal{F}_{Z \mathbb{Q}}) \ar[d]\\
R \Gamma_{\textrm{rig},c}(U,  E^{\dag}_{K_{\vert U}}) \ar[r] \ar[d]_{1-\phi_{U}} & R \Gamma_{\textrm{rig},c}(X, E^{\dag}_{K}) \ar[r] \ar[d]_{1-\phi} & R \Gamma_{\textrm{rig},c}(Z, E_{K_{\vert Z}}^{\dag}) \ar[d]^{1-\phi_{Z}} \\
R \Gamma_{\textrm{rig},c}(U,  E^{\dag}_{K_{\vert U}}) \ar[r]  & R \Gamma_{\textrm{rig},c}(X, E^{\dag}_{K}) \ar[r]  & R \Gamma_{\textrm{rig},c}(Z, E_{K_{\vert Z}}^{\dag}) ; 
}
\end{array}
\leqno{(3.3.13.4)}
$$

\noindent par les axiomes des cat\'egories triangul\'ees [H 1, I, $\S$\ 1] on peut compl\'eter par un morphisme $R \Gamma_{\textrm{\'et},c} (X, \mathcal{F}_{\mathbb{Q}}) \rightarrow R \Gamma_{\textrm{rig},c}(X, E^{\dag}_{K})$.
\noindent Les suites exactes courtes (3.3.13.3) pour $U$ et $Z$ fournissent alors l'analogue pour $X$, d'o\`u le (2) du th\'eor\`eme (3.3.13).\\

\noindent \textit{Autre d\'emonstration du (2)}. Une autre m\'ethode consiste \`a appliquer le foncteur $\mathcal{C}$ \`a la suite exacte du th\'eor\`eme (3.3.12)\\

\noindent (3.3.13.5) $ 0 \rightarrow j_{!}\  \mathcal{F}_{n} \rightarrow j_{!}\ (\mathcal{F}_{n} \otimes_{\mathcal{V}^{\sigma}_{n}} \mathcal{O}^{m-\textrm{cris}}_{n,X}) \displaystyle \mathop{\rightarrow}_{1-\phi}  j_{!}\ (\mathcal{F}_{n} \otimes_{\mathcal{V}^{\sigma}_{n}} \mathcal{O}^{m-\textrm{cris}}_{n,X}) \rightarrow 0$

\noindent o\`u l'on remarque que $E^{m-\textrm{cris}}_{n} =  \mathcal{F}_{n} \otimes_{\mathcal{V}^{\sigma}_{n}} \mathcal{O}^{m-\textrm{cris}}_{n,X}$. Par le (1) du th\'eor\`eme le triangle distingu\'e ainsi obtenu s'identifie au triangle distingu\'e

$$
R \Gamma_{\textrm{\'et},c} (X, \mathcal{F}_{\mathbb{Q}}) \longrightarrow R \Gamma_{\textrm{rig},c}(X, E^{\dag}_{K}) \displaystyle \mathop{\rightarrow}_{1-\phi} R \Gamma_{\textrm{rig},c}(X, E^{\dag}_{K}).
$$

\noindent La suite exacte longue de cohomologie se scinde alors en suites exactes courtes par le m\^eme argument que ci-dessus. $\square$

\vskip 3mm
\noindent \textbf{Remarque (3.3.14)}. En supposant seulement que $k$ contient $\mathbb{F}_{q}$ [cf (3.0)] et que $\mathcal{F}_{n}$ est un $\mathcal{V}^{\sigma}_{n}$-module localement trivial sur SYNT$(X)$, on pose encore

$$
E^{m-\textrm{cris}}_{n} = T^{(m)}(\mathcal{F}_{n})^{\textrm{cris}} = \mathcal{F}_{n} \otimes_{\mathcal{V}^{\sigma}_{n}} \mathcal{O}^{m-\textrm{cris}}_{n,X}.
$$

\noindent En appliquant le foncteur $R \Gamma(\overline{X}_{\textrm{SYNT}, -})$ \`a la suite exacte (3.3.13.5), on obtient un triangle distingu\'e

$$
R \Gamma_{\textrm{\'et},c} (X, \mathcal{F}_{n}) \longrightarrow R \Gamma_{\textrm{synt},c}(X, E^{m-\textrm{cris}}_{n}) \displaystyle \mathop{\longrightarrow}_{1-\phi} R \Gamma_{\textrm{synt},c}(X, E^{m-\textrm{cris}}_{n}).
$$


\cleardoublepage


\vskip 10mm
\chapter*{VI.  Fonctions $L$}
\markboth{\sc j.-y. etesse}{\sc VI.  Fonctions $L$}
 Nous allons d\'efinir dans ce VI les fonctions $L$ attach\'ees \`a des vari\'et\'es sur des corps finis, et \`a coefficients dans des $F$-modules ou des $F$-(iso)cristaux : nous utiliserons deux m\'ethodes, l'une par les rel\`evements de Teichm¬\"uller, l'autre par voie cohomologique, et nous montrerons comment elles se rejoignent.\\

Sauf mention du contraire, on suppose dans ce VI que $k$ est un corps fini, $k = \mathbb{F}_{q}, q = p^a, \mathcal{V}$ est un anneau de valuation discr\`ete complet, d'id\'eal maximal $\mathfrak{m}$ et corps r\'esiduel $k = \mathbb{F}_{q}$. On suppose le corps des fractions $K$ de $\mathcal{V}$ de caract\'eristique 0,  on fixe une uniformisante $\pi$ et on note $e$ l'indice de ramification.\\

On rel\`eve la puissance $q$ sur $k$ en un automorphisme $\sigma$ de $\mathcal{V}$, tel que $\sigma(\pi) = \pi$, suivant la m\'ethode [Et 5, I 1.1], cf [II, 0] : on note encore $\sigma$ son extension \`a $K$.\\

\section*{1. Fonctions $L$ des $F$-modules convergents ou surconvergents}

\subsection*{1.1. Rel\`evement de Teichm¬\"uller}

On reprend les notations de [III, 3.1] en supposant cette fois que $k$ est un corps fini, $k = \mathbb{F}_{q}, q = p^a$.

\vskip 3mm
\subsection*{1.2. $F$-modules convergents}

Avec les notations du 1.1 on d\'esigne par $\mathbf{F^a\mbox{-}\textrm{\bf Mod}(\hat{A})}$ (resp. $\mathbf{F^a\mbox{-}\textrm{\bf Mod}(\hat{A}_{K}})$ la cat\'egorie des $\hat{A}$ (resp. $\hat{A}_{K})$-modules projectifs de type fini $\mathcal{M}$ munis d'un morphisme de Frobenius (non n\'ecessairement un isomorphisme) 

$$\phi_{\mathcal{M}} : \mathcal{M}^{\sigma} :=  F^{\ast}_{\mathcal{A}}(\mathcal{M}) \rightarrow \mathcal{M}
$$

\noindent avec $\mathcal{A} = \hat{A}$ (resp. $\hat{A}_{K}$). Un tel $\mathcal{M}$ est appel\'e $F$-\textbf{module convergent}.\\

Soit $\mathcal{M} \in F^a\mbox{-}\textrm{Mod}(\hat{A}_{K})$.
La fibre $\mathcal{M}_{x}$ de $\mathcal{M}$ en $x \in \vert X \vert$ est par d\'efinition\\

\noindent (1.2.1) $\qquad \qquad \qquad \mathcal{M}_{x} := \hat{\tau}_{K}(x)^{\ast} (\mathcal{M})$ , \\

\noindent et $\phi_{\mathcal{M}}$ induit\\

\noindent (1.2.2) $\qquad \qquad \phi_{x} = \phi_{\mathcal{M}} \otimes_{\hat{A}_{K}} K(x) : \sigma^{\ast}_{K(x)}(\mathcal{M}_{x}) \rightarrow \mathcal{M}_{x} $\ , \\

\noindent d'apr\`es la commutativit\'e du diagramme (1.1.1).\\

L'it\'er\'e deg $x$ fois de $\phi_{x}$ est un endomorphisme $K(x)$-lin\'eaire du $K(x)$-espace vectoriel de dimension finie $\mathcal{M}_{x}$

$$
\phi^{\textrm{deg}\  x}_{x} : \mathcal{M}_{x} \rightarrow \mathcal{M}_{x}\ .
$$

\noindent Notons\\

\noindent (1.2.3) $\qquad \qquad \textrm{det} (\mathcal{M}_{x},T) = \textrm{det} (1-T\ \phi^{\textrm{deg}\  x}_{x}, \mathcal{M}_{x})$\\

\noindent le ``polyn\^ome caract\'eristique'' de $\phi^{\textrm{deg}\  x}_{x}$.\\

Pour $\mathcal{M}' \in F^a\mbox{-}\textrm{Mod} (\hat{A})$ on d\'efinit de m\^eme\\

\noindent (1.2.4)$ \qquad \qquad \textrm{det}(\mathcal{M}'_{x}, T) = \textrm{det}(1-T\ \phi^{\textrm{deg}\  x}_{x}, \mathcal{M}'_{x}).$\\

\noindent \textbf{Lemme (1.2.5)}. \textit{Avec les notations pr\'ec\'edentes, on a :}

\begin{itemize}
\item[(i)] \textit{Si  $\mathcal{M} \in F^a\mbox{-}Mod(\hat{A}_K)$ , alors det $(\mathcal{M}_{x},T) \in K[T]  $ .}
\item[(ii)] \textit{Si $\mathcal{M}' \in F^a\mbox{-}Mod(\hat{A})$, alors det  $(\mathcal{M}'_{x},T)  \in \mathcal{V}[T] $.}
\end{itemize}

\vskip 3mm
\noindent \textit{D\'emonstration}. Pour (i), soient $(e_{i})_{i=1,...,r}$ (resp. $(e_{i} \otimes 1)_{i=1,...,r})$ une base locale de $\mathcal{M}$ (resp. de $\mathcal{M}^{\sigma})$, et $C(\underline{X})$ la matrice de $\phi_{\mathcal{M}}$ dans ces bases respectives. Alors 
$$
C(\underline{X}) = \displaystyle \mathop{\Sigma}_{\underline{u} \in \mathbb{N}^n } a_{\underline{u}}\ \underline{X}^{\underline{u}},
$$

\noindent avec $a_{\underline{u}}\ \in \pi^{\alpha}\ M_{r}(\mathcal{V})$ pour un $\alpha \in \mathbb{Z}$ et 

$$
det(\mathcal{M}_{x}, T) = det \{ 1-T\ C(\underline{t}(x)^{q^{\textrm{deg}\  x-1}}) \times ... \times C(\underline{t}(x)^{q}) \times C(\underline{t}(x)) \}
$$

\noindent o\`u $\underline{t}(x)^{\beta} = t_{1}(x)^{\beta} \times ... \times t_{n}(x)^{\beta}$ et les $t_{j}(x)$ sont les coordonn\'ees de $\hat{\tau}_{K}(x)$. Il est clair que det $(\mathcal{M}_{x},T)$ a des coefficients invariants par l'action du Frobenius $\sigma_{K(x)}$, car $\sigma_{K(x)}$ envoie $\underline{t}(x)$ sur $\underline{t}(x)^q$ ; d'o\`u le (i).\\

Le cas (ii) est analogue. $\square$\\

Si l'on note $\tilde{\mathcal{M}}_{x}$ l'espace vectoriel $\mathcal{M}_{x}$ vu comme $K$-espace vectoriel et $\tilde{\phi}_{x}$ son endomorphisme de Frobenius on a

$$
det(\mathcal{M}_{x}, T) = det(1-T\ \phi^{\textrm{deg}\ x}_{x}) = det_{K}(1-T\ \tilde{\phi}^{\textrm{deg}\ x}_{x})^{-1/\textrm{deg}\ x}\ .
$$

\vskip 3mm
\noindent \textbf{D\'efinition (1.2.6)}. \textit{La fonction $L$ de $(\mathcal{M}, \phi_{\mathcal{M}}) \in F^a\mbox{-}Mod(\hat{A}_{K})$ est d\'efinie par}

$$
L(Spec\ A_{0}, \mathcal{M}, t) = \displaystyle \mathop{\prod}_{x \in \vert Spec\ A_{0} \vert} det(1-t^{\textrm{deg}\ x}\  \phi^{\textrm{deg}\ x}_{x} \mid \mathcal{M}_{x})^{-1} \in K[[t]]
$$
$\qquad \qquad \qquad \qquad \qquad = \displaystyle \mathop{\prod}_{x \in \vert Spec\ A_{0} \vert} 
det(1-t^{\textrm{deg}\ x}\  \tilde{\phi}^{\textrm{deg}\ x}_{x} \mid \tilde{\mathcal{M}}_{x})^{-1/\textrm{deg}\ x}\ .$\\

\noindent \textit{Si $(\mathcal{M'}, \phi_{\mathcal{M'}}) \in F^a\mbox{-}\textrm{Mod}(\hat{A})$ on d\'efinit de m\^eme $L(Spec\ A_{0}, \mathcal{M'}, t)$ et alors}

$$
L(Spec\ A_{0}, \mathcal{M'}, t) = L(Spec\ A_{0}, \mathcal{M}'_{K},t) \in \mathcal{V}[[t])\ ,
$$

\noindent \textit{o\`u l'on a pos\'e}

$$
(\mathcal{M}'_{K}, \phi_{\mathcal{M}'_{K}}) := (\mathcal{M'}, \phi_{\mathcal{M'}}) \otimes_{\hat{A}} \hat{A}_{K}\ .
$$

\vskip 3mm
\noindent \textbf{Lemme (1.2.7)}. \textit{Pour \'etablir la m\'eromorphie $p$-adique de $L(Spec\ A_{0}, \mathcal{M}, t)$ pour $(\mathcal{M}, \phi_{\mathcal{M}}) \in F^a\mbox{-}Mod(\hat{A} _{K})$ on peut supposer qu'il existe $(\mathcal{M'}, \phi_{\mathcal{M'}}) \in F^a\mbox{-}Mod(\hat{A})$ tel que $\mathcal{M'}$ est libre et $(\mathcal{M}, \phi_{\mathcal{M}}) = (\mathcal{M'}, \pi^{\alpha} \phi_{\mathcal{M'}})  \otimes_{\hat{A}} \hat{A}_{K}$ pour un $\alpha \in \mathbb{Z}$. On a alors}

$$
L(Spec\ A_{0}, \mathcal{M}, t) = L(Spec\ A_{0}, \mathcal{M'}, \pi^{\alpha} t).
$$

\vskip 3mm
\noindent \textit{D\'emonstration}. Puisque $\mathcal{M}$ est projectif de type fini sur $\hat{A}_{K}$, et qu'un ouvert de $Spec\ \hat{A}_{K}$ est intersection d'un ouvert de $Spec\ \hat{A}$ avec $Spec\ \hat{A}_{K}$, il existe un recouvrement fini de $Spec\ \hat{A}_{K}$ par des ouverts $U_{K} = Spec\ B_{K}$, o\`u $B = \hat{A}\ [1/g], g \in \hat{A}$, tels que

$$
\mathcal{M} \otimes_{\hat{A}_{K}} B_{K}  \simeq \displaystyle \mathop{\oplus}_{i=1}^r B_{K}\ e_{i}.
$$

\noindent Relevons $g\  \textrm{mod}\  \pi =: g_{0}$ en $f \in A$ ; comme on a un isomorphisme $\hat{B} \simeq \widehat{A[1/f]}$ [Et 4, cor 1 du th\'eo 4] on en d\'eduit que

$$
\mathcal{N}  := \mathcal{M} \otimes_{\hat{A}_{K}} \hat{B}_{K} \simeq \displaystyle \mathop{\oplus}_{i=1}^r \widehat{A[1/f]_{K}}\  e_{i}
 $$

\noindent et que la matrice $C(\underline{X})$ du Frobenius $\phi_{\mathcal{N}}$ est \`a coefficients dans $\pi^{\alpha} \widehat{A[1/f]}$, avec $\alpha \in \mathbb{Z}$.\\

La fonction $L(\mathcal{M}, t)$ est d\'efinie par un produit eul\'erien sur les points ferm\'es de $Spec\ A_{0}$ et ceux-ci sont en bijection avec les points ferm\'es de $Spf\ \hat{A}$ : un recouvrement ouvert fini de $Spf\ \hat{A}$ est fourni par des $Spf\ \widehat{A[1/f]}$, $f \in A$ comme ci-dessus ; on peut donc prendre

$$
\mathcal{M'} = \displaystyle \mathop{\oplus}_{i=1}^r \widehat{A[1/f]}\  e_{i}
$$

\noindent avec pour matrice du Frobenius $\phi_{\mathcal{M'}}$ la matrice $\pi^{- \alpha}\  C(\underline{X})$, \`a coefficients dans $\widehat{A[1/f]}$ :

$$
\xymatrix{
\mathcal{M'}^{\sigma} \ar[rr]^{\phi_{\mathcal{M'}}}  \ar@{^{(}->}[d] && \mathcal{M'} \ar@{^{(}->}[d] \\
\mathcal{N}^{\sigma} = \mathcal{M}^{\prime\sigma}_{K} \ar[rr]_{\pi^{- \alpha}{\phi_{\mathcal{N}}}}&& \mathcal{N} = \mathcal{M}'_{K} .
}
$$

Le couple $(\mathcal{M'}, \phi_{\mathcal{M'}})$ est ce que Wan appelle une $\sigma$-module convergent [W 2], [W 3]. $\square$\\

\newpage
\subsection*{1.3.  $F$-modules surconvergents}

Avec les notations de 1.1 et par analogie avec 1.2 on d\'esigne par $\mathbf{F^a\mbox{-}\textrm{\bf Mod}(A^{\dag})}$ (resp. $\mathbf{F^a\mbox{-}\textrm{\bf Mod}(A^{\dag}_{K}})$ la cat\'egorie des $\textbf{F-\textrm{modules surconvergents}}$, i.e. la cat\'egorie des $A^{\dag}$ (resp. $A^{\dag}_{K})$-modules projectifs de type fini $M$ muni d'un morphisme de Frobenius (non n\'ecessairement un isomorphisme)

$$
\phi_{M} : M^{\sigma} = F^{\ast}_{\mathcal{A}}(M) \rightarrow M
$$

\noindent  avec $\mathcal{A} = A^{\dag}$ (resp. $A^{\dag}_{K}$).\\

Soit $M \in F^a\mbox{-}\textrm{Mod}(A^{\dag}_{K})$. La fibre $M_{x}$ de $M$ en $x \in \vert X \vert $ est par d\'efinition\\

\noindent (1.3.1) $\qquad \qquad \qquad \qquad M_{x} := \tau^{\dag}_{K}(x)^{\ast}(M)$,\\

\noindent et $\phi_{M}$ induit\\

\noindent (1.3.2) $\qquad \qquad \qquad \phi_{x} = \phi_{M} \otimes_{A^{\dag}_{K}} K(x) : \sigma^{\ast}_{K(x)}(M_{x}) \rightarrow M_{x}\ ,$\\

\noindent d'apr\`es la commutativit\'e du diagramme (1.1.1).\\

L'it\'er\'e deg $x$ fois de $\phi_{x}$ est un endomorphisme $K(x)$-lin\'eaire du $K(x)$-espace vectoriel de dimension finie $M_{x}$

$$
\phi^{\textrm{deg}\ x}_{x} : M_{x} \rightarrow M_{x}\ ;
$$

\noindent on notera $\tilde{M}_{x}$ l'espace vectoriel $M_{x}$ vu comme $K$-espace vectoriel et $\tilde{\phi}_{x}$ son morphisme de Frobenius. \\

Notons\\

\noindent {(1.3.3)}  $\qquad \qquad \qquad  \textrm{det}(M_{x}, T) =  \textrm{det}(1 - T \phi^{\textrm{deg}\ x}_{x}, M_{x})\ ;$\\

\noindent pour $M' \in F^a\mbox{-}\textrm{Mod}(A^{\dag})$ on pose de m\^eme\\

\noindent (1.3.4) $\qquad \qquad \qquad  \textrm{det}(M'_{x},T)\   =  \textrm{det}(1 - T\ \phi^{\textrm{deg}\ x}_{x}, M'_{x}).$\\

On d\'emontre le lemme suivant comme (1.2.5).

\vskip 3mm
\noindent \textbf{Lemme (1.3.5)}. \textit{Avec les notations pr\'ec\'edentes, on a : }

\begin{itemize}
\item[(i)] \textit{Si $M \in F^a\mbox{-}\textrm{Mod} (A^{\dag}_{K})$, alors $det(M_{x}, T) \in K[T]\  $}.
\item[(ii)] \textit{Si $M' \in F^a\mbox{-}\textrm{Mod} (A^{\dag})$, alors $det(M'_{x}, T) \in \mathcal{V}[T]\ $.}
\end{itemize}

\vskip 3mm
De m\^eme on a :

\vskip 3mm
\noindent \textbf{D\'efinition et proposition (1.3.6)}. \textit{Soit $(M, \phi_{M}) \in F^a \mbox{-}\textrm{Mod}(A^{\dag}_{K})$ et \\
$(\mathcal{M}, \phi_{\mathcal{M}}) \in F^a\mbox{-}\textrm{Mod}(\hat{A}_{K})$ son image canonique par l'extension des scalaires de $A^{\dag}_{K}$ \`a  $\hat{A}_{K}$. La fonction $L$ de $(M, \phi_{M})$ est d\'efinie par }

$$
L(Spec\ A_{0}, M, T) = \displaystyle \mathop{\prod}_{x \in \vert Spec\ A_{0} \vert}\ det(1-t^{\textrm{deg}\  x}\  \phi^{\textrm{deg}\ x}_{x}, M_{x})^{-1} \in K[[t]]
$$
$\qquad \qquad \qquad \qquad \qquad = \displaystyle \mathop{\prod}_{x \in \vert Spec\ A_{0} \vert}\ det(1-t^{\textrm{deg}\  x}\ \tilde{\phi}^{\textrm{deg}\ x}_{x}, \tilde{M}_{x})^{-1/\textrm{deg}\ x}\ ,
$

\noindent \textit{et on a}

$$L(Spec\ A_{0}, M, t) = L(Spec\ A_{0}, \mathcal{M}, t).$$

\vskip 3mm

\textit{Si $(M', \phi_{M'}) \in F^a\mbox{-}\textrm{Mod} (A^{\dag})$ on d\'efinit de m\^eme $L(Spec\ A_{0}, M', t)$ et alors}

$$
L(Spec\ A_{0}, M', t) = L(Spec\ A_{0}, M'_{K}, t) \in \mathcal{V}[[t]] ,
$$

\noindent \textit{o\`u l'on a pos\'e}

$$
(M'_{K}, \phi_{M'_{K}}) := (M', \phi_{M'}) \otimes_{A^{\dag}}\ A^{\dag}_{K}\ .
$$

\vskip 3mm
Comme (1.2.7) on montre :

\vskip 3mm
\noindent \textbf{Lemme (1.3.7)}. \textit{Pour \'etablir la m\'eromorphie $p$-adique de $L(Spec\ A_{0}, M, t)$ pour $(M, \phi_{M}) \in F^a\mbox{-}\textrm{Mod}(A^{\dag}_{K})$ on peut supposer qu'il existe $(M', \phi_{M'}) \in F^a\mbox{-}\textrm{Mod}(A^{\dag})$ tel que $M'$ est libre et $(M, \phi_{M}) = (M', \pi^{\alpha} \phi_{M'}) \otimes_{A^{\dag}} A^{\dag}_{K}$ pour un $\alpha \in \mathbb{Z}$. On a alors}

$$
L(Spec\ A_{0}, M, t) = L(Spec\ A_{0}, M', \pi^{\alpha} t)\ .
$$

\vskip 3mm
\noindent \textbf{Th\'eor\`eme (1.3.8)}. \textit{Soient $X = Spec\ A_{0}$ un $k$-sch\'ema lisse, $A$ une $\mathcal{V}$-alg\`ebre lisse relevant $A_{0}$ et $(M, \phi_{M}) \in F^a\mbox{-}\textrm{Mod}(A^{\dag}_{K})$.\\
Alors $L(X, M, t)$ est $p$-adiquement m\'eromorphe.}

\vskip 3mm
\noindent \textit{D\'emonstration}. Par le lemme (1.3.7), on peut remplacer $M$ par un $A^{\dag}$-module libre $M' $: $(M', \phi_{M'})$ est alors ce que Wan appelle un $\sigma$-module surconvergent [W 2] [W 3]. L'extension au cas affine et lisse de la formule des traces de Monsky-Washnitzer \'etablie par Wan [W 3, appendix] prouve alors la m\'eromorphie de $L(M', t)$, donc celle de $L(M, t)$ : on peut aussi se ramener classiquement au cas de l'espace affine [W 3, \S\ 9] et utiliser la formule des traces de Monsky [Dw 2, 7(a)]. $\square$\\

\subsection*{1.4. La conjecture de Dwork pour les $F$-modules surconvergents}

	Avec les notations de 1.1, 1.2, 1.3 d\'ecomposons le ''polyn\^ome caract\'eristique'' de $(\mathcal{M, \phi_{\mathcal{M}}}) \in F^a\mbox{-}\textrm{Mod}(\hat{A}_{K})$ au point $x \in \vert X \vert $ en
	
$$
det(\mathcal{M}_{x},t) := det(1-t\ \phi^{\textrm{deg}\ x}_{x} \vert \mathcal{M}_{x}) \in K[t]
$$
$
\qquad \qquad \qquad \qquad \qquad \qquad \quad = \displaystyle \mathop{\pi}_{j}\  (1 - a_{j,x}\  t) ,
$

\noindent o\`u les $a_{j,x}$ sont dans une cl\^oture alg\'ebrique $K^{\textrm{alg}}$ de $K$ : plus exactement les $a_{j,x}$ sont dans une extension finie (\'eventuellement ramifi\'ee) $K'(x) \subset K ^{\textrm{alg}}$ de $K(x)$ ; soient $\pi'(x)$ une uniformisante de $K'(x)$ et $\sigma_{K'(x)}$ un rel\`evement \`a $K'(x)$ de la puissance $p^a$ de $k(x)$ tel que $\sigma_{K'(x)}(\pi'(x)) = \pi'(x)$  [Et 5, 1.1]. Notons $\pi_{x} = \pi^{\textrm{deg}\ x}$ et $\textrm{ord}_{\pi_{x}}$ la valuation de $K'(x)$ normalis\'ee par 
$\textrm{ord}_{\pi_{x}}(\pi_{x})  = 1$.\\

Pour tout nombre rationnel $\alpha \in \mathbb{Q}$ on d\'efinit la partie de pente $\alpha$ du ``polyn\^ome caract\'eristique''  $det(\mathcal{M}_{x}, t)$ par le produit\\

\noindent (1.4.1) $\qquad \qquad \qquad det_{\alpha}(\mathcal{M}_{x}, t) := \displaystyle \mathop{\prod}_{\textrm{ord}_{\pi_{x}}(a_{j,x}) = \alpha}\ (1 - a_{j,x}\ t)$ .\\

\vskip 2mm
\noindent Si $a_{j,x}$ est l'inverse d'une racine de $det(\mathcal{M}_{x}, t)$, alors $\sigma_{K'(x)} (a_{j,x})$ en est une aussi par le m\^eme argument que pour la d\'emonstration du lemme (1.2.5), et puisque $\sigma_{K'(x)} : K'(x) \rightarrow K'(x)$ est une extension isom\'etrique on a $\textrm{ord}_{\pi_{x}}(a_{j,x}) = \textrm{ord}_{\pi_{x}}(\sigma_{K'(x)}(a_{jx}))$ ; par cons\'equent $det_{\alpha}(\mathcal{M}_{x}, t)$ est \`a coefficients dans $K$. Comme $\mathcal{M}$ est un $\hat{A}_{K}$-module projectif de type fini, il existe $\alpha_{0} \in \mathbb{Q}$ tel que pour tout $x \in \vert Spec\ A_{0} \vert$ et tout $\alpha \in \mathbb{Q}$, $\alpha < \alpha_{0}$, on ait $det_{\alpha}(\mathcal{M}_{x}, t) = 1$. D'autre part $det(\mathcal{M}_{x}, t)$ s'exprime par un produit fini\\

\noindent (1.4.2) $\qquad \qquad \qquad det(1 - t\ \phi^{\textrm{deg}\ x}_{x}, \mathcal{M}_{x}) = \displaystyle \mathop{\prod}_{\alpha \in \mathbb{Q}} det_{\alpha}(\mathcal{M}_{x}, t)$ . \\

La partie de pente $\alpha$ de la fonction $L(X, \mathcal{M}, t)$ est d\'efinie par l'expression\\

\noindent (1.4.3) $\qquad \qquad \qquad L_{\alpha}(X, \mathcal{M}, t) := \displaystyle \mathop{\prod}_{x \in \vert X \vert}  det_{\alpha}(\mathcal{M}_{x}, t^{\textrm{deg}\ x})^{-1} \in K[[t]]$ .\\

\noindent D'apr\`es (1.4.2) on a la relation\\

\noindent (1.4.4) $\qquad \qquad \qquad L(X, \mathcal{M}, t) = \displaystyle \mathop{\prod}_{\alpha \in \mathbb{Q}} L_{\alpha}(X, \mathcal{M}, t)$ ;\\

\noindent en fait le th\'eor\`eme de sp\'ecialisation de Grothendieck prouve que le produit (1.4.4) est fini [W 3, fin du \S\ 2].\\

Pour $r \in \mathbb{N}^{\ast}$ et $\alpha \in \mathbb{Q}$, on d\'efinit plus g\'en\'eralement\\

\noindent (1.4.5) $\qquad \qquad \qquad det^{(r)}(\mathcal{M}_{x}, t) := \displaystyle \mathop{\prod}_{j}\  (1 - (a_{j,x})^r\ t),$\\

\vskip 3mm
\noindent (1.4.6) $\qquad \qquad \qquad det^{(r)}_{a}(\mathcal{M}_{x}, t) := \displaystyle \mathop{\prod}_{\textrm{ord}_{\pi_{x}}(a_{j,x}) = \alpha}\ (1 - (a_{j,x})^r\ t) $ ,\\

\vskip 3mm
\noindent (1.4.7) $\qquad \quad   L^{(r)}(X, \mathcal{M}, t) := \displaystyle \mathop{\prod}_{x \in \vert X \vert}\  det^{(r)}  (\mathcal{M}_{x}, t^{\textrm{deg}\ x})^{-1} \in K[[t]]$ , \\

\vskip 3mm
\noindent (1.4.8) $\qquad   \quad L^{(r)}_{\alpha}(X, \mathcal{M}, t) := \displaystyle \mathop{\prod}_{x \in \vert X \vert}\  det^{(r)}_{\alpha}  (\mathcal{M}_{x}, t^{\textrm{deg}\ x})^{-1} \in K[[t]]$ .\\

\noindent On a encore :\\

\noindent (1.4.9) $\qquad   \quad L^{(r)} (X, \mathcal{M}, t) = \displaystyle \mathop{\prod}_{\alpha \in \mathbb{Q}}\ L^{(r)}_{\alpha}  (X, \mathcal{M}, t)$ .\\

\noindent Plus pr\'ecis\'ement, gr\^ace \`a l'expression \'etablie par Wan [W 2, lemma 4.4], on a en fait\\

\noindent (1.4.10) $\qquad   \quad L^{(r)} (X, \mathcal{M}, t) = \displaystyle \mathop{\prod}_{i \geqslant 1}\ L(X, \textrm{Sym}^{r-i}\ \mathcal{M} \ \otimes \displaystyle \mathop{\Lambda}^i\ \mathcal{M}, t)^{i \times (-1)^{i-1}}$ ,\\

\vskip 3mm
\noindent (1.4.11) $\qquad   \quad L^{(r)}_{\alpha} (X, \mathcal{M}, t) = \displaystyle \mathop{\prod}_{i \geqslant 1}\ L_{\alpha} (X, \textrm{Sym}^{r-i}\ \mathcal{M} \ \otimes \displaystyle \mathop{\Lambda}^i\ \mathcal{M}, t)^{i \times (-1)^{i-1}}$ .\\

\vskip 3 mm
Les d\'efinitions pr\'ec\'edentes s'\'etendent \`a $M \in F^a\mbox{-}\textrm{Mod}(A^{\dag}_{K})$.\\

Nous pouvons \`a pr\'esent \'enoncer la conjecture de Dwork g\'en\'eralis\'ee pour les $F^a$-modules surconvergents.

\vskip 3mm
\noindent (1.4.12) \textit{\textbf{Conjecture (Dwork)}. Soit $M \in F^a\mbox{-}\textrm{Mod}(A^{\dag}_{K})$. Alors, pour tout nombre rationnel $\alpha \in \mathbb{Q}$ et tout entier $r \in \mathbb{N}^{\ast}$, la fonction $L^{(r)}_{\alpha}(X, M, t)$ est $p$-adiquement m\'eromorphe.}\\

 Gr\^ace aux travaux de Wan, on d\'emontre ci-dessous (facilement) la conjecture :

\vskip 3mm
\noindent \textbf{Th\'eor\`eme (1.4.13)}. \textit{Soient $A$ une $\mathcal{V}$-alg\`ebre lisse, $X = Spec(A/\mathfrak{m} A)$ et $M \in F^a\mbox{-}\textrm{Mod}(A^{\dag}_{K})$. Alors, pour tout  $\alpha \in \mathbb{Q}$ et tout $r \in \mathbb{N}^{\ast}$, les fonctions $L^{(r)}_{\alpha}(X, M, t)$, et $L^{(r)}(X, M, t)$ sont $p$-adiquement m\'eromorphes.}

\vskip 3mm
\noindent \textit{D\'emonstration}. Par stabilit\'e de la cat\'egorie $F^a\mbox{-}\textrm{Mod}(A^{\dag}_{K})$ par puissances sym\'etriques et ext\'erieures [W 2, \S\ 3] la conjecture (1.4.12) se ram\`ene au cas $r = 1$. Par la m\^eme d\'emonstration que celle du lemme (1.3.7) on peut supposer qu'il existe un $F^a$-module libre surconvergent $(M', \phi_{M'}) \in F^a\mbox{-}\textrm{Mod}(A^{\dag})$ tel que

$$
(M, \phi_{M}) = (M', \pi^{\beta} \phi_{M'}) \otimes_{A^{\dag}} A^{\dag}_{K}\  \textrm{pour un}\  \beta \in \mathbb{Z}\ ,
$$

\noindent et $ \qquad \qquad \quad L(M,t) = L(M', \pi^{\beta} t)$ .\\

\noindent On est ramen\'e \`a montrer que $L_{\alpha}(Spec\ A_{0}, M', t)$ est $p$-adiquement m\'eromorphe, ce que Wan a \'etabli [W 3, theo 1.1]. $\square$\\

\section*{2. Fonctions $L$ des $F$-isocristaux convergents (resp. des $F$-cristaux)}.

\subsection*{2.1.  $F$-isocristaux convergents}

\textbf{2.1.1.}  Soient $S$ un $k$-sch\'ema s\'epar\'e de type fini, $\mathcal{E} \in F^a\mbox{-}\textrm{Isoc}(S/K)$ et $\phi_{\mathcal{E}} : F^{\ast}_{\sigma} \ \mathcal{E} \displaystyle \mathop{\rightarrow}^{\sim} \mathcal{E}$ son isomorphisme de Frobenius [B 3, 2]. Comme la cat\'egorie $F^a\mbox{-}\textrm{Isoc}(S/K)$ ne d\'epend que du sous-sch\'ema r\'eduit sous-jacent  \`a $S$ et que la fonction $L(S, \mathcal{E}, t)$, d\'efinie ci-dessous, est un produit eul\'erien sur l'ensemble $\vert S \vert$ des points ferm\'es de $S$, on peut supposer $S$ r\'eduit pour d\'efinir $L(S, \mathcal{E}, t)$.\\

D\'ecomposons alors $S$ en strates affines et lisses, $X = Spec\ A_{0}$, sur $k = \mathbb{F}_{q}$ ; pour $A$ une $\mathcal{V}$-alg\`ebre lisse relevant $A_{0}$ on reprend les notations du \S\ 1. Pour $x \in \vert X \vert $ le diagramme commutatif

$$
\xymatrix{
Spec\ k(x) \ar@{^{(}->}[r] \ar@{^{(}->}[d]^{i_{x}} & Spec\ \mathcal{V}(x)\ar[d] \\
S \ar[r] & Spec\ \mathcal{V}
}
$$

[resp.

$$
\xymatrix{
Spec\ k(x) \ar@{^{(}->}[r] \ar@{^{(}->}[d]^{i_{x}} & Spec\ \mathcal{V}(x)\ar@{=}[d] \\
S \ar[r] & Spec\ \mathcal{V}\ ]
}
$$

\noindent fournit un foncteur image inverse analogue \`a [B 3, (2.3.6), (2.3.7)]\\

$$
\begin{array}{c}
\xymatrix{
i^{\ast}_{x} : F^a\mbox{-}\textrm{Isoc}(S/K) \ar[r] & F^a\mbox{-}\textrm{Isoc}(Spec(k(x))/K(x)) \ar@{=}[d]\\
 & F^a\mbox{-}\textrm{Isoc}^{\dag} (Spec(k(x))/K(x))
}
\end{array}
\leqno{(2.1.1.1)}
$$

\noindent [resp.

$$
\begin{array}{c}
\xymatrix{
\tilde{i}^{\ast}_{x} : F^a\mbox{-}\textrm{Isoc}(S/K) \ar[r] & F^a\mbox{-}\textrm{Isoc}(Spec(k(x))/K) \ar@{=}[d]\\
 & F^a\mbox{-}\textrm{Isoc}^{\dag} (Spec(k(x))/K)\ ].
}
\end{array}
\leqno{(2.1.1.2)}
$$\\

Alors 

$$
\begin{array}{c}
\xymatrix{
\qquad \quad  \mathcal{E}_{x}  := H^0_{\textrm{conv}}(Spec(k(x))/K(x), i^{\ast}_{x} \mathcal{(E)})\  [ cf\ O3]\\
\simeq H^0_{\textrm{rig}} (Spec(k(x))/K(x), i^{\ast}_{x} \mathcal{(E)})
}
\end{array}
\leqno{(2.1.1.3)}
$$

\noindent [resp.

$$
\begin{array}{c}
\xymatrix{
\tilde{\mathcal{E}}_{x}  := H^0_{\textrm{conv}}(Spec(k(x))/K, \tilde{i}^{\ast}_{x} \mathcal{(E)})\\
\quad \simeq H^0_{\textrm{rig}}(Spec(k(x))/K, \tilde{i}^{\ast}_{x} \mathcal{(E)})\ ]
}
\end{array}
\leqno{(2.1.1.4)}
$$

\noindent est un $K(x)$-espace vectoriel [resp. $K$-espace vectoriel] de dimension finie muni d'un isomorphisme de Frobenius\\

\noindent (2.1.1.5)  $ \qquad \qquad \qquad \qquad \phi'_{x} : = H^0_{\textrm{rig}}(i^{\ast}_{x} (\phi))$\\

\noindent (2.1.1.6)  $ \quad  [\textrm{resp.} \qquad \qquad \tilde{\phi'}_{x} : = H^0_{\textrm{rig}}(\tilde{i}^{\ast}_{x} (\phi))\ ]\ .$\\

\textbf{2.1.2.} Il y a une deuxi\`eme fa\c{c}on, explicit\'ee ci-apr\`es en (2.1.2.3), de d\'efinir la fibre de $\mathcal{E}$ en $x \in \vert X \vert $ et son Frobenius, et ceci via les $F$- modules convergents : nous montrons que les deux fa\c{c}ons co\"{\i}ncident.\\

Pour $x \in \vert X \vert = \vert Spec\ A_{0} \vert $ le rel\`evement de Teichm¬\"uller de $x$, vu en 1.1,

$$
\xymatrix{
\hat{\tau}(x) : \hat{A} \ar@{->>}[r]  & \mathcal{V}(x)
}
$$

\noindent induit un morphisme

$$
u : P_{1} := Spf(\mathcal{V}(x)) \rightarrow P := Spf(\hat{A}) 
$$

\noindent dont la r\'eduction mod $\pi$ est not\'ee

$$
v : X_{1} := Spec(k(x)) \rightarrow Y := P_{k}
$$

\noindent et $v$ se factorise par

$$
i_{x} : X_{1} \hookrightarrow X = Spec\ A_{0}\ .
$$

On dispose ainsi d'un diagramme commutatif

$$
\begin{array}{c}
\xymatrix{
X_{1} = Spec(k(x)) \ar@{^{(}->}[r] \ar@{^{(}->}[d]_{i_{x}} & Y_{1} = X_{1} \ar[r] \ar[d]^{v} & P_{1} \ar@{=}[r] \ar[d]^{u} & Spf(\mathcal{V}(x)) = : \mathcal{S}_{1}\\
X   \ar@{^{(}->}[r]  & Y = P_{k}  \ar[r] & P \ar[r] & Spf\ \mathcal{V} =: \mathcal{S}\ ;
}
\end{array}
\leqno{(2.1.2.1)}
$$

\noindent de la m\^eme mani\`ere que [B 3, (2.3.2) (iv)] $u^{\ast}_{K}$ est donc une ``r\'ealisation''\ de $i^{\ast}_{x}$. De plus $u_{K}$ induit un morphisme [B 3, (0.2.7)]\\

\noindent (2.1.2.2) $\qquad \qquad \tau_{x} : \ ]X_{1}[_{P_{1}} = : \ ]x[_{P_{1}} = Spm(K(x)) \hookrightarrow ]X[_{P}\  = Spm(\hat{A}_{K})$ \\

\noindent commutant aux actions de Frobenius, et la restriction $\tau^{\ast}_{x}$ de $u^{\ast}_{K}$ est exacte (m\^eme argument que pour [B 3, (2.3.3) (iv)]) : \'evidemment $\tau_{x}$ est aussi induit par $\hat{\tau}_{K}(x) : \hat{A}_{K} \rightarrow K(x)$. Ainsi les fl\`eches canoniques
$$
\begin{array}{c}
\xymatrix{
\Gamma(]X[_{P}, \mathcal{E}_{\vert X})  = :  \mathcal{M} \ar[r] & \Gamma(]X[_{P}, R\ \tau_{x^{\ast}}\ \tau^{\ast}_{x}(\mathcal{E}))\ar[d]^{\simeq}&{}&{}& \\
{}& \Gamma(]x[_{P_{1}}, \tau^{\ast}_{x}(\mathcal{E}))\ar[d]^{\simeq}\ar[r]^{ \simeq}& \hat{\tau}_{K}(x)^{\ast} (\mathcal{M}) = \mathcal{M}_{x} \\
 {}& H^0_{\textrm{rig}}(Spec(k(x))/K(x), i^{\ast}_{x}(\mathcal{E})) = \mathcal{E}_{x}
}
\end{array} \leqno{(2.1.2.3)}
$$

\noindent sont compatibles aux Frobenius : par suite $\phi'_{x} : = H^0_{\textrm{rig}}(i^{\ast}_{x}(\phi))$ co\"{\i}ncide avec le morphisme

$$
\phi_{x} : \sigma^{\ast}_{K(x)}(\mathcal{M}_{x}) \rightarrow \mathcal{M}_{x}\ \textrm{du}\  (1.2.2)\  ,
$$

\vskip 3mm
\noindent (2.1.2.4) $\qquad \qquad \qquad \qquad \qquad \phi_{x} = \phi'_{x} \ ; $\\

\noindent de m\^eme pour $\tilde{\phi}_{x}$.\\

L'it\'er\'e deg $x$ fois de $\phi_{x}$ est un automorphisme $K(x)$-lin\'eaire du $K(x)$-espace vectoriel de dimension finie $\mathcal{E}_{x}$

$$
\phi^{\textrm{deg}\ x}_{x} : \mathcal{E}_{x} \displaystyle \mathop{\rightarrow}^{\sim} \mathcal{E}_{x}\ .
$$

\noindent Le $K$-espace vectoriel $\tilde{\mathcal{E}}_{x}$ d\'efini en (2.1.1.4) n'est autre que $\mathcal{E}_{x}$ vu comme $K$-espace vectoriel et $\tilde{\phi}_{x}$ son morphisme de Frobenius.\\

Gr\^ace au lemme (1.2.5) on obtient

\vskip 3mm
\noindent \textbf{Lemme (2.1.2.5)}. \textit{Avec les notations de (2.1.1) et (2.1.2) on a, pour $x \in \vert S \vert $,}

$$
det (1 - T\  \phi^{\textrm{deg}\ x}_{x} \vert \mathcal{E}_{x})\\
= det(1 - T\  \tilde{\phi}^{\textrm{deg}\ x}_{x} \vert \mathcal{E}_{x})^{1/\textrm{deg}\ x} \in K[T]\ .
$$
\vskip 2mm
On d\'eduit de (2.1.2.5) :

\vskip 3mm
\noindent \textbf{D\'efinition et proposition (2.1.3)}. \textit{Soient $S$ un $k$-sch\'ema s\'epar\'e de type fini et $(\mathcal{E}, \phi_{\mathcal{E}}) \in F^a\mbox{-}Isoc(S/K)$. La fonction $L$ de $\mathcal{E}$ est d\'efinie par}

$$
L(S, \mathcal{E}, t) = \displaystyle \mathop{\prod}_{x \in \vert S \vert}\ det(1 - t^{\textrm{deg}\ x}\ \phi^{\textrm{deg}\ x}_{x} \vert \mathcal{E}_{x})^{-1} \in K[[t]]
$$
$ \qquad \qquad \qquad \qquad \qquad = \displaystyle \mathop{\prod}_{x \in \vert S \vert}\ det(1 - t^{\textrm{deg}\ x}\ \tilde{\phi}^{\textrm{deg}\ x}_{x} \vert \tilde{\mathcal{E}}_{x})^{-1/ \textrm{deg}\ x} .$\\

Dans le cas d'une famille de vari\'et\'es on a le r\'esultat suivant :

\vskip 3mm
\noindent \textbf{Th\'eor\`eme (2.1.4)}. \textit{Supposons $e \leqslant p-1$. Soient $S$ un $k$-sch\'ema s\'epar\'e de type fini et $f : X \rightarrow S$ un $k$-morphisme propre et lisse. Alors, pout tout entier $i \geqslant 0$, la fonction $L(S, R^i f_{\textrm{conv}^{\ast}}(\mathcal{O}_{X/K}), t)$ est rationnelle :}

$$
L(S, R^i f_{\textrm{conv}^{\ast}}(\mathcal{O}_{X/K}), t) \in K(t)\ .
$$

\vskip 3mm
\noindent \textit{D\'emonstration}. On peut supposer $S$ lisse sur $k$. D'apr\`es le [III, th\'eor\`eme (3.2.1)], $R^i f_{\textrm{conv}^{\ast}}(\mathcal{O}_{X/K}) \in F^a\mbox{-}\textrm{Isoc}(S/K)$ et sa formation commute au passage aux fibres [III, prop (3.3.4)] : pour $s \in \vert S \vert $ on a $(R^i f_{\textrm{conv}^{\ast}}(\mathcal{O}_{X/K}))_{s} = H^i_{\textrm{cris}}(X_{s}/\mathcal{V}(s)) \otimes_{\mathcal{V}(s)}\ K(s)$.\\

La suite de la d\'emonstration est identique \`a celle du [th\'eor\`eme 3 (1)] de [Et 3] en rempla\c{c}ant le th\'eor\`eme de comparaison de Katz-Messing [K-M] du cas projectif lisse par celui de Chiarellotto-Le Stum [C-LS 1, cor 1.3] pour le cas propre et lisse. $\square$\\

\subsection*{2.2. $F$-cristaux}

Supposons, jusqu'\`a la fin de ce \S\ 2, que l'indice de ramification $e$ de $\mathcal{V}$ v\'erifie $e \leqslant p-1$, et notons $F^a\mbox{-}\textrm{Cris}(S/\mathcal{V})$ la cat\'egorie des $F^a$-cristaux localement libres de type fini [B 1] et [Et 2, II, \S\ 2]. Soit $E \in F^{a}\mbox{-}\textrm{Cris}(S/\mathcal{V})$ : rappelons la d\'efinition de la fonction $L(S, E, t)$ [Et 2]. Pour $x \in \vert S \vert $, notons $i_{x} : Spec\ k(x) \hookrightarrow S$ l'immersion canonique et $E_{x} : = i^{\ast}_{x}(E) (Spec(k(x)), Spec(\mathcal{V}(x)))$. L'it\'er\'e deg $x$ fois de $\phi : F^{\ast}_{S}(E) \rightarrow E$ induit un endomorphisme $\mathcal{V}(x)$-lin\'eaire, $F_{x}$, du $\mathcal{V}(x)$-module libre de rang fini $E_{x}$, et, par le m\^eme argument qu'en (1.2), le ``polyn\^ome caract\'eristique'' det$(1 - t\ F_{x})$ est en fait \`a coefficients dans $\mathcal{V}$. Notons $\tilde{E}_{x}$ le module $E_{x}$ vu comme $\mathcal{V}$-module et $\tilde{F}_{x}$ l'endomorphisme de $\tilde{E}_{x}$ d\'eduit de $F_{x}$. La fonction $L$ de $E$ est d\'efinie par [Et 2, II, \S\ 2, p 50-51]

$$\begin{array}{c}
\xymatrix{
L(S, E, t) : = \displaystyle \mathop{\prod}_{x \in \vert S \vert} det(1 - t^{\textrm{deg}\ x}\  \tilde{F}_{x}\vert \tilde{E}_{x})^{-1/\textrm{deg}\ x }\\
\qquad \qquad \qquad = \displaystyle \mathop{\prod}_{x \in \vert S \vert} det(1 - t^{\textrm{deg}\ x}\  F_{x}\vert E_{x})^{-1} \in \mathcal{V} [[t]]\ .
}
\end{array}
\leqno{(2.2.1)
}
$$

\noindent Notons $\mathcal{E }: = E^{an}$ le $F$-isocristal convergent associ\'e \`a $E$ [B 3, (2.4.2)]. Il est clair d'apr\`es les d\'efinitions que l'on a l'\'egalit\'e \\

\noindent (2.2.2) $\qquad \qquad \qquad L(S, E, t) = L(S, \mathcal{E}, t)\ .$\\

\section*{3. Fonctions $L$ des $F$-isocristaux Dwork -surconvergents}

\subsection*{3.1. $F$-isocristaux Dwork-surconvergents}

En 3.1 on suppose simplement que $k$ est un corps de caract\'eristique $p > 0$.

\vskip 2mm
Soit $X$ un $k$-sch\'ema s\'epar\'e de type fini. On suppose qu'il existe une stratification de $X$ par des $k$-sch\'emas lisses $X_{i}$, $X = \displaystyle \mathop{\coprod}_{i} X_{i}$ : ceci est possible par exemple si $k$ est parfait et $X$ r\'eduit [EGA IV, 17.15.13]. Fixons un recouvrement de $X_{i}$ par des $k$-sch\'emas affines et lisses $X_{ij}, X_{i} = \displaystyle \mathop{\cup}_{j}\ X_{ij}$, et des $\mathcal{V}$-sch\'emas formels lisses $\mathcal{X}_{ij}$ relevant $X_{ij}$, $\mathcal{X}_{ij} = Spf(\hat{A}_{ij})$ o\`u $A_{ij}$ est une $\mathcal{V}$-alg\`ebre lisse. Notons $\mathcal{X}^{\dag}_{ij} = Spff(A^{\dag}_{ij})$ le sch\'ema faiblement formel associ\'e \`a $A^{\dag}_{ij}$ [Mer], $F_{A^{\dag}_{ij}}$ un rel\`evement du Frobenius de $X_{ij}$ comme en 1.1, et $F_{\hat{A}_{ij}} = F_{A^{\dag}_{ij}}  \otimes_{A^{\dag}_{ij}} \hat{A}_{ij}$.\\

Consid\'erons maintenant $\mathcal{E} \in F^a\mbox{-}\textrm{Isoc}(X/K)$ et pour chaque $(i,j)$ soit $\mathcal{E}_{ij}$ un $\hat{A}_{ij\ K}$-module projectif de type fini correspondant \`a la restriction de $\mathcal{E}$ \`a $X_{ij}$, dont l'existence est assur\'ee par le corollaire (1.2.3) du II : par cette \'equivalence la structure de Frobenius de $\mathcal{E}$ fournit pour tout $(i,j)$ un isomorphisme

$$
\hat{\phi}_{ij} : F^{\ast}_{\hat{A}_{ij}}(\mathcal{E}_{ij}) \displaystyle \mathop{\rightarrow}^{\sim}
\mathcal{E}_{ij}$$

\noindent compatible aux connexions.\\

Puisque $(A^{\dag}_{ij}, A^{\dag}_{ij}/(\pi))$ est un couple hens\'elien [Et 4, th\'eo 3] il existe d'apr\`es Elkik [E$\ell$, cor 1 du th\'eo 3, p 573] un $A^{\dag}_{ij\ K}$-module projectif de type fini $E_{ij}$ tel que

$$
\mathcal{E}_{ij} \simeq E_{ij}\  \otimes_{A^{\dag}_{ij\ K}} \hat{A}_{ij\ K}\ .
$$

\vskip 3mm
\noindent \textbf{D\'efinition (3.1.1)}. \textit{Avec les notations ci-dessus, on dira que $\mathcal{E}$ est \textbf{Dwork-surconvergent} pour la stratification pr\'ec\'edente si $\hat{\phi}_{ij}$ provient, par l'extension des scalaires $A^{\dag}_{ij\ K} \hookrightarrow \hat{A}_{ij\ K}$, d'un isomorphisme} 

$$
\phi^{\dag}_{ij} : F^{\ast}_{A^{\dag}_{ij}}(E_{ij}) \displaystyle \mathop{\rightarrow}^{\sim} E_{ij}\ .
$$

\textit{Les morphismes entre de tels $\mathcal{E}$ se d\'efinissent de la mani\`ere naturelle.}\\

\textit{On notera}

$$
\mathbf{Dw^{\dag}(X, (X_{ij}), (\mathcal{X}^{\dag}_{ij}), (F_{\mathcal{X}^{\dag}_{ij}}))}
$$

\noindent \textit{la cat\'egorie ainsi obtenue et}

$$
\begin{array}{c}
\xymatrix{
\mathbf{\mathcal{D} : Dw^{\dag}(X, (X_{ij}), (\mathcal{X}^{\dag}_{ij}), (F_{\mathcal{X}^{\dag}_{ij}})} \ar[r]  & \mathbf{F^a\mbox{-}Isoc(X/K)}\\
\mathbf{\mathcal{E}^{\dag} : = (\mathcal{E}, E_{ij}, \phi^{\dag}_{ij})} \ar@{|->}[r] & \mathbf{\mathcal{E}}
}
\end{array}
\leqno{(3.1.2)}
$$

\noindent \textit{le foncteur naturel ``oubli des stratifications''.}

\vskip 3mm
\noindent \textbf{Remarque (3.1.3)}. \textit{La notion ``Dwork-surconvergent'' d\'epend de la stratification : un $F$-isocristal convergent $\mathcal{E}$ peut tr\`es bien \^etre Dwork-surconvergent pour une stratification et pas pour une autre. On retrouve le m\^eme probl\`eme que Wan avec les $\sigma$-modules surconvergents : [cf W 2, p 885].}\\

Les propri\'et\'es suivantes sont claires :

\vskip 3mm
\noindent \textbf{Proposition (3.1.4)}. \textit{Le foncteur $\mathcal{D}$ (3.1.2) est fid\`ele.}

\vskip 3mm
\noindent \textbf{Proposition (3.1.5)}. \textit{La cat\'egorie $Dw^{\dag}(X, (X_{ij}), (\mathcal{X}^{\dag}_{ij}), (F_{\mathcal{X}^{\dag}_{ij}}))$ est stable par puissances tensorielles, puissances sym\'etriques, puissances ext\'erieures et passage au dual}.\\

\newpage

\subsection*{3.2. Fonction $L$ des $F$-isocristaux Dwork-surconvergents}

\noindent \textbf{D\'efinition (3.2.1)}. \textit{Avec les notations de (3.1), et $k$ fini, on d\'efinit la fonction $L$ de $\mathcal{E}^{\dag} \in Dw^{\dag} (X, (X_{ij}), (\mathcal{X}^{\dag}_{ij}), (F_{\mathcal{X}^{\dag}_{ij}}))$ par
$$
L(X, \mathcal{E}^{\dag}, t) = L(X, \mathcal{E}, t)
$$
\noindent o\`u $\mathcal{E}$ est l'image de $\mathcal{E}^{\dag}$ par le foncteur d'oubli $\mathcal{D}$ (3.1.2).}\\

Pour un $k$-sch\'ema s\'epar\'e de type fini $X$ et $\mathcal{E} \in F^a\mbox{-}\textrm{Isoc}(X/K)$ on a vu en 2.1 que la d\'efinition de la fonction $L$ de $\mathcal{E}$ via la cohomologie rigide co\"{\i}ncide avec la  d\'efinition via les $F$-modules convergents : ainsi pour $\alpha \in \mathbb{Q}$ et $r \in \mathbb{N}^{\ast}$ on d\'efinit,  comme en (1.4.8),\\

\noindent (3.2.2) $\qquad \qquad L^{(r)}_{\alpha}(X, \mathcal{E}, t) := \displaystyle \mathop{\prod}_{x \in \vert X \vert} det^{(r)}_{\alpha}(\mathcal{E}_{x},  t^{\textrm{deg}\ x})^{-1} \in K [[t]].$\\

On a encore :\\

\noindent (3.2.3) $\qquad \qquad \qquad \qquad  L^{(r)}(X, \mathcal{E}, t) := \displaystyle \mathop{\prod}_{\alpha \in \mathbb{Q}} L^{(r)}_{\alpha}(X, \mathcal{E}, t)\ ;$\\

\noindent plus pr\'ecis\'ement, on a\\

\noindent (3.2.4) $\qquad \qquad \qquad   L^{(r)}(X, \mathcal{E}, t) = \displaystyle \mathop{\prod}_{i \geqslant 1} L(X, \textrm{Sym}^{r-i}\ \mathcal{E} \otimes \displaystyle \mathop{\Lambda}^i \mathcal{E}, t)^{i\  \times\ (-1)^{i-1}}\ ,$\\

\vskip 2mm
\noindent (3.2.5) $\qquad \qquad \qquad   L^{(r)}_{\alpha}(X, \mathcal{E}, t) = \displaystyle \mathop{\prod}_{i \geqslant 1} L_{\alpha}(X, \textrm{Sym}^{r-i}\ \mathcal{E} \otimes \displaystyle \mathop{\Lambda}^i \mathcal{E}, t)^{i\  \times\ (-1)^{i-1}}\ .$\\

Ces d\'efinitions s'\'etendent au cas des $F$-isocristaux Dwork-surconvergents $\mathcal{E}^{\dag}$ en posant, pour $\mathcal{E} = \mathcal{D}(\mathcal{E}^{\dag})$ :\\

\noindent  (3.2.6) $\qquad \qquad \qquad  \qquad \qquad  L^{(r)}_{\alpha}(X, \mathcal{E}^{\dag}, t) := L^{(r)}_{\alpha} (X, \mathcal{E}, t)\ ,$\\

\vskip 2mm
\noindent (3.2.7) $\qquad \qquad \qquad \qquad \qquad L^{(r)}(X, \mathcal{E}^{\dag}, t) = L^{(r)} (X, \mathcal{E}, t)\ .$\\

Si $E$ est un $F$-cristal, $E \in F^a\mbox{-}\textrm{Cris}(X/\mathcal{V})$, avec $e \leqslant p-1$ et $\mathcal{E} := E^{an}$ est le $F$-isocristal convergent associ\'e on pose \'egalement\\

\noindent (3.2.8)  $\qquad \qquad \qquad \qquad \qquad L^{(r)}_{\alpha}(X, E, t) := L^{(r)}_{\alpha} (X, \mathcal{E}, t)\ $\\

\vskip 2mm
\noindent (3.2.9)  $\qquad \qquad \qquad \qquad \qquad L^{(r)}(X, E, t) = L^{(r)} (X, \mathcal{E}, t)\ .$\\

La conjecture de Dwork est vraie pour les $F$-cristaux Dwork-surconvergents,  \`a savoir que l'on a le th\'eor\`eme suivant, par application directe du th\'eor\`eme (1.4.13) :

\vskip 3mm
\noindent \textbf{Th\'eor\`eme (3.2.10)}. \textit{Soient $X$ un $k$-sch\'ema s\'epar\'e de type fini muni d'une stratification comme en (3.1) et  $\mathcal{E}^{\dag} \in Dw^{\dag}(X, (X_{ij}), (\mathcal{X}_{ij}^{\dag}), (F_{\mathcal{X}^{\dag}_{ij}}))$. Alors, pour tout nombre rationnel $\alpha \in \mathbb{Q}$ et tout entier $r \in \mathbb{N}^{\ast}$, les fonctions $L^{(r)}_{\alpha}(X, \mathcal{E}^{\dag}, t)$ et $L^{(r)}(X, \mathcal{E}^{\dag}, t)$ sont $p$-adiquement m\'eromorphes.}

\vskip 3mm

\section*{4. Fonctions $L$ des $F$-isocristaux surconvergents}

\section*{4.1. D\'efinitions} Soient $S$ un $k$-sch\'ema s\'epar\'e de type fini,  $E \in F^a\mbox{-}\textrm{Isoc}^{\dag} (S/K)$ et  $\phi_{E} : F^{\ast}_{\sigma}\  E \displaystyle \mathop{\rightarrow}^{\sim} E$ son isomorphisme de Frobenius [B 3, (2)]. Avec les notations de (2.1) et  $x \in \vert Spec\ A_{0} \vert \subset \vert S \vert $ on a aussi des foncteurs images inverses\\

\noindent (4.1.1) $\qquad \qquad i^{\ast}_{x} : F^a\mbox{-}\textrm{Isoc}^{\dag}(S/K) \rightarrow F^a\mbox{-}\textrm{Isoc}^{\dag} (Spec(k(x))/K(x))$

\vskip 3mm
\noindent [resp.\\
\noindent (4.1.2) $\qquad \qquad \tilde{i}^{\ast}_{x} : F^a\mbox{-}\textrm{Isoc}^{\dag}(S/K) \rightarrow F^a\mbox{-}\textrm{Isoc}^{\dag} (Spec(k(x))/K)]\ ,$\\

\noindent et\\

\noindent (4.1.3) $\qquad \qquad E_{x} : = H^{0}_{\textrm{rig}} (Spec(k(x))/K(x), i^{\ast}_{x}(E))$\

\vskip 3mm
\noindent [resp.\\
\noindent (4.1.4) $\qquad \qquad \tilde{E}_{x} : = H^{0}_{\textrm{rig}} (Spec\ k(x)/K,\tilde{i}^{\ast}_{x}(E))]$\\

\noindent est un $K(x)$-espace vectoriel [resp. $K$-espace vectoriel] de dimension finie muni d'un isomorphisme de Frobenius\\

\noindent (4.1.5) $\qquad \qquad \qquad \qquad \phi'_{x} = H^o_{\textrm{rig}}(i^{\ast}_{x}(\phi))$

\vskip 3mm
\noindent [resp.\\
\noindent (4.1.6) $\qquad \qquad \qquad \qquad \tilde{\phi}'_{x} = H^o_{\textrm{rig}}(\tilde{i}^{\ast}_{x}(\phi))]\ .$\\

Posons\\

\noindent (4.1.7)  $\qquad \qquad \qquad M : = \Gamma(\mathcal{X}^{an}_{K}, E_{\vert X}), \phi_{M} : = \Gamma({\mathcal{X}}^{an}_{K}, \phi_{E_{\vert X}})$\\

\noindent o\`u $X = Spec\ A_{0}, \mathcal{X} = Spec\ A$ ; alors $M \in F^a\mbox{-}\textrm{Mod}(A^{\dag}_{K})$, et la fibre $M_{x}$ de $M$ en $x \in \vert X \vert $ est par d\'efinition (cf. (1.3.1))\\

\noindent (4.1.8) $\qquad \qquad \qquad \qquad \qquad M_{x} : = \tau^{\dag}_{K}(x)^{\ast} (M)$\\

\noindent et

$$
\phi_{M} :  M^{\sigma} : = F^{\ast}_{A^{\dag}_{K}} (M) \rightarrow M
$$

\noindent induit (cf. (1.3.2))\\

\noindent (4.1.9) $\qquad \qquad \qquad \qquad   \phi_{x} :  \sigma^{\ast}_{K(x)} (M_{x}) \rightarrow M_{x }\ .$\\

On montre, comme en (2.1.2), que l'on a une identification\\

\noindent (4.1.10) $\qquad \qquad \qquad \qquad (E_{x}, \phi'_{x}) \equiv (M_{x}, \phi_{x})\ .$\\

On d\'eduit de (1.3.5) le lemme suivant :

\vskip 3mm
\noindent \textbf{Lemme (4.1.11)}. \textit{Avec les notations de (4.1) on a, pour $x \in \vert S \vert$,}

$$
det (1 - T\ \phi^{\textrm{deg}\ x}_{x} \vert E_{x}) = det(1 - T\ \tilde{\phi}^{\textrm{deg}\ x}_{x} \vert \tilde{E}_{x})^{1/\textrm{deg}\ x} \in K[T]\ .
$$
\vskip 3mm
La d\'efinition de la fonction $L(S, E, t)$ donn\'ee ci-dessous ne d\'epend que du sous-sch\'ema r\'eduit sous-jacent \`a $S$ : si l'on suppose $S$ r\'eduit on peut alors le d\'ecomposer en strates $S_{ij}$-affines et lisses sur $k$ \`a la mani\`ere de 3.1 ; on dispose alors d'un foncteur naturel\\

\noindent (4.1.12) $\qquad \quad \mathcal{D}^{\dag} : F^a\mbox{-}\textrm{Isoc}^{\dag}(S/K) \rightarrow Dw^{\dag}(S, (S_{ij}), (\mathcal{S}^{\dag}_{ij}), (F_{\mathcal{S}^{\dag}_{ij}}))\ $

$$
E \longmapsto \mathcal{E}^{\dag} = (\mathcal{E}, E_{ij}, \phi^{\dag}_{ij})
$$

\noindent o\`u $\mathcal{E}$ est le $F$-isocristal convergent associ\'e \`a $E$, et $\phi^{\dag}_{ij} : F^{\ast}_{A^{\dag}_{ij\ K}}(E_{ij}) \simeq E_{ij}$ est l'isomorphisme de Frobenius provenant de l'\'equivalence de Berthelot [B 3, cor (2.5.8)]. De mani\`ere imag\'ee on dira qu'un $F$-isocristal surconvergent est toujours Dwork-surconvergent.\\

\newpage
Comme cons\'equence de (4.1.11) on a :\\

\noindent \textbf{D\'efinitions et proposition (4.1.13)}. \textit{Soient $S$ un $k$-sch\'ema s\'epar\'e de type fini, $(E, \phi_{E}) \in F^a\mbox{-}\textrm{Isoc}^{\dag}(S/K)$ et $(\mathcal{E}, \phi_{\mathcal{E}}) \in  F^a\mbox{-}\textrm{Isoc}(S/K)$ son image canonique par le foncteur d'oubli $F^a\mbox{-}\textrm{Isoc}^{\dag}(S/K) \rightarrow F^a\mbox{-}\textrm{Isoc}(S/K)$. La fonction $L$ de $E$ est d\'efinie par }\\

$\qquad \qquad L(S, E, t) = \displaystyle \mathop{\prod}_{x \in \vert S \vert} \ det(1 - t^{\textrm{deg}\ x}\ \phi^{\textrm{deg}\ x}_{x} \vert E_{x})^{-1} \in K[[t]]$\\

$\qquad \qquad \qquad \qquad  Ê  =  \displaystyle \mathop{\prod}_{x \in \vert S \vert} \ det(1 - t^{\textrm{deg}\ x}\ \tilde{\phi}^{\textrm{deg}\ x}_{x} \vert \tilde{E}_{x})^{- 1/\textrm{deg}\ x}\ ,$\\

\noindent \textit{et on a}

$$
L(S, E, t) = L(S, \mathcal{E}, t)\ .
$$\\

\noindent \textit{On pose, pour $\alpha \in \mathbb{Q}$ et $r \in \mathbb{N}^{\ast}$}

$$
L^{(r)}_{\alpha}(S, E, t) = L^{(r)}_{\alpha}(S, \mathcal{E}, t)
$$

$$
L^{(r)}(S, E, t) = L^{(r)}(S, \mathcal{E}, t)\ .
$$

\vskip 3mm
\noindent \textbf{Remarques (4.1.14)}. 

\begin{enumerate}
\item[(i)]
La deuxi\`eme expression d\'efinissant $L(S, E, t)$ \`a l'aide de $\tilde{E}_{x}$ est celle qui avait \'et\'e utilis\'ee dans la d\'efinition de [E-LS 1] : la nouvelle d\'efinition donn\'ee ici avec l'appui du lemme (1.3.5) co¬\"{\i}ncide donc avec celle de [E-LS 1].
\item[(ii)]
La fonction $L(S, E, t)$ ci-dessus est $p$-adiquement m\'eromorphe : ceci se d\'eduit ou bien du th\'eor\`eme (1.3.8), ou bien du th\'eor\`eme 6.3 de [E-LS 1] qui fournit l'expression cohomologique
$$
L(S, E, t) = \displaystyle \mathop{\prod}^{2\ \textrm{dim}\ S}_{i=0}\ det (1 - t\ F \vert H^i_{\textrm{rig}, c}(S/K,E))^{(-1)^{i+1}}\ .
$$
\noindent De plus cette expression cohomologique fournit la rationalit\'e, $L(S, E, t) \in K(t)$, depuis que Kedlaya a prouv\'e la finitude de la cohomologie rigide \`a supports compacts et \`a coefficients dans $E$ [Ked 1].
\end{enumerate}

\vskip 3mm
\subsection*{4.2. La conjecture de Dwork}

\noindent {\textbf{Th\'eor\`eme (4.2.1)}. \textit{Soient $S$ un $k$-sch\'ema s\'epar\'e de type fini et $(E, \phi_{E}) \in F^a\mbox{-}\textrm{Isoc}^{\dag}(S/K)$. Alors, pour tout rationnel $\alpha$ et tout entier $r \in \mathbb{N}^{\ast}$, on a :}

\begin{enumerate}
\item[(i)] \textit{La fonction $L^{(r)}_{\alpha}(S, E, t)$ est $p$-adiquement m\'eromorphe}
\item[(ii)] \textit{La fonction $L^{(r)}(S, E, t)$ est rationnelle}
\end{enumerate}

$$
L^{(r)}(S, E, t) \in K(t).
$$
\vskip 3mm
\noindent \textit{D\'emonstration}. \\

Le (i) r\'esulte de (3.2.10).

Le (ii) r\'esulte de (4.1.14) (ii) et de la relation (3.2.4) :

$$
L^{(r)}(S, E, t) = \displaystyle \mathop{\prod}_{i\geqslant 1} L(S, \textrm{Sym}^{r-i}\ E \otimes \displaystyle \mathop{\Lambda}^i\ E, t)^{i \times (-1)^{i-1}}\ . \qquad \square
$$

\vskip 3mm
\subsection*{4.3. La conjecture de Katz}

Dans le th\'eor\`eme suivant nous levons l'hypoth\`ese de prolongement \`a une compactification qui avait \'et\'e faite dans [E-LS 2, th\'eo 6.7] pour r\'esoudre la conjecture de Katz.\\

On note $\overline{k}$ une cl\^oture alg\'ebrique de $k = \mathbb{F}_{q}$.

\vskip 3mm
\noindent \textbf{Th\'eor\`eme (4.3.1)}. \textit{Soient $X$ un sch\'ema s\'epar\'e de type fini sur $\mathbb{F}_{q}$, $X_{\overline{k}} = X \times_{k} \overline{k}$ et $\mathcal{F}_{\mathbb{Q}}$ un $K^{\sigma}$-faisceau lisse sur $X$. On suppose que le $F$-isocristal convergent unit\'e $\mathcal{E}$ associ\'e \`a $\mathcal{F}_{\mathbb{Q}}$ provient d'un $F$-isocristal surconvergent $E$, $\mathcal{E} = \hat{E}$. Alors}\\

\noindent (4.3.1.1) $\qquad \qquad L(X, \mathcal{F}_{\mathbb{Q}}, t) = L(X, \mathcal{E}, t) \in K(t).$\\

\noindent (4.3.1.2)  \textit{Le quotient
$$
L(X, \mathcal{F}_{\mathbb{Q}}, t)/\displaystyle \mathop{\prod}_{i} det(1 - t\ F | H^i_{\textrm{\'et}, c}(X_{\overline{k}}, \mathcal{F}_{\mathbb{Q}}))^{(-1)^{i+1}}
$$
n'a ni z\'ero ni p\^ole sur la couronne unit\'e $\vert t \vert_{p} = 1.$}

 \vskip 3mm
\noindent  \textit{D\'emonstration}. La premi\`ere assertion est l\`a pour m\'emoire. Pour la deuxi\`eme il suffit d'utiliser la suite exacte (2) du th\'eor\`eme (3.3.13) du IV et [E-LS 2, prop 6.5]. $\square$\\

\noindent \textbf{Remarque (4.3.1.3)}. En utilisant l'expression cohomologique [E-LS 1, 6.3] de $L(X, \mathcal{F}_{\mathcal{Q}}, t) = L(X, E, t)$ rappel\'ee en (4.1.13) (ii), c'est chaque terme du produit altern\'e figurant dans l'expression (4.3.1.2) qui peut \^etre pr\'ecis\'e, c'est-\`a-dire \`a $i$ fix\'e, gr\^ace \`a [E-LS 1, 6.5] et au [IV, (3.3.16)(2)].

Ceci apporte une pr\'ecision \`a la d\'emonstration de la conjecture de Katz faite par Emerton et Kisin dans [E-K 1].\\

On g\'en\'eralise (4.3.1.1) de la fa\c{c}on suivante.

\vskip 3mm
\noindent \textbf{Th\'eor\`eme (4.3.2)}. \textit{Soient $S$ un $\mathbb{F}_{q}$-sch\'ema lisse et s\'epar\'e, $f : X \rightarrow S$ un $k$-morphisme propre et lisse satisfaisant aux hypoth\`eses de [IV, (3.1)] ou [IV, (5.2)] et $E$, $\mathcal{F}_{\mathbb{Q}}$ comme en (4.3.1). Alors, pour $(i, r, \alpha) \in \mathbb{N} \times \mathbb{N}^{\ast} \times \mathbb{Q}$, on a :}\\

\noindent (4.3.2.1) $\qquad \qquad  L(X, \mathcal{F}_{\mathbb{Q}}, t) = \displaystyle \mathop{\prod}_{i} L(S, R^i f_{\textrm{rig}\ast}(E), t)^{(-1)^{i}} \in{K(t)},$\\
\noindent\textit{c'est-\`a-dire}\\
$$
\displaystyle \mathop{\prod}_{j} det(1 - t\ F | H^j_{\textrm{rig}, c}(X/K, E))^{(-1)^{j+1}}= \displaystyle \mathop{\prod}_{i,j} det(1 - t\ F | H^j_{\textrm{rig}, c}(S/K, R^{i}f_{\textrm{rig}\ast}(E)))^{(-1)^{i+j+1}}\in K(t).
$$

\vskip 2mm
\noindent (4.3.2.2) $\qquad \qquad  L^{(r)}(S, R^i f_{\textrm{rig}\ast}(E), t) \in K(t).$

\vskip 3mm
\noindent (4.3.2.3) $\qquad \qquad  L^{(r)}_{\alpha}(S, R^i f_{\textrm{rig}\ast }(E), t)\  \textit{est}\  p\textit{-adiquement m\'eromorphe}.$

\vskip 3mm
\noindent (4.3.2.4) $\qquad \qquad  L_{0}(S, R^i f_{\textrm{rig}\ast }(E), t) = L(S, R^i f_{\textrm{et}\ast }(\mathcal{F}_{\mathbb{Q}}), t)$.

\vskip 3mm
\noindent (4.3.2.5)  \textit{Si de plus $\mathcal{F}$ est un $\mathcal{V}^{\sigma}$-faisceau lisse tel que $\mathcal{F}_{\mathbb{Q}} = \mathcal{F} \otimes \mathbb{Q}$, alors le quotient}

$$
L(S, R^i f_{\textrm{rig}\ast }(E), t)/\displaystyle \mathop{\prod}_{j} det(1 - t\ F \vert H^j_{\textrm{et},c}(S_{\overline{k}}, R^i f_{\textrm{et}\ast }(\mathcal{F}_{\mathcal{\mathbb{Q}}})))^{(-1)^{j+1}}
$$

\noindent \textit{n'a ni z\'ero ni p\^ole dans la couronne unit\'e $\vert t \vert_{p} = 1$ (ni m\^eme dans le disque ferm\'e $\vert t \vert_{p} \leqslant 1)$.}\\

\noindent \textit{D\'emonstration}. Par d\'efinition on a

$$
L(X, \mathcal{F}_{\mathbb{Q}}, t) = L(X, E, t),
$$

\noindent et il suffit d'appliquer [IV, (3.1) ou (5.2)], (4.1.13) et (4.2.1) ; d'o\`u les trois premi\`eres assertions.\\

Pour (4.3.2.4) on utilise la suite exacte de [V, (3.3.13.3)] et (E-LS 2, (6.2)].\\

Pour (4.3.2.5) on pose $E^i = R^i f_{\textrm{rig}\ast }(E), L_{0}(t) = L_{0}(S, E^i, t), L_{>0}(t) = L(S, E^i, t) / L_{0}(t)$.\\

\noindent Par d\'efinition on a

$$
L(S, E^i, t) = \displaystyle \mathop{\prod}_{s \in \vert S \vert}(det(1 - t^{\textrm{deg} \ s}\ F_{s} \vert H^i_{\textrm{rig}\ast }(X_{s}/K(s), E_{\vert X_{s}})))^{-1}
$$
$
 \qquad \qquad \qquad \qquad \qquad = :  \displaystyle \mathop{\prod}_{s \in \vert S \vert}(P_{s}(t))^{-1}\ .
$\\

\noindent Pour chaque $s \in \vert S \vert$ l'injection $k(s) \hookrightarrow \overline{k}$ s'\'etend en une injection

$$
\mathcal{V}(s) : = W(k(s)) \otimes_{W(k)} \mathcal{V} \hookrightarrow \overline{\mathcal{V}} : = W(\overline{k}) \otimes_{W(k)} \mathcal{V},
$$

\noindent et on a

$$
P_{s}(t) \in \mathcal{V}(s) [t] \subset \overline{\mathcal{V}}[t]\ ;
$$

\noindent en effet,  comme $E$ provient d'un $\mathcal{V}^{\sigma}$-faisceau, le polyn\^ome $P_{s}$ est \`a coefficients dans $\mathcal{V}(s)$ et pas seulement dans $K(s) = Frac(\mathcal{V}(s))$. Donc  $L(S, E^i, t), L_{0}(t)$ ainsi que leurs inverses sont \'el\'ements de $\overline{\mathcal{V}}[[t]]$ : par suite ils convergent tous les quatre dans le disque ouvert $\vert t \vert_{p} < 1$ ; de plus, par d\'efinition de $L_{0}(t)$ on en d\'eduit que $L_{>0}(t)$ et $1/L_{>0}(t)$ convergent dans le disque ferm\'e $\vert t \vert_{p} \leqslant 1$. Or d'apr\`es [E-K 1], le quotient

$$
L(S, R^i f_{\textrm{et}\ast }(\mathcal{F}_{\mathbb{Q}}), t) / \displaystyle \mathop{\prod}_{j} det(1 - t\ F \vert H^j_{\textrm{et},c}(S_{\overline{k}}, R^i f_{\textrm{et}\ast }(\mathcal{F}_{\mathbb{Q}})))^{(-1)^{j+1}}
$$

\noindent converge ainsi que son inverse dans le disque $\vert t \vert_{p} \leqslant 1$ ; d'o\`u (4.3.2.5) via (4.3.2.4). $\square$

\newpage
\section*{5. Sch\'emas ab\'eliens ordinaires}

On suppose $e \leqslant p-1$. 

Dans ce paragraphe 5 nous allons expliciter l'expression $L^{(r)}_{\alpha}(S, R^i f_{\textrm{cris}\ast}(\mathcal{O}_{X/W}), t)$ lorsque $f : X \rightarrow S$ est un sch\'ema ab\'elien ordinaire.\\

Apr\`es des rappels en 5.1 sur les $F$-cristaux ordinaires, et en 5.2 sur les sch\'emas ordinaires, nous caract\'erisons en (5.2.9) les sch\'emas ab\'eliens ordinaires. Nous revenons ensuite en 5.3, th\'eor\`eme  (5.3.1), \`a l'expression de $L^{(r)}_{\alpha}$.\\

\subsection*{5.1. $F$-cristaux ordinaires}

Soient $k$ un corps de caract\'eristique $p > 0$, $S$ un $k$-sch\'ema, $E \in F^a\mbox{-}\textrm{Cris}(S/\mathcal{V})$ un $F^a$-cristal localement libre de type fini sur $(S/\mathcal{V})$\ . \\

A la suite de Katz et Deligne [K 2, II, 2.4, Rks p 148], [De$\ell$ 2], on dit que \textbf{$E$ est un $F^a$-cristal ordinaire de niveau $n$} s'il existe une filtration de $E$ par des sous $F^a$-cristaux localement libres de type fini\\

\noindent (5.1.1) $\qquad \qquad \qquad 0 \subset U_{0}Ê\subset U_{1} \subset ... \subset U_{i-1} \subset U_{i} \subset ... \subset U_{n} = E$\\

\noindent tels que $U_{i}/U_{i-1}$ soit de la forme $V_{i}(-i) := (E_{i}, \nabla_{i}, \pi^{i}\  F_{i})$, o\`u $V_{i} := (E_{i}, \nabla_{i}, F_{i})$ est un $F^a$-cristal unit\'e : $V_{i}$ sera not\'e $E_{i}$ dans la suite. On dit que $E \in F^a\mbox{-}\textrm{Cris}(S/\mathcal{V})$ est $\textbf{ordinaire}$ s'il existe $n \in \mathbb{N}$ tel que $E$ soit ordinaire de niveau $n$.\\

Pour les d\'efinitions et propri\'et\'es des polygones de Hodge et de Newton de $E$ nous renvoyons le lecteur \`a [K 2, I, 1.2, 1.3 et 2.3 p 142].\\

Si $k$ est parfait et $E \in F^a\mbox{-}\textrm{Cris}(S/\mathcal{V})$ est ordinaire, il est clair qu'en tout point $s \in S$ les polygones de Hodge et de Newton de $E$ co¬\"{\i}ncident et qu'ils sont constants, i.e. ind\'ependants de $s \in S$ [loc. cit.]. La r\'eciproque est vraie si $S$ est un sch\'ema sur un corps parfait $k$ de caract\'eristique $p > 0$ et $S$ de l'un des deux types suivants [K 2, II, 2.4 Rks p 148], [W 2, lemma 3.6], [W 3, theo 7.2, cor 7.3] :\\

\begin{enumerate}
\item[(i)] $S$ est lisse sur $k$\ ,
\item[(ii)] $S = Spec\ k[[X_{1},...,X_{d}]]$.
\end{enumerate}

Lorsque $k = \mathbb{F}_{q}$ et que $E$ est ordinaire de niveau $n$, on a :\\

\noindent (5.1.2) $\qquad \qquad \qquad L(S, E, t) = \displaystyle \mathop{\prod}_{i=0}^n L(S, E_{i}, \pi^i t).$

\vskip 3mm
\noindent \textbf{Proposition (5.1.3)}. \textit{Si $A_{0}$ est une $\mathbb{F}_{q}$-alg\`ebre lisse, la cat\'egorie des objets ordinaires de $F^a\mbox{-}\textrm{Cris}(Spec\ A_{0}/\mathcal{V})$ est stable par puissances tensorielles, puissances sym\'etriques et puissances ext\'erieures.}

\vskip 3mm
\noindent \textit{D\'emonstration}. Il suffit de raisonner avec des bases locales de $E \in F^a\mbox{-}\textrm{Cris}(Spec\ A_{0}/\mathcal{V})$ et d'appliquer [K 2] ou [W 2, discussion apr\`es def 3.3]. $\square$\\

\subsection*{5.2. Caract\'erisation des sch\'emas ab\'eliens ordinaires}

\noindent \textbf{D\'efinition (5.2.1) [I$\ell$ 2]}. \textit{Soient $S$ un sch\'ema de caract\'eristique $p > 0$ et $f : X \rightarrow S$ un morphisme propre et lisse. On dit que $X$ est ordinaire sur $S$ si}

$$
R^j f_{\ast} B  \Omega^i_{X/S} = 0
$$

\noindent \textit{pour tout $i$ et tout $j$ ($B \Omega^i_{X/S}$ d\'esigne le sous-faisceau\  $d\ \Omega^{i-1}_{X/S}$ de $\Omega^i_{X/S}$, et les $R^j f_{\ast}$ sont calcul\'es pour la topologie de Zariski).}\\

On a la caract\'erisation suivante :

\vskip 3mm
\noindent \textbf{Proposition (5.2.2)}. \textit{Soient $S$ un sch\'ema de caract\'eristique $p > 0$ et $f : X \rightarrow S$ un morphisme propre et lisse. Alors les conditions suivantes sont \'equivalentes :}

\begin{enumerate} 
\item[ \textit{(1)}]  \textit{$X$ est ordinaire sur $S$.
\item[(2)] Pour tout $s \in S$, $X_{s}$ est ordinaire sur $Spec\ k(s)$. 
\item[(3)] Pour tout point ferm\'e $s \in \vert S \vert$, $X_{s}$ est ordinaire sur $Spec\ k(s)$.
\item[(4)] Pour tout $s \in S$, localit\'e d'un point g\'eom\'etrique $\overline{s} = Spec\ \overline{k(s)}$ (o\`u $\overline{k(s)}$ est une cl\^oture alg\'ebrique de $k(s))$, $X_{\overline{s}}$ est ordinaire sur $Spec\ \overline{k(s)}$.}
\end{enumerate}

\vskip 3mm

\noindent \textit{D\'emonstration}. 

 L'\'equivalence de \textit{(1)} et \textit{(2)} r\'esulte de [I$\ell$ 2, 1.2(b)].\\
 
 Montrons que \textit{(3)} $\Rightarrow$ \textit{(1)}. Notons $\vert S \vert$ l'ensemble des points ferm\'es de $S$ : on suppose donc que, pour tout $s \in \vert S \vert$, $X_{s}$ est ordinaire sur $k(s) $. D'apr\`es [Bour, AC II, \S\ 3, n$\circ$\ 3, th\'eo 1] il s'agit de v\'erifier que pour tout $i$, $j$ et tout $s \in \vert S \vert$ on a $(R^i f_{\ast} B \Omega^i_{X/S}) \otimes_{\mathcal{O}_{S}} O_{S,s} = 0$ ; par platitude de $\mathcal{O}_{S,s}$ sur $\mathcal{O}_{S}$ c'est \'equivalent \`a montrer que $(\mathbb{R} f_{\ast} B \Omega^i_{X/S}) \displaystyle \mathop{\otimes}^{\mathbb{L}} {_{_{\mathcal{O}_{S}}}}  \mathcal{O}_{S,s} = 0$, ce qui r\'esulte de l'implication (i) $\Rightarrow$ (iii) de [I$\ell$2, 1.2 (b)].\\
  
 L'implication \textit{(2)} $\Rightarrow$ \textit{(4)} s'obtient par le changement de base $k(s) \rightarrow \overline{k(s)}$ [I$\ell$ 2, 1.2 (a)].\\
 
 R\'eciproquement, montrons \textit{(4)} $\Rightarrow$ \textit{(2)}. Soit $s \in S$, alors $X_{\overline{s}}$ est ordinaire par hypoth\`ese, d'o\`u par d\'efinition [I$\ell$ 2, 1.1]
 
 $$
 \mathbb{R} g_{\ast} B \Omega^i_{X_{\overline{s}/ \overline{s}}} = 0 
 $$
 
 \noindent o\`u $g : X_{\overline{s}} \rightarrow Spec\ \overline{k(s)}$ est le morphisme structural, d\'eduit de $f : X_{s} \rightarrow Spec\ k(s)$. Or d'apr\`es [I$\ell$ 2, ($\ast$) p 379] on a
 
 $$
 0 = \mathbb{R}  g_{\ast} B \Omega^i_{X_{\overline{s}/ \overline{s}}} \displaystyle \mathop{\longleftarrow}^{\sim}( \mathbb{R} f_{\ast}B \Omega^i_{X_{s/s}}) \displaystyle \mathop{\otimes}^{\mathbb{L}} {_k{_{(s)}}}\  \overline{k(s)}\ ;
 $$
 
 \noindent donc par fid\`ele platitude de $\overline{k(s)}$ sur $k(s)$ on obtient
 $$
 \mathbb{R}  f_{\ast} B \Omega^i_{X_{s/s}} = 0\ , 
 $$
 
 \noindent  ce qui est l'ordinarit\'e de $X_{s}$ sur $k(s)$. $\square$\\
 
 Nous adopterons la d\'efinition suivante due \`a Raynaud  [R 3, Rq 4.2.2] :
 
 \vskip 3mm
 \noindent \textbf{D\'efinition (5.2.3)}. \textit{Un groupe $p$-divisible $G$ est dit ordinaire s'il est extension d'un groupe $p$-divisible \'etale par un groupe $p$-divisible de type multiplicatif.}\\
 
 Pour un groupe $p$-divisible $G$ sur un sch\'ema $S$ on note, pour tout entier $n$, $G(n) : = \textrm{Ker} \{ p^n : G \rightarrow G \}$. Le dual de Cartier $G^{\ast}$ de $G$ est d\'efini par
 
 $$
G^{\ast} = \displaystyle \mathop{\lim}_{\rightarrow \atop{n}} (G(n))^{\ast}, \textrm{o\`u}\   (G(n))^{\ast} : = Hom(G(n), \mathbb{G}_{m})\ ;
 $$
 
 \noindent le dual de Pontryagin $G^{\vee}$ de $G$ est d\'efini par $G^{\vee} : = \displaystyle \mathop{\lim}_{\rightarrow \atop{n}}\ (G(n))^{\vee}$, o\`u $(G(n))^{\vee} : = Hom(G(n), \mathbb{Q}_{p}/\mathbb{Z}_{p})$.

 \vskip 2mm
 \noindent Si $S$ est localement annul\'e par une puissance de $p$ on note $\mathbb{D}(G)$ son cristal de Dieudonn\'e [B-M 1], [B-B-M]. Rappelons que $G$ est infinit\'esimal sur $S$ si, pour tout entier $n$, le morphisme structural $G(n) \rightarrow S$ est radiciel [Me, I, 3.2].

 \vskip 3mm
 \noindent \textbf{Lemme (5.2.4)}. \textit{Soient $S$ un sch\'ema de caract\'eristique $p > 0$ et $G$ un groupe $p$-divisible de type multiplicatif. Alors $G$ est infinitis\'emal sur $S$ et son dual de Cartier $G^{\ast}$ est \`a fibres unipotentes.}
 
 \vskip 3mm
 \noindent \textit{D\'emonstration}. Comme le caract\`ere radiciel se v\'erifie fibre \`a fibre [EGA I, 3.7.4] nous supposerons que $S$ est spectre d'un corps. Pour tout $n$, le Verschiebung $V$ de $G(n)$ est un isomorphisme et le n$^{\textrm{i\`eme}}$ it\'er\'e de son Frobenius $F$, donn\'e par $F^n = p^n(V^{-1})^n$, est nul : donc $G(n)$ est infinit\'esimal et son dual de Cartier $(G(n))^{\ast}$ est unipotent [D-G, IV, $\S\ 3$, n$\circ$\ 5.3]. $\square$\\
 
 Pour un sch\'ema $S$ localement annul\'e par une puissance de $p$ et $\Sigma = Spec\ \mathbb{Z}_{p}$, on note [B-M 1, \S\ 6], [B-M 2, 2.4] $C_{S}$ la cat\'egorie des cristaux localement de pr\'esentation finie sur CRIS(S/$\Sigma$), munis d'homomorphismes
 
 $$
 F : E^{\sigma} \rightarrow E \qquad , \qquad V : E \rightarrow E^{\sigma}
 $$
 
 \noindent tels que $F \circ V = p$, $ V \circ F = p$, et $C^{\textrm{\'et}}_{S}$ (resp. $C^{tm}_{S})$ la sous-cat\'egorie pleine de $C_{S}$ form\'ee des cristaux localement annul\'es par une puissance de $p$, et pour lesquels $F$ (resp. $V$) est un isomorphisme.
 
 \vskip 3mm
 \noindent \textbf{Proposition (5.2.5)}. \textit{Soient $S$ un sch\'ema de caract\'eristique $p > 0$ et $H$ un groupe $p$-divisible de type multiplicatif sur $S$ (i.e. $V$ est un isomorphisme). Alors on a un isomorphisme de $F$-cristaux}
 
 $$
\mathbb{D}(H) \simeq \mathbb{D}(H^{\ast \vee})(-1)\ ,
 $$

\noindent \textit{o\`u ( )$^{\ast}$ est la dualit\'e de Cartier, ( )$^{\vee}$ la dualit\'e de Pontryagin des groupes $p$-divisibles \'etales et $(-1)$ le twist \`a la Tate}.

\vskip 3mm
\noindent \textit{D\'emonstration}. Avec les notations de [B-M 2, 2.2.3.1] on  a

$$
\mathbb{D}(H) \simeq \underline{H}^{\ast} \otimes_{\mathbb{Z}_{p}} \mathcal{O}_{S/\Sigma} \qquad \textrm{et} \qquad F_{\mathbb{D}(H)} = p(F_{\underline{H}^{\ast}})^{-1} \otimes 1\ ;
$$

\noindent de m\^eme [B-M 2, (2.1.3.1) et (2.1.4) (ii)] on a

$$
\mathbb{D}(H^{\ast {\vee}}) \simeq \underline{H}^{\ast} \otimes_{\mathbb{Z}_{p}} \mathcal{O}_{S/\Sigma} \qquad \textrm{et} \qquad F_{\mathbb{D}(H ^{\ast {\vee}})} = (F_{\underline{H}^{\ast}})^{-1} \otimes 1\ ,
$$

\noindent d'o\`u la proposition. $\square$

\vskip 3mm
\noindent \textbf{Proposition (5.2.6)}. \textit{Soient $S$ un sch\'ema de caract\'eristique $p > 0$ et $E$ un \'el\'ement de $C_{S}$. Les assertions suivantes sont \'equivalentes :}

\begin{enumerate}
\item[ \textit{(1)}] \textit{$E$ est un $F$-cristal ordinaire de niveau 1 [cf. (5.1)].
\item[(2)] Il existe un groupe $p$-divisible \'etale $G^{\textrm{\'et}}$ sur $S$ et un groupe $p$-divisible de type multiplicatif  $G^{tm}$ sur $S$ tels que $E$ soit extension de $\mathbb{D}(G^{tm})$ par  $\mathbb{D}(G^{\textrm{\'et}})$.}
\end{enumerate}

\vskip 3mm
\noindent \textit{D\'emonstration}. Ceci r\'esulte de la d\'efinition (5.1), de la proposition (5.2.5) pr\'ec\'edente et des \'equivalences de cat\'egories de [B-M 1, corollaire du th\'eor\`eme 6]. $\square$

\vskip 3mm
\noindent \textbf{Proposition (5.2.7)}. \textit{Soient $S$ un sch\'ema de caract\'eristique $p > 0$, normal, localement irr\'eductible et poss\'edant localement une $p$-base et $G$ un groupe $p$-divisible sur $S$. Alors les assertions suivantes sont \'equivalentes :}

\begin{enumerate}
\item[ \textit{(1)}] \textit{$G$ est un groupe $p$-divisible ordinaire.
\item[(2)] $\mathbb{D}(G)$ est un $F$-cristal ordinaire de niveau 1}.
\end{enumerate}

\vskip 3mm
\noindent \textit{D\'emonstration}. Si $G$ est ordinaire, on a une suite exacte

$$
0 \longrightarrow G^{tm} \longrightarrow G \longrightarrow G^{\textrm{\'et}} \longrightarrow 0
$$

\noindent o\`u $G^{\textrm{\'et}}$ est un groupe $p$-divisible \'etale et $G^{tm}$ un groupe $p$-divisible de type multiplicatif, ce qui fournit une suite exacte [B-M 1, \S\ 4.1 a)]

$$
0 \longrightarrow \mathbb{D}(G^{\textrm{\'et}}) \longrightarrow \mathbb{D}(G) \longrightarrow \mathbb{D}(G^{tm}) \longrightarrow 0\ ,
$$

\noindent donc $\mathbb{D}(G)$ est ordinaire de niveau 1 par la proposition (5.2.6) pr\'ec\'edente.\\

R\'eciproquement, supposons que $\mathbb{D}(G)$ est un $F$- cristal ordinaire de niveau 1, i.e. par la proposition (5.2.6) (dont nous adoptons les notations), que $\mathbb{D}(G)$ est extension de $\mathbb{D}(G^{tm})$ par $\mathbb{D}(G^{\textrm{\'et}})$. Sous nos hypoth\`eses, le foncteur $\mathbb{D}$ est pleinement fid\`ele [B-M 2, th\'eo (4.1.1)], donc la suite exacte

$$
0 \longrightarrow \mathbb{D}(G^{\textrm{\'et}}) \longrightarrow \mathbb{D}(G) \longrightarrow \mathbb{D}(G^{tm}) \longrightarrow 0\ ,
$$

\noindent fournit une suite exacte

$$
0 \longrightarrow G^{tm} \longrightarrow G \longrightarrow G^{\textrm{\'et}} \longrightarrow 0,
$$

\noindent et $G$ est ordinaire. $\square$

\vskip 3mm
\noindent \textbf{Remarques (5.2.8) sur la proposition (5.2.7)}.

 \begin{enumerate}
\item[(1)] Les hypoth\`eses sur $S$ sont satisfaites si $S$ est un sch\'ema sur un corps $k$ de caract\'eristique $p > 0$ tel que $S$ est de l'un des deux types suivants [B-M 1, \S\ 6] :
\begin{itemize}
\item[(i)] $S$ est lisse sur $k$.
\item[(ii)] $S = Spec\ k[[X_{1},..., X_{n}]]$, avec $[k : k^p] < + \infty$.
\end{itemize}
\item[(2)] Le (1) de la proposition implique le (2) en supposant seulement que $S$ est un $k$-sch\'ema ; ce n'est que pour la r\'eciproque que nos hypoth\`eses sur $S$ sont utiles.
\end{enumerate}

\vskip 3mm
\noindent \textbf{Th\'eor\`eme (5.2.9)}. \textit{Soient $S$ un sch\'ema de caract\'eristique $p$, $f : X \rightarrow S$ un sch\'ema ab\'elien de dimension relative $d$, $G$ le groupe $p$-divisible associ\'e \`a $X/S$ et $\mathbb{D}(G)$ son cristal de Dieudonn\'e [B-B-M]. Alors les propri\'et\'es (1) \`a (5) ci-dessous sont \'equivalentes :
\begin{enumerate}
\item[(1)] $X$ est ordinaire sur $S$.
\item[(2)] Pour tout point g\'eom\'etrique $\overline{s}$ de $S$ le polygone de Newton $Nwt_{n}$ d\'efini par les pentes de Frobenius agissant sur le groupe de cohomologie cristalline $H^n_{\textrm{cris}}(X_{\overline{s}}/W(k(\overline{s})))$ co¬\"{\i}ncide, pour tout $n$, avec le polygone de Hodge $Hdg_{n}$ d\'efini par les nombres de Hodge
$$
h^{i,n-i} = dim_{k(\overline{s})}\ H^{n-i}(X_{\overline{s}}, \Omega^i_{X_{\overline{s}/k(\overline{s})}})
$$
(i.e ayant pour pente $i$ avec la multiplicit\'e $h^{i,n-i})$.
\item[(3)] Pour tout point g\'eom\'etrique $\overline{s}$ de $S$, le $p$-rang de $X_{\overline{s}}$ est $d$ (par d\'efinition le $p$-rang est \'egal \`a $dim_{\mathbb{Z}/p \mathbb{Z}}(H^1_{\textrm{\'et}}(X_{\overline{s}}, \mathbb{Z}/p \mathbb{Z})))$.
\item[(4)] $G$ est ordinaire.
\item[(5)] Le $F$-cristal de Dieudonn\'e $\mathbb{D}(G) \simeq \mathbb{D}(X) : = R^1 f_{\textrm{cris} \ast}(\mathcal{O}_{X/\Sigma})$ est ordinaire de niveau 1.
\end{enumerate}}

\noindent \textit{Si de plus $S$ est un sch\'ema sur un corps parfait $k$ de caract\'eristique $p > 0$ et si $S$ est de l'un des deux types suivants :}
\textit{
\begin{itemize}
\item[(i)] $S$ est lisse sur $k$,
\item[(ii)] $S = Spec\ k[[t_{1} ,..., t_{n}]]$,
\end{itemize}}
\noindent \textit{alors ces propri\'et\'es \'equivalent aussi \`a :
\begin{itemize}
\item[(6)] Pour tout $i$, $0 \leqslant i \leqslant 2d$, le $F$-cristal $R^i  f_{\textrm{cris}^{\ast}}(\mathcal{O}_{X/\Sigma}) = R^i  f_{\textrm{cris}^{\ast}}(\mathcal{O}_{X/W(k)})$ est ordinaire de niveau $i$.
\end{itemize}}

\vskip 3mm
\noindent \textit{D\'emonstration}. L'\'equivalence de  \textit{(1)},  \textit{(2)} et  \textit{(3)} r\'esulte du (4) de la proposition (5.2.2) de [I$\ell$ 2, \S\ 1] et [B-K, 7.2, 7.3, 7.4].\\

Montrons que  \textit{(1)} $\Rightarrow$  \textit{(4)}. Soient $s = Spec\ k(s)$ un point de $S$ et $\overline{s} = Spec\ \overline{k(s)}$, o\`u $\overline{k(s)}$ est une cl\^oture alg\'ebrique de $k(s)$. Comme $X$ est ordinaire sur $S$, $X_{\overline{s}}$ est une vari\'et\'e ab\'elienne ordinaire [I$\ell$ 2, 1.2 (a)] au sens usuel ([I$\ell$ 2, 1.1] et [B-K, 7.2 et 7.4]), donc le $p$-rang de $X_{\overline{s}}$ [Mu, p. 147] est $d$ [B-K, 7.4] : par suite [Mu, p. 147], sur le corps $\overline{k(s)}$, le groupe $p$-divisible $G_{\overline{s}}$ est le produit d'un groupe $p$-divisible \'etale $G^{\textrm{\'et}}_{\overline{s}}$ par un groupe $p$-divisible de type multiplicatif $G^{tm}_{\overline{s}}$ ; par dualit\'e de Cartier on a aussi $G^{\ast}_{\overline{s}}
 \simeq G^{tm\ \ast}_{\overline{s}} \times G^{\textrm{\'et}\ \ast}_{\overline{s}}$. Or le rang s\'eparable de $G(1)_{s}$ [EGA I, 6.5.9] est aussi celui de $G(1)_{\overline{s}}$ [EGA I, 6.5.11], et comme $G^{tm}$ est infinit\'esimal [lemme (5.2.4)] ce rang s\'eparable de $G(1)_{\overline{s}}$ est le rang de $G(1)^{\textrm{\'et}}_{\overline{s}}$, qui est \'egal \`a $p^d$ puisque $X_{\overline{s}}$ est ordinaire [Mu, p. 147]. Ainsi la fonction $s \mapsto$ rang s\'eparable de $G(1)_{s}$ (resp. de $G(1)^{\ast}_{s})$ est constante : 
par [Me, II, prop 4.9] $G$ (resp. $G^{\ast}$) est extension d'un groupe $p$-divisible \'etale $G''$  
(resp. $(G^{\ast})''$) par un groupe $p$-divisible infinit\'esimal $G'$  (resp. $(G^{\ast})'$).
La dualit\'e de Cartier pr\'eservant les suites exactes, on a les deux suites exactes :\\

 \noindent (S1) $\qquad \qquad \qquad 0  \longrightarrow ((G^{\ast})'')^{\ast}  \displaystyle \mathop{\longrightarrow}^i  G \longrightarrow ((G^{\ast})')^{\ast} \longrightarrow 0$\\
 
  \noindent (S2) $\qquad \qquad \qquad \qquad  0  \longrightarrow G'  \longrightarrow G \displaystyle \mathop{\longrightarrow}^j  {G''}  \longrightarrow  0.$\\
  
 \noindent Comme $(G'')^{\ast}$ est de type multiplicatif, il est infinit\'esimal [lemme (5.2.4] donc $G''$ est \`a fibres unipotentes [D-G, IV, \S\ 3, n$\circ$\ 5.3] : il r\'esulte alors de [SGA 3, XVII, \S\ 2, lemme 2.5] que le morphisme compos\'e $j \circ i$ est nul. Par suite, l'identit\'e de $G$ induit un morphisme de la suite exacte (S1) dans la suite exacte (S2) ; soit $f : ((G^{\ast})')^{\ast} \rightarrow G''$ le morphisme induit. Pour tout entier $n$, le morphisme $G(n) \rightarrow G''(n)$ est un \'epimorphisme de $S$-sch\'emas en groupes, donc il est fid\`element plat : en effet, par le crit\`ere de platitude fibre par fibre [EGA IV, 11.3.11] on est ramen\'e au cas o\`u $S$ est le spectre d'un corps, et alors le r\'esultat est standard [D-G, III, \S\ 3, 7.4]. De la m\^eme fa\c{c}on $G(n) \rightarrow ((G^{\ast}(n))')^{\ast}$ est fid\`element plat, donc $f$ est fid\`element plat [EGA IV, 2.2.13] :  montrons que $f$ est un isomorphisme, i.e. que pour tout $n$, $f_{n} : ((G^{\ast}(n))')^{\ast} \rightarrow G''(n)$ est un isomorphisme. Pour cela on peut supposer que $S$ est le spectre $s$ d'un corps $k$  [EGA IV, 17.9.5], et m\^eme que $k = \overline{k}$ est alg\'ebriquement clos, par fid\`ele platitude de $\overline{k}$ sur $k$ [EGA IV, 2.7.1, (viii)].  Mais alors $(G^{\ast}_{\overline{s}})' \simeq G^{\textrm{\'et} \ast}_{\overline{s}}$, d'o\`u $((G^{\ast}_{\overline{s}})')^{\ast} \simeq G^{\textrm{\'et}}_{\overline{s}}$ qui est de hauteur $ht\  (G^{\textrm{\'et}}_{\overline{s}}) = d$ car $X_{\overline{s}}$ est ordinaire ; de m\^eme $G''_{\overline{s}} \simeq G^{\textrm{\'et}}_{\overline{s}}$ et $f_{n}$ induit un morphisme fid\`element plat $f_{n, \overline{s}} : G^{\textrm{\'et}}_{\overline{s}}(n) \rightarrow G^{\textrm{\'et}}_{\overline{s}}(n)$ : c'est donc un isomorphisme puisqu'au niveau des anneaux il fournit une injection entre $\overline{k}$-espaces vectoriels de m\^eme dimension $p^{nd}$. Ainsi $f$ est un isomorphisme et l'identit\'e de $G$ induit un isomorphisme

$$
((G^{\ast})'')^{\ast} \tilde{\longrightarrow}\  G'\ ;
$$

\noindent puisque $(G^{\ast})''$ est \'etale, $G'$ est de type multiplicatif et donc $G$ est ordinaire.\\

L'implication  \textit{(4)} $\Rightarrow$  \textit{(5)} provient de la remarque (2) de (5.2.8).\\

Montrons que  \textit{(5)} $\Rightarrow$  \textit{(3)}. Supposant  \textit{(5)} on a une suite exacte de cristaux sur $S$ [prop (5.2.6)]

$$
0 \longrightarrow \mathbb{D} (G^{\textrm{\'et}}) \longrightarrow \mathbb{D}(G) \longrightarrow \mathbb{D} (G^{tm}) \longrightarrow 0\ ,
$$

\noindent qui, pour tout point g\'eom\'etrique $\overline{s}$ de $S$, reste exacte par image inverse sur $Spec\ k(\overline{s})$ et s'identifie [B-B-M, 1.3.3] \`a

$$
0 \longrightarrow \mathbb{D} (G^{\textrm{\'et}}_{\overline{s}}) \longrightarrow \mathbb{D}(G_{\overline{s}}) \longrightarrow \mathbb{D} (G^{tm}_{\overline{s}}) \longrightarrow 0\ .
$$

\noindent Le corps $k(\overline{s})$ v\'erifiant les hypoth\`eses de la proposition (5.2.7) on en d\'eduit une suite exacte

$$
0 \longrightarrow G^{tm}_{\overline{s}} \longrightarrow G_{\overline{s}} \longrightarrow G^{\textrm{\'et}}_{\overline{s}} \longrightarrow 0\ , 
$$

\noindent qui est scind\'ee puisque $k(\overline{s})$ est parfait. Ainsi la composante infinit\'esimale unipotente de  $G_{\overline{s}}$ est nulle [Dem,  p.39] ; gr\^ace \`a [Mu, p. 147] on en d\'eduit que le $p$-rang de $X_{\overline{s}}$ est $d$, et ceci pour tout point g\'eom\'etrique $\overline{s}$, d'o\`u le  \textit{(3)}.\\

Il reste \`a \'etablir l'implication  \textit{(5)} $\Rightarrow$  \textit{(6)}. Sous les conditions de l'\'enonc\'e l'ordinarit\'e se v\'erifie sur les polygones de Hodge et de Newton de $R^i f_{\textrm{cris} \ast} (\mathcal{O}_{S/\Sigma})$ en chaque point $s \in S$ [cf. (5.1)] : la co¬\"{\i}ncidence de ces polygones pour un $i$, $0 \leqslant i \leqslant 2d$, r\'esulte de leur co¬\"{\i}ncidence pour $i = 1$ [K 2, I, 1.2, 1.3], car on a un isomorphisme [B-B-M, (2.5.5.1)]

$$
R^i f_{\textrm{cris}\ast}(\mathcal{O}_{S/\Sigma}) \simeq \displaystyle \mathop{\Lambda}^i R^1 f_{\textrm{cris} \ast}(\mathcal{O}_{S/\Sigma})\ .
$$
 
 \noindent L'assertion sur le niveau est claire. $ \square$
 
 \vskip 3mm
\subsection*{5.3. Explicitation des fonctions $L^{(r)}_{\alpha}$ et conjecture de Dwork}
 
 On reprend les hypoth\`eses et notations du \S\ 4. Ainsi $k = \mathbb{F}_{q}$ et $S$ est un $k$-sch\'ema s\'epar\'e de type fini. On supposera dor\'enavant que $f : X \longrightarrow S$ est un sch\'ema ab\'elien ordinaire de dimension relative $g$, de groupe $p$-divisible $G$ : par le th\'eor\`eme (5.2.9), $G$ est extension d'un groupe $p$-divisible \'etale $G^{\textrm{\'et}}$ par un groupe $p$-divisible de type multiplicatif $G^{tm}$.\\
 
 Nous omettrons dor\'enavant $S$ et $t$ dans l'\'ecriture des fonctions $L$, par exemple $L_{\alpha}(E) : = L_{\alpha}(S, E, t)$ pour $\alpha \in \mathbb{Q}$.

 \vskip 3mm
 \noindent \textbf{Th\'eor\`eme (5.3.1)}. \textit{Sous les hypoth\`eses 5.3 on a, pour tout $i$, $0 \leqslant i \leqslant 2g$ :} 
\begin{enumerate}
\item[(1)] 
	\begin{itemize}
	\item[(1.1)] $L_{0}(R^i f_{\textrm{cris}\ast}(\mathcal{O}_{X/W})) = L(R^i f_{\textrm{\'et}\ast} (\mathbb		{Z}_{p}) \otimes_{\mathbb{Z}_{p}} \mathcal{O}_{X/W})$\  ,  \textit{o\`u}
	$$R^i f_{\textrm{\'et}\ast} (\mathbb{Z}_{p}) \otimes_{\mathbb{Z}_{p}} \mathcal{O}_{X/W} \simeq 		\displaystyle \mathop{\Lambda}^i \mathbb{D}(G^{\textrm{\'et}})$$ 
	\textit{est le sous F-cristal unit\'e de} $ R^i f_{\textrm{cris}\ast}(\mathcal{O}_{X/W}).$
	\item[(1.2)] $L_{i}(R^i f_{\textrm{cris}\ast}(\mathcal{O}_{X/W})) = L(\displaystyle \mathop{\Lambda}^i 	\mathbb{D}(G^{tm})).$
	\item[(1.3)] $L_{\alpha} (R^i f_{\textrm{cris}\ast}(\mathcal{O}_{X/W})) = 1\ si\  \alpha < 0\  ou\  \alpha 		> i, ou\  \alpha \in [0,i] \setminus \mathbb{N}$,\\
	 $$\qquad\qquad\qquad\qquad\qquad = L((\displaystyle \mathop{\Lambda}^{i - \alpha} \mathbb{D}(G^{\textrm		{\'et}})) \otimes (\displaystyle \mathop{\Lambda}^{\alpha} \mathbb{D} (G^{tm})))\ si\ \alpha \in [0, i] 		\cap \mathbb{N}\ .$$
	\end{itemize}
\item[(2)] \textit{Pour tout $\alpha \in \mathbb{Q}$ et $r \in \mathbb{N}^{\ast}$ on a:}\\

	\begin{itemize}
	\item[(2.1)] $L^{(r)}_{\alpha} (R^i f_{\textrm{cris}\ast}(\mathcal{O}_{X/W}))$ \textit{est $p$-adiquement m\'eromorphe}.\\
	
	\item[(2.2)] $L^{(r)}_{\alpha} (R^i f_{\textrm{cris}\ast}(\mathcal{O}_{X/W})) = 1\  si\  \alpha < 0, ou\  		\alpha > i, ou\   \alpha \in [0, i] \setminus \mathbb{N}.$\\
	
	\item[(2.3)] $L^{(r)}_{\alpha} (R^i f_{\textrm{cris}\ast}(\mathcal{O}_{X/W})) =\\
	 \displaystyle \mathop{\prod}_{j \geqslant 1} L \{Sym^{r-j}[(\displaystyle \mathop{\Lambda}^{i- \alpha}  		\mathbb{D}(G^{\textrm{\'et}})) \otimes (\displaystyle \mathop{\Lambda}^{\alpha} \mathbb{D}(G^		{tm}))]  \otimes
	\displaystyle \mathop{\Lambda}^{j} [(\displaystyle \mathop{\Lambda}^{i- \alpha}  \mathbb{D}(G^		{\textrm{\'et}})) \otimes (\displaystyle \mathop{\Lambda}^{\alpha} \mathbb{D}(G^{tm}))] \} ^{j \times		(-1)^{j-1}}\\
 	si\  \alpha \in [0, i] \capÊ\mathbb{N}\ .$
	\end{itemize}
\end{enumerate}

\vskip 3mm
\noindent \textit{D\'emonstration}. Pour (1.1) on rappelle l'isomorphisme [B-B-M]

$$
R^i f_{\textrm{cris}\ast}(\mathcal{O}_{X/W})) \simeq \displaystyle \mathop{\Lambda}^{i} R^1 f_{\textrm{cris}\ast} (\mathcal{O}_{X/W})) \simeq \displaystyle \mathop{\Lambda}^i \mathbb{D}(G)\ ,
$$

\noindent et il r\'esulte de [E-LS 2, prop 6.5] que

$$
L_{0}(R^i f_{\textrm{cris}\ast}(\mathcal{O}_{X/W})) = L(R^i f_{\textrm{\'et}\ast} (\mathbb{Z}_{p}))\ ;
$$

\noindent pour $i = 1$, $R^1 f_{\textrm{\'et}\ast} (\mathbb{Z}_{p})$ et $\mathbb{D}(G^{\textrm{\'et}})$ ont m\^eme fonction $L$, ce qui prouve (1.1).\\

Pour le (1.2) et (1.3) il suffit de remarquer que

$$
L(R^1 f_{\textrm{cris}\ast}(\mathcal{O}_{X/W})) = L(\mathbb{D}(G)) = L(\mathbb{D}(G^{\textrm{\'et}})) \times L(\mathbb{D}(G^{tm}
))\ .$$

\noindent Compte tenu de [Et 5, th\'eo 8] le (2.1) est prouv\'e en (4.2.1)(i) ; (2.2) r\'esulte de (1.3).\\

Pour le (2.3) on utilise l'expression (1.3) du th\'eor\`eme et la formule (3.2.5). $\square$


\cleardoublepage
 \bibliographystyle{plain}


\end{document}